\documentclass[11pt]{amsart}

\usepackage{amsmath}
\usepackage{amssymb}
\usepackage{amsthm}
\usepackage{mathtools}

\usepackage{mathrsfs}
\usepackage{bm}

\usepackage{enumitem}
\usepackage[protrusion=true,expansion=false]{microtype}

\usepackage{graphicx}
\usepackage{booktabs}

\usepackage{xcolor}
\usepackage{listings}
\usepackage{needspace}

\definecolor{codebg}{HTML}{F6F7F9}
\definecolor{codeframe}{HTML}{CBD2D9}

\lstdefinelanguage{Isabelle}{
  morekeywords={theory,imports,begin,end,type_synonym,record,definition,where,theorem,shows,locale,for,fixes,assumes},
  sensitive=true,
  morecomment=[s]{(*}{*)},
  morestring=[b]"
}

\usepackage[
    textsize=small
]{todonotes}

\usepackage[
    colorlinks=true,
    linkcolor=blue,
    citecolor=blue,
    urlcolor=blue
]{hyperref}

\usepackage[nameinlink,capitalise,noabbrev]{cleveref}

\theoremstyle{plain}

\theoremstyle{definition}

\theoremstyle{remark}

\newcommand{\R}{\mathbb{R}}

\newcommand{\DN}{\Lambda}

\newcommand{\reviewlabel}[1]{\par\smallskip\noindent\textsf{\bfseries #1}\par\smallskip}
\newcommand{\dossierentry}[1]{\par\medskip\noindent{\bfseries #1}\par\smallskip}
\newcommand{\dossiersubhead}[1]{\par\bigskip\noindent{\itshape #1}\par\smallskip}
\newcommand{\dossiersection}[1]{\par\bigskip\noindent{\large\bfseries #1}\par\medskip}

\title{Relative formalization in Isabelle/HOL of a result in inverse problems}

\author{C\u{a}t\u{a}lin I. C\^{a}rstea}

\address{Department of Applied Mathematics,
National Yang Ming Chiao Tung University,
Hsinchu 30050, Taiwan}

\email{catalin.carstea@gmail.com}
\date{\today}

\begin{document}

\begin{abstract}
This reports on an experiment in autoformalization of \\ arxiv:2606.15977, an inverse problems result for piecewise polynomial anisotropic conductivities, using Isabelle/HOL. The proof of the result is relative to a number of external results which were deemed to be well known and generally accepted to be true. The formalization files are made available at a GitHub repository. Translation issues are discussed.
\end{abstract}

\maketitle

%

\section{Introduction}

It has been observed \cite{Tao2026AgeAI,PetrovEtAl2026,MaEtAl2026} that AI generated mathematical prose may, at times, be harder to read and to check than more traditionally written books and papers. This leads to the justified concern that false results may now enter the published record at higher rates than in the past \cite{Tao2026AgeAI}. One of the solutions that are beginning to arise in response is the use of formalization, which is now more accessible to more researchers, also thanks to the availability of AI tools \cite{Tao2026AgeAI,WengEtAl2025}. Most recent reports use Lean, including formalizations of modern PDE results \cite{ArmstrongKempe2026,Miller2026}, AI-assisted work on Erd\H{o}s problems \cite{Sothanaphan2026,ZhangEtAlLeanMarathon2026}, and collections of formalized research papers in other areas \cite{Garg2026}; see also the survey \cite{WengEtAl2025}. Isabelle/HOL has also been used \cite{WuEtAl2022,ZhangValentinoFreitas2025,BryantEtAl2026,KappelmannEtAl2026}, as well as other systems such as Coq/Rocq \cite{CunninghamEtAl2022,ZhangEtAlKELPS2025}.

As an experiment and learning exercise I have attempted to formalize one of my recent preprints \cite{Carstea2026Polynomial}. Anecdotally, formalization in Lean appears to be somewhat expensive in terms of both time and money. For this experiment I have opted for Isabelle/HOL as the formalization language, and have used Codex and ChatGPT 5.6 Sol with a \$200/month subscription to produce the formalization. The final result was produced in under a week, which would put the monetary cost of the exercise at about \$50. At these costs formalizing a paper becomes more feasible as a routine additional verification step.

This paper and the experiment it reports on is meant as a proof of concept. It undoubtedly has many weaknesses. I gladly welcome any constructive criticism and suggestions.

\section{What was done}

The formalization translates the main theorem, with its hypotheses and relevant definitions into Isabelle/HOL, then reconstructs the proof as a sequence of machine-checked lemmas. One key decision taken was to assume, rather than prove, certain results taken from literature. Without this measure the formalization work would have probably taken considerably more time. Each cited result is also translated into Isabelle and appears within an Isabelle ``locale''. These assumptions are  visible in the dependency structure. The final theorem records that the claimed conclusion follows from Isabelle’s logical foundations and libraries, together with this explicitly identified literature base. The formalization was considered complete when Isabelle checked the exact intended theorem (the Isabelle language version!) from beginning to end, with no omitted proofs or provisional assumptions, and when every dependency belonged either to the checked Isabelle libraries, to a proved component of the project, or to one of the declared literature interfaces. 

The workflow used treats the formalization as a single, extended, long running project. Codex provided the general ability to read the mathematical source, search the Isabelle libraries, write and check formal proofs, and use available tools. I supplemented this with a project-specific harness, consisting mainly of a library of reusable instructions and a collection of project files specifying rules of operation, available tools, templates, and the mathematical context that was to be preserved. A useful feature of the setup is that it can be given relatively long tasks and be allowed to work on them with little supervision.  The surrounding harness provided continuity and discipline for these longer runs, including instructions about what information was to be treated as fixed and when a question should be left for human review. This made it possible to carry out a substantial amount of formalization and checking without continuous interaction, while leaving the resulting Isabelle files and a record of the intermediate work available for later inspection.

The formalization code can be found at \href{https://github.com/catalin-carstea/arxiv-2606.15977-formalization}{\nolinkurl{https://github.com/catalin-carstea/arxiv-2606.15977-formalization}}.

\section{The translation problem}

Besides the possibility of bugs in the Isabelle kernel, the approach taken here has a clear weakness: translation. The total size of the Isabelle code for the theorem proven, the definitions and imported literature results is about 50KB, or about 1,400 lines of code. This would be hard to check for a human reader in a reasonable amount of time, even with perfect understanding of Isabelle language. I have instead resorted to a number of other approaches which may not lead to perfect confidence in the accuracy of what was done, but may at least suggest that the effort has some significance.

It may be worth noting that a small mistake in translating the main theorem, unless it reduces it to a known result or a trivial one, might be less serious if it only means that the English language theorem is slightly misstated. A small mistake in translating a known result from a monograph, on the other hand, could in principle invalidate the whole project, as such results are often optimal or nearly optimal.

\subsection*{Translation testing approaches}\mbox{ }
\bigskip

\paragraph{1.} All translations have been checked and approved by both ChatGPT 5.6 Sol Pro (in the chat, outside the Codex context where they were produced) and by Claude Opus 5.
\bigskip

\paragraph{2.} The translation of the main theorem was done first, before any of the other formalization work, and was kept fixed throughout. This measure was intended to prevent the LLM from independently altering the statement until it got something that would allow it to more quickly finish its task \cite{ZhangEtAl2026BeyondCompilation,CornishEtAl2026FaithformBench}. Translations of cited results were also frozen after being made, for the same reason. Incidentally, in this particular experiment, it became necessary to revise the translations once during the formalization process. The version of the paper that was given to Codex was imprecise regarding what is meant by a ``Lipschitz domain'', and the initial translation took it in a sense that was different from what one of the cited results used. 
\bigskip

\paragraph{3.} In order to check the translations of the literature results, I have asked ChatGPT 5.6 Sol in Codex (effort: extra high) to try for 24 hours to derive a contradiction between the translated statements and the Isabelle libraries\footnote{I would like to thank Matti Lassas for this idea.}. It found none.
\bigskip

\paragraph{4.} Finally, I have had Codex translate the Isabelle code versions of the definitions and statements of the main theorem and imported literature results back into English, using subagents that were not given access to the original manuscripts. Though not perfect, this method does should give a degree of confidence in the accuracy of translations, as two translation errors canceling each other should be a unlikely  To make the approach more robust, this was done twice, independently. The results of this effort are included as Appendices~\ref{app:translation-a} and~\ref{app:translation-b}, in the form of interlinear translations with Isabelle code and its English translation for definitions, and English original, Isabelle code, and translation back to English for theorem statements.

\section{Provisional conclusions}

In a strict mathematical sense, without perfect confidence in the correspondence between natural language mathematical statements and their corresponding Isabelle code translations this experiment proves little. It is my personal impression however that, together with a sufficient number of varied translation tests, the existence of a formalization that checks does raise the level of confidence in the result somewhat. To be sure, even a traditionally written and peer reviewed paper may on occasion turn out to contain errors. A question one may ask, even though it is probably ill-posed, is then: can formalization raise the confidence level at least to match what would usually be assigned to traditional work? This seems to require more reflection and experimentation.

On the practical side, the experiment does show that relative formalization in Isabelle/HOL is achievable at least in some cases. With more future attempts we may end up in a better position to answer the more philosophical questions above.

\bigskip
\paragraph*{{\bf Acknowledgements}}
I am grateful to Lauri Oksanen for many constructive discussions. This work was supported by NSTC grant 113-2115-M-A49-018-MY3.

\clearpage

\appendix

\clearpage
\section{Independent back-translation A}\label{app:translation-a}
\begingroup
\small
\emergencystretch=3em
\sloppy
\dossiersection{Definitions}

\dossierentry{\texttt{\detokenize{cgu_point}}}
\reviewlabel{Isabelle code}
\begin{lstlisting}
type_synonym 'n cgu_point = "real ^ 'n"
\end{lstlisting}
\reviewlabel{English mathematical translation}
\noindent\textit{Mathematical role:} Euclidean vector type \(\mathbb R^n\)\par\smallskip
For any index type \texttt{'n}, \texttt{cgu\_\allowbreak{}point} is the type of real-valued vectors indexed by \texttt{'n}, namely \texttt{real\allowbreak\ \textasciicircum{}\allowbreak\ 'n}.

\dossierentry{\texttt{\detokenize{cgu_monomial}}}
\reviewlabel{Isabelle code}
\begin{lstlisting}
type_synonym 'n cgu_monomial = "'n \<Rightarrow>\<^sub>0 nat"
\end{lstlisting}
\reviewlabel{English mathematical translation}
\noindent\textit{Mathematical role:} finitely supported multi-index \(\alpha\)\par\smallskip
For any type \texttt{'n}, \texttt{cgu\_\allowbreak{}monomial} is the type of finitely supported mappings from \texttt{'n} to natural numbers.

\dossierentry{\texttt{\detokenize{cgu_polynomial}}}
\reviewlabel{Isabelle code}
\begin{lstlisting}
type_synonym 'n cgu_polynomial = "'n cgu_monomial \<Rightarrow>\<^sub>0 real"
\end{lstlisting}
\reviewlabel{English mathematical translation}
\noindent\textit{Mathematical role:} finitely supported scalar coefficient family \((p_\alpha)\)\par\smallskip
For any type \texttt{'n}, \texttt{cgu\_\allowbreak{}polynomial} is the type of finitely supported mappings from objects of type \texttt{cgu\_\allowbreak{}monomial} to real numbers.

\dossierentry{\texttt{\detokenize{cgu_matrix_polynomial}}}
\reviewlabel{Isabelle code}
\begin{lstlisting}
type_synonym 'n cgu_matrix_polynomial = "'n cgu_polynomial ^ 'n ^ 'n"
\end{lstlisting}
\reviewlabel{English mathematical translation}
\noindent\textit{Mathematical role:} coefficient-represented matrix polynomial \(P=(p_{ij})\)\par\smallskip
For any type \texttt{'n}, \texttt{cgu\_\allowbreak{}matrix\_\allowbreak{}polynomial} is an \texttt{'n}-by-\texttt{'n} vector-of-vectors whose entries have type \texttt{cgu\_\allowbreak{}polynomial}.

\dossierentry{\texttt{\detokenize{cgu_monomial_total_degree}}}
\reviewlabel{Isabelle code}
\begin{lstlisting}
definition cgu_monomial_total_degree :: "'n cgu_monomial \<Rightarrow> nat"
where
  "cgu_monomial_total_degree m =
    (\<Sum>i\<in>Poly_Mapping.keys m. Poly_Mapping.lookup m i)"
\end{lstlisting}
\reviewlabel{English mathematical translation}
\noindent\textit{Mathematical role:} total degree of a multi-index\par\smallskip
For \texttt{m\allowbreak\ :\allowbreak\ cgu\_\allowbreak{}monomial}, \texttt{cgu\_\allowbreak{}monomial\_\allowbreak{}total\_\allowbreak{}degree\allowbreak\ m} is the natural-number sum, over all keys in the finite support of \texttt{m}, of the value looked up at that key.

\dossierentry{\texttt{\detokenize{cgu_poly_degree_le}}}
\reviewlabel{Isabelle code}
\begin{lstlisting}
definition cgu_poly_degree_le :: "nat \<Rightarrow> 'n cgu_polynomial \<Rightarrow> bool"
where
  "cgu_poly_degree_le d p \<longleftrightarrow>
    (\<forall>m\<in>Poly_Mapping.keys p. cgu_monomial_total_degree m \<le> d)"
\end{lstlisting}
\reviewlabel{English mathematical translation}
\noindent\textit{Mathematical role:} stored degree predicate \(\deg_*p\le d\)\par\smallskip
\texttt{cgu\_\allowbreak{}poly\_\allowbreak{}degree\_\allowbreak{}le\allowbreak\ d\allowbreak\ p} holds exactly when every \texttt{m} in the finite support of \texttt{p\allowbreak\ :\allowbreak\ cgu\_\allowbreak{}polynomial} satisfies \texttt{cgu\_\allowbreak{}monomial\_\allowbreak{}total\_\allowbreak{}degree\allowbreak\ m\allowbreak\ <=\allowbreak\ d}. The universal condition is vacuous when the support of \texttt{p} is empty.

\dossierentry{\texttt{\detokenize{cgu_monomial_eval}}}
\reviewlabel{Isabelle code}
\begin{lstlisting}
definition cgu_monomial_eval ::
  "'n cgu_monomial \<Rightarrow> 'n::finite cgu_point \<Rightarrow> real"
where
  "cgu_monomial_eval m x =
    (\<Prod>i\<in>Poly_Mapping.keys m. x $ i ^ Poly_Mapping.lookup m i)"
\end{lstlisting}
\reviewlabel{English mathematical translation}
\noindent\textit{Mathematical role:} monomial \(x^\alpha\)\par\smallskip
For a finite index type \texttt{'n}, \texttt{cgu\_\allowbreak{}monomial\_\allowbreak{}eval\allowbreak\ m\allowbreak\ x} is the product over all indices \texttt{i} in the finite support of \texttt{m} of \texttt{(x\_\allowbreak{}i)\textasciicircum{}(m(i))}.

\dossierentry{\texttt{\detokenize{cgu_poly_eval}}}
\reviewlabel{Isabelle code}
\begin{lstlisting}
definition cgu_poly_eval ::
  "'n cgu_polynomial \<Rightarrow> 'n::finite cgu_point \<Rightarrow> real"
where
  "cgu_poly_eval p x =
    (\<Sum>m\<in>Poly_Mapping.keys p.
      Poly_Mapping.lookup p m * cgu_monomial_eval m x)"
\end{lstlisting}
\reviewlabel{English mathematical translation}
\noindent\textit{Mathematical role:} polynomial evaluation \(\widehat p(x)\)\par\smallskip
For finite \texttt{'n}, \texttt{cgu\_\allowbreak{}poly\_\allowbreak{}eval\allowbreak\ p\allowbreak\ x} is the finite sum over every \texttt{m} in the support of \texttt{p} of the coefficient \texttt{p(m)} times \texttt{cgu\_\allowbreak{}monomial\_\allowbreak{}eval\allowbreak\ m\allowbreak\ x}.

\dossierentry{\texttt{\detokenize{cgu_matrix_poly_eval}}}
\reviewlabel{Isabelle code}
\begin{lstlisting}
definition cgu_matrix_poly_eval ::
  "'n::finite cgu_matrix_polynomial \<Rightarrow> 'n cgu_point \<Rightarrow>
    real ^ 'n ^ 'n"
where
  "cgu_matrix_poly_eval P x = (\<chi> i j. cgu_poly_eval (P $ i $ j) x)"
\end{lstlisting}
\reviewlabel{English mathematical translation}
\noindent\textit{Mathematical role:} matrix evaluation \(\widehat P(x)\)\par\smallskip
For finite \texttt{'n}, \texttt{cgu\_\allowbreak{}matrix\_\allowbreak{}poly\_\allowbreak{}eval\allowbreak\ P\allowbreak\ x} is the real \texttt{'n}-by-\texttt{'n} matrix whose \texttt{(i,j)} entry is \texttt{cgu\_\allowbreak{}poly\_\allowbreak{}eval\allowbreak\ (P[i,j])\allowbreak\ x}.

\dossierentry{\texttt{\detokenize{cgu_matrix_poly_degree_le}}}
\reviewlabel{Isabelle code}
\begin{lstlisting}
definition cgu_matrix_poly_degree_le ::
  "nat \<Rightarrow> 'n::finite cgu_matrix_polynomial \<Rightarrow> bool"
where
  "cgu_matrix_poly_degree_le d P \<longleftrightarrow>
    (\<forall>i j. cgu_poly_degree_le d (P $ i $ j))"
\end{lstlisting}
\reviewlabel{English mathematical translation}
\noindent\textit{Mathematical role:} entrywise stored degree bound\par\smallskip
For finite \texttt{'n}, \texttt{cgu\_\allowbreak{}matrix\_\allowbreak{}poly\_\allowbreak{}degree\_\allowbreak{}le\allowbreak\ d\allowbreak\ P} holds exactly when, for every pair of indices \texttt{i,j}, every supported exponent mapping of the entry \texttt{P[i,j]} has total \texttt{cgu\_\allowbreak{}monomial\_\allowbreak{}total\_\allowbreak{}degree} value at most \texttt{d}.

\dossierentry{\texttt{\detokenize{cgu_symmetric_matrix_polynomial}}}
\reviewlabel{Isabelle code}
\begin{lstlisting}
definition cgu_symmetric_matrix_polynomial ::
  "'n::finite cgu_matrix_polynomial \<Rightarrow> bool"
where
  "cgu_symmetric_matrix_polynomial P \<longleftrightarrow>
    (\<forall>i j. P $ i $ j = P $ j $ i)"
\end{lstlisting}
\reviewlabel{English mathematical translation}
\noindent\textit{Mathematical role:} coefficient symmetry\par\smallskip
For finite \texttt{'n}, \texttt{cgu\_\allowbreak{}symmetric\_\allowbreak{}matrix\_\allowbreak{}polynomial\allowbreak\ P} holds exactly when \texttt{P[i,j]\allowbreak\ =\allowbreak\ P[j,i]} for every \texttt{i,j}.

\dossierentry{\texttt{\detokenize{cgu_matrix_bilinear}}}
\reviewlabel{Isabelle code}
\begin{lstlisting}
definition cgu_matrix_bilinear ::
  "(real ^ 'n::finite ^ 'n) \<Rightarrow> real ^ 'n \<Rightarrow>
    real ^ 'n \<Rightarrow> real"
where
  "cgu_matrix_bilinear A x y = inner x (A *v y)"
\end{lstlisting}
\reviewlabel{English mathematical translation}
\noindent\textit{Mathematical role:} bilinear value \(x\cdot Ay\)\par\smallskip
For a real square matrix \texttt{A} and real vectors \texttt{x,y} of the same finite coordinate type, \texttt{cgu\_\allowbreak{}matrix\_\allowbreak{}bilinear\allowbreak\ A\allowbreak\ x\allowbreak\ y} is the inner product of \texttt{x} with the matrix-vector product \texttt{A*y}.

\dossierentry{\texttt{\detokenize{cgu_symmetric_positive_definite_matrix}}}
\reviewlabel{Isabelle code}
\begin{lstlisting}
definition cgu_symmetric_positive_definite_matrix ::
  "(real ^ 'n::finite ^ 'n) \<Rightarrow> bool"
where
  "cgu_symmetric_positive_definite_matrix A \<longleftrightarrow>
    transpose A = A \<and>
    (\<forall>x. x \<noteq> 0 \<longrightarrow> 0 < cgu_matrix_bilinear A x x)"
\end{lstlisting}
\reviewlabel{English mathematical translation}
\noindent\textit{Mathematical role:} symmetric positive-definite matrix\par\smallskip
\texttt{cgu\_\allowbreak{}symmetric\_\allowbreak{}positive\_\allowbreak{}definite\_\allowbreak{}matrix\allowbreak\ A} holds exactly when both:

\begin{itemize}
\item \texttt{transpose\allowbreak\ A\allowbreak\ =\allowbreak\ A}; and
\item for every vector \texttt{x}, if \texttt{x\allowbreak\ !=\allowbreak\ 0}, then \texttt{0\allowbreak\ <\allowbreak\ cgu\_\allowbreak{}matrix\_\allowbreak{}bilinear\allowbreak\ A\allowbreak\ x\allowbreak\ x}.
\end{itemize}

Thus the second conjunct is strict positivity only on nonzero vectors.

\dossierentry{\texttt{\detokenize{cgu_uniformly_positive_definite_on}}}
\reviewlabel{Isabelle code}
\begin{lstlisting}
definition cgu_uniformly_positive_definite_on ::
  "'n::finite cgu_point set \<Rightarrow> 'n cgu_matrix_polynomial \<Rightarrow> bool"
where
  "cgu_uniformly_positive_definite_on D P \<longleftrightarrow>
    (\<exists>k>0. \<forall>x\<in>D. \<forall>xi.
      k * norm xi ^ 2 \<le>
        cgu_matrix_bilinear (cgu_matrix_poly_eval P x) xi xi)"
\end{lstlisting}
\reviewlabel{English mathematical translation}
\noindent\textit{Mathematical role:} uniform ellipticity on a set\par\smallskip
For finite \texttt{'n}, \texttt{cgu\_\allowbreak{}uniformly\_\allowbreak{}positive\_\allowbreak{}definite\_\allowbreak{}on\allowbreak\ D\allowbreak\ P} holds exactly when there exists one real \texttt{k\allowbreak\ >\allowbreak\ 0} such that, for every \texttt{x\allowbreak\ in\allowbreak\ D} and every vector \texttt{xi}, \texttt{k\allowbreak\ *\allowbreak\ ||xi||\textasciicircum{}2\allowbreak\ <=\allowbreak\ cgu\_\allowbreak{}matrix\_\allowbreak{}bilinear\allowbreak\ (cgu\_\allowbreak{}matrix\_\allowbreak{}poly\_\allowbreak{}eval\allowbreak\ P\allowbreak\ x)\allowbreak\ xi\allowbreak\ xi}. The same \texttt{k} must work for all \texttt{x} and \texttt{xi}.

\dossierentry{\texttt{\detokenize{cgu_model}}}
\reviewlabel{Isabelle code}
\begin{lstlisting}
record ('n, 'c, 'h) cgu_model =
  cgu_domain :: "'n cgu_point set"
  cgu_measured :: "'n cgu_point set"
  cgu_cell :: "'c \<Rightarrow> 'n cgu_point set"
  cgu_degree :: "'c \<Rightarrow> nat"
  cgu_order :: "nat \<Rightarrow> 'c"
  cgu_faces :: "'c \<Rightarrow> 'h set"
  cgu_normal :: "'c \<Rightarrow> 'h \<Rightarrow> 'n cgu_point"
  cgu_offset :: "'c \<Rightarrow> 'h \<Rightarrow> real"
  cgu_face_patch :: "'c \<Rightarrow> 'h \<Rightarrow> 'n cgu_point set"
  cgu_witness_faces :: "'c \<Rightarrow> 'h \<Rightarrow> 'h set"
  cgu_triple_region :: "'c \<Rightarrow> 'h \<Rightarrow> 'h \<Rightarrow> 'h \<Rightarrow>
    'n cgu_point set"
  cgu_step_interface_pieces :: "'c \<Rightarrow> 'h set"
  cgu_step_interface_patch :: "'c \<Rightarrow> 'h \<Rightarrow> 'n cgu_point set"
  cgu_step_interface_normal :: "'c \<Rightarrow> 'h \<Rightarrow> 'n cgu_point"
  cgu_step_interface_offset :: "'c \<Rightarrow> 'h \<Rightarrow> real"
  cgu_step_outer_patch :: "'c \<Rightarrow> 'n cgu_point set"
  cgu_step_adjacent_cell :: "'c \<Rightarrow> 'c"
  cgu_step_center :: "'c \<Rightarrow> 'n cgu_point"
  cgu_step_interior_normal :: "'c \<Rightarrow> 'n cgu_point"
  cgu_step_radius :: "'c \<Rightarrow> real"
\end{lstlisting}
\reviewlabel{English mathematical translation}
\noindent\textit{Mathematical role:} geometric record \(M\); its 20 fields map as follows: \texttt{cgu\_\allowbreak{}domain} \(\Omega\), \texttt{cgu\_\allowbreak{}measured} \(\Gamma\), \texttt{cgu\_\allowbreak{}cell} \(\Omega_c\), \texttt{cgu\_\allowbreak{}degree} \(N_c\), \texttt{cgu\_\allowbreak{}order} \(\sigma\), \texttt{cgu\_\allowbreak{}faces} \(I_c\), \texttt{cgu\_\allowbreak{}normal} \(\nu_{c,i}\), \texttt{cgu\_\allowbreak{}offset} \(b_{c,i}\), \texttt{cgu\_\allowbreak{}face\_\allowbreak{}patch} \(\gamma_{c,i}\), \texttt{cgu\_\allowbreak{}witness\_\allowbreak{}faces} \(J_{c,i}\), \texttt{cgu\_\allowbreak{}triple\_\allowbreak{}region} \(W_{c,i,j,k}\), \texttt{cgu\_\allowbreak{}step\_\allowbreak{}interface\_\allowbreak{}pieces} auxiliary label set, \texttt{cgu\_\allowbreak{}step\_\allowbreak{}interface\_\allowbreak{}patch} auxiliary patch, \texttt{cgu\_\allowbreak{}step\_\allowbreak{}interface\_\allowbreak{}normal} auxiliary normal, \texttt{cgu\_\allowbreak{}step\_\allowbreak{}interface\_\allowbreak{}offset} auxiliary offset, \texttt{cgu\_\allowbreak{}step\_\allowbreak{}outer\_\allowbreak{}patch} \(\Lambda_c\), \texttt{cgu\_\allowbreak{}step\_\allowbreak{}adjacent\_\allowbreak{}cell} \(a(c)\), \texttt{cgu\_\allowbreak{}step\_\allowbreak{}center} \(q_c\), \texttt{cgu\_\allowbreak{}step\_\allowbreak{}interior\_\allowbreak{}normal} \(\mu_c\), \texttt{cgu\_\allowbreak{}step\_\allowbreak{}radius} \(r_c\)\par\smallskip
\texttt{cgu\_\allowbreak{}model} is a record with type parameters \texttt{'n}, \texttt{'c}, and \texttt{'h} and the following fields, with no class restrictions imposed by the record declaration itself:

\begin{itemize}
\item \texttt{cgu\_\allowbreak{}domain}: a set of \texttt{'n}-indexed real vectors;
\item \texttt{cgu\_\allowbreak{}measured}: another such set;
\item \texttt{cgu\_\allowbreak{}cell}: from each \texttt{'c} to such a set;
\item \texttt{cgu\_\allowbreak{}degree}: from each \texttt{'c} to a natural number;
\item \texttt{cgu\_\allowbreak{}order}: from natural numbers to \texttt{'c};
\item \texttt{cgu\_\allowbreak{}faces}: from each \texttt{'c} to a set of \texttt{'h} values;
\item \texttt{cgu\_\allowbreak{}normal}: from \texttt{'c} and \texttt{'h} to an \texttt{'n}-indexed real vector;
\item \texttt{cgu\_\allowbreak{}offset}: from \texttt{'c} and \texttt{'h} to a real number;
\item \texttt{cgu\_\allowbreak{}face\_\allowbreak{}patch}: from \texttt{'c} and \texttt{'h} to a set of \texttt{'n}-indexed real vectors;
\item \texttt{cgu\_\allowbreak{}witness\_\allowbreak{}faces}: from \texttt{'c} and \texttt{'h} to a set of \texttt{'h} values;
\item \texttt{cgu\_\allowbreak{}triple\_\allowbreak{}region}: from \texttt{'c} and three \texttt{'h} values to a set of \texttt{'n}-indexed real vectors;
\item \texttt{cgu\_\allowbreak{}step\_\allowbreak{}interface\_\allowbreak{}pieces}: from each \texttt{'c} to a set of \texttt{'h} values;
\item \texttt{cgu\_\allowbreak{}step\_\allowbreak{}interface\_\allowbreak{}patch}: from \texttt{'c} and \texttt{'h} to a set of \texttt{'n}-indexed real vectors;
\item \texttt{cgu\_\allowbreak{}step\_\allowbreak{}interface\_\allowbreak{}normal}: from \texttt{'c} and \texttt{'h} to an \texttt{'n}-indexed real vector;
\item \texttt{cgu\_\allowbreak{}step\_\allowbreak{}interface\_\allowbreak{}offset}: from \texttt{'c} and \texttt{'h} to a real number;
\item \texttt{cgu\_\allowbreak{}step\_\allowbreak{}outer\_\allowbreak{}patch}: from each \texttt{'c} to a set of \texttt{'n}-indexed real vectors;
\item \texttt{cgu\_\allowbreak{}step\_\allowbreak{}adjacent\_\allowbreak{}cell}: from \texttt{'c} to \texttt{'c};
\item \texttt{cgu\_\allowbreak{}step\_\allowbreak{}center}: from \texttt{'c} to an \texttt{'n}-indexed real vector;
\item \texttt{cgu\_\allowbreak{}step\_\allowbreak{}interior\_\allowbreak{}normal}: from \texttt{'c} to an \texttt{'n}-indexed real vector; and
\item \texttt{cgu\_\allowbreak{}step\_\allowbreak{}radius}: from \texttt{'c} to a real number.
\end{itemize}

\dossierentry{\texttt{\detokenize{cgu_piecewise_polynomial_class}}}
\reviewlabel{Isabelle code}
\begin{lstlisting}
definition cgu_piecewise_polynomial_class ::
  "('n::finite, 'c::finite, 'h) cgu_model \<Rightarrow>
    ('c \<Rightarrow> 'n cgu_matrix_polynomial) \<Rightarrow> bool"
where
  "cgu_piecewise_polynomial_class M P \<longleftrightarrow>
    (\<forall>c.
      cgu_matrix_poly_degree_le (cgu_degree M c) (P c) \<and>
      cgu_symmetric_matrix_polynomial (P c) \<and>
      cgu_uniformly_positive_definite_on (cgu_cell M c) (P c))"
\end{lstlisting}
\reviewlabel{English mathematical translation}
\noindent\textit{Mathematical role:} admissible coefficient family\par\smallskip
For finite \texttt{'n} and finite \texttt{'c}, \texttt{cgu\_\allowbreak{}piecewise\_\allowbreak{}polynomial\_\allowbreak{}class\allowbreak\ M\allowbreak\ P} holds exactly when, for every \texttt{c}, all three conditions hold: \texttt{P\allowbreak\ c} satisfies \texttt{cgu\_\allowbreak{}matrix\_\allowbreak{}poly\_\allowbreak{}degree\_\allowbreak{}le} with bound \texttt{cgu\_\allowbreak{}degree\allowbreak\ M\allowbreak\ c}; \texttt{P\allowbreak\ c} satisfies \texttt{cgu\_\allowbreak{}symmetric\_\allowbreak{}matrix\_\allowbreak{}polynomial}; and \texttt{P\allowbreak\ c} satisfies \texttt{cgu\_\allowbreak{}uniformly\_\allowbreak{}positive\_\allowbreak{}definite\_\allowbreak{}on} on \texttt{cgu\_\allowbreak{}cell\allowbreak\ M\allowbreak\ c}.

\dossierentry{\texttt{\detokenize{cgu_affine_hyperplane}}}
\reviewlabel{Isabelle code}
\begin{lstlisting}
definition cgu_affine_hyperplane ::
  "'n::finite cgu_point \<Rightarrow> real \<Rightarrow> 'n cgu_point set"
where
  "cgu_affine_hyperplane normal offset = {x. inner normal x = offset}"
\end{lstlisting}
\reviewlabel{English mathematical translation}
\noindent\textit{Mathematical role:} hyperplane \(H(\nu,b)\)\par\smallskip
For finite \texttt{'n}, \texttt{cgu\_\allowbreak{}affine\_\allowbreak{}hyperplane\allowbreak\ normal\allowbreak\ offset} is the set of all vectors \texttt{x} satisfying \texttt{inner\allowbreak\ normal\allowbreak\ x\allowbreak\ =\allowbreak\ offset}.

\dossierentry{\texttt{\detokenize{cgu_coordinate_hyperplane_projection_v5}}}
\reviewlabel{Isabelle code}
\begin{lstlisting}
definition cgu_coordinate_hyperplane_projection_v5 ::
  "'i::finite \<Rightarrow> 'i cgu_point \<Rightarrow> 'i cgu_point"
where
  "cgu_coordinate_hyperplane_projection_v5 i y =
    (\<chi> j. if j = i then 0 else y $ j)"
\end{lstlisting}
\reviewlabel{English mathematical translation}
\noindent\textit{Mathematical role:} coordinate-deleted vector \(y^{(i)}\)\par\smallskip
For finite \texttt{'i}, \texttt{cgu\_\allowbreak{}coordinate\_\allowbreak{}hyperplane\_\allowbreak{}projection\_\allowbreak{}v5\allowbreak\ i\allowbreak\ y} is the vector whose \texttt{j}th coordinate is zero if \texttt{j=i}, and is \texttt{y\_\allowbreak{}j} otherwise.

\dossierentry{\texttt{\detokenize{cgu_local_lipschitz_boundary_at}}}
\reviewlabel{Isabelle code}
\begin{lstlisting}
definition cgu_local_lipschitz_boundary_at ::
  "'n::finite cgu_point set \<Rightarrow> 'n cgu_point \<Rightarrow> bool"
where
  "cgu_local_lipschitz_boundary_at D x \<longleftrightarrow>
    (\<exists>(U :: 'n cgu_point set)
       (f :: 'n cgu_point \<Rightarrow> 'n cgu_point)
       (g :: 'n cgu_point \<Rightarrow> 'n cgu_point) (i :: 'n) L M.
      open U \<and> x \<in> U \<and> 0 < L \<and> 0 < M \<and>
      homeomorphism U (f ` U) f g \<and>
      lipschitz_on L U f \<and> lipschitz_on M (f ` U) g \<and>
      f ` (D \<inter> U) = (f ` U) \<inter> {y. 0 < y $ i})"
\end{lstlisting}
\reviewlabel{English mathematical translation}
\noindent\textit{Mathematical role:} local bi-Lipschitz half-space chart\par\smallskip
For finite \texttt{'n}, \texttt{cgu\_\allowbreak{}local\_\allowbreak{}lipschitz\_\allowbreak{}boundary\_\allowbreak{}at\allowbreak\ D\allowbreak\ x} holds exactly when there exist a set \texttt{U}, functions \texttt{f,g}, an index \texttt{i}, and real constants \texttt{L,M} such that:

\begin{itemize}
\item \texttt{U} is open and contains \texttt{x};
\item \texttt{L>0} and \texttt{M>0};
\item \texttt{f} and \texttt{g} give a homeomorphism between \texttt{U} and \texttt{f(U)};
\item \texttt{f} is \texttt{L}-Lipschitz on \texttt{U}, and \texttt{g} is \texttt{M}-Lipschitz on \texttt{f(U)}; and
\item \texttt{f(D\allowbreak\ intersect\allowbreak\ U)} is exactly \texttt{f(U)\allowbreak\ intersect\allowbreak\ \{y\allowbreak\ |\allowbreak\ 0\allowbreak\ <\allowbreak\ y\_\allowbreak{}i\}}.
\end{itemize}

All witnesses are inside one existential scope.

\dossierentry{\texttt{\detokenize{cgu_lipschitz_domain}}}
\reviewlabel{Isabelle code}
\begin{lstlisting}
definition cgu_lipschitz_domain ::
  "'n::finite cgu_point set \<Rightarrow> bool"
where
  "cgu_lipschitz_domain D \<longleftrightarrow>
    D \<noteq> {} \<and> open D \<and> connected D \<and> bounded D \<and>
    (\<forall>x\<in>frontier D. cgu_local_lipschitz_boundary_at D x)"
\end{lstlisting}
\reviewlabel{English mathematical translation}
\noindent\textit{Mathematical role:} bi-Lipschitz half-space domain\par\smallskip
\texttt{cgu\_\allowbreak{}lipschitz\_\allowbreak{}domain\allowbreak\ D} holds exactly when \texttt{D} is nonempty, open, connected, and bounded, and every \texttt{x} in \texttt{frontier\allowbreak\ D} satisfies \texttt{cgu\_\allowbreak{}local\_\allowbreak{}lipschitz\_\allowbreak{}boundary\_\allowbreak{}at\allowbreak\ D\allowbreak\ x}.

\dossierentry{\texttt{\detokenize{cgu_rigid_motion_v5}}}
\reviewlabel{Isabelle code}
\begin{lstlisting}
definition cgu_rigid_motion_v5 ::
  "('n::finite cgu_point \<Rightarrow> 'n cgu_point) \<Rightarrow> bool"
where
  "cgu_rigid_motion_v5 R \<longleftrightarrow>
    (\<exists>A b.
      linear A \<and>
      surj A \<and>
      (\<forall>u v. inner (A u) (A v) = inner u v) \<and>
      R = (\<lambda>x. A x + b))"
\end{lstlisting}
\reviewlabel{English mathematical translation}
\noindent\textit{Mathematical role:} rigid motion\par\smallskip
\texttt{cgu\_\allowbreak{}rigid\_\allowbreak{}motion\_\allowbreak{}v5\allowbreak\ R} holds exactly when there exist a function \texttt{A} and vector \texttt{b} such that \texttt{A} is linear and surjective, preserves all inner products (\texttt{inner\allowbreak\ (A\allowbreak\ u)\allowbreak\ (A\allowbreak\ v)\allowbreak\ =\allowbreak\ inner\allowbreak\ u\allowbreak\ v} for every \texttt{u,v}), and \texttt{R\allowbreak\ x\allowbreak\ =\allowbreak\ A\allowbreak\ x\allowbreak\ +\allowbreak\ b} for every \texttt{x}.

\dossierentry{\texttt{\detokenize{cgu_local_graph_boundary_at_v5}}}
\reviewlabel{Isabelle code}
\begin{lstlisting}
definition cgu_local_graph_boundary_at_v5 ::
  "'n::finite cgu_point set \<Rightarrow> 'n cgu_point \<Rightarrow> bool"
where
  "cgu_local_graph_boundary_at_v5 D x \<longleftrightarrow>
    (\<exists>U R i zeta L.
      open U \<and>
      x \<in> U \<and>
      cgu_rigid_motion_v5 R \<and>
      0 < L \<and>
      lipschitz_on L UNIV zeta \<and>
      R ` (D \<inter> U) =
        R ` U \<inter>
          {y. y $ i < zeta (cgu_coordinate_hyperplane_projection_v5 i y)} \<and>
      R ` (frontier D \<inter> U) =
        R ` U \<inter>
          {y. y $ i = zeta (cgu_coordinate_hyperplane_projection_v5 i y)})"
\end{lstlisting}
\reviewlabel{English mathematical translation}
\noindent\textit{Mathematical role:} local Lipschitz graph condition\par\smallskip
For finite \texttt{'n}, \texttt{cgu\_\allowbreak{}local\_\allowbreak{}graph\_\allowbreak{}boundary\_\allowbreak{}at\_\allowbreak{}v5\allowbreak\ D\allowbreak\ x} holds exactly when there exist \texttt{U,R,i,zeta,L} such that:

\begin{itemize}
\item \texttt{U} is open and contains \texttt{x};
\item \texttt{cgu\_\allowbreak{}rigid\_\allowbreak{}motion\_\allowbreak{}v5\allowbreak\ R} holds;
\item \texttt{L>0}, and \texttt{zeta} is \texttt{L}-Lipschitz on the whole vector space;
\item \texttt{R(D\allowbreak\ intersect\allowbreak\ U)\allowbreak\ =\allowbreak\ R(U)\allowbreak\ intersect\allowbreak\ \{y\allowbreak\ |\allowbreak\ y\_\allowbreak{}i\allowbreak\ <\allowbreak\ zeta(cgu\_\allowbreak{}coordinate\_\allowbreak{}hyperplane\_\allowbreak{}projection\_\allowbreak{}v5\allowbreak\ i\allowbreak\ y)\}}; and
\item \texttt{R(frontier\allowbreak\ D\allowbreak\ intersect\allowbreak\ U)\allowbreak\ =\allowbreak\ R(U)\allowbreak\ intersect\allowbreak\ \{y\allowbreak\ |\allowbreak\ y\_\allowbreak{}i\allowbreak\ =\allowbreak\ zeta(cgu\_\allowbreak{}coordinate\_\allowbreak{}hyperplane\_\allowbreak{}projection\_\allowbreak{}v5\allowbreak\ i\allowbreak\ y)\}}.
\end{itemize}

The same witnesses serve both displayed equalities.

\dossierentry{\texttt{\detokenize{cgu_graph_lipschitz_domain_v5}}}
\reviewlabel{Isabelle code}
\begin{lstlisting}
definition cgu_graph_lipschitz_domain_v5 ::
  "'n::finite cgu_point set \<Rightarrow> bool"
where
  "cgu_graph_lipschitz_domain_v5 D \<longleftrightarrow>
    cgu_lipschitz_domain D \<and>
    (\<forall>x\<in>frontier D. cgu_local_graph_boundary_at_v5 D x)"
\end{lstlisting}
\reviewlabel{English mathematical translation}
\noindent\textit{Mathematical role:} Lipschitz-graph domain\par\smallskip
\texttt{cgu\_\allowbreak{}graph\_\allowbreak{}lipschitz\_\allowbreak{}domain\_\allowbreak{}v5\allowbreak\ D} holds exactly when \texttt{cgu\_\allowbreak{}lipschitz\_\allowbreak{}domain\allowbreak\ D} holds and every point \texttt{x} of \texttt{frontier\allowbreak\ D} satisfies \texttt{cgu\_\allowbreak{}local\_\allowbreak{}graph\_\allowbreak{}boundary\_\allowbreak{}at\_\allowbreak{}v5\allowbreak\ D\allowbreak\ x}.

\dossierentry{\texttt{\detokenize{cgu_locally_flat_boundary_point}}}
\reviewlabel{Isabelle code}
\begin{lstlisting}
definition cgu_locally_flat_boundary_point ::
  "'n::finite cgu_point set \<Rightarrow> 'n cgu_point \<Rightarrow> bool"
where
  "cgu_locally_flat_boundary_point D x \<longleftrightarrow>
    x \<in> frontier D \<and>
    (\<exists>normal offset U.
      normal \<noteq> 0 \<and> open U \<and> x \<in> U \<and>
      U \<inter> frontier D =
        U \<inter> cgu_affine_hyperplane normal offset)"
\end{lstlisting}
\reviewlabel{English mathematical translation}
\noindent\textit{Mathematical role:} locally flat boundary point\par\smallskip
\texttt{cgu\_\allowbreak{}locally\_\allowbreak{}flat\_\allowbreak{}boundary\_\allowbreak{}point\allowbreak\ D\allowbreak\ x} holds exactly when \texttt{x} belongs to \texttt{frontier\allowbreak\ D} and there exist \texttt{normal}, \texttt{offset}, and an open set \texttt{U} containing \texttt{x}, with \texttt{normal\allowbreak\ !=\allowbreak\ 0}, such that \texttt{U\allowbreak\ intersect\allowbreak\ frontier\allowbreak\ D\allowbreak\ =\allowbreak\ U\allowbreak\ intersect\allowbreak\ cgu\_\allowbreak{}affine\_\allowbreak{}hyperplane\allowbreak\ normal\allowbreak\ offset}.

\dossierentry{\texttt{\detokenize{cgu_relative_frontier}}}
\reviewlabel{Isabelle code}
\begin{lstlisting}
definition cgu_relative_frontier ::
  "'a::topological_space set \<Rightarrow> 'a set \<Rightarrow> 'a set"
where
  "cgu_relative_frontier X A = closure A \<inter> closure (X - A)"
\end{lstlisting}
\reviewlabel{English mathematical translation}
\noindent\textit{Mathematical role:} \(\overline A\cap\overline{X\setminus A}\)\par\smallskip
In any topological space, \texttt{cgu\_\allowbreak{}relative\_\allowbreak{}frontier\allowbreak\ X\allowbreak\ A\allowbreak\ =\allowbreak\ closure\allowbreak\ A\allowbreak\ intersect\allowbreak\ closure\allowbreak\ (X\allowbreak\ -\allowbreak{}\allowbreak\ A)}. Both closures are the ambient closures appearing in the code.

\dossierentry{\texttt{\detokenize{cgu_subdivision_edge_corner_set}}}
\reviewlabel{Isabelle code}
\begin{lstlisting}
definition cgu_subdivision_edge_corner_set ::
  "('n::finite, 'c::finite, 'h) cgu_model \<Rightarrow> 'n cgu_point set"
where
  "cgu_subdivision_edge_corner_set M =
    (\<Union>c.
      let X = frontier (cgu_cell M c);
          F = {x\<in>X. cgu_locally_flat_boundary_point (cgu_cell M c) x}
      in cgu_relative_frontier X F)"
\end{lstlisting}
\reviewlabel{English mathematical translation}
\noindent\textit{Mathematical role:} singular flatness set \(\Sigma_M\)\par\smallskip
For each \texttt{c}, let \texttt{X\allowbreak\ =\allowbreak\ frontier\allowbreak\ (cgu\_\allowbreak{}cell\allowbreak\ M\allowbreak\ c)} and let \texttt{F} be the points \texttt{x\allowbreak\ in\allowbreak\ X} satisfying \texttt{cgu\_\allowbreak{}locally\_\allowbreak{}flat\_\allowbreak{}boundary\_\allowbreak{}point\allowbreak\ (cgu\_\allowbreak{}cell\allowbreak\ M\allowbreak\ c)\allowbreak\ x}. Then \texttt{cgu\_\allowbreak{}subdivision\_\allowbreak{}edge\_\allowbreak{}corner\_\allowbreak{}set\allowbreak\ M} is the union over all \texttt{c} of \texttt{closure\allowbreak\ F\allowbreak\ intersect\allowbreak\ closure\allowbreak\ (X-\allowbreak{}F)}.

\dossierentry{\texttt{\detokenize{cgu_graph_subdivision_core_v7}}}
\reviewlabel{Isabelle code}
\begin{lstlisting}
definition cgu_graph_subdivision_core_v7 ::
  "('n::finite, 'c::finite, 'h::finite) cgu_model \<Rightarrow> bool"
where
  "cgu_graph_subdivision_core_v7 M \<longleftrightarrow>
    cgu_lipschitz_domain (cgu_domain M) \<and>
    cgu_measured M \<noteq> {} \<and>
    openin (top_of_set (frontier (cgu_domain M))) (cgu_measured M) \<and>
    cgu_measured M \<subseteq> frontier (cgu_domain M) \<and>
    (\<forall>c. cgu_cell M c \<subseteq> cgu_domain M \<and>
      cgu_lipschitz_domain (cgu_cell M c)) \<and>
    (\<forall>c d. c \<noteq> d \<longrightarrow> cgu_cell M c \<inter> cgu_cell M d = {}) \<and>
    closure (cgu_domain M) = (\<Union>c. closure (cgu_cell M c)) \<and>
    cgu_graph_lipschitz_domain_v5 (cgu_domain M) \<and>
    (\<forall>c. cgu_graph_lipschitz_domain_v5 (cgu_cell M c))"
\end{lstlisting}
\reviewlabel{English mathematical translation}
\noindent\textit{Mathematical role:} regular cell decomposition\par\smallskip
For finite \texttt{'n}, \texttt{'c}, and \texttt{'h}, \texttt{cgu\_\allowbreak{}graph\_\allowbreak{}subdivision\_\allowbreak{}core\_\allowbreak{}v7\allowbreak\ M} holds exactly when all of the following conjuncts hold:

\begin{itemize}
\item \texttt{cgu\_\allowbreak{}lipschitz\_\allowbreak{}domain\allowbreak\ (cgu\_\allowbreak{}domain\allowbreak\ M)};
\item \texttt{cgu\_\allowbreak{}measured\allowbreak\ M} is nonempty, is open in the subspace topology on \texttt{frontier\allowbreak\ (cgu\_\allowbreak{}domain\allowbreak\ M)}, and is a subset of that frontier;
\item for every \texttt{c}, \texttt{cgu\_\allowbreak{}cell\allowbreak\ M\allowbreak\ c} is a subset of \texttt{cgu\_\allowbreak{}domain\allowbreak\ M} and satisfies \texttt{cgu\_\allowbreak{}lipschitz\_\allowbreak{}domain};
\item for every distinct \texttt{c,d}, the sets \texttt{cgu\_\allowbreak{}cell\allowbreak\ M\allowbreak\ c} and \texttt{cgu\_\allowbreak{}cell\allowbreak\ M\allowbreak\ d} are disjoint;
\item \texttt{closure\allowbreak\ (cgu\_\allowbreak{}domain\allowbreak\ M)} equals the union, over all \texttt{c}, of \texttt{closure\allowbreak\ (cgu\_\allowbreak{}cell\allowbreak\ M\allowbreak\ c)};
\item \texttt{cgu\_\allowbreak{}graph\_\allowbreak{}lipschitz\_\allowbreak{}domain\_\allowbreak{}v5\allowbreak\ (cgu\_\allowbreak{}domain\allowbreak\ M)}; and
\item for every \texttt{c}, \texttt{cgu\_\allowbreak{}graph\_\allowbreak{}lipschitz\_\allowbreak{}domain\_\allowbreak{}v5\allowbreak\ (cgu\_\allowbreak{}cell\allowbreak\ M\allowbreak\ c)}.
\end{itemize}

\dossierentry{\texttt{\detokenize{cgu_graph_boundary_regularities_imply_lipschitz_domain}}}
\reviewlabel{Isabelle code}
\begin{lstlisting}
theorem cgu_graph_boundary_regularities_imply_lipschitz_domain:
  shows \<open>D \<noteq> {} \<and> open D \<and> connected D \<and> bounded D \<and>
    (\<forall>x\<in>frontier D. cgu_local_graph_boundary_at_v5 D x)
    \<longrightarrow> cgu_lipschitz_domain D\<close>
\end{lstlisting}
\reviewlabel{English mathematical translation}
\noindent\textit{Mathematical role:} graph-condition-to-half-space-domain proposition\par\smallskip
\texttt{D} is an implicitly generalized schematic term parameter whose finite-dimensional vector-set type is forced by the referenced constants. The single object-level proposition shown by the theorem is:

If \texttt{D} is nonempty, open, connected, and bounded, and every \texttt{x\allowbreak\ in\allowbreak\ frontier\allowbreak\ D} satisfies \texttt{cgu\_\allowbreak{}local\_\allowbreak{}graph\_\allowbreak{}boundary\_\allowbreak{}at\_\allowbreak{}v5\allowbreak\ D\allowbreak\ x}, then \texttt{cgu\_\allowbreak{}lipschitz\_\allowbreak{}domain\allowbreak\ D} holds.

The entire five-part conjunction is the antecedent of that object-level implication. The displayed proposition itself contains no object-level universal quantifier over \texttt{D}.

\dossierentry{\texttt{\detokenize{cgu_share_flat_face_in}}}
\reviewlabel{Isabelle code}
\begin{lstlisting}
definition cgu_share_flat_face_in ::
  "'n::finite cgu_point set \<Rightarrow> 'n cgu_point set \<Rightarrow>
    'n cgu_point set \<Rightarrow> bool"
where
  "cgu_share_flat_face_in G E F \<longleftrightarrow>
    (\<exists>normal offset U.
      normal \<noteq> 0 \<and> open U \<and>
      U \<inter> cgu_affine_hyperplane normal offset \<noteq> {} \<and>
      U \<inter> cgu_affine_hyperplane normal offset
        \<subseteq> G \<inter> frontier E \<inter> frontier F)"
\end{lstlisting}
\reviewlabel{English mathematical translation}
\noindent\textit{Mathematical role:} flat common facet visible in \(G\)\par\smallskip
\texttt{cgu\_\allowbreak{}share\_\allowbreak{}flat\_\allowbreak{}face\_\allowbreak{}in\allowbreak\ G\allowbreak\ E\allowbreak\ F} holds exactly when there exist \texttt{normal}, \texttt{offset}, and an open set \texttt{U} such that \texttt{normal\allowbreak\ !=\allowbreak\ 0}, the set \texttt{U\allowbreak\ intersect\allowbreak\ cgu\_\allowbreak{}affine\_\allowbreak{}hyperplane\allowbreak\ normal\allowbreak\ offset} is nonempty, and this set is contained in \texttt{G\allowbreak\ intersect\allowbreak\ frontier\allowbreak\ E\allowbreak\ intersect\allowbreak\ frontier\allowbreak\ F}.

\dossierentry{\texttt{\detokenize{cgu_linearly_independent3}}}
\reviewlabel{Isabelle code}
\begin{lstlisting}
definition cgu_linearly_independent3 ::
  "real ^ 'n::finite \<Rightarrow> real ^ 'n \<Rightarrow> real ^ 'n \<Rightarrow> bool"
where
  "cgu_linearly_independent3 a b c \<longleftrightarrow>
    (\<forall>x y z::real.
      x *\<^sub>R a + y *\<^sub>R b + z *\<^sub>R c = 0
      \<longrightarrow> x = 0 \<and> y = 0 \<and> z = 0)"
\end{lstlisting}
\reviewlabel{English mathematical translation}
\noindent\textit{Mathematical role:} linear independence of three normals\par\smallskip
\texttt{cgu\_\allowbreak{}linearly\_\allowbreak{}independent3\allowbreak\ a\allowbreak\ b\allowbreak\ c} holds exactly when, for every real \texttt{x,y,z}, the equality \texttt{x*a\allowbreak\ +\allowbreak\ y*b\allowbreak\ +\allowbreak\ z*c\allowbreak\ =\allowbreak\ 0} implies \texttt{x=0}, \texttt{y=0}, and \texttt{z=0}. The conclusion is a three-way conjunction.

\dossierentry{\texttt{\detokenize{cgu_n_generic_on}}}
\reviewlabel{Isabelle code}
\begin{lstlisting}
definition cgu_n_generic_on ::
  "nat \<Rightarrow> 'h::finite set \<Rightarrow>
    ('h \<Rightarrow> real ^ 'n::finite) \<Rightarrow> ('h \<Rightarrow> real) \<Rightarrow>
    ('h \<Rightarrow> 'h set) \<Rightarrow> bool"
where
  "cgu_n_generic_on N I normal offset J \<longleftrightarrow>
    inj_on (\<lambda>i. cgu_affine_hyperplane (normal i) (offset i)) I \<and>
    (\<forall>i\<in>I.
      J i \<subseteq> I - {i} \<and> card (J i) = N + 2 \<and>
      (\<forall>j\<in>J i. \<forall>k\<in>J i. j \<noteq> k \<longrightarrow>
        cgu_linearly_independent3 (normal i) (normal j) (normal k) \<and>
        aff_dim
          (cgu_affine_hyperplane (normal i) (offset i) \<inter>
           cgu_affine_hyperplane (normal j) (offset j) \<inter>
           cgu_affine_hyperplane (normal k) (offset k))
        = int (CARD('n)) - 3) \<and>
      (\<forall>j\<in>J i. \<forall>k\<in>J i. \<forall>l\<in>J i.
        j \<noteq> k \<and> l \<noteq> j \<and> l \<noteq> k \<longrightarrow>
        \<not> cgu_affine_hyperplane (normal i) (offset i) \<inter>
             cgu_affine_hyperplane (normal j) (offset j) \<inter>
             cgu_affine_hyperplane (normal k) (offset k)
           \<subseteq> cgu_affine_hyperplane (normal l) (offset l)))"
\end{lstlisting}
\reviewlabel{English mathematical translation}
\noindent\textit{Mathematical role:} hyperplane general-position system \((I,J,N)\)\par\smallskip
For finite \texttt{'h} and finite \texttt{'n}, \texttt{cgu\_\allowbreak{}n\_\allowbreak{}generic\_\allowbreak{}on\allowbreak\ N\allowbreak\ I\allowbreak\ normal\allowbreak\ offset\allowbreak\ J} holds exactly when:

\begin{itemize}
\item the map \texttt{i\allowbreak\ |-\allowbreak{}>\allowbreak\ cgu\_\allowbreak{}affine\_\allowbreak{}hyperplane\allowbreak\ (normal\allowbreak\ i)\allowbreak\ (offset\allowbreak\ i)} is injective on \texttt{I}; and
\item for every \texttt{i\allowbreak\ in\allowbreak\ I}:
\item \texttt{J\allowbreak\ i} is a subset of \texttt{I-\allowbreak{}\{i\}} and has cardinality exactly \texttt{N+2};
\item for every \texttt{j,k\allowbreak\ in\allowbreak\ J\allowbreak\ i} with \texttt{j\allowbreak\ !=\allowbreak\ k}, the three vectors \texttt{normal\allowbreak\ i}, \texttt{normal\allowbreak\ j}, and \texttt{normal\allowbreak\ k} satisfy \texttt{cgu\_\allowbreak{}linearly\_\allowbreak{}independent3}, and the affine dimension of the intersection of their three corresponding \texttt{cgu\_\allowbreak{}affine\_\allowbreak{}hyperplane} sets is \texttt{int(CARD('n))-\allowbreak{}3}; and
\item for every \texttt{j,k,l\allowbreak\ in\allowbreak\ J\allowbreak\ i} with \texttt{j\allowbreak\ !=\allowbreak\ k}, \texttt{l\allowbreak\ !=\allowbreak\ j}, and \texttt{l\allowbreak\ !=\allowbreak\ k}, the intersection of the \texttt{i}, \texttt{j}, and \texttt{k} corresponding \texttt{cgu\_\allowbreak{}affine\_\allowbreak{}hyperplane} sets is not a subset of the \texttt{l} corresponding \texttt{cgu\_\allowbreak{}affine\_\allowbreak{}hyperplane} set.
\end{itemize}

\dossierentry{\texttt{\detokenize{cgu_triple_region_geometry}}}
\reviewlabel{Isabelle code}
\begin{lstlisting}
definition cgu_triple_region_geometry ::
  "'n::finite cgu_point set \<Rightarrow>
    'n cgu_point \<Rightarrow> real \<Rightarrow> 'n cgu_point \<Rightarrow> real \<Rightarrow>
    'n cgu_point \<Rightarrow> real \<Rightarrow> bool"
where
  "cgu_triple_region_geometry W ni oi nj oj nk ok \<longleftrightarrow>
    W \<noteq> {} \<and>
    openin
      (top_of_set
        (cgu_affine_hyperplane ni oi \<inter>
         cgu_affine_hyperplane nj oj \<inter>
         cgu_affine_hyperplane nk ok)) W \<and>
    W \<subseteq>
      cgu_affine_hyperplane ni oi \<inter>
      cgu_affine_hyperplane nj oj \<inter>
      cgu_affine_hyperplane nk ok"
\end{lstlisting}
\reviewlabel{English mathematical translation}
\noindent\textit{Mathematical role:} relatively open triple-intersection window\par\smallskip
\texttt{cgu\_\allowbreak{}triple\_\allowbreak{}region\_\allowbreak{}geometry\allowbreak\ W\allowbreak\ ni\allowbreak\ oi\allowbreak\ nj\allowbreak\ oj\allowbreak\ nk\allowbreak\ ok} holds exactly when \texttt{W} is nonempty, is open in the subspace topology of the intersection of the three displayed \texttt{cgu\_\allowbreak{}affine\_\allowbreak{}hyperplane} sets, and is a subset of that triple intersection.

\dossierentry{\texttt{\detokenize{cgu_polynomial_positive_on}}}
\reviewlabel{Isabelle code}
\begin{lstlisting}
definition cgu_polynomial_positive_on ::
  "'n::finite cgu_matrix_polynomial \<Rightarrow> 'n cgu_point set \<Rightarrow> bool"
where
  "cgu_polynomial_positive_on P W \<longleftrightarrow>
    (\<forall>x\<in>W.
      cgu_symmetric_positive_definite_matrix (cgu_matrix_poly_eval P x))"
\end{lstlisting}
\reviewlabel{English mathematical translation}
\noindent\textit{Mathematical role:} pointwise positive definiteness on a window\par\smallskip
\texttt{cgu\_\allowbreak{}polynomial\_\allowbreak{}positive\_\allowbreak{}on\allowbreak\ P\allowbreak\ W} holds exactly when, for every \texttt{x\allowbreak\ in\allowbreak\ W}, the evaluated matrix \texttt{cgu\_\allowbreak{}matrix\_\allowbreak{}poly\_\allowbreak{}eval\allowbreak\ P\allowbreak\ x} satisfies \texttt{cgu\_\allowbreak{}symmetric\_\allowbreak{}positive\_\allowbreak{}definite\_\allowbreak{}matrix}.

\dossierentry{\texttt{\detokenize{cgu_order_valid}}}
\reviewlabel{Isabelle code}
\begin{lstlisting}
definition cgu_order_valid ::
  "('n, 'c::finite, 'h) cgu_model \<Rightarrow> bool"
where
  "cgu_order_valid M \<longleftrightarrow>
    bij_betw (cgu_order M) {..<CARD('c)} UNIV"
\end{lstlisting}
\reviewlabel{English mathematical translation}
\noindent\textit{Mathematical role:} bijective recovery enumeration \(\sigma\)\par\smallskip
For finite \texttt{'c}, \texttt{cgu\_\allowbreak{}order\_\allowbreak{}valid\allowbreak\ M} holds exactly when \texttt{cgu\_\allowbreak{}order\allowbreak\ M} is a bijection from the natural-number interval \texttt{\{0,.\allowbreak{}.\allowbreak{}.\allowbreak{},CARD('c)-\allowbreak{}1\}} onto all values of type \texttt{'c}.

\dossierentry{\texttt{\detokenize{cgu_before_indices}}}
\reviewlabel{Isabelle code}
\begin{lstlisting}
definition cgu_before_indices ::
  "('n, 'c, 'h) cgu_model \<Rightarrow> nat \<Rightarrow> 'c set"
where
  "cgu_before_indices M r = cgu_order M ` {..<r}"
\end{lstlisting}
\reviewlabel{English mathematical translation}
\noindent\textit{Mathematical role:} recovered-label set \(\mathcal C_s^-\)\par\smallskip
\texttt{cgu\_\allowbreak{}before\_\allowbreak{}indices\allowbreak\ M\allowbreak\ r} is the image under \texttt{cgu\_\allowbreak{}order\allowbreak\ M} of the natural numbers strictly less than \texttt{r}.

\dossierentry{\texttt{\detokenize{cgu_after_indices}}}
\reviewlabel{Isabelle code}
\begin{lstlisting}
definition cgu_after_indices ::
  "('n, 'c::finite, 'h) cgu_model \<Rightarrow> nat \<Rightarrow> 'c set"
where
  "cgu_after_indices M r = cgu_order M ` {r..<CARD('c)}"
\end{lstlisting}
\reviewlabel{English mathematical translation}
\noindent\textit{Mathematical role:} remaining-label set \(\mathcal C_s^+\)\par\smallskip
For finite \texttt{'c}, \texttt{cgu\_\allowbreak{}after\_\allowbreak{}indices\allowbreak\ M\allowbreak\ r} is the image under \texttt{cgu\_\allowbreak{}order\allowbreak\ M} of the natural numbers \texttt{q} satisfying \texttt{r\allowbreak\ <=\allowbreak\ q\allowbreak\ <\allowbreak\ CARD('c)}.

\dossierentry{\texttt{\detokenize{cgu_remaining_domain}}}
\reviewlabel{Isabelle code}
\begin{lstlisting}
definition cgu_remaining_domain ::
  "('n::finite, 'c::finite, 'h) cgu_model \<Rightarrow> nat \<Rightarrow>
    'n cgu_point set"
where
  "cgu_remaining_domain M r =
    interior (\<Union>c\<in>cgu_after_indices M r. closure (cgu_cell M c))"
\end{lstlisting}
\reviewlabel{English mathematical translation}
\noindent\textit{Mathematical role:} remaining region \(R_s\)\par\smallskip
\texttt{cgu\_\allowbreak{}remaining\_\allowbreak{}domain\allowbreak\ M\allowbreak\ r} is the interior of the union, over \texttt{c\allowbreak\ in\allowbreak\ cgu\_\allowbreak{}after\_\allowbreak{}indices\allowbreak\ M\allowbreak\ r}, of \texttt{closure\allowbreak\ (cgu\_\allowbreak{}cell\allowbreak\ M\allowbreak\ c)}.

\dossierentry{\texttt{\detokenize{cgu_recovered_region}}}
\reviewlabel{Isabelle code}
\begin{lstlisting}
definition cgu_recovered_region ::
  "('n::finite, 'c::finite, 'h) cgu_model \<Rightarrow> nat \<Rightarrow>
    'n cgu_point set"
where
  "cgu_recovered_region M r =
    cgu_domain M - closure (cgu_remaining_domain M r)"
\end{lstlisting}
\reviewlabel{English mathematical translation}
\noindent\textit{Mathematical role:} recovered region \(K_s\)\par\smallskip
\texttt{cgu\_\allowbreak{}recovered\_\allowbreak{}region\allowbreak\ M\allowbreak\ r\allowbreak\ =\allowbreak\ cgu\_\allowbreak{}domain\allowbreak\ M\allowbreak\ -\allowbreak{}\allowbreak\ closure\allowbreak\ (cgu\_\allowbreak{}remaining\_\allowbreak{}domain\allowbreak\ M\allowbreak\ r)}.

\dossierentry{\texttt{\detokenize{cgu_face_connected}}}
\reviewlabel{Isabelle code}
\begin{lstlisting}
definition cgu_face_connected ::
  "'n::finite cgu_point set \<Rightarrow> 'c set \<Rightarrow>
    ('c \<Rightarrow> 'n cgu_point set) \<Rightarrow> bool"
where
  "cgu_face_connected G I cell \<longleftrightarrow>
    (\<forall>i\<in>I. \<forall>j\<in>I.
      rtranclp
        (\<lambda>a b. a \<in> I \<and> b \<in> I \<and>
          cgu_share_flat_face_in G (cell a) (cell b)) i j)"
\end{lstlisting}
\reviewlabel{English mathematical translation}
\noindent\textit{Mathematical role:} connectivity through flat common facets\par\smallskip
\texttt{cgu\_\allowbreak{}face\_\allowbreak{}connected\allowbreak\ G\allowbreak\ I\allowbreak\ cell} holds exactly when, for every \texttt{i,j\allowbreak\ in\allowbreak\ I}, \texttt{j} is reachable from \texttt{i} by the reflexive-transitive closure of the relation that sends \texttt{a} to \texttt{b} when \texttt{a\allowbreak\ in\allowbreak\ I}, \texttt{b\allowbreak\ in\allowbreak\ I}, and \texttt{cgu\_\allowbreak{}share\_\allowbreak{}flat\_\allowbreak{}face\_\allowbreak{}in\allowbreak\ G\allowbreak\ (cell\allowbreak\ a)\allowbreak\ (cell\allowbreak\ b)} holds.

\dossierentry{\texttt{\detokenize{cgu_flat_face_patch}}}
\reviewlabel{Isabelle code}
\begin{lstlisting}
definition cgu_flat_face_patch ::
  "'n::finite cgu_point set \<Rightarrow> 'n cgu_point \<Rightarrow> real \<Rightarrow>
    'n cgu_point set \<Rightarrow> bool"
where
  "cgu_flat_face_patch D normal offset gamma \<longleftrightarrow>
    norm normal = 1 \<and>
    gamma \<noteq> {} \<and>
    openin (top_of_set (cgu_affine_hyperplane normal offset)) gamma \<and>
    gamma \<subseteq> frontier D \<inter> cgu_affine_hyperplane normal offset \<and>
    (\<exists>U. open U \<and> gamma \<subseteq> U \<and>
      U \<inter> frontier D = U \<inter> cgu_affine_hyperplane normal offset \<and>
      (D \<inter> U \<subseteq> {x. inner normal x < offset} \<or>
       D \<inter> U \<subseteq> {x. offset < inner normal x}))"
\end{lstlisting}
\reviewlabel{English mathematical translation}
\noindent\textit{Mathematical role:} one-sided planar boundary patch\par\smallskip
\texttt{cgu\_\allowbreak{}flat\_\allowbreak{}face\_\allowbreak{}patch\allowbreak\ D\allowbreak\ normal\allowbreak\ offset\allowbreak\ gamma} holds exactly when:

\begin{itemize}
\item \texttt{||normal||=1};
\item \texttt{gamma} is nonempty, is open in the subspace topology on \texttt{cgu\_\allowbreak{}affine\_\allowbreak{}hyperplane\allowbreak\ normal\allowbreak\ offset}, and is contained in \texttt{frontier\allowbreak\ D\allowbreak\ intersect\allowbreak\ cgu\_\allowbreak{}affine\_\allowbreak{}hyperplane\allowbreak\ normal\allowbreak\ offset}; and
\item there exists an open \texttt{U} containing \texttt{gamma} such that \texttt{U\allowbreak\ intersect\allowbreak\ frontier\allowbreak\ D\allowbreak\ =\allowbreak\ U\allowbreak\ intersect\allowbreak\ cgu\_\allowbreak{}affine\_\allowbreak{}hyperplane\allowbreak\ normal\allowbreak\ offset}, and either all of \texttt{D\allowbreak\ intersect\allowbreak\ U} satisfies \texttt{inner\allowbreak\ normal\allowbreak\ x\allowbreak\ <\allowbreak\ offset}, or all of it satisfies \texttt{offset\allowbreak\ <\allowbreak\ inner\allowbreak\ normal\allowbreak\ x}.
\end{itemize}

The final disjunction is inside the scope of the one existential \texttt{U}.

\dossierentry{\texttt{\detokenize{cgu_clipped_recovery_interface_v6}}}
\reviewlabel{Isabelle code}
\begin{lstlisting}
definition cgu_clipped_recovery_interface_v6 ::
  "('n::finite, 'c::finite, 'h) cgu_model \<Rightarrow> nat \<Rightarrow> 'n cgu_point set"
where
  "cgu_clipped_recovery_interface_v6 M r =
    cgu_domain M \<inter>
      frontier (cgu_recovered_region M r) \<inter>
      frontier (cgu_remaining_domain M r)"
\end{lstlisting}
\reviewlabel{English mathematical translation}
\noindent\textit{Mathematical role:} stored interface \(\Omega\cap\partial K_s\cap\partial R_s\)\par\smallskip
\texttt{cgu\_\allowbreak{}clipped\_\allowbreak{}recovery\_\allowbreak{}interface\_\allowbreak{}v6\allowbreak\ M\allowbreak\ r} is \texttt{cgu\_\allowbreak{}domain\allowbreak\ M\allowbreak\ intersect\allowbreak\ frontier\allowbreak\ (cgu\_\allowbreak{}recovered\_\allowbreak{}region\allowbreak\ M\allowbreak\ r)\allowbreak\ intersect\allowbreak\ frontier\allowbreak\ (cgu\_\allowbreak{}remaining\_\allowbreak{}domain\allowbreak\ M\allowbreak\ r)}.

\dossierentry{\texttt{\detokenize{cgu_active_pde_face_route_v4}}}
\reviewlabel{Isabelle code}
\begin{lstlisting}
definition cgu_active_pde_face_route_v4 ::
  "('n::finite, 'c::finite, 'h::finite) cgu_model \<Rightarrow>
    nat \<Rightarrow> 'c \<Rightarrow> 'h \<Rightarrow> bool"
where
  "cgu_active_pde_face_route_v4 M r c i \<longleftrightarrow>
    (let remaining = cgu_remaining_domain M r;
         recovered = cgu_recovered_region M r;
         interface = frontier recovered \<inter> frontier remaining;
         patch = cgu_face_patch M c i;
         normal = cgu_normal M c i;
         offset = cgu_offset M c i
     in if r = 0 then
          patch \<subseteq> cgu_measured M \<and>
          cgu_flat_face_patch (cgu_domain M) normal offset patch
        else
          (patch \<subseteq> cgu_measured M \<and>
           cgu_flat_face_patch (cgu_domain M) normal offset patch) \<or>
          (patch \<subseteq> cgu_domain M \<inter> interface \<and>
           cgu_flat_face_patch remaining normal offset patch))"
\end{lstlisting}
\reviewlabel{English mathematical translation}
\noindent\textit{Mathematical role:} stage source-patch eligibility\par\smallskip
Let \texttt{remaining=cgu\_\allowbreak{}remaining\_\allowbreak{}domain\allowbreak\ M\allowbreak\ r}, \texttt{recovered=cgu\_\allowbreak{}recovered\_\allowbreak{}region\allowbreak\ M\allowbreak\ r}, \texttt{interface=frontier\allowbreak\ recovered\allowbreak\ intersect\allowbreak\ frontier\allowbreak\ remaining}, \texttt{patch=cgu\_\allowbreak{}face\_\allowbreak{}patch\allowbreak\ M\allowbreak\ c\allowbreak\ i}, \texttt{normal=cgu\_\allowbreak{}normal\allowbreak\ M\allowbreak\ c\allowbreak\ i}, and \texttt{offset=cgu\_\allowbreak{}offset\allowbreak\ M\allowbreak\ c\allowbreak\ i}. Then:

\begin{itemize}
\item if \texttt{r=0}, \texttt{cgu\_\allowbreak{}active\_\allowbreak{}pde\_\allowbreak{}face\_\allowbreak{}route\_\allowbreak{}v4\allowbreak\ M\allowbreak\ r\allowbreak\ c\allowbreak\ i} holds exactly when \texttt{patch} is a subset of \texttt{cgu\_\allowbreak{}measured\allowbreak\ M} and \texttt{cgu\_\allowbreak{}flat\_\allowbreak{}face\_\allowbreak{}patch\allowbreak\ (cgu\_\allowbreak{}domain\allowbreak\ M)\allowbreak\ normal\allowbreak\ offset\allowbreak\ patch} holds;
\item if \texttt{r!=0}, it holds exactly when either that same two-part condition holds, or \texttt{patch} is a subset of \texttt{cgu\_\allowbreak{}domain\allowbreak\ M\allowbreak\ intersect\allowbreak\ interface} and \texttt{cgu\_\allowbreak{}flat\_\allowbreak{}face\_\allowbreak{}patch\allowbreak\ remaining\allowbreak\ normal\allowbreak\ offset\allowbreak\ patch} holds.
\end{itemize}

\dossierentry{\texttt{\detokenize{cgu_normal_height}}}
\reviewlabel{Isabelle code}
\begin{lstlisting}
definition cgu_normal_height ::
  "'n::finite cgu_point \<Rightarrow> 'n cgu_point \<Rightarrow>
    'n cgu_point \<Rightarrow> real"
where
  "cgu_normal_height center normal x = inner normal (x - center)"
\end{lstlisting}
\reviewlabel{English mathematical translation}
\noindent\textit{Mathematical role:} signed normal coordinate \(\tau_c\)\par\smallskip
\texttt{cgu\_\allowbreak{}normal\_\allowbreak{}height\allowbreak\ center\allowbreak\ normal\allowbreak\ x} is \texttt{inner\allowbreak\ normal\allowbreak\ (x-\allowbreak{}center)}.

\dossierentry{\texttt{\detokenize{cgu_tangent_part}}}
\reviewlabel{Isabelle code}
\begin{lstlisting}
definition cgu_tangent_part ::
  "'n::finite cgu_point \<Rightarrow> 'n cgu_point \<Rightarrow>
    'n cgu_point \<Rightarrow> 'n cgu_point"
where
  "cgu_tangent_part center normal x =
    (x - center) - cgu_normal_height center normal x *\<^sub>R normal"
\end{lstlisting}
\reviewlabel{English mathematical translation}
\noindent\textit{Mathematical role:} tangential displacement \(\pi_c\)\par\smallskip
\texttt{cgu\_\allowbreak{}tangent\_\allowbreak{}part\allowbreak\ center\allowbreak\ normal\allowbreak\ x} is \texttt{(x-\allowbreak{}center)\allowbreak\ -\allowbreak{}\allowbreak\ cgu\_\allowbreak{}normal\_\allowbreak{}height\allowbreak\ center\allowbreak\ normal\allowbreak\ x\allowbreak\ *\allowbreak\ normal}, where the multiplication of the real scalar by \texttt{normal} is scalar multiplication.

\dossierentry{\texttt{\detokenize{cgu_flat_disc}}}
\reviewlabel{Isabelle code}
\begin{lstlisting}
definition cgu_flat_disc ::
  "'n::finite cgu_point \<Rightarrow> 'n cgu_point \<Rightarrow> real \<Rightarrow>
    'n cgu_point set"
where
  "cgu_flat_disc center normal r =
    {x. cgu_normal_height center normal x = 0 \<and>
      norm (cgu_tangent_part center normal x) < r}"
\end{lstlisting}
\reviewlabel{English mathematical translation}
\noindent\textit{Mathematical role:} planar disk \(\Delta_c\)\par\smallskip
\texttt{cgu\_\allowbreak{}flat\_\allowbreak{}disc\allowbreak\ center\allowbreak\ normal\allowbreak\ r} is the set of all \texttt{x} such that \texttt{cgu\_\allowbreak{}normal\_\allowbreak{}height\allowbreak\ center\allowbreak\ normal\allowbreak\ x\allowbreak\ =\allowbreak\ 0} and \texttt{||cgu\_\allowbreak{}tangent\_\allowbreak{}part\allowbreak\ center\allowbreak\ normal\allowbreak\ x||\allowbreak\ <\allowbreak\ r}.

\dossierentry{\texttt{\detokenize{cgu_exterior_half_ball}}}
\reviewlabel{Isabelle code}
\begin{lstlisting}
definition cgu_exterior_half_ball ::
  "'n::finite cgu_point \<Rightarrow> 'n cgu_point \<Rightarrow> real \<Rightarrow>
    'n cgu_point set"
where
  "cgu_exterior_half_ball center normal r =
    ball center r \<inter> {x. cgu_normal_height center normal x < 0}"
\end{lstlisting}
\reviewlabel{English mathematical translation}
\noindent\textit{Mathematical role:} exterior half-ball \(B_c^-\)\par\smallskip
\texttt{cgu\_\allowbreak{}exterior\_\allowbreak{}half\_\allowbreak{}ball\allowbreak\ center\allowbreak\ normal\allowbreak\ r} is the intersection of the open ball of radius \texttt{r} about \texttt{center} with \texttt{\{x\allowbreak\ |\allowbreak\ cgu\_\allowbreak{}normal\_\allowbreak{}height\allowbreak\ center\allowbreak\ normal\allowbreak\ x\allowbreak\ <\allowbreak\ 0\}}.

\dossierentry{\texttt{\detokenize{cgu_auxiliary_half_ball_geometry_step}}}
\reviewlabel{Isabelle code}
\begin{lstlisting}
definition cgu_auxiliary_half_ball_geometry_step ::
  "('n::finite, 'c::finite, 'h) cgu_model \<Rightarrow> nat \<Rightarrow> bool"
where
  "cgu_auxiliary_half_ball_geometry_step M r \<longleftrightarrow>
    (let c = cgu_order M r;
         recovered = cgu_recovered_region M r;
         center = cgu_step_center M c;
         normal = cgu_step_interior_normal M c;
         rad = cgu_step_radius M c;
         patch = cgu_step_outer_patch M c;
         adjacent = cgu_step_adjacent_cell M c
     in norm normal = 1 \<and> 0 < rad \<and>
        patch \<noteq> {} \<and>
        openin (top_of_set (frontier (cgu_domain M))) patch \<and>
        patch \<subseteq> cgu_measured M \<inter> frontier recovered \<and>
        closure (cgu_flat_disc center normal rad)
          \<subseteq> patch - cgu_subdivision_edge_corner_set M \<and>
        cgu_domain M \<inter> ball center rad =
          ball center rad \<inter> {x. 0 < cgu_normal_height center normal x} \<and>
        adjacent \<in> cgu_before_indices M r \<and>
        cgu_flat_disc center normal rad
          \<subseteq> frontier (cgu_domain M) \<inter> frontier (cgu_cell M adjacent))"
\end{lstlisting}
\reviewlabel{English mathematical translation}
\noindent\textit{Mathematical role:} exterior visibility configuration\par\smallskip
Let \texttt{c=cgu\_\allowbreak{}order\allowbreak\ M\allowbreak\ r}, \texttt{recovered=cgu\_\allowbreak{}recovered\_\allowbreak{}region\allowbreak\ M\allowbreak\ r}, \texttt{center=cgu\_\allowbreak{}step\_\allowbreak{}center\allowbreak\ M\allowbreak\ c}, \texttt{normal=cgu\_\allowbreak{}step\_\allowbreak{}interior\_\allowbreak{}normal\allowbreak\ M\allowbreak\ c}, \texttt{rad=cgu\_\allowbreak{}step\_\allowbreak{}radius\allowbreak\ M\allowbreak\ c}, \texttt{patch=cgu\_\allowbreak{}step\_\allowbreak{}outer\_\allowbreak{}patch\allowbreak\ M\allowbreak\ c}, and \texttt{adjacent=cgu\_\allowbreak{}step\_\allowbreak{}adjacent\_\allowbreak{}cell\allowbreak\ M\allowbreak\ c}. Then \texttt{cgu\_\allowbreak{}auxiliary\_\allowbreak{}half\_\allowbreak{}ball\_\allowbreak{}geometry\_\allowbreak{}step\allowbreak\ M\allowbreak\ r} holds exactly when all of the following hold:

\begin{itemize}
\item \texttt{||normal||=1} and \texttt{rad>0};
\item \texttt{patch} is nonempty, is open in the subspace topology on \texttt{frontier\allowbreak\ (cgu\_\allowbreak{}domain\allowbreak\ M)}, and is contained in \texttt{cgu\_\allowbreak{}measured\allowbreak\ M\allowbreak\ intersect\allowbreak\ frontier\allowbreak\ recovered};
\item \texttt{closure\allowbreak\ (cgu\_\allowbreak{}flat\_\allowbreak{}disc\allowbreak\ center\allowbreak\ normal\allowbreak\ rad)} is contained in \texttt{patch\allowbreak\ -\allowbreak{}\allowbreak\ cgu\_\allowbreak{}subdivision\_\allowbreak{}edge\_\allowbreak{}corner\_\allowbreak{}set\allowbreak\ M};
\item \texttt{cgu\_\allowbreak{}domain\allowbreak\ M\allowbreak\ intersect\allowbreak\ ball\allowbreak\ center\allowbreak\ rad} equals \texttt{ball\allowbreak\ center\allowbreak\ rad\allowbreak\ intersect\allowbreak\ \{x\allowbreak\ |\allowbreak\ 0\allowbreak\ <\allowbreak\ cgu\_\allowbreak{}normal\_\allowbreak{}height\allowbreak\ center\allowbreak\ normal\allowbreak\ x\}};
\item \texttt{adjacent\allowbreak\ in\allowbreak\ cgu\_\allowbreak{}before\_\allowbreak{}indices\allowbreak\ M\allowbreak\ r}; and
\item \texttt{cgu\_\allowbreak{}flat\_\allowbreak{}disc\allowbreak\ center\allowbreak\ normal\allowbreak\ rad} is contained in \texttt{frontier\allowbreak\ (cgu\_\allowbreak{}domain\allowbreak\ M)\allowbreak\ intersect\allowbreak\ frontier\allowbreak\ (cgu\_\allowbreak{}cell\allowbreak\ M\allowbreak\ adjacent)}.
\end{itemize}

\dossierentry{\texttt{\detokenize{cgu_recovery_geometry_operational_v6}}}
\reviewlabel{Isabelle code}
\begin{lstlisting}
definition cgu_recovery_geometry_operational_v6 ::
  "('n::finite, 'c::finite, 'h::finite) cgu_model \<Rightarrow> bool"
where
  "cgu_recovery_geometry_operational_v6 M \<longleftrightarrow>
    cgu_order_valid M \<and>
    (\<forall>r<CARD('c).
      let c = cgu_order M r;
          remaining = cgu_remaining_domain M r;
          recovered = cgu_recovered_region M r;
          I = cgu_faces M c
      in cgu_lipschitz_domain remaining \<and>
         (r = 0 \<or>
          (recovered \<noteq> {} \<and>
           cgu_lipschitz_domain recovered \<and>
           cgu_face_connected recovered (cgu_before_indices M r) (cgu_cell M) \<and>
           cgu_auxiliary_half_ball_geometry_step M r)) \<and>
         (\<forall>i\<in>I.
           cgu_flat_face_patch (cgu_cell M c)
             (cgu_normal M c i) (cgu_offset M c i) (cgu_face_patch M c i) \<and>
           cgu_face_patch M c i \<subseteq>
             - cgu_subdivision_edge_corner_set M \<and>
           cgu_face_patch M c i \<subseteq> frontier remaining \<and>
           cgu_active_pde_face_route_v4 M r c i) \<and>
         cgu_degree M c + 3 \<le> card I \<and>
         cgu_n_generic_on (cgu_degree M c) I
           (cgu_normal M c) (cgu_offset M c) (cgu_witness_faces M c) \<and>
         (\<forall>i\<in>I. \<forall>j\<in>cgu_witness_faces M c i.
           \<forall>k\<in>cgu_witness_faces M c i. j \<noteq> k \<longrightarrow>
           cgu_triple_region_geometry (cgu_triple_region M c i j k)
             (cgu_normal M c i) (cgu_offset M c i)
             (cgu_normal M c j) (cgu_offset M c j)
             (cgu_normal M c k) (cgu_offset M c k))) \<and>
    (\<forall>r<CARD('c).
      let remaining = cgu_remaining_domain M r;
          recovered = cgu_recovered_region M r
      in cgu_graph_lipschitz_domain_v5 remaining \<and>
         (r = 0 \<or> cgu_graph_lipschitz_domain_v5 recovered))"
\end{lstlisting}
\reviewlabel{English mathematical translation}
\noindent\textit{Mathematical role:} full recovery geometry\par\smallskip
For finite \texttt{'n}, \texttt{'c}, and \texttt{'h}, \texttt{cgu\_\allowbreak{}recovery\_\allowbreak{}geometry\_\allowbreak{}operational\_\allowbreak{}v6\allowbreak\ M} is the conjunction of \texttt{cgu\_\allowbreak{}order\_\allowbreak{}valid\allowbreak\ M} and the following two universally quantified stage conditions.

First, for every natural \texttt{r\allowbreak\ <\allowbreak\ CARD('c)}, set \texttt{c=cgu\_\allowbreak{}order\allowbreak\ M\allowbreak\ r}, \texttt{remaining=cgu\_\allowbreak{}remaining\_\allowbreak{}domain\allowbreak\ M\allowbreak\ r}, \texttt{recovered=cgu\_\allowbreak{}recovered\_\allowbreak{}region\allowbreak\ M\allowbreak\ r}, and \texttt{I=cgu\_\allowbreak{}faces\allowbreak\ M\allowbreak\ c}. Then all of these must hold:

\begin{itemize}
\item \texttt{cgu\_\allowbreak{}lipschitz\_\allowbreak{}domain\allowbreak\ remaining};
\item either \texttt{r=0}, or \texttt{recovered} is nonempty, satisfies \texttt{cgu\_\allowbreak{}lipschitz\_\allowbreak{}domain}, satisfies \texttt{cgu\_\allowbreak{}face\_\allowbreak{}connected\allowbreak\ recovered\allowbreak\ (cgu\_\allowbreak{}before\_\allowbreak{}indices\allowbreak\ M\allowbreak\ r)\allowbreak\ (cgu\_\allowbreak{}cell\allowbreak\ M)}, and \texttt{cgu\_\allowbreak{}auxiliary\_\allowbreak{}half\_\allowbreak{}ball\_\allowbreak{}geometry\_\allowbreak{}step\allowbreak\ M\allowbreak\ r} holds;
\item every \texttt{i\allowbreak\ in\allowbreak\ I} has \texttt{cgu\_\allowbreak{}flat\_\allowbreak{}face\_\allowbreak{}patch\allowbreak\ (cgu\_\allowbreak{}cell\allowbreak\ M\allowbreak\ c)\allowbreak\ (cgu\_\allowbreak{}normal\allowbreak\ M\allowbreak\ c\allowbreak\ i)\allowbreak\ (cgu\_\allowbreak{}offset\allowbreak\ M\allowbreak\ c\allowbreak\ i)\allowbreak\ (cgu\_\allowbreak{}face\_\allowbreak{}patch\allowbreak\ M\allowbreak\ c\allowbreak\ i)}, has \texttt{cgu\_\allowbreak{}face\_\allowbreak{}patch\allowbreak\ M\allowbreak\ c\allowbreak\ i} contained in the complement of \texttt{cgu\_\allowbreak{}subdivision\_\allowbreak{}edge\_\allowbreak{}corner\_\allowbreak{}set\allowbreak\ M}, has it contained in \texttt{frontier\allowbreak\ remaining}, and satisfies \texttt{cgu\_\allowbreak{}active\_\allowbreak{}pde\_\allowbreak{}face\_\allowbreak{}route\_\allowbreak{}v4\allowbreak\ M\allowbreak\ r\allowbreak\ c\allowbreak\ i};
\item \texttt{cgu\_\allowbreak{}degree\allowbreak\ M\allowbreak\ c\allowbreak\ +\allowbreak\ 3\allowbreak\ <=\allowbreak\ card\allowbreak\ I};
\item \texttt{cgu\_\allowbreak{}n\_\allowbreak{}generic\_\allowbreak{}on\allowbreak\ (cgu\_\allowbreak{}degree\allowbreak\ M\allowbreak\ c)\allowbreak\ I\allowbreak\ (cgu\_\allowbreak{}normal\allowbreak\ M\allowbreak\ c)\allowbreak\ (cgu\_\allowbreak{}offset\allowbreak\ M\allowbreak\ c)\allowbreak\ (cgu\_\allowbreak{}witness\_\allowbreak{}faces\allowbreak\ M\allowbreak\ c)}; and
\item for every \texttt{i\allowbreak\ in\allowbreak\ I}, every \texttt{j,k\allowbreak\ in\allowbreak\ cgu\_\allowbreak{}witness\_\allowbreak{}faces\allowbreak\ M\allowbreak\ c\allowbreak\ i} with \texttt{j\allowbreak\ !=\allowbreak\ k}, the set \texttt{cgu\_\allowbreak{}triple\_\allowbreak{}region\allowbreak\ M\allowbreak\ c\allowbreak\ i\allowbreak\ j\allowbreak\ k} satisfies \texttt{cgu\_\allowbreak{}triple\_\allowbreak{}region\_\allowbreak{}geometry} for the three normal/offset pairs belonging to \texttt{i,j,k}.
\end{itemize}

Second, for every \texttt{r\allowbreak\ <\allowbreak\ CARD('c)}, with the same definitions of \texttt{remaining} and \texttt{recovered}, \texttt{cgu\_\allowbreak{}graph\_\allowbreak{}lipschitz\_\allowbreak{}domain\_\allowbreak{}v5\allowbreak\ remaining} holds and either \texttt{r=0} or \texttt{cgu\_\allowbreak{}graph\_\allowbreak{}lipschitz\_\allowbreak{}domain\_\allowbreak{}v5\allowbreak\ recovered} holds.

\dossierentry{\texttt{\detokenize{cgu_auxiliary_extension_elliptic}}}
\reviewlabel{Isabelle code}
\begin{lstlisting}
definition cgu_auxiliary_extension_elliptic ::
  "('n::finite, 'c::finite, 'h) cgu_model \<Rightarrow>
    ('c \<Rightarrow> 'n cgu_matrix_polynomial) \<Rightarrow> nat \<Rightarrow> bool"
where
  "cgu_auxiliary_extension_elliptic M P r \<longleftrightarrow>
    (let c = cgu_order M r;
         adjacent = cgu_step_adjacent_cell M c
     in cgu_uniformly_positive_definite_on
          (cgu_exterior_half_ball (cgu_step_center M c)
            (cgu_step_interior_normal M c) (cgu_step_radius M c))
          (P adjacent))"
\end{lstlisting}
\reviewlabel{English mathematical translation}
\noindent\textit{Mathematical role:} exterior-half-ball ellipticity condition\par\smallskip
Let \texttt{c=cgu\_\allowbreak{}order\allowbreak\ M\allowbreak\ r} and \texttt{adjacent=cgu\_\allowbreak{}step\_\allowbreak{}adjacent\_\allowbreak{}cell\allowbreak\ M\allowbreak\ c}. Then \texttt{cgu\_\allowbreak{}auxiliary\_\allowbreak{}extension\_\allowbreak{}elliptic\allowbreak\ M\allowbreak\ P\allowbreak\ r} is \texttt{cgu\_\allowbreak{}uniformly\_\allowbreak{}positive\_\allowbreak{}definite\_\allowbreak{}on\allowbreak\ (cgu\_\allowbreak{}exterior\_\allowbreak{}half\_\allowbreak{}ball\allowbreak\ (cgu\_\allowbreak{}step\_\allowbreak{}center\allowbreak\ M\allowbreak\ c)\allowbreak\ (cgu\_\allowbreak{}step\_\allowbreak{}interior\_\allowbreak{}normal\allowbreak\ M\allowbreak\ c)\allowbreak\ (cgu\_\allowbreak{}step\_\allowbreak{}radius\allowbreak\ M\allowbreak\ c))\allowbreak\ (P\allowbreak\ adjacent)}.

\dossierentry{\texttt{\detokenize{cgu_recovery_extension_admissible}}}
\reviewlabel{Isabelle code}
\begin{lstlisting}
definition cgu_recovery_extension_admissible ::
  "('n::finite, 'c::finite, 'h::finite) cgu_model \<Rightarrow>
    ('c \<Rightarrow> 'n cgu_matrix_polynomial) \<Rightarrow> bool"
where
  "cgu_recovery_extension_admissible M P \<longleftrightarrow>
    (\<forall>r<CARD('c).
      let c = cgu_order M r;
          I = cgu_faces M c
      in (r = 0 \<or> cgu_auxiliary_extension_elliptic M P r) \<and>
         (\<forall>i\<in>I. \<forall>j\<in>cgu_witness_faces M c i.
           \<forall>k\<in>cgu_witness_faces M c i. j \<noteq> k \<longrightarrow>
             cgu_polynomial_positive_on (P c)
              (cgu_triple_region M c i j k)))"
\end{lstlisting}
\reviewlabel{English mathematical translation}
\noindent\textit{Mathematical role:} coefficient compatibility with recovery\par\smallskip
For every natural \texttt{r\allowbreak\ <\allowbreak\ CARD('c)}, let \texttt{c=cgu\_\allowbreak{}order\allowbreak\ M\allowbreak\ r} and \texttt{I=cgu\_\allowbreak{}faces\allowbreak\ M\allowbreak\ c}. Then \texttt{cgu\_\allowbreak{}recovery\_\allowbreak{}extension\_\allowbreak{}admissible\allowbreak\ M\allowbreak\ P} requires both:

\begin{itemize}
\item \texttt{r=0} or \texttt{cgu\_\allowbreak{}auxiliary\_\allowbreak{}extension\_\allowbreak{}elliptic\allowbreak\ M\allowbreak\ P\allowbreak\ r}; and
\item for every \texttt{i\allowbreak\ in\allowbreak\ I} and every \texttt{j,k\allowbreak\ in\allowbreak\ cgu\_\allowbreak{}witness\_\allowbreak{}faces\allowbreak\ M\allowbreak\ c\allowbreak\ i} with \texttt{j\allowbreak\ !=\allowbreak\ k}, \texttt{cgu\_\allowbreak{}polynomial\_\allowbreak{}positive\_\allowbreak{}on\allowbreak\ (P\allowbreak\ c)\allowbreak\ (cgu\_\allowbreak{}triple\_\allowbreak{}region\allowbreak\ M\allowbreak\ c\allowbreak\ i\allowbreak\ j\allowbreak\ k)}.
\end{itemize}

The two requirements are within the scope of the universal quantifier over \texttt{r}.

\dossierentry{\texttt{\detokenize{cgu_potential}}}
\reviewlabel{Isabelle code}
\begin{lstlisting}
type_synonym 'n cgu_potential = "'n cgu_point \<Rightarrow> real"
\end{lstlisting}
\reviewlabel{English mathematical translation}
\noindent\textit{Mathematical role:} scalar functions on \(\mathbb R^n\)\par\smallskip
For any \texttt{'n}, \texttt{cgu\_\allowbreak{}potential} is the type of real-valued functions on \texttt{cgu\_\allowbreak{}point} vectors.

\dossierentry{\texttt{\detokenize{cgu_gradient}}}
\reviewlabel{Isabelle code}
\begin{lstlisting}
type_synonym 'n cgu_gradient = "'n cgu_point \<Rightarrow> 'n cgu_point"
\end{lstlisting}
\reviewlabel{English mathematical translation}
\noindent\textit{Mathematical role:} vector functions on \(\mathbb R^n\)\par\smallskip
For any \texttt{'n}, \texttt{cgu\_\allowbreak{}gradient} is the type of functions from \texttt{cgu\_\allowbreak{}point} vectors to \texttt{cgu\_\allowbreak{}point} vectors.

\dossierentry{\texttt{\detokenize{cgu_partial_derivative}}}
\reviewlabel{Isabelle code}
\begin{lstlisting}
definition cgu_partial_derivative ::
  "('n::finite cgu_potential) \<Rightarrow> 'n \<Rightarrow> 'n cgu_point \<Rightarrow> real"
where
  "cgu_partial_derivative phi i x =
    frechet_derivative phi (at x) (axis i 1)"
\end{lstlisting}
\reviewlabel{English mathematical translation}
\noindent\textit{Mathematical role:} coordinate derivative operator \(\partial_i\phi\)\par\smallskip
For finite \texttt{'n}, \texttt{cgu\_\allowbreak{}partial\_\allowbreak{}derivative\allowbreak\ phi\allowbreak\ i\allowbreak\ x} is the Fréchet derivative selected by Isabelle's \texttt{frechet\_\allowbreak{}derivative} operator for \texttt{phi} at \texttt{x}, applied to the coordinate-axis vector \texttt{axis\allowbreak\ i\allowbreak\ 1}. The definition itself does not add a differentiability premise.

\dossierentry{\texttt{\detokenize{cgu_classical_gradient}}}
\reviewlabel{Isabelle code}
\begin{lstlisting}
definition cgu_classical_gradient ::
  "'n::finite cgu_potential \<Rightarrow> 'n cgu_gradient"
where
  "cgu_classical_gradient phi x =
    (\<chi> i. cgu_partial_derivative phi i x)"
\end{lstlisting}
\reviewlabel{English mathematical translation}
\noindent\textit{Mathematical role:} gradient vector \(\nabla\phi\)\par\smallskip
For finite \texttt{'n}, \texttt{cgu\_\allowbreak{}classical\_\allowbreak{}gradient\allowbreak\ phi\allowbreak\ x} is the vector whose coordinate \texttt{i} is \texttt{cgu\_\allowbreak{}partial\_\allowbreak{}derivative\allowbreak\ phi\allowbreak\ i\allowbreak\ x}.

\dossierentry{\texttt{\detokenize{cgu_test_function_on}}}
\reviewlabel{Isabelle code}
\begin{lstlisting}
definition cgu_test_function_on ::
  "'n::finite cgu_point set \<Rightarrow> 'n cgu_potential \<Rightarrow> bool"
where
  "cgu_test_function_on U phi \<longleftrightarrow>
    smooth_on UNIV phi \<and>
    compact (closure {x. phi x \<noteq> 0}) \<and>
    closure {x. phi x \<noteq> 0} \<subseteq> U"
\end{lstlisting}
\reviewlabel{English mathematical translation}
\noindent\textit{Mathematical role:} test class \(\mathscr D(U)\)\par\smallskip
\texttt{cgu\_\allowbreak{}test\_\allowbreak{}function\_\allowbreak{}on\allowbreak\ U\allowbreak\ phi} holds exactly when \texttt{phi} is smooth on the entire vector space, the closure of \texttt{\{x\allowbreak\ |\allowbreak\ phi\allowbreak\ x\allowbreak\ !=\allowbreak\ 0\}} is compact, and that closure is contained in \texttt{U}.

\dossierentry{\texttt{\detokenize{cgu_weak_gradient_on}}}
\reviewlabel{Isabelle code}
\begin{lstlisting}
definition cgu_weak_gradient_on ::
  "'n::finite cgu_point set \<Rightarrow> 'n cgu_potential \<Rightarrow>
    'n cgu_gradient \<Rightarrow> bool"
where
  "cgu_weak_gradient_on U u Du \<longleftrightarrow>
    (\<forall>phi. cgu_test_function_on U phi \<longrightarrow>
      (\<forall>i.
        set_integrable lborel U
          (\<lambda>x. u x * cgu_partial_derivative phi i x) \<and>
        set_integrable lborel U (\<lambda>x. Du x $ i * phi x) \<and>
        set_lebesgue_integral lborel U
          (\<lambda>x. u x * cgu_partial_derivative phi i x) =
        - set_lebesgue_integral lborel U (\<lambda>x. Du x $ i * phi x)))"
\end{lstlisting}
\reviewlabel{English mathematical translation}
\noindent\textit{Mathematical role:} weak-gradient integration-by-parts identity\par\smallskip
\texttt{cgu\_\allowbreak{}weak\_\allowbreak{}gradient\_\allowbreak{}on\allowbreak\ U\allowbreak\ u\allowbreak\ Du} holds exactly when, for every \texttt{phi} satisfying \texttt{cgu\_\allowbreak{}test\_\allowbreak{}function\_\allowbreak{}on\allowbreak\ U\allowbreak\ phi}, and for every coordinate \texttt{i}:

\begin{itemize}
\item \texttt{x\allowbreak\ |-\allowbreak{}>\allowbreak\ u\allowbreak\ x\allowbreak\ *\allowbreak\ cgu\_\allowbreak{}partial\_\allowbreak{}derivative\allowbreak\ phi\allowbreak\ i\allowbreak\ x} is Lebesgue-integrable on \texttt{U};
\item \texttt{x\allowbreak\ |-\allowbreak{}>\allowbreak\ (Du\allowbreak\ x)\_\allowbreak{}i\allowbreak\ *\allowbreak\ phi\allowbreak\ x} is Lebesgue-integrable on \texttt{U}; and
\item the integral over \texttt{U} of the first function is the negative of the integral over \texttt{U} of the second.
\end{itemize}

The quantifier over \texttt{i} lies inside the implication whose antecedent is \texttt{cgu\_\allowbreak{}test\_\allowbreak{}function\_\allowbreak{}on\allowbreak\ U\allowbreak\ phi}.

\dossierentry{\texttt{\detokenize{cgu_h1_pair_on}}}
\reviewlabel{Isabelle code}
\begin{lstlisting}
definition cgu_h1_pair_on ::
  "'n::finite cgu_point set \<Rightarrow> 'n cgu_potential \<Rightarrow>
    'n cgu_gradient \<Rightarrow> bool"
where
  "cgu_h1_pair_on U u Du \<longleftrightarrow>
    u \<in> borel_measurable (restrict_space lborel U) \<and>
    Du \<in> borel_measurable (restrict_space lborel U) \<and>
    cgu_weak_gradient_on U u Du \<and>
    set_integrable lborel U (\<lambda>x. u x ^ 2) \<and>
    set_integrable lborel U (\<lambda>x. norm (Du x) ^ 2)"
\end{lstlisting}
\reviewlabel{English mathematical translation}
\noindent\textit{Mathematical role:} weak \(H^1\) pair \(\mathcal H^1(U)\)\par\smallskip
\texttt{cgu\_\allowbreak{}h1\_\allowbreak{}pair\_\allowbreak{}on\allowbreak\ U\allowbreak\ u\allowbreak\ Du} holds exactly when \texttt{u} and \texttt{Du} are Borel-measurable with respect to Lebesgue measure restricted to \texttt{U}, \texttt{cgu\_\allowbreak{}weak\_\allowbreak{}gradient\_\allowbreak{}on\allowbreak\ U\allowbreak\ u\allowbreak\ Du} holds, and the functions \texttt{x\allowbreak\ |-\allowbreak{}>\allowbreak\ (u\allowbreak\ x)\textasciicircum{}2} and \texttt{x\allowbreak\ |-\allowbreak{}>\allowbreak\ ||Du\allowbreak\ x||\textasciicircum{}2} are each Lebesgue-integrable on \texttt{U}.

\dossierentry{\texttt{\detokenize{cgu_h1_squared_distance}}}
\reviewlabel{Isabelle code}
\begin{lstlisting}
definition cgu_h1_squared_distance ::
  "'n::finite cgu_point set \<Rightarrow>
    'n cgu_potential \<Rightarrow> 'n cgu_gradient \<Rightarrow>
    'n cgu_potential \<Rightarrow> 'n cgu_gradient \<Rightarrow> real"
where
  "cgu_h1_squared_distance U u Du v Dv =
    set_lebesgue_integral lborel U
      (\<lambda>x. (u x - v x) ^ 2 + norm (Du x - Dv x) ^ 2)"
\end{lstlisting}
\reviewlabel{English mathematical translation}
\noindent\textit{Mathematical role:} squared pair-distance \(d_U^2\)\par\smallskip
\texttt{cgu\_\allowbreak{}h1\_\allowbreak{}squared\_\allowbreak{}distance\allowbreak\ U\allowbreak\ u\allowbreak\ Du\allowbreak\ v\allowbreak\ Dv} is the Lebesgue integral over \texttt{U} of \texttt{(u\allowbreak\ x-\allowbreak{}v\allowbreak\ x)\textasciicircum{}2\allowbreak\ +\allowbreak\ ||Du\allowbreak\ x-\allowbreak{}Dv\allowbreak\ x||\textasciicircum{}2}.

\dossierentry{\texttt{\detokenize{cgu_h1_zero_pair_on}}}
\reviewlabel{Isabelle code}
\begin{lstlisting}
definition cgu_h1_zero_pair_on ::
  "(real ^ 'n::finite) set \<Rightarrow>
    ((real ^ 'n) \<Rightarrow> real) \<Rightarrow>
    ((real ^ 'n) \<Rightarrow> (real ^ 'n)) \<Rightarrow> bool"
where
  "cgu_h1_zero_pair_on U u Du \<longleftrightarrow>
    cgu_h1_pair_on U u Du \<and>
    (\<exists>(seq :: nat \<Rightarrow> 'n cgu_potential).
      (\<forall>m. cgu_test_function_on U (seq m) \<and>
        cgu_h1_pair_on U (seq m) (cgu_classical_gradient (seq m))) \<and>
      (\<forall>epsilon>0. \<exists>N. \<forall>m\<ge>N.
        cgu_h1_squared_distance U (seq m)
          (cgu_classical_gradient (seq m)) u Du < epsilon))"
\end{lstlisting}
\reviewlabel{English mathematical translation}
\noindent\textit{Mathematical role:} closure-defined \(\mathcal H_0^1(U)\)\par\smallskip
\texttt{cgu\_\allowbreak{}h1\_\allowbreak{}zero\_\allowbreak{}pair\_\allowbreak{}on\allowbreak\ U\allowbreak\ u\allowbreak\ Du} holds exactly when \texttt{cgu\_\allowbreak{}h1\_\allowbreak{}pair\_\allowbreak{}on\allowbreak\ U\allowbreak\ u\allowbreak\ Du} holds and there exists a sequence \texttt{seq} of real-valued functions such that:

\begin{itemize}
\item for every natural \texttt{m}, both \texttt{cgu\_\allowbreak{}test\_\allowbreak{}function\_\allowbreak{}on\allowbreak\ U\allowbreak\ (seq\allowbreak\ m)} and \texttt{cgu\_\allowbreak{}h1\_\allowbreak{}pair\_\allowbreak{}on\allowbreak\ U\allowbreak\ (seq\allowbreak\ m)\allowbreak\ (cgu\_\allowbreak{}classical\_\allowbreak{}gradient\allowbreak\ (seq\allowbreak\ m))} hold; and
\item for every real \texttt{epsilon>0}, there exists \texttt{N} such that for every \texttt{m>=N}, \texttt{cgu\_\allowbreak{}h1\_\allowbreak{}squared\_\allowbreak{}distance\allowbreak\ U\allowbreak\ (seq\allowbreak\ m)\allowbreak\ (cgu\_\allowbreak{}classical\_\allowbreak{}gradient\allowbreak\ (seq\allowbreak\ m))\allowbreak\ u\allowbreak\ Du\allowbreak\ <\allowbreak\ epsilon}.
\end{itemize}

The same sequence must satisfy both conditions.

\dossierentry{\texttt{\detokenize{cgu_cellwise_coefficient}}}
\reviewlabel{Isabelle code}
\begin{lstlisting}
definition cgu_cellwise_coefficient ::
  "('n::finite, 'c::finite, 'h) cgu_model \<Rightarrow>
    ('c \<Rightarrow> 'n cgu_matrix_polynomial) \<Rightarrow>
    'n cgu_point \<Rightarrow> real ^ 'n ^ 'n"
where
  "cgu_cellwise_coefficient M P x =
    (if \<exists>c. x \<in> cgu_cell M c
     then cgu_matrix_poly_eval (P (SOME c. x \<in> cgu_cell M c)) x
     else 0)"
\end{lstlisting}
\reviewlabel{English mathematical translation}
\noindent\textit{Mathematical role:} assembled coefficient field \(A_{M,P}\)\par\smallskip
\texttt{cgu\_\allowbreak{}cellwise\_\allowbreak{}coefficient\allowbreak\ M\allowbreak\ P\allowbreak\ x} is defined by cases. If there exists a \texttt{c} with \texttt{x\allowbreak\ in\allowbreak\ cgu\_\allowbreak{}cell\allowbreak\ M\allowbreak\ c}, it is \texttt{cgu\_\allowbreak{}matrix\_\allowbreak{}poly\_\allowbreak{}eval\allowbreak\ (P\allowbreak\ c*)\allowbreak\ x}, where \texttt{c*} is an Isabelle Hilbert-choice witness selected by \texttt{SOME\allowbreak\ c.\allowbreak{}\allowbreak\ x\allowbreak\ in\allowbreak\ cgu\_\allowbreak{}cell\allowbreak\ M\allowbreak\ c}. If there is no such \texttt{c}, it is the zero matrix. The definition itself neither asserts uniqueness of the witness nor fixes which witness is chosen when several exist.

\dossierentry{\texttt{\detokenize{cgu_boundary_test_on}}}
\reviewlabel{Isabelle code}
\begin{lstlisting}
definition cgu_boundary_test_on ::
  "('n::finite, 'c::finite, 'h::finite) cgu_model \<Rightarrow>
    'n cgu_point set \<Rightarrow> 'n cgu_point set \<Rightarrow>
    'n cgu_potential \<Rightarrow> bool"
where
  "cgu_boundary_test_on M U gamma f \<longleftrightarrow>
    smooth_on UNIV f \<and>
    compact (closure {x. f x \<noteq> 0}) \<and>
    closure {x. x \<in> frontier U \<and> f x \<noteq> 0}
      \<subseteq> gamma - cgu_subdivision_edge_corner_set M \<and>
    cgu_h1_pair_on U f (cgu_classical_gradient f)"
\end{lstlisting}
\reviewlabel{English mathematical translation}
\noindent\textit{Mathematical role:} boundary-supported smooth datum\par\smallskip
\texttt{cgu\_\allowbreak{}boundary\_\allowbreak{}test\_\allowbreak{}on\allowbreak\ M\allowbreak\ U\allowbreak\ gamma\allowbreak\ f} holds exactly when:

\begin{itemize}
\item \texttt{f} is smooth on the entire vector space;
\item \texttt{closure\allowbreak\ \{x\allowbreak\ |\allowbreak\ f\allowbreak\ x\allowbreak\ !=\allowbreak\ 0\}} is compact;
\item \texttt{closure\allowbreak\ \{x\allowbreak\ |\allowbreak\ x\allowbreak\ in\allowbreak\ frontier\allowbreak\ U\allowbreak\ and\allowbreak\ f\allowbreak\ x\allowbreak\ !=\allowbreak\ 0\}} is contained in \texttt{gamma\allowbreak\ -\allowbreak{}\allowbreak\ cgu\_\allowbreak{}subdivision\_\allowbreak{}edge\_\allowbreak{}corner\_\allowbreak{}set\allowbreak\ M}; and
\item \texttt{cgu\_\allowbreak{}h1\_\allowbreak{}pair\_\allowbreak{}on\allowbreak\ U\allowbreak\ f\allowbreak\ (cgu\_\allowbreak{}classical\_\allowbreak{}gradient\allowbreak\ f)} holds.
\end{itemize}

\dossierentry{\texttt{\detokenize{cgu_local_trace_pair_on}}}
\reviewlabel{Isabelle code}
\begin{lstlisting}
definition cgu_local_trace_pair_on ::
  "('n::finite, 'c::finite, 'h::finite) cgu_model \<Rightarrow>
    'n cgu_point set \<Rightarrow> 'n cgu_point set \<Rightarrow>
    'n cgu_potential \<Rightarrow> 'n cgu_gradient \<Rightarrow> bool"
where
  "cgu_local_trace_pair_on M U gamma f Df \<longleftrightarrow>
    cgu_h1_pair_on U f Df \<and>
    (\<exists>(seq :: nat \<Rightarrow> 'n cgu_potential).
      (\<forall>m. cgu_boundary_test_on M U gamma (seq m)) \<and>
      (\<forall>epsilon>0. \<exists>N. \<forall>m\<ge>N.
        \<exists>z Dz. cgu_h1_zero_pair_on U z Dz \<and>
          cgu_h1_squared_distance U (seq m)
            (cgu_classical_gradient (seq m))
            (\<lambda>x. f x + z x) (\<lambda>x. Df x + Dz x) < epsilon))"
\end{lstlisting}
\reviewlabel{English mathematical translation}
\noindent\textit{Mathematical role:} boundary-data closure predicate \(\mathcal T_M(U,\gamma)\)\par\smallskip
\texttt{cgu\_\allowbreak{}local\_\allowbreak{}trace\_\allowbreak{}pair\_\allowbreak{}on\allowbreak\ M\allowbreak\ U\allowbreak\ gamma\allowbreak\ f\allowbreak\ Df} holds exactly when \texttt{cgu\_\allowbreak{}h1\_\allowbreak{}pair\_\allowbreak{}on\allowbreak\ U\allowbreak\ f\allowbreak\ Df} holds and there exists a sequence \texttt{seq} such that:

\begin{itemize}
\item every \texttt{seq\allowbreak\ m} satisfies \texttt{cgu\_\allowbreak{}boundary\_\allowbreak{}test\_\allowbreak{}on\allowbreak\ M\allowbreak\ U\allowbreak\ gamma}; and
\item for every \texttt{epsilon>0}, there is \texttt{N} such that for every \texttt{m>=N} there exist \texttt{z,Dz} satisfying \texttt{cgu\_\allowbreak{}h1\_\allowbreak{}zero\_\allowbreak{}pair\_\allowbreak{}on\allowbreak\ U\allowbreak\ z\allowbreak\ Dz} and \texttt{cgu\_\allowbreak{}h1\_\allowbreak{}squared\_\allowbreak{}distance\allowbreak\ U\allowbreak\ (seq\allowbreak\ m)\allowbreak\ (cgu\_\allowbreak{}classical\_\allowbreak{}gradient\allowbreak\ (seq\allowbreak\ m))\allowbreak\ (x\allowbreak\ |-\allowbreak{}>\allowbreak\ f\allowbreak\ x+z\allowbreak\ x)\allowbreak\ (x\allowbreak\ |-\allowbreak{}>\allowbreak\ Df\allowbreak\ x+Dz\allowbreak\ x)\allowbreak\ <\allowbreak\ epsilon}.
\end{itemize}

The witnesses \texttt{z,Dz} may depend on \texttt{epsilon} and \texttt{m} because their existential quantifier is inside both eventual quantifiers.

\dossierentry{\texttt{\detokenize{cgu_weak_solution_on}}}
\reviewlabel{Isabelle code}
\begin{lstlisting}
definition cgu_weak_solution_on ::
  "('n::finite, 'c::finite, 'h) cgu_model \<Rightarrow>
    ('c \<Rightarrow> 'n cgu_matrix_polynomial) \<Rightarrow>
    'n cgu_point set \<Rightarrow> 'n cgu_potential \<Rightarrow>
    'n cgu_gradient \<Rightarrow> bool"
where
  "cgu_weak_solution_on M P U u Du \<longleftrightarrow>
    cgu_h1_pair_on U u Du \<and>
    (\<forall>phi. cgu_test_function_on U phi \<longrightarrow>
      set_integrable lborel U
        (\<lambda>x. inner (cgu_cellwise_coefficient M P x *v Du x)
          (cgu_classical_gradient phi x)) \<and>
      set_lebesgue_integral lborel U
        (\<lambda>x. inner (cgu_cellwise_coefficient M P x *v Du x)
          (cgu_classical_gradient phi x)) = 0)"
\end{lstlisting}
\reviewlabel{English mathematical translation}
\noindent\textit{Mathematical role:} \(A_{M,P}\)-harmonic weak pair\par\smallskip
\texttt{cgu\_\allowbreak{}weak\_\allowbreak{}solution\_\allowbreak{}on\allowbreak\ M\allowbreak\ P\allowbreak\ U\allowbreak\ u\allowbreak\ Du} holds exactly when \texttt{cgu\_\allowbreak{}h1\_\allowbreak{}pair\_\allowbreak{}on\allowbreak\ U\allowbreak\ u\allowbreak\ Du} holds and, for every \texttt{phi} satisfying \texttt{cgu\_\allowbreak{}test\_\allowbreak{}function\_\allowbreak{}on\allowbreak\ U\allowbreak\ phi}, the function \texttt{x\allowbreak\ |-\allowbreak{}>\allowbreak\ inner\allowbreak\ (cgu\_\allowbreak{}cellwise\_\allowbreak{}coefficient\allowbreak\ M\allowbreak\ P\allowbreak\ x\allowbreak\ *\allowbreak\ Du\allowbreak\ x)\allowbreak\ (cgu\_\allowbreak{}classical\_\allowbreak{}gradient\allowbreak\ phi\allowbreak\ x)} is Lebesgue-integrable on \texttt{U} and has integral zero there.

\dossierentry{\texttt{\detokenize{cgu_dirichlet_solution_for}}}
\reviewlabel{Isabelle code}
\begin{lstlisting}
definition cgu_dirichlet_solution_for ::
  "('n::finite, 'c::finite, 'h) cgu_model \<Rightarrow>
    ('c \<Rightarrow> 'n cgu_matrix_polynomial) \<Rightarrow>
    'n cgu_point set \<Rightarrow> 'n cgu_potential \<Rightarrow>
    'n cgu_gradient \<Rightarrow> 'n cgu_potential \<Rightarrow>
    'n cgu_gradient \<Rightarrow> bool"
where
  "cgu_dirichlet_solution_for M P U f Df u Du \<longleftrightarrow>
    cgu_h1_pair_on U f Df \<and>
    cgu_weak_solution_on M P U u Du \<and>
    cgu_h1_zero_pair_on U (\<lambda>x. u x - f x)
      (\<lambda>x. Du x - Df x)"
\end{lstlisting}
\reviewlabel{English mathematical translation}
\noindent\textit{Mathematical role:} variational Dirichlet solution\par\smallskip
\texttt{cgu\_\allowbreak{}dirichlet\_\allowbreak{}solution\_\allowbreak{}for\allowbreak\ M\allowbreak\ P\allowbreak\ U\allowbreak\ f\allowbreak\ Df\allowbreak\ u\allowbreak\ Du} holds exactly when \texttt{cgu\_\allowbreak{}h1\_\allowbreak{}pair\_\allowbreak{}on\allowbreak\ U\allowbreak\ f\allowbreak\ Df}, \texttt{cgu\_\allowbreak{}weak\_\allowbreak{}solution\_\allowbreak{}on\allowbreak\ M\allowbreak\ P\allowbreak\ U\allowbreak\ u\allowbreak\ Du}, and \texttt{cgu\_\allowbreak{}h1\_\allowbreak{}zero\_\allowbreak{}pair\_\allowbreak{}on\allowbreak\ U\allowbreak\ (u-\allowbreak{}f)\allowbreak\ (Du-\allowbreak{}Df)} all hold, where both subtractions are pointwise.

\dossierentry{\texttt{\detokenize{cgu_dn_energy}}}
\reviewlabel{Isabelle code}
\begin{lstlisting}
definition cgu_dn_energy ::
  "('n::finite, 'c::finite, 'h) cgu_model \<Rightarrow>
    ('c \<Rightarrow> 'n cgu_matrix_polynomial) \<Rightarrow>
    'n cgu_point set \<Rightarrow> 'n cgu_gradient \<Rightarrow>
    'n cgu_gradient \<Rightarrow> real"
where
  "cgu_dn_energy M P U Du Dh =
    set_lebesgue_integral lborel U
      (\<lambda>x. inner (cgu_cellwise_coefficient M P x *v Du x)
        (Dh x))"
\end{lstlisting}
\reviewlabel{English mathematical translation}
\noindent\textit{Mathematical role:} boundary-energy pairing value\par\smallskip
\texttt{cgu\_\allowbreak{}dn\_\allowbreak{}energy\allowbreak\ M\allowbreak\ P\allowbreak\ U\allowbreak\ Du\allowbreak\ Dh} is the Lebesgue integral over \texttt{U} of \texttt{inner\allowbreak\ (cgu\_\allowbreak{}cellwise\_\allowbreak{}coefficient\allowbreak\ M\allowbreak\ P\allowbreak\ x\allowbreak\ *\allowbreak\ Du\allowbreak\ x)\allowbreak\ (Dh\allowbreak\ x)}.

\dossierentry{\texttt{\detokenize{cgu_dn_energy_integrable}}}
\reviewlabel{Isabelle code}
\begin{lstlisting}
definition cgu_dn_energy_integrable ::
  "('n::finite, 'c::finite, 'h) cgu_model \<Rightarrow>
    ('c \<Rightarrow> 'n cgu_matrix_polynomial) \<Rightarrow>
    'n cgu_point set \<Rightarrow> 'n cgu_gradient \<Rightarrow>
    'n cgu_gradient \<Rightarrow> bool"
where
  "cgu_dn_energy_integrable M P U Du Dh \<longleftrightarrow>
    set_integrable lborel U
      (\<lambda>x. inner (cgu_cellwise_coefficient M P x *v Du x)
        (Dh x))"
\end{lstlisting}
\reviewlabel{English mathematical translation}
\noindent\textit{Mathematical role:} boundary-energy integrability predicate\par\smallskip
\texttt{cgu\_\allowbreak{}dn\_\allowbreak{}energy\_\allowbreak{}integrable\allowbreak\ M\allowbreak\ P\allowbreak\ U\allowbreak\ Du\allowbreak\ Dh} holds exactly when the integrand used in \texttt{cgu\_\allowbreak{}dn\_\allowbreak{}energy\allowbreak\ M\allowbreak\ P\allowbreak\ U\allowbreak\ Du\allowbreak\ Dh} is Lebesgue-integrable on \texttt{U}.

\dossierentry{\texttt{\detokenize{cgu_local_dn_equal}}}
\reviewlabel{Isabelle code}
\begin{lstlisting}
definition cgu_local_dn_equal ::
  "('n::finite, 'c::finite, 'h::finite) cgu_model \<Rightarrow>
    ('c \<Rightarrow> 'n cgu_matrix_polynomial) \<Rightarrow>
    ('c \<Rightarrow> 'n cgu_matrix_polynomial) \<Rightarrow>
    'n cgu_point set \<Rightarrow> 'n cgu_point set \<Rightarrow> bool"
where
  "cgu_local_dn_equal M P1 P2 U gamma \<longleftrightarrow>
    (\<forall>f Df. cgu_local_trace_pair_on M U gamma f Df \<longrightarrow>
      (\<exists>u Du. cgu_dirichlet_solution_for M P1 U f Df u Du) \<and>
      (\<exists>u Du. cgu_dirichlet_solution_for M P2 U f Df u Du)) \<and>
    (\<forall>f Df h Dh u1 Du1 u2 Du2.
      cgu_local_trace_pair_on M U gamma f Df \<and>
      cgu_local_trace_pair_on M U gamma h Dh \<and>
      cgu_dirichlet_solution_for M P1 U f Df u1 Du1 \<and>
      cgu_dirichlet_solution_for M P2 U f Df u2 Du2
      \<longrightarrow>
      cgu_dn_energy_integrable M P1 U Du1 Dh \<and>
      cgu_dn_energy_integrable M P2 U Du2 Dh \<and>
      cgu_dn_energy M P1 U Du1 Dh = cgu_dn_energy M P2 U Du2 Dh)"
\end{lstlisting}
\reviewlabel{English mathematical translation}
\noindent\textit{Mathematical role:} equality of local boundary-energy data\par\smallskip
\texttt{cgu\_\allowbreak{}local\_\allowbreak{}dn\_\allowbreak{}equal\allowbreak\ M\allowbreak\ P1\allowbreak\ P2\allowbreak\ U\allowbreak\ gamma} is the conjunction of two clauses.

First, for every \texttt{f,Df}, if \texttt{cgu\_\allowbreak{}local\_\allowbreak{}trace\_\allowbreak{}pair\_\allowbreak{}on\allowbreak\ M\allowbreak\ U\allowbreak\ gamma\allowbreak\ f\allowbreak\ Df}, then there exist \texttt{u,Du} satisfying \texttt{cgu\_\allowbreak{}dirichlet\_\allowbreak{}solution\_\allowbreak{}for\allowbreak\ M\allowbreak\ P1\allowbreak\ U\allowbreak\ f\allowbreak\ Df\allowbreak\ u\allowbreak\ Du}, and there also exist (not necessarily the same) \texttt{u,Du} satisfying the analogous predicate for \texttt{P2}.

Second, for every \texttt{f,Df,h,Dh,u1,Du1,u2,Du2}, if \texttt{f,Df} and \texttt{h,Dh} each satisfy \texttt{cgu\_\allowbreak{}local\_\allowbreak{}trace\_\allowbreak{}pair\_\allowbreak{}on\allowbreak\ M\allowbreak\ U\allowbreak\ gamma}, and \texttt{u1,Du1} satisfies \texttt{cgu\_\allowbreak{}dirichlet\_\allowbreak{}solution\_\allowbreak{}for} for \texttt{P1,f,Df} while \texttt{u2,Du2} satisfies it for \texttt{P2,f,Df}, then both integrability predicates \texttt{cgu\_\allowbreak{}dn\_\allowbreak{}energy\_\allowbreak{}integrable\allowbreak\ M\allowbreak\ P1\allowbreak\ U\allowbreak\ Du1\allowbreak\ Dh} and \texttt{cgu\_\allowbreak{}dn\_\allowbreak{}energy\_\allowbreak{}integrable\allowbreak\ M\allowbreak\ P2\allowbreak\ U\allowbreak\ Du2\allowbreak\ Dh} hold and the two corresponding values of \texttt{cgu\_\allowbreak{}dn\_\allowbreak{}energy} are equal.

\dossierentry{\texttt{\detokenize{lu_metric}}}
\reviewlabel{Isabelle code}
\begin{lstlisting}
type_synonym 'n lu_metric = "'n cgu_point \<Rightarrow> real ^ 'n ^ 'n"
\end{lstlisting}
\reviewlabel{English mathematical translation}
\noindent\textit{Mathematical role:} matrix-field function type\par\smallskip
For any \texttt{'n}, \texttt{lu\_\allowbreak{}metric} is the type of functions from \texttt{cgu\_\allowbreak{}point} vectors to real \texttt{'n}-by-\texttt{'n} matrices.

\dossierentry{\texttt{\detokenize{lu_density}}}
\reviewlabel{Isabelle code}
\begin{lstlisting}
type_synonym 'n lu_density = "'n cgu_point \<Rightarrow> real"
\end{lstlisting}
\reviewlabel{English mathematical translation}
\noindent\textit{Mathematical role:} scalar-field function type\par\smallskip
For any \texttt{'n}, \texttt{lu\_\allowbreak{}density} is the type of real-valued functions on \texttt{cgu\_\allowbreak{}point} vectors.

\dossierentry{\texttt{\detokenize{lu_smooth_metric_density_data_on}}}
\reviewlabel{Isabelle code}
\begin{lstlisting}
definition lu_smooth_metric_density_data_on ::
  "'n::finite cgu_point set \<Rightarrow> 'n lu_metric \<Rightarrow>
    'n lu_metric \<Rightarrow> 'n lu_density \<Rightarrow> bool"
where
  "lu_smooth_metric_density_data_on V g a rho \<longleftrightarrow>
    open V \<and>
    (\<forall>i j. smooth_on V (\<lambda>x. g x $ i $ j)) \<and>
    (\<forall>i j. smooth_on V (\<lambda>x. a x $ i $ j)) \<and>
    smooth_on V rho \<and>
    (\<forall>x\<in>V.
      cgu_symmetric_positive_definite_matrix (g x) \<and>
      cgu_symmetric_positive_definite_matrix (a x) \<and>
      0 < rho x \<and>
      rho x ^ 2 = det (g x) \<and>
      (\<forall>v. a x *v (g x *v v) = rho x *\<^sub>R v))"
\end{lstlisting}
\reviewlabel{English mathematical translation}
\noindent\textit{Mathematical role:} metric/conductivity representation on \(V\)\par\smallskip
\texttt{lu\_\allowbreak{}smooth\_\allowbreak{}metric\_\allowbreak{}density\_\allowbreak{}data\_\allowbreak{}on\allowbreak\ V\allowbreak\ g\allowbreak\ a\allowbreak\ rho} holds exactly when:

\begin{itemize}
\item \texttt{V} is open;
\item for every \texttt{i,j}, the scalar functions \texttt{x\allowbreak\ |-\allowbreak{}>\allowbreak\ g(x)\_\allowbreak{}\{ij\}} and \texttt{x\allowbreak\ |-\allowbreak{}>\allowbreak\ a(x)\_\allowbreak{}\{ij\}} are smooth on \texttt{V};
\item \texttt{rho} is smooth on \texttt{V}; and
\item for every \texttt{x\allowbreak\ in\allowbreak\ V}, both \texttt{g\allowbreak\ x} and \texttt{a\allowbreak\ x} satisfy \texttt{cgu\_\allowbreak{}symmetric\_\allowbreak{}positive\_\allowbreak{}definite\_\allowbreak{}matrix}, \texttt{rho\allowbreak\ x>0}, \texttt{(rho\allowbreak\ x)\textasciicircum{}2\allowbreak\ =\allowbreak\ det(g\allowbreak\ x)}, and for every vector \texttt{v}, \texttt{a\allowbreak\ x\allowbreak\ *\allowbreak\ (g\allowbreak\ x\allowbreak\ *\allowbreak\ v)\allowbreak\ =\allowbreak\ rho\allowbreak\ x\allowbreak\ *\allowbreak\ v}.
\end{itemize}

\dossierentry{\texttt{\detokenize{lu_smooth_metric_extension_near}}}
\reviewlabel{Isabelle code}
\begin{lstlisting}
definition lu_smooth_metric_extension_near ::
  "('n::finite, 'c::finite, 'h) cgu_model \<Rightarrow>
    ('c \<Rightarrow> 'n cgu_matrix_polynomial) \<Rightarrow>
    'n cgu_point set \<Rightarrow> 'n cgu_point set \<Rightarrow>
    'n lu_metric \<Rightarrow> 'n lu_metric \<Rightarrow>
    'n lu_density \<Rightarrow> bool"
where
  "lu_smooth_metric_extension_near M P U gamma g a rho \<longleftrightarrow>
    (\<exists>V. gamma \<subseteq> V \<and>
      lu_smooth_metric_density_data_on V g a rho \<and>
      (\<forall>x\<in>U \<inter> V. a x = cgu_cellwise_coefficient M P x))"
\end{lstlisting}
\reviewlabel{English mathematical translation}
\noindent\textit{Mathematical role:} local representation of \(A_{M,P}\)\par\smallskip
\texttt{lu\_\allowbreak{}smooth\_\allowbreak{}metric\_\allowbreak{}extension\_\allowbreak{}near\allowbreak\ M\allowbreak\ P\allowbreak\ U\allowbreak\ gamma\allowbreak\ g\allowbreak\ a\allowbreak\ rho} holds exactly when there exists a set \texttt{V} containing \texttt{gamma} such that \texttt{lu\_\allowbreak{}smooth\_\allowbreak{}metric\_\allowbreak{}density\_\allowbreak{}data\_\allowbreak{}on\allowbreak\ V\allowbreak\ g\allowbreak\ a\allowbreak\ rho} holds and, for every \texttt{x\allowbreak\ in\allowbreak\ U\allowbreak\ intersect\allowbreak\ V}, \texttt{a\allowbreak\ x\allowbreak\ =\allowbreak\ cgu\_\allowbreak{}cellwise\_\allowbreak{}coefficient\allowbreak\ M\allowbreak\ P\allowbreak\ x}. Openness of \texttt{V} follows only through \texttt{lu\_\allowbreak{}smooth\_\allowbreak{}metric\_\allowbreak{}density\_\allowbreak{}data\_\allowbreak{}on}; no separate relation between \texttt{g} and \texttt{cgu\_\allowbreak{}cellwise\_\allowbreak{}coefficient} is asserted beyond the equations inside \texttt{lu\_\allowbreak{}smooth\_\allowbreak{}metric\_\allowbreak{}density\_\allowbreak{}data\_\allowbreak{}on}.

\dossierentry{\texttt{\detokenize{lu_induced_boundary_metrics_equal_on}}}
\reviewlabel{Isabelle code}
\begin{lstlisting}
definition lu_induced_boundary_metrics_equal_on ::
  "'n::finite cgu_point set \<Rightarrow> 'n cgu_point \<Rightarrow>
    'n lu_metric \<Rightarrow> 'n lu_metric \<Rightarrow> bool"
where
  "lu_induced_boundary_metrics_equal_on gamma normal g1 g2 \<longleftrightarrow>
    (\<forall>x\<in>gamma. \<forall>xi eta.
      inner normal xi = 0 \<and> inner normal eta = 0
      \<longrightarrow>
      inner xi (g1 x *v eta) = inner xi (g2 x *v eta))"
\end{lstlisting}
\reviewlabel{English mathematical translation}
\noindent\textit{Mathematical role:} tangential bilinear-pairing agreement\par\smallskip
\texttt{lu\_\allowbreak{}induced\_\allowbreak{}boundary\_\allowbreak{}metrics\_\allowbreak{}equal\_\allowbreak{}on\allowbreak\ gamma\allowbreak\ normal\allowbreak\ g1\allowbreak\ g2} holds exactly when, for every \texttt{x\allowbreak\ in\allowbreak\ gamma} and all vectors \texttt{xi,eta}, if both are orthogonal to \texttt{normal}, then \texttt{inner\allowbreak\ xi\allowbreak\ (g1\allowbreak\ x\allowbreak\ *\allowbreak\ eta)\allowbreak\ =\allowbreak\ inner\allowbreak\ xi\allowbreak\ (g2\allowbreak\ x\allowbreak\ *\allowbreak\ eta)}.

\dossierentry{\texttt{\detokenize{mclean_coefficient}}}
\reviewlabel{Isabelle code}
\begin{lstlisting}
type_synonym 'n mclean_coefficient =
  "'n cgu_point \<Rightarrow> real ^ 'n ^ 'n"
\end{lstlisting}
\reviewlabel{English mathematical translation}
\noindent\textit{Mathematical role:} general matrix coefficient field type\par\smallskip
For any \texttt{'n}, \texttt{mclean\_\allowbreak{}coefficient} is the same function shape as \texttt{lu\_\allowbreak{}metric}: functions from \texttt{cgu\_\allowbreak{}point} vectors to real \texttt{'n}-by-\texttt{'n} matrices.

\dossierentry{\texttt{\detokenize{mclean_source}}}
\reviewlabel{Isabelle code}
\begin{lstlisting}
type_synonym 'n mclean_source =
  "'n cgu_potential \<Rightarrow> 'n cgu_gradient \<Rightarrow> real"
\end{lstlisting}
\reviewlabel{English mathematical translation}
\noindent\textit{Mathematical role:} two-argument functional type \(\mathcal F(u,G)\)\par\smallskip
For any \texttt{'n}, \texttt{mclean\_\allowbreak{}source} is the type of two-argument real-valued functionals taking first an \texttt{cgu\_\allowbreak{}potential} function and then an \texttt{cgu\_\allowbreak{}gradient} function.

\dossierentry{\texttt{\detokenize{mclean_h1_norm}}}
\reviewlabel{Isabelle code}
\begin{lstlisting}
definition mclean_h1_norm ::
  "'n::finite cgu_point set \<Rightarrow> 'n cgu_potential \<Rightarrow>
    'n cgu_gradient \<Rightarrow> real"
where
  "mclean_h1_norm U u Du =
    sqrt (cgu_h1_squared_distance U u Du (\<lambda>_. 0) (\<lambda>_. 0))"
\end{lstlisting}
\reviewlabel{English mathematical translation}
\noindent\textit{Mathematical role:} pair norm \(\lVert(u,G)\rVert_U\)\par\smallskip
\texttt{mclean\_\allowbreak{}h1\_\allowbreak{}norm\allowbreak\ U\allowbreak\ u\allowbreak\ Du} is the square root of \texttt{cgu\_\allowbreak{}h1\_\allowbreak{}squared\_\allowbreak{}distance\allowbreak\ U\allowbreak\ u\allowbreak\ Du\allowbreak\ 0\allowbreak\ 0}, where both zero arguments are the pointwise zero functions.

\dossierentry{\texttt{\detokenize{mclean_same_trace_on}}}
\reviewlabel{Isabelle code}
\begin{lstlisting}
definition mclean_same_trace_on ::
  "'n::finite cgu_point set \<Rightarrow>
    'n cgu_potential \<Rightarrow> 'n cgu_gradient \<Rightarrow>
    'n cgu_potential \<Rightarrow> 'n cgu_gradient \<Rightarrow> bool"
where
  "mclean_same_trace_on U f Df g Dg \<longleftrightarrow>
    cgu_h1_pair_on U f Df \<and>
    cgu_h1_pair_on U g Dg \<and>
    cgu_h1_zero_pair_on U (\<lambda>x. f x - g x) (\<lambda>x. Df x - Dg x)"
\end{lstlisting}
\reviewlabel{English mathematical translation}
\noindent\textit{Mathematical role:} one-sided zero-boundary-difference relation \(R_U\)\par\smallskip
\texttt{mclean\_\allowbreak{}same\_\allowbreak{}trace\_\allowbreak{}on\allowbreak\ U\allowbreak\ f\allowbreak\ Df\allowbreak\ g\allowbreak\ Dg} holds exactly when both \texttt{(f,Df)} and \texttt{(g,Dg)} satisfy \texttt{cgu\_\allowbreak{}h1\_\allowbreak{}pair\_\allowbreak{}on\allowbreak\ U}, and their pointwise difference \texttt{(f-\allowbreak{}g,Df-\allowbreak{}Dg)} satisfies \texttt{cgu\_\allowbreak{}h1\_\allowbreak{}zero\_\allowbreak{}pair\_\allowbreak{}on\allowbreak\ U}.

\dossierentry{\texttt{\detokenize{mclean_trace_norm}}}
\reviewlabel{Isabelle code}
\begin{lstlisting}
definition mclean_trace_norm ::
  "'n::finite cgu_point set \<Rightarrow> 'n cgu_potential \<Rightarrow>
    'n cgu_gradient \<Rightarrow> real"
where
  "mclean_trace_norm U f Df =
    Inf {r. \<exists>g Dg. mclean_same_trace_on U f Df g Dg \<and>
      r = mclean_h1_norm U g Dg}"
\end{lstlisting}
\reviewlabel{English mathematical translation}
\noindent\textit{Mathematical role:} raw trace-size infimum \(q_U\)\par\smallskip
\texttt{mclean\_\allowbreak{}trace\_\allowbreak{}norm\allowbreak\ U\allowbreak\ f\allowbreak\ Df} is the real infimum of the set of all \texttt{r} for which there exist \texttt{g,Dg} with \texttt{mclean\_\allowbreak{}same\_\allowbreak{}trace\_\allowbreak{}on\allowbreak\ U\allowbreak\ f\allowbreak\ Df\allowbreak\ g\allowbreak\ Dg} and \texttt{r=mclean\_\allowbreak{}h1\_\allowbreak{}norm\allowbreak\ U\allowbreak\ g\allowbreak\ Dg}. The definition does not separately assert that this set is nonempty or bounded below.

\dossierentry{\texttt{\detokenize{mclean_energy}}}
\reviewlabel{Isabelle code}
\begin{lstlisting}
definition mclean_energy ::
  "'n::finite cgu_point set \<Rightarrow> 'n mclean_coefficient \<Rightarrow>
    'n cgu_gradient \<Rightarrow> 'n cgu_gradient \<Rightarrow> real"
where
  "mclean_energy U a Du Dv =
    set_lebesgue_integral lborel U (\<lambda>x. inner (a x *v Du x) (Dv x))"
\end{lstlisting}
\reviewlabel{English mathematical translation}
\noindent\textit{Mathematical role:} bilinear energy \(B_{U,a}\)\par\smallskip
\texttt{mclean\_\allowbreak{}energy\allowbreak\ U\allowbreak\ a\allowbreak\ Du\allowbreak\ Dv} is the Lebesgue integral over \texttt{U} of \texttt{inner\allowbreak\ (a\allowbreak\ x\allowbreak\ *\allowbreak\ Du\allowbreak\ x)\allowbreak\ (Dv\allowbreak\ x)}.

\dossierentry{\texttt{\detokenize{mclean_energy_integrable}}}
\reviewlabel{Isabelle code}
\begin{lstlisting}
definition mclean_energy_integrable ::
  "'n::finite cgu_point set \<Rightarrow> 'n mclean_coefficient \<Rightarrow>
    'n cgu_gradient \<Rightarrow> 'n cgu_gradient \<Rightarrow> bool"
where
  "mclean_energy_integrable U a Du Dv \<longleftrightarrow>
    set_integrable lborel U (\<lambda>x. inner (a x *v Du x) (Dv x))"
\end{lstlisting}
\reviewlabel{English mathematical translation}
\noindent\textit{Mathematical role:} energy-integrability predicate\par\smallskip
\texttt{mclean\_\allowbreak{}energy\_\allowbreak{}integrable\allowbreak\ U\allowbreak\ a\allowbreak\ Du\allowbreak\ Dv} holds exactly when the integrand used by \texttt{mclean\_\allowbreak{}energy\allowbreak\ U\allowbreak\ a\allowbreak\ Du\allowbreak\ Dv} is Lebesgue-integrable on \texttt{U}.

\dossierentry{\texttt{\detokenize{mclean_variational_form_on}}}
\reviewlabel{Isabelle code}
\begin{lstlisting}
definition mclean_variational_form_on ::
  "'n::finite cgu_point set \<Rightarrow> 'n mclean_coefficient \<Rightarrow> bool"
where
  "mclean_variational_form_on U a \<longleftrightarrow>
    (\<exists>coercive continuous. 0 < coercive \<and> 0 \<le> continuous \<and>
      (\<forall>u Du v Dv.
        cgu_h1_pair_on U u Du \<and> cgu_h1_pair_on U v Dv
        \<longrightarrow>
        mclean_energy_integrable U a Du Dv \<and>
        abs (mclean_energy U a Du Dv) \<le>
          continuous * mclean_h1_norm U u Du * mclean_h1_norm U v Dv) \<and>
      (\<forall>z Dz. cgu_h1_zero_pair_on U z Dz \<longrightarrow>
        coercive * mclean_h1_norm U z Dz ^ 2 \<le>
          mclean_energy U a Dz Dz))"
\end{lstlisting}
\reviewlabel{English mathematical translation}
\noindent\textit{Mathematical role:} stable bounded/coercive form\par\smallskip
\texttt{mclean\_\allowbreak{}variational\_\allowbreak{}form\_\allowbreak{}on\allowbreak\ U\allowbreak\ a} holds exactly when there exist real constants \texttt{coercive,continuous} with \texttt{coercive>0} and \texttt{continuous>=0} such that both clauses below hold:

\begin{itemize}
\item for every \texttt{u,Du,v,Dv}, if both pairs satisfy \texttt{cgu\_\allowbreak{}h1\_\allowbreak{}pair\_\allowbreak{}on\allowbreak\ U}, then \texttt{mclean\_\allowbreak{}energy\_\allowbreak{}integrable\allowbreak\ U\allowbreak\ a\allowbreak\ Du\allowbreak\ Dv} and \texttt{abs(mclean\_\allowbreak{}energy\allowbreak\ U\allowbreak\ a\allowbreak\ Du\allowbreak\ Dv)\allowbreak\ <=\allowbreak\ continuous\allowbreak\ *\allowbreak\ mclean\_\allowbreak{}h1\_\allowbreak{}norm\allowbreak\ U\allowbreak\ u\allowbreak\ Du\allowbreak\ *\allowbreak\ mclean\_\allowbreak{}h1\_\allowbreak{}norm\allowbreak\ U\allowbreak\ v\allowbreak\ Dv};
\item for every \texttt{z,Dz} satisfying \texttt{cgu\_\allowbreak{}h1\_\allowbreak{}zero\_\allowbreak{}pair\_\allowbreak{}on\allowbreak\ U}, \texttt{coercive\allowbreak\ *\allowbreak\ (mclean\_\allowbreak{}h1\_\allowbreak{}norm\allowbreak\ U\allowbreak\ z\allowbreak\ Dz)\textasciicircum{}2\allowbreak\ <=\allowbreak\ mclean\_\allowbreak{}energy\allowbreak\ U\allowbreak\ a\allowbreak\ Dz\allowbreak\ Dz}.
\end{itemize}

The same two constants work for all quantified functions.

\dossierentry{\texttt{\detokenize{mclean_weak_solution_on}}}
\reviewlabel{Isabelle code}
\begin{lstlisting}
definition mclean_weak_solution_on ::
  "'n::finite cgu_point set \<Rightarrow> 'n mclean_coefficient \<Rightarrow>
    'n cgu_potential \<Rightarrow> 'n cgu_gradient \<Rightarrow> bool"
where
  "mclean_weak_solution_on U a u Du \<longleftrightarrow>
    cgu_h1_pair_on U u Du \<and>
    (\<forall>phi. cgu_test_function_on U phi \<longrightarrow>
      mclean_energy_integrable U a Du (cgu_classical_gradient phi) \<and>
      mclean_energy U a Du (cgu_classical_gradient phi) = 0)"
\end{lstlisting}
\reviewlabel{English mathematical translation}
\noindent\textit{Mathematical role:} weak \(a\)-harmonic pair\par\smallskip
\texttt{mclean\_\allowbreak{}weak\_\allowbreak{}solution\_\allowbreak{}on\allowbreak\ U\allowbreak\ a\allowbreak\ u\allowbreak\ Du} holds exactly when \texttt{cgu\_\allowbreak{}h1\_\allowbreak{}pair\_\allowbreak{}on\allowbreak\ U\allowbreak\ u\allowbreak\ Du} and, for every \texttt{phi} satisfying \texttt{cgu\_\allowbreak{}test\_\allowbreak{}function\_\allowbreak{}on\allowbreak\ U\allowbreak\ phi}, both \texttt{mclean\_\allowbreak{}energy\_\allowbreak{}integrable\allowbreak\ U\allowbreak\ a\allowbreak\ Du\allowbreak\ (cgu\_\allowbreak{}classical\_\allowbreak{}gradient\allowbreak\ phi)} and \texttt{mclean\_\allowbreak{}energy\allowbreak\ U\allowbreak\ a\allowbreak\ Du\allowbreak\ (cgu\_\allowbreak{}classical\_\allowbreak{}gradient\allowbreak\ phi)=0} hold.

\dossierentry{\texttt{\detokenize{mclean_dirichlet_solution_for}}}
\reviewlabel{Isabelle code}
\begin{lstlisting}
definition mclean_dirichlet_solution_for ::
  "'n::finite cgu_point set \<Rightarrow> 'n mclean_coefficient \<Rightarrow>
    'n cgu_potential \<Rightarrow> 'n cgu_gradient \<Rightarrow>
    'n cgu_potential \<Rightarrow> 'n cgu_gradient \<Rightarrow> bool"
where
  "mclean_dirichlet_solution_for U a f Df u Du \<longleftrightarrow>
    cgu_h1_pair_on U f Df \<and>
    mclean_weak_solution_on U a u Du \<and>
    mclean_same_trace_on U f Df u Du"
\end{lstlisting}
\reviewlabel{English mathematical translation}
\noindent\textit{Mathematical role:} abstract Dirichlet solution\par\smallskip
\texttt{mclean\_\allowbreak{}dirichlet\_\allowbreak{}solution\_\allowbreak{}for\allowbreak\ U\allowbreak\ a\allowbreak\ f\allowbreak\ Df\allowbreak\ u\allowbreak\ Du} holds exactly when \texttt{cgu\_\allowbreak{}h1\_\allowbreak{}pair\_\allowbreak{}on\allowbreak\ U\allowbreak\ f\allowbreak\ Df}, \texttt{mclean\_\allowbreak{}weak\_\allowbreak{}solution\_\allowbreak{}on\allowbreak\ U\allowbreak\ a\allowbreak\ u\allowbreak\ Du}, and \texttt{mclean\_\allowbreak{}same\_\allowbreak{}trace\_\allowbreak{}on\allowbreak\ U\allowbreak\ f\allowbreak\ Df\allowbreak\ u\allowbreak\ Du} all hold.

\dossierentry{\texttt{\detokenize{mclean_green_solution_for}}}
\reviewlabel{Isabelle code}
\begin{lstlisting}
definition mclean_green_solution_for ::
  "'n::finite cgu_point set \<Rightarrow> 'n mclean_coefficient \<Rightarrow>
    'n mclean_source \<Rightarrow> 'n cgu_potential \<Rightarrow>
    'n cgu_gradient \<Rightarrow> bool"
where
  "mclean_green_solution_for U a F u Du \<longleftrightarrow>
    cgu_h1_zero_pair_on U u Du \<and>
    (\<forall>v Dv. cgu_h1_zero_pair_on U v Dv \<longrightarrow>
      mclean_energy_integrable U a Du Dv \<and>
      mclean_energy U a Du Dv = F v Dv)"
\end{lstlisting}
\reviewlabel{English mathematical translation}
\noindent\textit{Mathematical role:} source solution\par\smallskip
\texttt{mclean\_\allowbreak{}green\_\allowbreak{}solution\_\allowbreak{}for\allowbreak\ U\allowbreak\ a\allowbreak\ F\allowbreak\ u\allowbreak\ Du} holds exactly when \texttt{cgu\_\allowbreak{}h1\_\allowbreak{}zero\_\allowbreak{}pair\_\allowbreak{}on\allowbreak\ U\allowbreak\ u\allowbreak\ Du} and, for every \texttt{v,Dv} satisfying \texttt{cgu\_\allowbreak{}h1\_\allowbreak{}zero\_\allowbreak{}pair\_\allowbreak{}on\allowbreak\ U}, both \texttt{mclean\_\allowbreak{}energy\_\allowbreak{}integrable\allowbreak\ U\allowbreak\ a\allowbreak\ Du\allowbreak\ Dv} and \texttt{mclean\_\allowbreak{}energy\allowbreak\ U\allowbreak\ a\allowbreak\ Du\allowbreak\ Dv\allowbreak\ =\allowbreak\ F\allowbreak\ v\allowbreak\ Dv} hold.

\dossierentry{\texttt{\detokenize{mclean_source_bound_on}}}
\reviewlabel{Isabelle code}
\begin{lstlisting}
definition mclean_source_bound_on ::
  "'n::finite cgu_point set \<Rightarrow> 'n mclean_source \<Rightarrow> real \<Rightarrow> bool"
where
  "mclean_source_bound_on U F K \<longleftrightarrow>
    0 \<le> K \<and> F (\<lambda>_. 0) (\<lambda>_. 0) = 0 \<and>
    (\<forall>u Du v Dv.
      cgu_h1_zero_pair_on U u Du \<and> cgu_h1_zero_pair_on U v Dv
      \<longrightarrow>
      F (\<lambda>x. u x + v x) (\<lambda>x. Du x + Dv x) = F u Du + F v Dv) \<and>
    (\<forall>r u Du. cgu_h1_zero_pair_on U u Du \<longrightarrow>
      F (\<lambda>x. r * u x) (\<lambda>x. r *\<^sub>R Du x) = r * F u Du) \<and>
    (\<forall>u Du. cgu_h1_zero_pair_on U u Du \<longrightarrow>
      abs (F u Du) \<le> K * mclean_h1_norm U u Du)"
\end{lstlisting}
\reviewlabel{English mathematical translation}
\noindent\textit{Mathematical role:} bounded linear functional with bound \(K\)\par\smallskip
\texttt{mclean\_\allowbreak{}source\_\allowbreak{}bound\_\allowbreak{}on\allowbreak\ U\allowbreak\ F\allowbreak\ K} holds exactly when:

\begin{itemize}
\item \texttt{K>=0} and \texttt{F\allowbreak\ 0\allowbreak\ 0\allowbreak\ =\allowbreak\ 0};
\item for all \texttt{(u,Du)} and \texttt{(v,Dv)} satisfying \texttt{cgu\_\allowbreak{}h1\_\allowbreak{}zero\_\allowbreak{}pair\_\allowbreak{}on\allowbreak\ U}, \texttt{F(u+v,Du+Dv)=F(u,Du)+F(v,Dv)};
\item for every real \texttt{r} and every \texttt{(u,Du)} satisfying \texttt{cgu\_\allowbreak{}h1\_\allowbreak{}zero\_\allowbreak{}pair\_\allowbreak{}on\allowbreak\ U}, \texttt{F(r*u,r*Du)=r*F(u,Du)}; and
\item for every \texttt{(u,Du)} satisfying \texttt{cgu\_\allowbreak{}h1\_\allowbreak{}zero\_\allowbreak{}pair\_\allowbreak{}on\allowbreak\ U}, \texttt{abs(F\allowbreak\ u\allowbreak\ Du)\allowbreak\ <=\allowbreak\ K\allowbreak\ *\allowbreak\ mclean\_\allowbreak{}h1\_\allowbreak{}norm\allowbreak\ U\allowbreak\ u\allowbreak\ Du}.
\end{itemize}

The equations constrain \texttt{F} only on arguments covered by the stated \texttt{cgu\_\allowbreak{}h1\_\allowbreak{}zero\_\allowbreak{}pair\_\allowbreak{}on} premises, apart from the explicit value at the zero pair.

\dossierentry{\texttt{\detokenize{mclean_flat_chart}}}
\reviewlabel{Isabelle code}
\begin{lstlisting}
definition mclean_flat_chart ::
  "('m::finite cgu_point \<Rightarrow> 'n::finite cgu_point) \<Rightarrow>
    'n cgu_point \<Rightarrow> bool"
where
  "mclean_flat_chart chart normal \<longleftrightarrow>
    norm normal = 1 \<and>
    (\<forall>x y. norm (chart x - chart y) = norm (x - y)) \<and>
    (\<forall>x y. inner normal (chart x - chart y) = 0)"
\end{lstlisting}
\reviewlabel{English mathematical translation}
\noindent\textit{Mathematical role:} isometric hyperplane chart condition\par\smallskip
\texttt{mclean\_\allowbreak{}flat\_\allowbreak{}chart\allowbreak\ chart\allowbreak\ normal} holds exactly when \texttt{||normal||=1}, \texttt{chart} preserves every pairwise distance, and every displacement \texttt{chart\allowbreak\ x-\allowbreak{}chart\allowbreak\ y} is orthogonal to \texttt{normal}.

\dossierentry{\texttt{\detokenize{mclean_flat_offset_action}}}
\reviewlabel{Isabelle code}
\begin{lstlisting}
definition mclean_flat_offset_action ::
  "('m::finite cgu_point \<Rightarrow> 'n::finite cgu_point) \<Rightarrow>
    'n cgu_point \<Rightarrow> real \<Rightarrow> 'm cgu_potential \<Rightarrow>
    'n cgu_potential \<Rightarrow> real"
where
  "mclean_flat_offset_action chart normal t g phi =
    integral\<^sup>L lborel
      (\<lambda>x. g x * phi (chart x + t *\<^sub>R normal))"
\end{lstlisting}
\reviewlabel{English mathematical translation}
\noindent\textit{Mathematical role:} moving-layer integral\par\smallskip
\texttt{mclean\_\allowbreak{}flat\_\allowbreak{}offset\_\allowbreak{}action\allowbreak\ chart\allowbreak\ normal\allowbreak\ t\allowbreak\ g\allowbreak\ phi} is the Lebesgue integral over the entire \texttt{'m}-coordinate space of \texttt{g\allowbreak\ x\allowbreak\ *\allowbreak\ phi(chart\allowbreak\ x\allowbreak\ +\allowbreak\ t*normal)}.

\dossierentry{\texttt{\detokenize{mclean_source_extends_flat_offset_action}}}
\reviewlabel{Isabelle code}
\begin{lstlisting}
definition mclean_source_extends_flat_offset_action ::
  "'n::finite cgu_point set \<Rightarrow>
    ('m::finite cgu_point \<Rightarrow> 'n cgu_point) \<Rightarrow>
    'n cgu_point \<Rightarrow> real \<Rightarrow> 'm cgu_potential \<Rightarrow>
    'n mclean_source \<Rightarrow> bool"
where
  "mclean_source_extends_flat_offset_action Omega chart normal t g F \<longleftrightarrow>
    (\<forall>phi. cgu_test_function_on Omega phi \<longrightarrow>
      F phi (cgu_classical_gradient phi) =
        mclean_flat_offset_action chart normal t g phi)"
\end{lstlisting}
\reviewlabel{English mathematical translation}
\noindent\textit{Mathematical role:} functional agreement with a layer source\par\smallskip
\texttt{mclean\_\allowbreak{}source\_\allowbreak{}extends\_\allowbreak{}flat\_\allowbreak{}offset\_\allowbreak{}action\allowbreak\ Omega\allowbreak\ chart\allowbreak\ normal\allowbreak\ t\allowbreak\ g\allowbreak\ F} holds exactly when, for every \texttt{phi} satisfying \texttt{cgu\_\allowbreak{}test\_\allowbreak{}function\_\allowbreak{}on\allowbreak\ Omega\allowbreak\ phi}, \texttt{F\allowbreak\ phi\allowbreak\ (cgu\_\allowbreak{}classical\_\allowbreak{}gradient\allowbreak\ phi)\allowbreak\ =\allowbreak\ mclean\_\allowbreak{}flat\_\allowbreak{}offset\_\allowbreak{}action\allowbreak\ chart\allowbreak\ normal\allowbreak\ t\allowbreak\ g\allowbreak\ phi}.

\dossierentry{\texttt{\detokenize{mclean_flat_offset_cylinder_inside}}}
\reviewlabel{Isabelle code}
\begin{lstlisting}
definition mclean_flat_offset_cylinder_inside ::
  "'n::finite cgu_point set \<Rightarrow>
    ('m::finite cgu_point \<Rightarrow> 'n cgu_point) \<Rightarrow>
    'n cgu_point \<Rightarrow> 'm cgu_potential \<Rightarrow> real \<Rightarrow> bool"
where
  "mclean_flat_offset_cylinder_inside Omega chart normal g delta \<longleftrightarrow>
    {chart x + t *\<^sub>R normal |x t.
      x \<in> closure {y. g y \<noteq> 0} \<and> 0 \<le> t \<and> t \<le> delta}
      \<subseteq> Omega"
\end{lstlisting}
\reviewlabel{English mathematical translation}
\noindent\textit{Mathematical role:} tubular support containment\par\smallskip
\texttt{mclean\_\allowbreak{}flat\_\allowbreak{}offset\_\allowbreak{}cylinder\_\allowbreak{}inside\allowbreak\ Omega\allowbreak\ chart\allowbreak\ normal\allowbreak\ g\allowbreak\ delta} holds exactly when every point \texttt{chart\allowbreak\ x\allowbreak\ +\allowbreak\ t*normal} with \texttt{x\allowbreak\ in\allowbreak\ closure\allowbreak\ \{y\allowbreak\ |\allowbreak\ g\allowbreak\ y\allowbreak\ !=\allowbreak\ 0\}} and \texttt{0<=t<=delta} belongs to \texttt{Omega}.

\dossierentry{\texttt{\detokenize{mclean_source_real_subspace}}}
\reviewlabel{Isabelle code}
\begin{lstlisting}
definition mclean_source_real_subspace ::
  "'n::finite mclean_source set \<Rightarrow> bool"
where
  "mclean_source_real_subspace E \<longleftrightarrow>
    0 \<in> E \<and>
    (\<forall>F\<in>E. \<forall>G\<in>E. F + G \<in> E) \<and>
    (\<forall>r F. F \<in> E \<longrightarrow> r *\<^sub>R F \<in> E)"
\end{lstlisting}
\reviewlabel{English mathematical translation}
\noindent\textit{Mathematical role:} vector-space closure of \(E\)\par\smallskip
\texttt{mclean\_\allowbreak{}source\_\allowbreak{}real\_\allowbreak{}subspace\allowbreak\ E} holds exactly when the zero functional belongs to \texttt{E}, \texttt{E} is closed under addition of any two of its members, and \texttt{E} is closed under multiplication of any member by any real scalar.

\dossierentry{\texttt{\detokenize{mclean_source_linear_on}}}
\reviewlabel{Isabelle code}
\begin{lstlisting}
definition mclean_source_linear_on ::
  "'n::finite mclean_source set \<Rightarrow>
    ('n mclean_source \<Rightarrow> real) \<Rightarrow> bool"
where
  "mclean_source_linear_on E g \<longleftrightarrow>
    g 0 = 0 \<and>
    (\<forall>F\<in>E. \<forall>G\<in>E. g (F + G) = g F + g G) \<and>
    (\<forall>r F. F \<in> E \<longrightarrow> g (r *\<^sub>R F) = r * g F)"
\end{lstlisting}
\reviewlabel{English mathematical translation}
\noindent\textit{Mathematical role:} linearity of \(g:E\to\mathbb R\)\par\smallskip
\texttt{mclean\_\allowbreak{}source\_\allowbreak{}linear\_\allowbreak{}on\allowbreak\ E\allowbreak\ g} holds exactly when \texttt{g\allowbreak\ 0=0}, \texttt{g(F+G)=g\allowbreak\ F+g\allowbreak\ G} for all \texttt{F,G\allowbreak\ in\allowbreak\ E}, and \texttt{g(r*F)=r*g\allowbreak\ F} for every real \texttt{r} and every \texttt{F\allowbreak\ in\allowbreak\ E}.

\dossierentry{\texttt{\detokenize{mclean_source_seminorm_on}}}
\reviewlabel{Isabelle code}
\begin{lstlisting}
definition mclean_source_seminorm_on ::
  "'n::finite mclean_source set \<Rightarrow>
    ('n mclean_source \<Rightarrow> real) \<Rightarrow> bool"
where
  "mclean_source_seminorm_on E p \<longleftrightarrow>
    p 0 = 0 \<and>
    (\<forall>F\<in>E. 0 \<le> p F) \<and>
    (\<forall>F\<in>E. \<forall>G\<in>E. p (F + G) \<le> p F + p G) \<and>
    (\<forall>r F. F \<in> E \<longrightarrow> p (r *\<^sub>R F) = \<bar>r\<bar> * p F)"
\end{lstlisting}
\reviewlabel{English mathematical translation}
\noindent\textit{Mathematical role:} seminorm axioms for \(p\)\par\smallskip
\texttt{mclean\_\allowbreak{}source\_\allowbreak{}seminorm\_\allowbreak{}on\allowbreak\ E\allowbreak\ p} holds exactly when \texttt{p\allowbreak\ 0=0}; \texttt{p\allowbreak\ F>=0} for every \texttt{F\allowbreak\ in\allowbreak\ E}; \texttt{p(F+G)<=p\allowbreak\ F+p\allowbreak\ G} for all \texttt{F,G\allowbreak\ in\allowbreak\ E}; and \texttt{p(r*F)=abs(r)*p\allowbreak\ F} for every real \texttt{r} and every \texttt{F\allowbreak\ in\allowbreak\ E}.

\dossierentry{\texttt{\detokenize{mclean_h1_zero_source_bidual_functional_on_v2}}}
\reviewlabel{Isabelle code}
\begin{lstlisting}
definition mclean_h1_zero_source_bidual_functional_on_v2 ::
  "'n::finite mclean_source set \<Rightarrow>
    ('n mclean_source \<Rightarrow> real) \<Rightarrow>
    ('n mclean_source \<Rightarrow> real) \<Rightarrow> bool"
where
  "mclean_h1_zero_source_bidual_functional_on_v2 E p g \<longleftrightarrow>
    mclean_source_linear_on E g \<and>
    (\<exists>C\<ge>0. \<forall>F\<in>E. \<bar>g F\<bar> \<le> C * p F)"
\end{lstlisting}
\reviewlabel{English mathematical translation}
\noindent\textit{Mathematical role:} \(p\)-bounded linear functional on \(E\)\par\smallskip
\texttt{mclean\_\allowbreak{}h1\_\allowbreak{}zero\_\allowbreak{}source\_\allowbreak{}bidual\_\allowbreak{}functional\_\allowbreak{}on\_\allowbreak{}v2\allowbreak\ E\allowbreak\ p\allowbreak\ g} holds exactly when \texttt{mclean\_\allowbreak{}source\_\allowbreak{}linear\_\allowbreak{}on\allowbreak\ E\allowbreak\ g} and there exists one real \texttt{C>=0} such that \texttt{abs(g\allowbreak\ F)\allowbreak\ <=\allowbreak\ C*p\allowbreak\ F} for every \texttt{F\allowbreak\ in\allowbreak\ E}.

\dossierentry{\texttt{\detokenize{mclean_h1_zero_bidual_realization_on_v2}}}
\reviewlabel{Isabelle code}
\begin{lstlisting}
definition mclean_h1_zero_bidual_realization_on_v2 ::
  "'n::finite cgu_point set \<Rightarrow> 'n mclean_source set \<Rightarrow>
    ('n mclean_source \<Rightarrow> real) \<Rightarrow> bool"
where
  "mclean_h1_zero_bidual_realization_on_v2 U E p \<longleftrightarrow>
    mclean_source_real_subspace E \<and>
    mclean_source_seminorm_on E p \<and>
    (\<forall>F. F \<in> E \<longleftrightarrow>
      (\<exists>K. mclean_source_bound_on U F K)) \<and>
    (\<forall>F\<in>E.
      mclean_source_bound_on U F (p F) \<and>
      (\<forall>K. mclean_source_bound_on U F K \<longrightarrow> p F \<le> K)) \<and>
    (\<forall>g.
      mclean_h1_zero_source_bidual_functional_on_v2 E p g
      \<longrightarrow>
      (\<exists>u Du.
        cgu_h1_zero_pair_on U u Du \<and>
        (\<forall>F\<in>E. g F = F u Du)))"
\end{lstlisting}
\reviewlabel{English mathematical translation}
\noindent\textit{Mathematical role:} total-functional representation package \((E,p)\)\par\smallskip
\texttt{mclean\_\allowbreak{}h1\_\allowbreak{}zero\_\allowbreak{}bidual\_\allowbreak{}realization\_\allowbreak{}on\_\allowbreak{}v2\allowbreak\ U\allowbreak\ E\allowbreak\ p} holds exactly when all of the following hold:

\begin{itemize}
\item \texttt{mclean\_\allowbreak{}source\_\allowbreak{}real\_\allowbreak{}subspace\allowbreak\ E} and \texttt{mclean\_\allowbreak{}source\_\allowbreak{}seminorm\_\allowbreak{}on\allowbreak\ E\allowbreak\ p};
\item for every functional \texttt{F}, \texttt{F\allowbreak\ in\allowbreak\ E} if and only if there exists a real \texttt{K} with \texttt{mclean\_\allowbreak{}source\_\allowbreak{}bound\_\allowbreak{}on\allowbreak\ U\allowbreak\ F\allowbreak\ K};
\item for every \texttt{F\allowbreak\ in\allowbreak\ E}, \texttt{mclean\_\allowbreak{}source\_\allowbreak{}bound\_\allowbreak{}on\allowbreak\ U\allowbreak\ F\allowbreak\ (p\allowbreak\ F)}, and \texttt{p\allowbreak\ F\allowbreak\ <=\allowbreak\ K} for every \texttt{K} satisfying \texttt{mclean\_\allowbreak{}source\_\allowbreak{}bound\_\allowbreak{}on\allowbreak\ U\allowbreak\ F\allowbreak\ K}; and
\item for every \texttt{g}, if \texttt{mclean\_\allowbreak{}h1\_\allowbreak{}zero\_\allowbreak{}source\_\allowbreak{}bidual\_\allowbreak{}functional\_\allowbreak{}on\_\allowbreak{}v2\allowbreak\ E\allowbreak\ p\allowbreak\ g}, then there exist \texttt{u,Du} satisfying \texttt{cgu\_\allowbreak{}h1\_\allowbreak{}zero\_\allowbreak{}pair\_\allowbreak{}on\allowbreak\ U\allowbreak\ u\allowbreak\ Du} such that \texttt{g\allowbreak\ F\allowbreak\ =\allowbreak\ F\allowbreak\ u\allowbreak\ Du} for every \texttt{F\allowbreak\ in\allowbreak\ E}.
\end{itemize}

The final representing pair is existential; uniqueness is not asserted.

\dossierentry{\texttt{\detokenize{mclean_zero_extension_on}}}
\reviewlabel{Isabelle code}
\begin{lstlisting}
definition mclean_zero_extension_on ::
  "'n::finite cgu_point set \<Rightarrow> ('n cgu_point \<Rightarrow> 'a::zero) \<Rightarrow>
    'n cgu_point \<Rightarrow> 'a"
where
  "mclean_zero_extension_on U f x = (if x \<in> U then f x else 0)"
\end{lstlisting}
\reviewlabel{English mathematical translation}
\noindent\textit{Mathematical role:} pointwise zero extension\par\smallskip
For a codomain with a zero, \texttt{mclean\_\allowbreak{}zero\_\allowbreak{}extension\_\allowbreak{}on\allowbreak\ U\allowbreak\ f\allowbreak\ x} equals \texttt{f\allowbreak\ x} when \texttt{x\allowbreak\ in\allowbreak\ U} and equals zero otherwise.

\dossierentry{\texttt{\detokenize{ucp_point}}}
\reviewlabel{Isabelle code}
\begin{lstlisting}
type_synonym 'n ucp_point = "real ^ 'n"
\end{lstlisting}
\reviewlabel{English mathematical translation}
\noindent\textit{Mathematical role:} second Euclidean-vector alias\par\smallskip
For any \texttt{'n}, \texttt{ucp\_\allowbreak{}point} is the type of real-valued vectors indexed by \texttt{'n}.

\dossierentry{\texttt{\detokenize{ucp_coefficient}}}
\reviewlabel{Isabelle code}
\begin{lstlisting}
type_synonym 'n ucp_coefficient =
  "'n ucp_point \<Rightarrow> (real ^ 'n ^ 'n)"
\end{lstlisting}
\reviewlabel{English mathematical translation}
\noindent\textit{Mathematical role:} local matrix-field function type\par\smallskip
For any \texttt{'n}, \texttt{ucp\_\allowbreak{}coefficient} is the type of functions from \texttt{ucp\_\allowbreak{}point} vectors to real \texttt{'n}-by-\texttt{'n} matrices.

\dossierentry{\texttt{\detokenize{ucp_symmetric_matrix}}}
\reviewlabel{Isabelle code}
\begin{lstlisting}
definition ucp_symmetric_matrix ::
  "(real ^ 'n::finite ^ 'n) \<Rightarrow> bool"
where
  "ucp_symmetric_matrix A \<longleftrightarrow> transpose A = A"
\end{lstlisting}
\reviewlabel{English mathematical translation}
\noindent\textit{Mathematical role:} matrix symmetry\par\smallskip
\texttt{ucp\_\allowbreak{}symmetric\_\allowbreak{}matrix\allowbreak\ A} holds exactly when \texttt{transpose\allowbreak\ A\allowbreak\ =\allowbreak\ A}.

\dossierentry{\texttt{\detokenize{ucp_matrix_polynomial}}}
\reviewlabel{Isabelle code}
\begin{lstlisting}
definition ucp_matrix_polynomial ::
  "'n::finite ucp_coefficient \<Rightarrow> bool"
where
  "ucp_matrix_polynomial P \<longleftrightarrow>
    (\<forall>i j. polynomial_function (\<lambda>x. P x $ i $ j))"
\end{lstlisting}
\reviewlabel{English mathematical translation}
\noindent\textit{Mathematical role:} entrywise polynomial field\par\smallskip
\texttt{ucp\_\allowbreak{}matrix\_\allowbreak{}polynomial\allowbreak\ P} holds exactly when, for every \texttt{i,j}, the scalar component function \texttt{x\allowbreak\ |-\allowbreak{}>\allowbreak\ P(x)\_\allowbreak{}\{ij\}} is a polynomial function.

\dossierentry{\texttt{\detokenize{ucp_uniformly_elliptic_near}}}
\reviewlabel{Isabelle code}
\begin{lstlisting}
definition ucp_uniformly_elliptic_near ::
  "'n::finite ucp_point set \<Rightarrow> 'n ucp_coefficient \<Rightarrow> bool"
where
  "ucp_uniformly_elliptic_near E P \<longleftrightarrow>
    (\<exists>V k. open V \<and> closure E \<subseteq> V \<and> 0 < k \<and>
      (\<forall>x\<in>V. \<forall>xi. k * norm xi ^ 2 \<le> inner xi (P x *v xi)))"
\end{lstlisting}
\reviewlabel{English mathematical translation}
\noindent\textit{Mathematical role:} ellipticity on a neighborhood of \(\overline X\)\par\smallskip
\texttt{ucp\_\allowbreak{}uniformly\_\allowbreak{}elliptic\_\allowbreak{}near\allowbreak\ E\allowbreak\ P} holds exactly when there exist an open set \texttt{V} and a real \texttt{k>0} such that \texttt{closure\allowbreak\ E} is contained in \texttt{V} and, for every \texttt{x\allowbreak\ in\allowbreak\ V} and every vector \texttt{xi}, \texttt{k*||xi||\textasciicircum{}2\allowbreak\ <=\allowbreak\ inner\allowbreak\ xi\allowbreak\ (P\allowbreak\ x\allowbreak\ *\allowbreak\ xi)}. The same \texttt{V,k} work uniformly for all \texttt{x,xi}.

\dossierentry{\texttt{\detokenize{ucp_partial_derivative}}}
\reviewlabel{Isabelle code}
\begin{lstlisting}
definition ucp_partial_derivative ::
  "('n::finite ucp_point \<Rightarrow> real) \<Rightarrow> 'n \<Rightarrow>
   'n ucp_point \<Rightarrow> real"
where
  "ucp_partial_derivative phi i x =
    frechet_derivative phi (at x) (axis i 1)"
\end{lstlisting}
\reviewlabel{English mathematical translation}
\noindent\textit{Mathematical role:} local coordinate derivative operator\par\smallskip
\texttt{ucp\_\allowbreak{}partial\_\allowbreak{}derivative\allowbreak\ phi\allowbreak\ i\allowbreak\ x} is the Fréchet derivative selected by Isabelle's \texttt{frechet\_\allowbreak{}derivative} operator for \texttt{phi} at \texttt{x}, applied to \texttt{axis\allowbreak\ i\allowbreak\ 1}. The definition itself supplies no differentiability premise.

\dossierentry{\texttt{\detokenize{ucp_test_function_on}}}
\reviewlabel{Isabelle code}
\begin{lstlisting}
definition ucp_test_function_on ::
  "'n::finite ucp_point set \<Rightarrow>
   ('n ucp_point \<Rightarrow> real) \<Rightarrow> bool"
where
  "ucp_test_function_on G phi \<longleftrightarrow>
    smooth_on UNIV phi \<and>
    compact (closure {x. phi x \<noteq> 0}) \<and>
    closure {x. phi x \<noteq> 0} \<subseteq> G"
\end{lstlisting}
\reviewlabel{English mathematical translation}
\noindent\textit{Mathematical role:} local test class \(\mathscr D(X)\)\par\smallskip
\texttt{ucp\_\allowbreak{}test\_\allowbreak{}function\_\allowbreak{}on\allowbreak\ G\allowbreak\ phi} holds exactly when \texttt{phi} is smooth on the entire vector space, \texttt{closure\allowbreak\ \{x\allowbreak\ |\allowbreak\ phi\allowbreak\ x\allowbreak\ !=\allowbreak\ 0\}} is compact, and that closure is contained in \texttt{G}.

\dossierentry{\texttt{\detokenize{ucp_locally_square_integrable_on}}}
\reviewlabel{Isabelle code}
\begin{lstlisting}
definition ucp_locally_square_integrable_on ::
  "'n::finite ucp_point set \<Rightarrow>
   ('n ucp_point \<Rightarrow> 'a::euclidean_space) \<Rightarrow> bool"
where
  "ucp_locally_square_integrable_on G f \<longleftrightarrow>
    (\<forall>K. compact K \<and> K \<subseteq> G \<longrightarrow>
      set_integrable lborel K (\<lambda>x. norm (f x) ^ 2))"
\end{lstlisting}
\reviewlabel{English mathematical translation}
\noindent\textit{Mathematical role:} local square integrability\par\smallskip
\texttt{ucp\_\allowbreak{}locally\_\allowbreak{}square\_\allowbreak{}integrable\_\allowbreak{}on\allowbreak\ G\allowbreak\ f} holds exactly when, for every compact set \texttt{K} contained in \texttt{G}, the function \texttt{x\allowbreak\ |-\allowbreak{}>\allowbreak\ ||f\allowbreak\ x||\textasciicircum{}2} is Lebesgue-integrable on \texttt{K}.

\dossierentry{\texttt{\detokenize{ucp_weak_gradient_on}}}
\reviewlabel{Isabelle code}
\begin{lstlisting}
definition ucp_weak_gradient_on ::
  "'n::finite ucp_point set \<Rightarrow>
   ('n ucp_point \<Rightarrow> real) \<Rightarrow>
   ('n ucp_point \<Rightarrow> 'n ucp_point) \<Rightarrow> bool"
where
  "ucp_weak_gradient_on G u Du \<longleftrightarrow>
    (\<forall>phi. ucp_test_function_on G phi \<longrightarrow>
      (\<forall>i.
        set_integrable lborel G
          (\<lambda>x. u x * ucp_partial_derivative phi i x) \<and>
        set_integrable lborel G (\<lambda>x. Du x $ i * phi x) \<and>
        set_lebesgue_integral lborel G
          (\<lambda>x. u x * ucp_partial_derivative phi i x) =
        - set_lebesgue_integral lborel G (\<lambda>x. Du x $ i * phi x)))"
\end{lstlisting}
\reviewlabel{English mathematical translation}
\noindent\textit{Mathematical role:} local weak-gradient identity\par\smallskip
\texttt{ucp\_\allowbreak{}weak\_\allowbreak{}gradient\_\allowbreak{}on\allowbreak\ G\allowbreak\ u\allowbreak\ Du} holds exactly when, for every \texttt{phi} satisfying \texttt{ucp\_\allowbreak{}test\_\allowbreak{}function\_\allowbreak{}on\allowbreak\ G\allowbreak\ phi} and every coordinate \texttt{i}, both \texttt{x\allowbreak\ |-\allowbreak{}>\allowbreak\ u\allowbreak\ x\allowbreak\ *\allowbreak\ ucp\_\allowbreak{}partial\_\allowbreak{}derivative\allowbreak\ phi\allowbreak\ i\allowbreak\ x} and \texttt{x\allowbreak\ |-\allowbreak{}>\allowbreak\ (Du\allowbreak\ x)\_\allowbreak{}i\allowbreak\ *\allowbreak\ phi\allowbreak\ x} are Lebesgue-integrable on \texttt{G}, and the integral of the first over \texttt{G} equals the negative of the integral of the second.

\dossierentry{\texttt{\detokenize{ucp_test_gradient}}}
\reviewlabel{Isabelle code}
\begin{lstlisting}
definition ucp_test_gradient ::
  "('n::finite ucp_point \<Rightarrow> real) \<Rightarrow>
   'n ucp_point \<Rightarrow> 'n ucp_point"
where
  "ucp_test_gradient phi x =
    (\<chi> i. ucp_partial_derivative phi i x)"
\end{lstlisting}
\reviewlabel{English mathematical translation}
\noindent\textit{Mathematical role:} local gradient vector\par\smallskip
\texttt{ucp\_\allowbreak{}test\_\allowbreak{}gradient\allowbreak\ phi\allowbreak\ x} is the vector whose coordinate \texttt{i} is \texttt{ucp\_\allowbreak{}partial\_\allowbreak{}derivative\allowbreak\ phi\allowbreak\ i\allowbreak\ x}.

\dossierentry{\texttt{\detokenize{ucp_weak_solution_on}}}
\reviewlabel{Isabelle code}
\begin{lstlisting}
definition ucp_weak_solution_on ::
  "'n::finite ucp_point set \<Rightarrow> 'n ucp_coefficient \<Rightarrow>
   ('n ucp_point \<Rightarrow> real) \<Rightarrow> bool"
where
  "ucp_weak_solution_on G a0 u \<longleftrightarrow>
    (\<exists>Du.
      ucp_locally_square_integrable_on G u \<and>
      ucp_locally_square_integrable_on G Du \<and>
      ucp_weak_gradient_on G u Du \<and>
      (\<forall>phi. ucp_test_function_on G phi \<longrightarrow>
        set_integrable lborel G
          (\<lambda>x. inner (a0 x *v Du x) (ucp_test_gradient phi x)) \<and>
        set_lebesgue_integral lborel G
          (\<lambda>x. inner (a0 x *v Du x) (ucp_test_gradient phi x)) = 0))"
\end{lstlisting}
\reviewlabel{English mathematical translation}
\noindent\textit{Mathematical role:} locally weak solution of \(\operatorname{div}(P\nabla u)=0\)\par\smallskip
\texttt{ucp\_\allowbreak{}weak\_\allowbreak{}solution\_\allowbreak{}on\allowbreak\ G\allowbreak\ a0\allowbreak\ u} holds exactly when there exists a vector-valued function \texttt{Du} such that \texttt{ucp\_\allowbreak{}locally\_\allowbreak{}square\_\allowbreak{}integrable\_\allowbreak{}on\allowbreak\ G\allowbreak\ u}, \texttt{ucp\_\allowbreak{}locally\_\allowbreak{}square\_\allowbreak{}integrable\_\allowbreak{}on\allowbreak\ G\allowbreak\ Du}, and \texttt{ucp\_\allowbreak{}weak\_\allowbreak{}gradient\_\allowbreak{}on\allowbreak\ G\allowbreak\ u\allowbreak\ Du} hold, and for every \texttt{phi} satisfying \texttt{ucp\_\allowbreak{}test\_\allowbreak{}function\_\allowbreak{}on\allowbreak\ G\allowbreak\ phi}, the function \texttt{x\allowbreak\ |-\allowbreak{}>\allowbreak\ inner\allowbreak\ (a0\allowbreak\ x\allowbreak\ *\allowbreak\ Du\allowbreak\ x)\allowbreak\ (ucp\_\allowbreak{}test\_\allowbreak{}gradient\allowbreak\ phi\allowbreak\ x)} is Lebesgue-integrable on \texttt{G} and has integral zero there.

\dossierentry{\texttt{\detokenize{hormander_regular_distribution_test_zero_on}}}
\reviewlabel{Isabelle code}
\begin{lstlisting}
definition hormander_regular_distribution_test_zero_on ::
  "'n::finite cgu_point set \<Rightarrow> 'n cgu_potential \<Rightarrow> bool"
where
  "hormander_regular_distribution_test_zero_on X u \<longleftrightarrow>
    set_integrable lborel X u \<and>
    (\<forall>phi. cgu_test_function_on X phi \<longrightarrow>
      set_integrable lborel X (\<lambda>x. u x * phi x) \<and>
      set_lebesgue_integral lborel X (\<lambda>x. u x * phi x) = 0)"
\end{lstlisting}
\reviewlabel{English mathematical translation}
\noindent\textit{Mathematical role:} integrable distributional-zero predicate\par\smallskip
\texttt{hormander\_\allowbreak{}regular\_\allowbreak{}distribution\_\allowbreak{}test\_\allowbreak{}zero\_\allowbreak{}on\allowbreak\ X\allowbreak\ u} holds exactly when \texttt{u} is Lebesgue-integrable on \texttt{X} and, for every \texttt{phi} satisfying \texttt{cgu\_\allowbreak{}test\_\allowbreak{}function\_\allowbreak{}on\allowbreak\ X\allowbreak\ phi}, the product \texttt{x\allowbreak\ |-\allowbreak{}>\allowbreak\ u\allowbreak\ x*phi\allowbreak\ x} is Lebesgue-integrable on \texttt{X} and has integral zero there.

\dossierentry{\texttt{\detokenize{evans_bounded_domain}}}
\reviewlabel{Isabelle code}
\begin{lstlisting}
definition evans_bounded_domain ::
  "'n::finite cgu_point set \<Rightarrow> bool"
where
  "evans_bounded_domain U \<longleftrightarrow> bounded U"
\end{lstlisting}
\reviewlabel{English mathematical translation}
\noindent\textit{Mathematical role:} boundedness of \(U\)\par\smallskip
\texttt{evans\_\allowbreak{}bounded\_\allowbreak{}domain\allowbreak\ U} holds exactly when \texttt{U} is bounded.

\dossierentry{\texttt{\detokenize{evans_elliptic_coefficient_on}}}
\reviewlabel{Isabelle code}
\begin{lstlisting}
definition evans_elliptic_coefficient_on ::
  "'n::finite cgu_point set \<Rightarrow> 'n mclean_coefficient \<Rightarrow> bool"
where
  "evans_elliptic_coefficient_on U a \<longleftrightarrow>
    a \<in> borel_measurable (restrict_space lborel U) \<and>
    (\<exists>B\<ge>0. \<forall>x\<in>U. norm (a x) \<le> B) \<and>
    (\<forall>x\<in>U. transpose (a x) = a x) \<and>
    (\<exists>theta>0. \<forall>x\<in>U. \<forall>xi.
      theta * norm xi ^ 2 \<le> inner (a x *v xi) xi)"
\end{lstlisting}
\reviewlabel{English mathematical translation}
\noindent\textit{Mathematical role:} measurable, bounded, symmetric, uniformly elliptic coefficient\par\smallskip
\texttt{evans\_\allowbreak{}elliptic\_\allowbreak{}coefficient\_\allowbreak{}on\allowbreak\ U\allowbreak\ a} holds exactly when:

\begin{itemize}
\item \texttt{a} is Borel-measurable with respect to Lebesgue measure restricted to \texttt{U};
\item there exists \texttt{B>=0} such that \texttt{||a\allowbreak\ x||<=B} for every \texttt{x\allowbreak\ in\allowbreak\ U};
\item \texttt{transpose(a\allowbreak\ x)=a\allowbreak\ x} for every \texttt{x\allowbreak\ in\allowbreak\ U}; and
\item there exists \texttt{theta>0} such that, for every \texttt{x\allowbreak\ in\allowbreak\ U} and every vector \texttt{xi}, \texttt{theta*||xi||\textasciicircum{}2\allowbreak\ <=\allowbreak\ inner\allowbreak\ (a\allowbreak\ x\allowbreak\ *\allowbreak\ xi)\allowbreak\ xi}.
\end{itemize}

\dossierentry{\texttt{\detokenize{evans_poincare_elliptic_green}}}
\reviewlabel{Isabelle code}
\begin{lstlisting}
locale evans_poincare_elliptic_green =
  fixes dimension_type :: "'n::finite itself"
  assumes evans_poincare_elliptic_green:
    "\<forall>U a. evans_poincare_elliptic_green_claim
      (U :: 'n cgu_point set) a"
\end{lstlisting}
\reviewlabel{English mathematical translation}
\noindent\textit{Mathematical role:} An atomic dimension-indexed universal assumption about \texttt{evans\_\allowbreak{}poincare\_\allowbreak{}elliptic\_\allowbreak{}green\_\allowbreak{}claim}.\par\smallskip
Fix \texttt{dimension\_\allowbreak{}type} of type \texttt{'n::finite\allowbreak\ itself}, selecting the finite-dimensional carrier type \texttt{'n} without specifying its dimension. Assume, under the name \texttt{evans\_\allowbreak{}poincare\_\allowbreak{}elliptic\_\allowbreak{}green}, that for every \texttt{U} and every \texttt{a},

\texttt{evans\_\allowbreak{}poincare\_\allowbreak{}elliptic\_\allowbreak{}green\_\allowbreak{}claim\allowbreak\ U\allowbreak\ a}

holds, where \texttt{U} is explicitly a set of objects of type \texttt{'n\allowbreak\ cgu\_\allowbreak{}point}. The quantifiers range over all well-typed \texttt{U} and \texttt{a}; no further restriction on either variable is stated.

\dossierentry{\texttt{\detokenize{mclean_smooth_test_h1_zero}}}
\reviewlabel{Isabelle code}
\begin{lstlisting}
locale mclean_smooth_test_h1_zero =
  fixes dimension_type :: "'n::finite itself"
  assumes mclean_smooth_test_h1_zero:
    "mclean_smooth_test_h1_zero_claim dimension_type"
\end{lstlisting}
\reviewlabel{English mathematical translation}
\noindent\textit{Mathematical role:} An atomic dimension-indexed assumption asserting \texttt{mclean\_\allowbreak{}smooth\_\allowbreak{}test\_\allowbreak{}h1\_\allowbreak{}zero\_\allowbreak{}claim}.\par\smallskip
Fix \texttt{dimension\_\allowbreak{}type} of type \texttt{'n::finite\allowbreak\ itself}, selecting the finite-dimensional carrier type \texttt{'n} without specifying its dimension. Assume, under the name \texttt{mclean\_\allowbreak{}smooth\_\allowbreak{}test\_\allowbreak{}h1\_\allowbreak{}zero}, the proposition

\texttt{mclean\_\allowbreak{}smooth\_\allowbreak{}test\_\allowbreak{}h1\_\allowbreak{}zero\_\allowbreak{}claim\allowbreak\ dimension\_\allowbreak{}type}.

\dossierentry{\texttt{\detokenize{hormander_analytic_elliptic_ucp}}}
\reviewlabel{Isabelle code}
\begin{lstlisting}
locale hormander_analytic_elliptic_ucp =
  fixes dimension_type :: "'n::finite itself"
  assumes hormander_analytic_elliptic_ucp:
    "\<forall>X P u.
      hormander_analytic_elliptic_ucp_claim
        (X :: 'n ucp_point set) P u"
\end{lstlisting}
\reviewlabel{English mathematical translation}
\noindent\textit{Mathematical role:} An atomic dimension-indexed universal assumption about \texttt{hormander\_\allowbreak{}analytic\_\allowbreak{}elliptic\_\allowbreak{}ucp\_\allowbreak{}claim}.\par\smallskip
Fix \texttt{dimension\_\allowbreak{}type} of type \texttt{'n::finite\allowbreak\ itself}, selecting the finite-dimensional carrier type \texttt{'n} without specifying its dimension. Assume, under the name \texttt{hormander\_\allowbreak{}analytic\_\allowbreak{}elliptic\_\allowbreak{}ucp}, that for every \texttt{X}, every \texttt{P}, and every \texttt{u},

\texttt{hormander\_\allowbreak{}analytic\_\allowbreak{}elliptic\_\allowbreak{}ucp\_\allowbreak{}claim\allowbreak\ X\allowbreak\ P\allowbreak\ u}

holds, where \texttt{X} is explicitly a set of objects of type \texttt{'n\allowbreak\ ucp\_\allowbreak{}point}. The argument order is \texttt{X}, then \texttt{P}, then \texttt{u}. No further restriction on these universally quantified variables is stated.

\dossierentry{\texttt{\detokenize{hormander_weak_gradient_locality}}}
\reviewlabel{Isabelle code}
\begin{lstlisting}
locale hormander_weak_gradient_locality =
  fixes dimension_type :: "'n::finite itself"
  assumes hormander_weak_gradient_locality:
    "\<forall>X u Du.
      hormander_weak_gradient_locality_claim
        (X :: 'n ucp_point set) u Du"
\end{lstlisting}
\reviewlabel{English mathematical translation}
\noindent\textit{Mathematical role:} An atomic dimension-indexed universal locality assumption about \texttt{hormander\_\allowbreak{}weak\_\allowbreak{}gradient\_\allowbreak{}locality\_\allowbreak{}claim}.\par\smallskip
Fix \texttt{dimension\_\allowbreak{}type} of type \texttt{'n::finite\allowbreak\ itself}, selecting the finite-dimensional carrier type \texttt{'n} without specifying its dimension. Assume, under the name \texttt{hormander\_\allowbreak{}weak\_\allowbreak{}gradient\_\allowbreak{}locality}, that for every \texttt{X}, every \texttt{u}, and every \texttt{Du},

\texttt{hormander\_\allowbreak{}weak\_\allowbreak{}gradient\_\allowbreak{}locality\_\allowbreak{}claim\allowbreak\ X\allowbreak\ u\allowbreak\ Du}

holds, where \texttt{X} is explicitly a set of objects of type \texttt{'n\allowbreak\ ucp\_\allowbreak{}point}. The argument order is \texttt{X}, then \texttt{u}, then \texttt{Du}. No further restriction on these universally quantified variables is stated.

\dossierentry{\texttt{\detokenize{mclean_h1_zero_extension_localization}}}
\reviewlabel{Isabelle code}
\begin{lstlisting}
locale mclean_h1_zero_extension_localization =
  fixes dimension_type :: "'n::finite itself"
  assumes mclean_h1_zero_extension_localization:
    "mclean_h1_zero_extension_localization_claim dimension_type"
\end{lstlisting}
\reviewlabel{English mathematical translation}
\noindent\textit{Mathematical role:} An atomic dimension-indexed assumption asserting \texttt{mclean\_\allowbreak{}h1\_\allowbreak{}zero\_\allowbreak{}extension\_\allowbreak{}localization\_\allowbreak{}claim}.\par\smallskip
Fix \texttt{dimension\_\allowbreak{}type} of type \texttt{'n::finite\allowbreak\ itself}, selecting the finite-dimensional carrier type \texttt{'n} without specifying its dimension. Assume, under the name \texttt{mclean\_\allowbreak{}h1\_\allowbreak{}zero\_\allowbreak{}extension\_\allowbreak{}localization}, the proposition

\texttt{mclean\_\allowbreak{}h1\_\allowbreak{}zero\_\allowbreak{}extension\_\allowbreak{}localization\_\allowbreak{}claim\allowbreak\ dimension\_\allowbreak{}type}.

\dossierentry{\texttt{\detokenize{mclean_trace_dirichlet_green_v3}}}
\reviewlabel{Isabelle code}
\begin{lstlisting}
locale mclean_trace_dirichlet_green_v3 =
  fixes dimension_type :: "'n::finite itself"
  assumes mclean_trace_dirichlet_green_v3:
    "\<forall>U a. mclean_trace_dirichlet_green_claim_v3
      (U :: 'n cgu_point set) a"
\end{lstlisting}
\reviewlabel{English mathematical translation}
\noindent\textit{Mathematical role:} An atomic dimension-indexed universal trace assumption about \texttt{mclean\_\allowbreak{}trace\_\allowbreak{}dirichlet\_\allowbreak{}green\_\allowbreak{}claim\_\allowbreak{}v3}.\par\smallskip
Fix \texttt{dimension\_\allowbreak{}type} of type \texttt{'n::finite\allowbreak\ itself}, selecting the finite-dimensional carrier type \texttt{'n} without specifying its dimension. Assume, under the name \texttt{mclean\_\allowbreak{}trace\_\allowbreak{}dirichlet\_\allowbreak{}green\_\allowbreak{}v3}, that for every \texttt{U} and every \texttt{a},

\texttt{mclean\_\allowbreak{}trace\_\allowbreak{}dirichlet\_\allowbreak{}green\_\allowbreak{}claim\_\allowbreak{}v3\allowbreak\ U\allowbreak\ a}

holds, where \texttt{U} is explicitly a set of objects of type \texttt{'n\allowbreak\ cgu\_\allowbreak{}point}. The argument order is \texttt{U}, then \texttt{a}. No further restriction on these universally quantified variables is stated.

\dossierentry{\texttt{\detokenize{hormander_regular_distribution_injectivity}}}
\reviewlabel{Isabelle code}
\begin{lstlisting}
locale hormander_regular_distribution_injectivity =
  fixes dimension_type :: "'n::finite itself"
  assumes hormander_regular_distribution_injectivity:
    "hormander_regular_distribution_injectivity_claim dimension_type"
\end{lstlisting}
\reviewlabel{English mathematical translation}
\noindent\textit{Mathematical role:} An atomic dimension-indexed injectivity assumption asserting \texttt{hormander\_\allowbreak{}regular\_\allowbreak{}distribution\_\allowbreak{}injectivity\_\allowbreak{}claim}.\par\smallskip
Fix \texttt{dimension\_\allowbreak{}type} of type \texttt{'n::finite\allowbreak\ itself}, selecting the finite-dimensional carrier type \texttt{'n} without specifying its dimension. Assume, under the name \texttt{hormander\_\allowbreak{}regular\_\allowbreak{}distribution\_\allowbreak{}injectivity}, the proposition

\texttt{hormander\_\allowbreak{}regular\_\allowbreak{}distribution\_\allowbreak{}injectivity\_\allowbreak{}claim\allowbreak\ dimension\_\allowbreak{}type}.

\dossierentry{\texttt{\detokenize{mclean_graph_supported_h1_zero}}}
\reviewlabel{Isabelle code}
\begin{lstlisting}
locale mclean_graph_supported_h1_zero =
  fixes dimension_type :: "'n::finite itself"
  assumes mclean_graph_supported_h1_zero:
    "mclean_graph_supported_h1_zero_claim dimension_type"
\end{lstlisting}
\reviewlabel{English mathematical translation}
\noindent\textit{Mathematical role:} An atomic dimension-indexed graph-support assumption asserting \texttt{mclean\_\allowbreak{}graph\_\allowbreak{}supported\_\allowbreak{}h1\_\allowbreak{}zero\_\allowbreak{}claim}.\par\smallskip
Fix \texttt{dimension\_\allowbreak{}type} of type \texttt{'n::finite\allowbreak\ itself}, selecting the finite-dimensional carrier type \texttt{'n} without specifying its dimension. Assume, under the name \texttt{mclean\_\allowbreak{}graph\_\allowbreak{}supported\_\allowbreak{}h1\_\allowbreak{}zero}, the proposition

\texttt{mclean\_\allowbreak{}graph\_\allowbreak{}supported\_\allowbreak{}h1\_\allowbreak{}zero\_\allowbreak{}claim\allowbreak\ dimension\_\allowbreak{}type}.

\dossierentry{\texttt{\detokenize{mclean_h1_zero_bidual_representation_v2}}}
\reviewlabel{Isabelle code}
\begin{lstlisting}
locale mclean_h1_zero_bidual_representation_v2 =
  fixes dimension_type :: "'n::finite itself"
  assumes mclean_h1_zero_bidual_representation_v2:
    "mclean_h1_zero_bidual_representation_claim_v2 dimension_type"
\end{lstlisting}
\reviewlabel{English mathematical translation}
\noindent\textit{Mathematical role:} An atomic dimension-indexed bidual assumption asserting \texttt{mclean\_\allowbreak{}h1\_\allowbreak{}zero\_\allowbreak{}bidual\_\allowbreak{}representation\_\allowbreak{}claim\_\allowbreak{}v2}.\par\smallskip
Fix \texttt{dimension\_\allowbreak{}type} of type \texttt{'n::finite\allowbreak\ itself}, selecting the finite-dimensional carrier type \texttt{'n} without specifying its dimension. Assume, under the name \texttt{mclean\_\allowbreak{}h1\_\allowbreak{}zero\_\allowbreak{}bidual\_\allowbreak{}representation\_\allowbreak{}v2}, the proposition

\texttt{mclean\_\allowbreak{}h1\_\allowbreak{}zero\_\allowbreak{}bidual\_\allowbreak{}representation\_\allowbreak{}claim\_\allowbreak{}v2\allowbreak\ dimension\_\allowbreak{}type}.

\dossierentry{\texttt{\detokenize{mclean_smooth_generator_ambient_lift_v2}}}
\reviewlabel{Isabelle code}
\begin{lstlisting}
locale mclean_smooth_generator_ambient_lift_v2 =
  assumes mclean_smooth_generator_ambient_lift_v2:
    "mclean_smooth_generator_ambient_lift_claim_v2
      TYPE('n::finite) TYPE('c::finite) TYPE('h::finite)"
\end{lstlisting}
\reviewlabel{English mathematical translation}
\noindent\textit{Mathematical role:} An atomic assumption involving three independently selected finite-dimensional carrier types.\par\smallskip
Assume, under the name \texttt{mclean\_\allowbreak{}smooth\_\allowbreak{}generator\_\allowbreak{}ambient\_\allowbreak{}lift\_\allowbreak{}v2}, the proposition

\texttt{mclean\_\allowbreak{}smooth\_\allowbreak{}generator\_\allowbreak{}ambient\_\allowbreak{}lift\_\allowbreak{}claim\_\allowbreak{}v2\allowbreak\ TYPE('n::finite)\allowbreak\ TYPE('c::finite)\allowbreak\ TYPE('h::finite)}.

Here the three arguments, in order, select the finite-dimensional carrier types \texttt{'n}, \texttt{'c}, and \texttt{'h}. No dimension is specified for any of them, and no equality or relationship among the three carrier types is assumed.

\dossierentry{\texttt{\detokenize{mclean_hyperplane_offset_source}}}
\reviewlabel{Isabelle code}
\begin{lstlisting}
locale mclean_hyperplane_offset_source =
  assumes mclean_hyperplane_offset_source:
    "mclean_hyperplane_offset_source_claim
      TYPE('m::finite) TYPE('n::finite)"
\end{lstlisting}
\reviewlabel{English mathematical translation}
\noindent\textit{Mathematical role:} An atomic assumption involving two independently selected finite-dimensional carrier types.\par\smallskip
Assume, under the name \texttt{mclean\_\allowbreak{}hyperplane\_\allowbreak{}offset\_\allowbreak{}source}, the proposition

\texttt{mclean\_\allowbreak{}hyperplane\_\allowbreak{}offset\_\allowbreak{}source\_\allowbreak{}claim\allowbreak\ TYPE('m::finite)\allowbreak\ TYPE('n::finite)}.

Here the two arguments, in order, select the finite-dimensional carrier types \texttt{'m} and \texttt{'n}. No dimension is specified for either, and no equality or relationship between the two carrier types is assumed.

\dossierentry{\texttt{\detokenize{kang_yun_ltu_local_boundary_metric_v6}}}
\reviewlabel{Isabelle code}
\begin{lstlisting}
locale kang_yun_ltu_local_boundary_metric_v6 =
  fixes dimension_type :: "'n::finite itself"
  assumes kang_yun_ltu_local_boundary_metric_v6:
    "\<forall>M P1 P2 U gamma normal offset g1 g2 a1 a2 rho1 rho2.
      kang_yun_ltu_local_boundary_metric_claim_v6
        (M :: ('n, 'c::finite, 'h::finite) cgu_model)
        P1 P2 U gamma normal offset g1 g2 a1 a2 rho1 rho2"
\end{lstlisting}
\reviewlabel{English mathematical translation}
\noindent\textit{Mathematical role:} An atomic dimension-indexed universal boundary-metric assumption about \texttt{kang\_\allowbreak{}yun\_\allowbreak{}ltu\_\allowbreak{}local\_\allowbreak{}boundary\_\allowbreak{}metric\_\allowbreak{}claim\_\allowbreak{}v6}.\par\smallskip
Fix \texttt{dimension\_\allowbreak{}type} of type \texttt{'n::finite\allowbreak\ itself}, selecting the finite-dimensional carrier type \texttt{'n} without specifying its dimension. Assume, under the name \texttt{kang\_\allowbreak{}yun\_\allowbreak{}ltu\_\allowbreak{}local\_\allowbreak{}boundary\_\allowbreak{}metric\_\allowbreak{}v6}, that for every

\texttt{M}, \texttt{P1}, \texttt{P2}, \texttt{U}, \texttt{gamma}, \texttt{normal}, \texttt{offset}, \texttt{g1}, \texttt{g2}, \texttt{a1}, \texttt{a2}, \texttt{rho1}, and \texttt{rho2},

the proposition

\texttt{kang\_\allowbreak{}yun\_\allowbreak{}ltu\_\allowbreak{}local\_\allowbreak{}boundary\_\allowbreak{}metric\_\allowbreak{}claim\_\allowbreak{}v6\allowbreak\ M\allowbreak\ P1\allowbreak\ P2\allowbreak\ U\allowbreak\ gamma\allowbreak\ normal\allowbreak\ offset\allowbreak\ g1\allowbreak\ g2\allowbreak\ a1\allowbreak\ a2\allowbreak\ rho1\allowbreak\ rho2}

holds. The variable \texttt{M} is explicitly of type \texttt{('n,\allowbreak\ 'c::finite,\allowbreak\ 'h::finite)\allowbreak\ cgu\_\allowbreak{}model}: its first type parameter is the finite carrier type \texttt{'n} selected by \texttt{dimension\_\allowbreak{}type}, while \texttt{'c} and \texttt{'h} are also finite carrier types; no dimension is specified for any of them. The arguments occur in exactly the displayed order. All thirteen variables are universally quantified, and no further restriction on them is stated.

\dossierentry{\texttt{\detokenize{piecewise_polynomial_ucp_literature}}}
\reviewlabel{Isabelle code}
\begin{lstlisting}
locale piecewise_polynomial_ucp_literature =
  analytic: hormander_analytic_elliptic_ucp dimension_type +
  locality: hormander_weak_gradient_locality dimension_type
  for dimension_type :: "'n::finite itself"
\end{lstlisting}
\reviewlabel{English mathematical translation}
\noindent\textit{Mathematical role:} A dimension-indexed package combining named analytic and locality packages.\par\smallskip
Fix \texttt{dimension\_\allowbreak{}type} of type \texttt{'n::finite\allowbreak\ itself}, selecting the finite-dimensional carrier type \texttt{'n} without specifying its dimension. The locale \texttt{piecewise\_\allowbreak{}polynomial\_\allowbreak{}ucp\_\allowbreak{}literature} assumes exactly these two parent packages:

\begin{enumerate}[start=1]
\item \texttt{hormander\_\allowbreak{}analytic\_\allowbreak{}elliptic\_\allowbreak{}ucp}, instantiated with \texttt{dimension\_\allowbreak{}type} and included as the named component \texttt{analytic}.
\item \texttt{hormander\_\allowbreak{}weak\_\allowbreak{}gradient\_\allowbreak{}locality}, instantiated with \texttt{dimension\_\allowbreak{}type} and included as the named component \texttt{locality}.
\end{enumerate}

No additional assumption is stated.

\dossierentry{\texttt{\detokenize{cgu_volume_annihilator_literature_base}}}
\reviewlabel{Isabelle code}
\begin{lstlisting}
locale cgu_volume_annihilator_literature_base =
  smooth: mclean_smooth_test_h1_zero dimension_type +
  injectivity: hormander_regular_distribution_injectivity dimension_type
  for dimension_type :: "'n::finite itself"
\end{lstlisting}
\reviewlabel{English mathematical translation}
\noindent\textit{Mathematical role:} A dimension-indexed volume package combining smoothness and injectivity packages.\par\smallskip
Fix \texttt{dimension\_\allowbreak{}type} of type \texttt{'n::finite\allowbreak\ itself}, selecting the finite-dimensional carrier type \texttt{'n} without specifying its dimension. The locale \texttt{cgu\_\allowbreak{}volume\_\allowbreak{}annihilator\_\allowbreak{}literature\_\allowbreak{}base} assumes exactly these two parent packages:

\begin{enumerate}[start=1]
\item \texttt{mclean\_\allowbreak{}smooth\_\allowbreak{}test\_\allowbreak{}h1\_\allowbreak{}zero}, instantiated with \texttt{dimension\_\allowbreak{}type} and included as the named component \texttt{smooth}.
\item \texttt{hormander\_\allowbreak{}regular\_\allowbreak{}distribution\_\allowbreak{}injectivity}, instantiated with \texttt{dimension\_\allowbreak{}type} and included as the named component \texttt{injectivity}.
\end{enumerate}

No additional assumption is stated.

\dossierentry{\texttt{\detokenize{cgu_volume_gradient_annihilator_literature_base}}}
\reviewlabel{Isabelle code}
\begin{lstlisting}
locale cgu_volume_gradient_annihilator_literature_base =
  volume: cgu_volume_annihilator_literature_base dimension_type +
  locality: hormander_weak_gradient_locality dimension_type
  for dimension_type :: "'n::finite itself"
\end{lstlisting}
\reviewlabel{English mathematical translation}
\noindent\textit{Mathematical role:} A dimension-indexed package combining named volume and locality packages.\par\smallskip
Fix \texttt{dimension\_\allowbreak{}type} of type \texttt{'n::finite\allowbreak\ itself}, selecting the finite-dimensional carrier type \texttt{'n} without specifying its dimension. The locale \texttt{cgu\_\allowbreak{}volume\_\allowbreak{}gradient\_\allowbreak{}annihilator\_\allowbreak{}literature\_\allowbreak{}base} assumes exactly these two parent packages:

\begin{enumerate}[start=1]
\item \texttt{cgu\_\allowbreak{}volume\_\allowbreak{}annihilator\_\allowbreak{}literature\_\allowbreak{}base}, instantiated with \texttt{dimension\_\allowbreak{}type} and included as the named component \texttt{volume}.
\item \texttt{hormander\_\allowbreak{}weak\_\allowbreak{}gradient\_\allowbreak{}locality}, instantiated with \texttt{dimension\_\allowbreak{}type} and included as the named component \texttt{locality}.
\end{enumerate}

No additional assumption is stated.

\dossierentry{\texttt{\detokenize{cgu_graph_interface_annihilator_literature_base}}}
\reviewlabel{Isabelle code}
\begin{lstlisting}
locale cgu_graph_interface_annihilator_literature_base =
  volume_gradient: cgu_volume_gradient_annihilator_literature_base dimension_type +
  localization: mclean_h1_zero_extension_localization dimension_type +
  graph_support: mclean_graph_supported_h1_zero dimension_type
  for dimension_type :: "'n::finite itself"
\end{lstlisting}
\reviewlabel{English mathematical translation}
\noindent\textit{Mathematical role:} A dimension-indexed package combining volume-gradient, localization, and graph-support packages.\par\smallskip
Fix \texttt{dimension\_\allowbreak{}type} of type \texttt{'n::finite\allowbreak\ itself}, selecting the finite-dimensional carrier type \texttt{'n} without specifying its dimension. The locale \texttt{cgu\_\allowbreak{}graph\_\allowbreak{}interface\_\allowbreak{}annihilator\_\allowbreak{}literature\_\allowbreak{}base} assumes exactly these three parent packages:

\begin{enumerate}[start=1]
\item \texttt{cgu\_\allowbreak{}volume\_\allowbreak{}gradient\_\allowbreak{}annihilator\_\allowbreak{}literature\_\allowbreak{}base}, instantiated with \texttt{dimension\_\allowbreak{}type} and included as the named component \texttt{volume\_\allowbreak{}gradient}.
\item \texttt{mclean\_\allowbreak{}h1\_\allowbreak{}zero\_\allowbreak{}extension\_\allowbreak{}localization}, instantiated with \texttt{dimension\_\allowbreak{}type} and included as the named component \texttt{localization}.
\item \texttt{mclean\_\allowbreak{}graph\_\allowbreak{}supported\_\allowbreak{}h1\_\allowbreak{}zero}, instantiated with \texttt{dimension\_\allowbreak{}type} and included as the named component \texttt{graph\_\allowbreak{}support}.
\end{enumerate}

No additional assumption is stated.

\dossierentry{\texttt{\detokenize{cgu_graph_interface_source_approximation_literature_base}}}
\reviewlabel{Isabelle code}
\begin{lstlisting}
locale cgu_graph_interface_source_approximation_literature_base =
  annihilator: cgu_graph_interface_annihilator_literature_base dimension_type +
  bidual: mclean_h1_zero_bidual_representation_v2 dimension_type
  for dimension_type :: "'n::finite itself"
\end{lstlisting}
\reviewlabel{English mathematical translation}
\noindent\textit{Mathematical role:} A dimension-indexed graph-source package combining named annihilator and bidual packages.\par\smallskip
Fix \texttt{dimension\_\allowbreak{}type} of type \texttt{'n::finite\allowbreak\ itself}, selecting the finite-dimensional carrier type \texttt{'n} without specifying its dimension. The locale \texttt{cgu\_\allowbreak{}graph\_\allowbreak{}interface\_\allowbreak{}source\_\allowbreak{}approximation\_\allowbreak{}literature\_\allowbreak{}base} assumes exactly these two parent packages:

\begin{enumerate}[start=1]
\item \texttt{cgu\_\allowbreak{}graph\_\allowbreak{}interface\_\allowbreak{}annihilator\_\allowbreak{}literature\_\allowbreak{}base}, instantiated with \texttt{dimension\_\allowbreak{}type} and included as the named component \texttt{annihilator}.
\item \texttt{mclean\_\allowbreak{}h1\_\allowbreak{}zero\_\allowbreak{}bidual\_\allowbreak{}representation\_\allowbreak{}v2}, instantiated with \texttt{dimension\_\allowbreak{}type} and included as the named component \texttt{bidual}.
\end{enumerate}

No additional assumption is stated.

\dossierentry{\texttt{\detokenize{cgu_runge_literature_base}}}
\reviewlabel{Isabelle code}
\begin{lstlisting}
locale cgu_runge_literature_base =
  evans: evans_poincare_elliptic_green dimension_type +
  smooth: mclean_smooth_test_h1_zero dimension_type +
  ucp: piecewise_polynomial_ucp_literature dimension_type +
  localization: mclean_h1_zero_extension_localization dimension_type
  for dimension_type :: "'n::finite itself"
\end{lstlisting}
\reviewlabel{English mathematical translation}
\noindent\textit{Mathematical role:} A dimension-indexed package assembling four named assumption packages.\par\smallskip
Fix \texttt{dimension\_\allowbreak{}type} of type \texttt{'n::finite\allowbreak\ itself}, selecting the finite-dimensional carrier type \texttt{'n} without specifying its dimension. The locale \texttt{cgu\_\allowbreak{}runge\_\allowbreak{}literature\_\allowbreak{}base} assumes exactly these four parent packages:

\begin{enumerate}[start=1]
\item \texttt{evans\_\allowbreak{}poincare\_\allowbreak{}elliptic\_\allowbreak{}green}, instantiated with \texttt{dimension\_\allowbreak{}type} and included as the named component \texttt{evans}.
\item \texttt{mclean\_\allowbreak{}smooth\_\allowbreak{}test\_\allowbreak{}h1\_\allowbreak{}zero}, instantiated with \texttt{dimension\_\allowbreak{}type} and included as the named component \texttt{smooth}.
\item \texttt{piecewise\_\allowbreak{}polynomial\_\allowbreak{}ucp\_\allowbreak{}literature}, instantiated with \texttt{dimension\_\allowbreak{}type} and included as the named component \texttt{ucp}.
\item \texttt{mclean\_\allowbreak{}h1\_\allowbreak{}zero\_\allowbreak{}extension\_\allowbreak{}localization}, instantiated with \texttt{dimension\_\allowbreak{}type} and included as the named component \texttt{localization}.
\end{enumerate}

No additional assumption is stated.

\dossierentry{\texttt{\detokenize{cgu_graph_dirichlet_literature_base_v3}}}
\reviewlabel{Isabelle code}
\begin{lstlisting}
locale cgu_graph_dirichlet_literature_base_v3 =
  cgu_runge_literature_base dimension_type +
  trace: mclean_trace_dirichlet_green_v3 dimension_type
  for dimension_type :: "'n::finite itself"
\end{lstlisting}
\reviewlabel{English mathematical translation}
\noindent\textit{Mathematical role:} A dimension-indexed package combining a base package with a named trace package.\par\smallskip
Fix \texttt{dimension\_\allowbreak{}type} of type \texttt{'n::finite\allowbreak\ itself}, selecting the finite-dimensional carrier type \texttt{'n} without specifying its dimension. The locale \texttt{cgu\_\allowbreak{}graph\_\allowbreak{}dirichlet\_\allowbreak{}literature\_\allowbreak{}base\_\allowbreak{}v3} assumes exactly these two parent packages:

\begin{enumerate}[start=1]
\item \texttt{cgu\_\allowbreak{}runge\_\allowbreak{}literature\_\allowbreak{}base}, instantiated with \texttt{dimension\_\allowbreak{}type}.
\item \texttt{mclean\_\allowbreak{}trace\_\allowbreak{}dirichlet\_\allowbreak{}green\_\allowbreak{}v3}, instantiated with \texttt{dimension\_\allowbreak{}type} and included as the named component \texttt{trace}.
\end{enumerate}

No additional assumption is stated.

\dossierentry{\texttt{\detokenize{cgu_graph_interface_conormal_context}}}
\reviewlabel{Isabelle code}
\begin{lstlisting}
locale cgu_graph_interface_conormal_context =
  cgu_graph_dirichlet_literature_base_v3 dimension_type +
  graph_source: cgu_graph_interface_source_approximation_literature_base dimension_type
  for dimension_type :: "'n::finite itself"
\end{lstlisting}
\reviewlabel{English mathematical translation}
\noindent\textit{Mathematical role:} A dimension-indexed package combining a base package with a named graph-source package.\par\smallskip
Fix \texttt{dimension\_\allowbreak{}type} of type \texttt{'n::finite\allowbreak\ itself}, selecting the finite-dimensional carrier type \texttt{'n} without specifying its dimension. The locale \texttt{cgu\_\allowbreak{}graph\_\allowbreak{}interface\_\allowbreak{}conormal\_\allowbreak{}context} assumes exactly these two parent packages:

\begin{enumerate}[start=1]
\item \texttt{cgu\_\allowbreak{}graph\_\allowbreak{}dirichlet\_\allowbreak{}literature\_\allowbreak{}base\_\allowbreak{}v3}, instantiated with \texttt{dimension\_\allowbreak{}type}.
\item \texttt{cgu\_\allowbreak{}graph\_\allowbreak{}interface\_\allowbreak{}source\_\allowbreak{}approximation\_\allowbreak{}literature\_\allowbreak{}base}, instantiated with \texttt{dimension\_\allowbreak{}type} and included as the named component \texttt{graph\_\allowbreak{}source}.
\end{enumerate}

No additional assumption is stated.

\dossierentry{\texttt{\detokenize{cgu_selected_face_dn_context}}}
\reviewlabel{Isabelle code}
\begin{lstlisting}
locale cgu_selected_face_dn_context =
  cgu_graph_interface_conormal_context dimension_type +
  mclean_smooth_generator_ambient_lift_v2
  for dimension_type :: "'n::finite itself"
\end{lstlisting}
\reviewlabel{English mathematical translation}
\noindent\textit{Mathematical role:} A dimension-indexed package combining a structural package with an additional unparameterized package.\par\smallskip
Fix \texttt{dimension\_\allowbreak{}type} of type \texttt{'n::finite\allowbreak\ itself}; that is, \texttt{dimension\_\allowbreak{}type} selects the finite-dimensional carrier type \texttt{'n}, without specifying its dimension. The locale \texttt{cgu\_\allowbreak{}selected\_\allowbreak{}face\_\allowbreak{}dn\_\allowbreak{}context} assumes exactly these two parent packages:

\begin{enumerate}[start=1]
\item \texttt{cgu\_\allowbreak{}graph\_\allowbreak{}interface\_\allowbreak{}conormal\_\allowbreak{}context}, instantiated with \texttt{dimension\_\allowbreak{}type}.
\item \texttt{mclean\_\allowbreak{}smooth\_\allowbreak{}generator\_\allowbreak{}ambient\_\allowbreak{}lift\_\allowbreak{}v2}.
\end{enumerate}

No additional assumption is stated.

\dossierentry{\texttt{\detokenize{cgu_terminal_recovery_context_v7}}}
\reviewlabel{Isabelle code}
\begin{lstlisting}
locale cgu_terminal_recovery_context_v7 =
  cgu_selected_face_dn_context "TYPE('n::finite)" +
  offset: mclean_hyperplane_offset_source +
  boundary_metric: kang_yun_ltu_local_boundary_metric_v6 "TYPE('n)"
\end{lstlisting}
\reviewlabel{English mathematical translation}
\noindent\textit{Mathematical role:} A top-level package combining three assumption packages.\par\smallskip
The locale \texttt{cgu\_\allowbreak{}terminal\_\allowbreak{}recovery\_\allowbreak{}context\_\allowbreak{}v7} assumes exactly the following three parent packages:

\begin{enumerate}[start=1]
\item \texttt{cgu\_\allowbreak{}selected\_\allowbreak{}face\_\allowbreak{}dn\_\allowbreak{}context}, instantiated at \texttt{TYPE('n::finite)}. Here \texttt{TYPE('n::finite)} selects the finite-dimensional carrier type \texttt{'n}; it does not specify or imply any particular dimension.
\item \texttt{mclean\_\allowbreak{}hyperplane\_\allowbreak{}offset\_\allowbreak{}source}, included as the named component \texttt{offset}.
\item \texttt{kang\_\allowbreak{}yun\_\allowbreak{}ltu\_\allowbreak{}local\_\allowbreak{}boundary\_\allowbreak{}metric\_\allowbreak{}v6}, instantiated at \texttt{TYPE('n)} and included as the named component \texttt{boundary\_\allowbreak{}metric}; this selects the same finite-dimensional carrier type \texttt{'n}, without specifying its dimension.
\end{enumerate}

No additional parameters or assumptions are stated in this locale.

\clearpage
\dossiersection{Statements}
\dossiersubhead{Manuscript conclusions}

\dossierentry{MAIN-001: \texttt{\detokenize{cgu_global_uniqueness_claim_v7}}}
\reviewlabel{English original and manuscript locator}
Theorem \texttt{thm:main} of \cite{Carstea2026Polynomial} states:

Let \(n\ge3\), let \(\Omega\subset\R^n\) be a bounded Lipschitz domain, and let
\(\mathcal D=\{D_\alpha\}_{\alpha\in A}\) be a known finite Lipschitz
subdivision of \(\Omega\) that is internally flat up to a skeleton.  Let
\(N_\alpha\ge0\).  Suppose that
\(\sigma^{(1)},\sigma^{(2)}\in\mathcal P(\mathcal D,\{N_\alpha\})\), and assume
that, at every recovery step and for \(j=1,2\), the polynomial extension of
\(\sigma^{(j)}|_{D_{\alpha_r}}\) is positive definite on the triple-intersection
sets \(W_{r,ijk}\) appearing in Definition \textup{[source reference: \texttt{def:stripping-order}]}.  We also
assume that, for \(j=1,2\), the polynomial extensions in the recovered cells adjacent to the
auxiliary exterior half-balls remain elliptic in those half-balls whenever the
half-balls are used.

Let \(\Sigma\subset\partial\Omega\) be a nonempty relatively open measured set.
Assume that the subdivision admits an admissible flat-face recovery order
relative to \(\Sigma\), in the sense of Definition \textup{[source reference: \texttt{def:stripping-order}]}.
If
\[
  \DN^\Omega_{\sigma^{(1)},\Sigma}
  =\DN^\Omega_{\sigma^{(2)},\Sigma},
\]
then
\[
  \sigma^{(1)}=\sigma^{(2)}\quad\text{in }\Omega.
\]
Equivalently, the polynomial matrices defining the two conductivities agree on
every cell of the subdivision.
\reviewlabel{Isabelle code}
\begin{lstlisting}
definition cgu_global_uniqueness_claim_v7 ::
  "('n::finite, 'c::finite, 'h::finite) cgu_model \<Rightarrow>
    ('c \<Rightarrow> 'n cgu_matrix_polynomial) \<Rightarrow>
    ('c \<Rightarrow> 'n cgu_matrix_polynomial) \<Rightarrow> bool"
where
  "cgu_global_uniqueness_claim_v7 M P1 P2 \<longleftrightarrow>
    (3 \<le> CARD('n) \<and>
     cgu_graph_subdivision_core_v7 M \<and>
     cgu_piecewise_polynomial_class M P1 \<and>
     cgu_piecewise_polynomial_class M P2 \<and>
     cgu_recovery_geometry_operational_v6 M \<and>
     cgu_recovery_extension_admissible M P1 \<and>
     cgu_recovery_extension_admissible M P2 \<and>
     cgu_local_dn_equal M P1 P2 (cgu_domain M) (cgu_measured M))
    \<longrightarrow> P1 = P2"
\end{lstlisting}
\reviewlabel{English mathematical translation of the Isabelle code}
Let \(n\ge3\). If \(M\) is a regular cell decomposition with the full recovery geometry of Section 2, if \(P_1,P_2\) are admissible coefficient families compatible with recovery, and if they have the same local boundary-energy data on \((\Omega,\Gamma)\), then

\[
 P_1=P_2.
\]

This is equality of the entire coefficient families as finite coefficient representations, not merely almost-everywhere or pointwise equality of their evaluated matrices.

\dossierentry{MAIN-002: \texttt{\detokenize{cgu_uniform_degree_admissible_injectivity_claim_v7}}}
\reviewlabel{English original and manuscript locator}
The corollary of \cite{Carstea2026Polynomial} states:

Suppose \(N_\alpha\le N\) for all cells.  If the subdivision admits an
admissible flat-face recovery order relative to \(\Sigma\) in which each cell is
accessible through at least \(N+3\) flat faces whose supporting hyperplanes
satisfy the corresponding genericity and triple-intersection geometry, and the
recovered regions satisfy the stated face-connectedness assumptions, then the
following holds.  For any
\(\sigma^{(1)},\sigma^{(2)}\in\mathcal P(\mathcal D,N)\) that satisfy the
coefficient-dependent triple-region positivity and auxiliary half-ball
ellipticity hypotheses of Theorem \textup{[source reference: \texttt{thm:main}]} for this recovery order,
equality of their local DN maps on \(\Sigma\) implies
\(\sigma^{(1)}=\sigma^{(2)}\).  Equivalently, the local DN map is injective on
the subclass of \(\mathcal P(\mathcal D,N)\) satisfying those
coefficient-dependent admissibility conditions.
\reviewlabel{Isabelle code}
\begin{lstlisting}
definition cgu_uniform_degree_admissible_injectivity_claim_v7 ::
  "nat \<Rightarrow> ('n::finite, 'c::finite, 'h::finite) cgu_model \<Rightarrow> bool"
where
  "cgu_uniform_degree_admissible_injectivity_claim_v7 N M \<longleftrightarrow>
    (3 \<le> CARD('n) \<and>
     cgu_graph_subdivision_core_v7 M \<and>
     (\<forall>c. cgu_degree M c = N) \<and>
     cgu_recovery_geometry_operational_v6 M)
    \<longrightarrow>
    (\<forall>P1 P2.
      cgu_piecewise_polynomial_class M P1 \<and>
      cgu_piecewise_polynomial_class M P2 \<and>
      cgu_recovery_extension_admissible M P1 \<and>
      cgu_recovery_extension_admissible M P2 \<and>
      cgu_local_dn_equal M P1 P2 (cgu_domain M) (cgu_measured M)
      \<longrightarrow> P1 = P2)"
\end{lstlisting}
\reviewlabel{English mathematical translation of the Isabelle code}
Fix \(N\in\mathbb N\). Let \(n\ge3\), let \(M\) be a regular cell decomposition with full recovery geometry, and assume \(N_c=N\) for every cell. Then for every \(P_1,P_2\), if both are admissible and recovery-compatible and have the same local boundary-energy data on \((\Omega,\Gamma)\), one has \(P_1=P_2\) as coefficient families.

\clearpage
\dossiersubhead{Cited results}

\dossierentry{\texttt{\detokenize{ANALYTIC-ELLIPTIC}}}
\noindent\textit{Formal identifier:} \texttt{\detokenize{hormander_analytic_elliptic_ucp_claim}}\par
\noindent\textit{Relationship to source:} derived.\par\smallskip
\reviewlabel{Source attribution and locator}
Hörmander \cite{Hormander1983}, Section 8.3, printed p. 271, equations (8.3.1)--(8.3.4) and Corollary 8.3.2.
\reviewlabel{Complete source theorem statement (faithful mathematical restatement)}
For \(P(x,D)=\sum_{|\alpha|\leq m}a_\alpha(x)D^\alpha\), the principal symbol is the homogeneous degree-\(m\) part \(p_m(x,\xi)=\sum_{|\alpha|=m}a_\alpha(x)\xi^\alpha\).  The characteristic set consists of the nonzero covectors \((x,\xi)\) for which \(p_m(x,\xi)=0\).  The operator is elliptic precisely when this set is empty, equivalently when \(p_m(x,\xi)\neq0\) for every \(\xi\neq0\).
\reviewlabel{Source attribution and locator}
Hörmander \cite{Hormander1983}, Theorem 8.6.1, printed p. 306.
\reviewlabel{Complete source theorem statement (faithful mathematical restatement)}
If the coefficients of a differential operator \(P\) are real analytic, then an analytic singularity of a distribution \(u\) which is not already an analytic singularity of \(Pu\) can occur only at a characteristic covector; in standard notation,

\[
  \operatorname{WF}_A(u)\subseteq
  \operatorname{WF}_A(Pu)\cup\operatorname{Char}P.
\]
\reviewlabel{Source attribution and locator}
Hörmander \cite{Hormander1983}, Theorem 8.6.5 and its stated elliptic consequence, printed p. 309.
\reviewlabel{Complete source theorem statement (faithful mathematical restatement)}
If \(u\) is a distribution on an open set \(X\) and \(Pu=0\), then every exterior conormal to \(\operatorname{supp}u\) is characteristic for \(P\).  Hence, if \(P\) is elliptic, the support of \(u\) has no boundary point in \(X\).  Consequently, when \(X\) is connected, a solution which vanishes in a neighborhood of one point vanishes throughout \(X\).
\reviewlabel{Source attribution and locator}
Hörmander \cite{Hormander1983}, Theorem 1.2.5, printed p. 15.
\reviewlabel{Complete source theorem statement (faithful mathematical restatement)}
If \(f,g\in L^1_{\mathrm{loc}}(X)\) and \(\int_X f\phi=\int_X g\phi\) for every \(\phi\in C_c^\infty(X)\), then \(f=g\) almost everywhere on \(X\).
\reviewlabel{Source attribution and locator}
Hörmander \cite{Hormander1983}, regular-distribution identification on printed p. 37 and Definition 3.1.1 on printed p. 55.
\reviewlabel{Complete source theorem statement (faithful mathematical restatement)}
Printed p. 37 identifies \(L^1_{\mathrm{loc}}(X)\) modulo almost-everywhere equality injectively with the regular distributions \(\phi\mapsto\int_X f\phi\).  Definition 3.1.1, printed p. 55, defines \((\partial_i u)(\phi)=-u(\partial_i\phi)\) and multiplication of a distribution by a smooth function.
\reviewlabel{Exact derived mathematical result formalized}
Let \(X\subset\mathbb R^n\), \(n\geq1\), be open and connected.  Let \(A(x)\) be a symmetric matrix polynomial which is uniformly positive definite on an open neighborhood of \(\overline X\).  If a real locally square-integrable regular distribution \(u\) is a weak solution of \(-\operatorname{div}(A\nabla u)=0\) on \(X\), and \(u=0\) almost everywhere on some nonempty open subset of \(X\), then \(u=0\) almost everywhere on \(X\).  This is the exact mathematical content of the formal interface.

For the notation correspondence, polynomial entries are real analytic and the principal symbol of \(-\operatorname{div}(A\nabla\cdot)\), in Hörmander's \(D_j=(1/i)\partial_j\) convention, is \(\xi^{T}A(x)\xi\).  Uniform positive definiteness makes it nonzero for \(\xi\neq0\).  The weak variational identity is the distribution equation, and the regular-distribution results convert distributional vanishing to almost-everywhere vanishing.  Polynomiality and uniform ellipticity near the closure are conservative strengthenings of local analyticity and ellipticity.
\reviewlabel{Isabelle code}
\begin{lstlisting}
definition hormander_analytic_elliptic_ucp_claim ::
  "'n::finite ucp_point set \<Rightarrow> 'n ucp_coefficient \<Rightarrow>
   ('n ucp_point \<Rightarrow> real) \<Rightarrow> bool"
where
  "hormander_analytic_elliptic_ucp_claim X P u \<longleftrightarrow>
    ((0 < CARD('n) \<and> open X \<and> connected X \<and>
      ucp_matrix_polynomial P \<and>
      (\<forall>x. ucp_symmetric_matrix (P x)) \<and>
      ucp_uniformly_elliptic_near X P \<and>
      ucp_weak_solution_on X P u \<and>
      (\<exists>V. V \<noteq> {} \<and> open V \<and> V \<subseteq> X \<and>
        (AE x in lborel. x \<in> V \<longrightarrow> u x = 0)))
     \<longrightarrow> (AE x in lborel. x \<in> X \<longrightarrow> u x = 0))"
\end{lstlisting}
\reviewlabel{English mathematical translation of the Isabelle code}
Let \(n>0\), let \(X\subset\mathbb R^n\) be open and connected, and let \(P:\mathbb R^n\to\mathbb R^{n\times n}\) have polynomial entries, be symmetric at every ambient point, and be uniformly elliptic on an open neighborhood of \(\overline X\). If \(u\) is a locally weak solution of \(\operatorname{div}(P\nabla u)=0\) in \(X\), and there exists a nonempty open \(V\subset X\) such that

\[
 u=0\quad\text{for Lebesgue-a.e. ambient }x\in V,
\]

then

\[
 u=0\quad\text{for Lebesgue-a.e. ambient }x\in X.
\]

\dossierentry{\texttt{\detokenize{WEAK-GRADIENT-LOCALITY}}}
\noindent\textit{Formal identifier:} \texttt{\detokenize{hormander_weak_gradient_locality_claim}}\par
\noindent\textit{Relationship to source:} derived.\par\smallskip
\reviewlabel{Source attribution and locator}
Hörmander \cite{Hormander1983}, Definitions 1.2.1--1.2.2, the following zero-extension observation, and Lemma 1.2.3, printed p. 14.
\reviewlabel{Complete source theorem statement (faithful mathematical restatement)}
The test space on an open set \(X\) is \(C_c^\infty(X)\), with support defined as the closure of the nonzero set; its members extend smoothly by zero to the ambient Euclidean space.  The lemma constructs a nonnegative compactly supported smooth bump positive at a prescribed point.  Translation, scaling, a finite compact subcover, and summation therefore give, for every compact \(K\Subset X\), a nonnegative test function which is strictly positive on \(K\).
\reviewlabel{Source attribution and locator}
Hörmander \cite{Hormander1983}, Theorem 1.2.5, printed p. 15.
\reviewlabel{Complete source theorem statement (faithful mathematical restatement)}
Two locally integrable functions with identical pairings against every member of \(C_c^\infty(X)\) are equal almost everywhere.
\reviewlabel{Source attribution and locator}
Hörmander \cite{Hormander1983}, regular-distribution identification on printed p. 37 and Definition 3.1.1 on printed p. 55.
\reviewlabel{Complete source theorem statement (faithful mathematical restatement)}
Printed p. 37 identifies locally integrable functions modulo almost-everywhere equality with their regular distributions.  Definition 3.1.1, printed p. 55, defines the distributional derivative by \((\partial_i u)(\phi)=-u(\partial_i\phi)\).
\reviewlabel{Exact derived mathematical result formalized}
Let \(X\) be open.  Suppose real \(u\) and real vector field \(D u\) satisfy the project's local-square-integrability and weak-gradient predicates on \(X\), and suppose \(u=0\) almost everywhere on \(X\).  Then \(D u=0\) almost everywhere on \(X\).

Coordinatewise, the weak identity gives \(\int_X u\,\partial_i\phi=-\int_X(Du)_i\phi\); the left side is zero. The weak-gradient product-integrability clause, combined with a positive localized test, proves \((Du)_i\in L^1_{\mathrm{loc}}(X)\).  Theorem 1.2.5 then makes each component zero almost everywhere, and finite dimensionality combines the component conclusions.  The formal real test class is a conservative specialization of Hörmander's complex convention: real and imaginary parts recover all complex pairings.
\reviewlabel{Isabelle code}
\begin{lstlisting}
definition hormander_weak_gradient_locality_claim ::
  "'n::finite ucp_point set \<Rightarrow>
   ('n ucp_point \<Rightarrow> real) \<Rightarrow>
   ('n ucp_point \<Rightarrow> 'n ucp_point) \<Rightarrow> bool"
where
  "hormander_weak_gradient_locality_claim X u Du \<longleftrightarrow>
    ((open X \<and>
      ucp_locally_square_integrable_on X u \<and>
      ucp_locally_square_integrable_on X Du \<and>
      ucp_weak_gradient_on X u Du \<and>
      (AE x in lborel. x \<in> X \<longrightarrow> u x = 0))
     \<longrightarrow> (AE x in lborel. x \<in> X \<longrightarrow> Du x = 0))"
\end{lstlisting}
\reviewlabel{English mathematical translation of the Isabelle code}
Let \(X\) be open and let \((u,G)\in H^1_{\mathrm{loc},\mathrm{pair}}(X)\). If \(u(x)=0\) for Lebesgue-a.e. ambient \(x\in X\), then \(G(x)=0\) for Lebesgue-a.e. ambient \(x\in X\).

\dossierentry{\texttt{\detokenize{HORMANDER-REGULAR-DISTRIBUTION-INJECTIVITY}}}
\noindent\textit{Formal identifier:} \texttt{\detokenize{hormander_regular_distribution_injectivity_claim}}\par
\noindent\textit{Relationship to source:} direct.\par\smallskip
\reviewlabel{Source attribution and locator}
Hörmander \cite{Hormander1983}, Definitions 1.2.1--1.2.2 and Lemma 1.2.3, printed p. 14.
\reviewlabel{Complete source theorem statement (faithful mathematical restatement)}
These results define the smooth compactly supported tests on an open set, their support, and the localized smooth bumps used to separate a nonzero locally integrable function from the zero distribution.
\reviewlabel{Source attribution and locator}
Hörmander \cite{Hormander1983}, Theorem 1.2.5, printed p. 15.
\reviewlabel{Complete source theorem statement (faithful mathematical restatement)}
For an open \(X\subset\mathbb R^n\), if \(f,g\in L^1_{\mathrm{loc}}(X)\) satisfy \(\int_X f\phi=\int_Xg\phi\) for every \(\phi\in C_c^\infty(X)\), then \(f=g\) almost everywhere on \(X\).
\reviewlabel{Source attribution and locator}
Hörmander \cite{Hormander1983}, regular-distribution identification on printed p. 37.
\reviewlabel{Complete source theorem statement (faithful mathematical restatement)}
The source expresses the preceding result as injectivity of the regular-distribution embedding of \(L^1_{\mathrm{loc}}(X)\) modulo almost-everywhere equality.

\textbf{Direct formal specialization.}

For every open \(X\) and real \(u\in L^1(X)\), if \(u\phi\) is integrable and \(\int_Xu\phi=0\) for every project smooth compact test \(\phi\), then \(u=0\) almost everywhere on \(X\).  Global \(L^1(X)\) is stronger than the source's local integrability, the comparison function is specialized to zero, and the project test predicate is the globally zero-extended realization of \(C_c^\infty(X)\).  The explicit product-integrability premise prevents use of Isabelle's totalized integral outside the ordinary integrable case.
\reviewlabel{Isabelle code}
\begin{lstlisting}
definition hormander_regular_distribution_injectivity_claim ::
  "'n::finite itself \<Rightarrow> bool"
where
  "hormander_regular_distribution_injectivity_claim dimension_type \<longleftrightarrow>
    (\<forall>X :: 'n cgu_point set. \<forall>u.
      open X \<and> hormander_regular_distribution_test_zero_on X u
      \<longrightarrow>
      (AE x in restrict_space lborel X. u x = 0))"
\end{lstlisting}
\reviewlabel{English mathematical translation of the Isabelle code}
In every finite dimension, for every open \(X\) and every Lebesgue-integrable \(u\) on \(X\), if for every \(\phi\in\mathscr D(X)\) the product \(u\phi\) is integrable and

\[
 \int_Xu\phi=0,
\]

then \(u=0\) almost everywhere with respect to Lebesgue measure restricted to \(X\).

\dossierentry{\texttt{\detokenize{EVANS-POINCARE-ELLIPTIC-GREEN}}}
\noindent\textit{Formal identifier:} \texttt{\detokenize{evans_poincare_elliptic_green_claim}}\par
\noindent\textit{Relationship to source:} derived.\par\smallskip
\reviewlabel{Source attribution and locator}
Evans \cite{Evans2010}, Section 5.6.1, Theorem 3 and its separately displayed particular consequence, printed pp. 279--280.
\reviewlabel{Complete source theorem statement (faithful mathematical restatement)}
If \(U\subset\mathbb R^n\) is bounded and open and \(u\in W^{1,p}_0(U)\), \(1\leq p<\infty\), then \(\|u\|_{L^p(U)}\leq C\|Du\|_{L^p(U)}\), with \(C\) depending on \(U,n,p\).  The interface uses only \(p=2\).
\reviewlabel{Source attribution and locator}
Evans \cite{Evans2010}, Section 5.9.1, Theorem 1, printed pp. 299--300.
\reviewlabel{Complete source theorem statement (faithful mathematical restatement)}
The negative Sobolev space \(H^{-1}(U)\) is the bounded dual of \(H^1_0(U)\).  Equivalently, every such functional has a representation by \(L^2\) functions \(f=f^0-\sum_{i=1}^n\partial_i f^i\), and conversely every such expression defines a bounded functional on \(H^1_0(U)\).
\reviewlabel{Source attribution and locator}
Evans \cite{Evans2010}, Sections 6.1.1--6.1.2, printed pp. 311--315.
\reviewlabel{Complete source theorem statement (faithful mathematical restatement)}
On a bounded open domain, a divergence-form operator with bounded measurable leading coefficients, uniformly elliptic almost everywhere (and symmetric in the symmetric case), gives the bilinear energy form on \(H^1_0(U)\).  A weak zero-Dirichlet solution is an element \(u\in H^1_0(U)\) satisfying that bilinear identity against every \(v\in H^1_0(U)\).
\reviewlabel{Source attribution and locator}
Evans \cite{Evans2010}, Section 6.2.1, Theorem 1 (Lax--Milgram), printed pp. 315--317.
\reviewlabel{Complete source theorem statement (faithful mathematical restatement)}
If \(H\) is a real Hilbert space and a bilinear form \(B:H\times H\to\mathbb R\) is bounded and coercive, \(|B[u,v]|\leq M\|u\|\|v\|\) and \(B[u,u]\geq\beta\|u\|^2\) with \(\beta>0\), then every \(f\in H^*\) has a unique \(u\in H\) satisfying \(B[u,v]=f(v)\) for all \(v\in H\), and \(\|u\|\leq\|f\|/\beta\).
\reviewlabel{Source attribution and locator}
Evans \cite{Evans2010}, Section 6.2.2, Theorem 2, printed pp. 317--318.
\reviewlabel{Complete source theorem statement (faithful mathematical restatement)}
Under the elliptic coefficient hypotheses, the weak bilinear form satisfies \(|B[u,v]|\leq C\|u\|_{H^1_0(U)}\|v\|_{H^1_0(U)}\) and a Gårding bound \(B[u,u]\geq\beta\|u\|_{H^1_0(U)}^2-\gamma\|u\|_{L^2(U)}^2\) for constants \(C,\beta>0\) and \(\gamma\geq0\).  Thus the leading elliptic part controls the gradient norm, while lower-order terms cost only an \(L^2\) term.
\reviewlabel{Source attribution and locator}
Evans \cite{Evans2010}, Section 6.2.2, Theorem 3, printed pp. 318--319.
\reviewlabel{Complete source theorem statement (faithful mathematical restatement)}
There is \(\gamma\geq0\) such that, for \(f\in L^2(U)\) and every \(\mu\geq\gamma\), the shifted zero-Dirichlet weak problem \(Lu+\mu u=f\) has a unique solution in \(H^1_0(U)\), with the corresponding energy estimate.  This theorem does not itself state the unshifted arbitrary-\(H^{-1}\) result below.
\reviewlabel{Exact derived mathematical result formalized}
Let \(U\neq\varnothing\) be bounded and open, and let \(a(x)\) be a chosen restricted-Borel measurable, pointwise bounded, pointwise symmetric matrix with one positive pointwise uniform ellipticity constant.  Then the pure leading energy form is bounded and strictly coercive on the project's \(H^1_0(U)\) quotient.  Moreover, one constant \(C\geq0\), depending only on \(U,a\), works for every bounded real linear source \(F\) with bound \(K\): there is a zero-boundary weak solution of \(-\operatorname{div}(a\nabla u)=F\), it satisfies \(\|u\|_{H^1(U)}\leq CK\), and it is unique modulo zero \(H^1\) distance.

The derivation takes \(H=H^1_0(U)\) and \(B[u,v]=\int_U(aDu)\cdot Dv\).  Coefficient boundedness gives continuity; ellipticity plus the \(p=2\) Poincaré inequality gives coercivity in the full project \(H^1\) norm; Section 5.9.1 identifies the source with an element of \(H^*\); Lax--Milgram gives the unshifted solution, uniqueness, and the uniform estimate.  The pointwise representative assumptions conservatively strengthen Evans's almost-everywhere coefficient assumptions.
\reviewlabel{Isabelle code}
\begin{lstlisting}
definition evans_poincare_elliptic_green_claim ::
  "'n::finite cgu_point set \<Rightarrow> 'n mclean_coefficient \<Rightarrow> bool"
where
  "evans_poincare_elliptic_green_claim U a \<longleftrightarrow>
    (0 < CARD('n) \<and> U \<noteq> {} \<and> open U \<and>
      evans_bounded_domain U \<and> evans_elliptic_coefficient_on U a)
    \<longrightarrow>
    (mclean_variational_form_on U a \<and>
      (\<exists>C. 0 \<le> C \<and>
        (\<forall>F K. mclean_source_bound_on U F K \<longrightarrow>
          (\<exists>u Du. mclean_green_solution_for U a F u Du \<and>
            mclean_h1_norm U u Du \<le> C * K \<and>
            (\<forall>v Dv. mclean_green_solution_for U a F v Dv
              \<longrightarrow>
              cgu_h1_squared_distance U u Du v Dv = 0)))))"
\end{lstlisting}
\reviewlabel{English mathematical translation of the Isabelle code}
Let \(n>0\), let \(U\ne\varnothing\) be open and bounded, and let \(a\) be Borel measurable on \(U\), pointwise bounded there by one constant \(B\ge0\), symmetric at every point of \(U\), and uniformly elliptic there: for some \(\theta>0\),

\[
 \theta\lVert\xi\rVert^2\le \xi\cdot a(x)\xi
 \quad(x\in U,\ \xi\in\mathbb R^n).
\]

Then the energy form \(B_{U,a}\) is stable, and there exists \(C\ge0\), independent of \(\mathcal F\) and \(K\), such that whenever \(\mathcal F\) is bounded linear on \(\mathcal H_0^1(U)\) with bound \(K\), there is a source solution \((u,G)\) satisfying

\[
 \lVert(u,G)\rVert_U\le CK.
\]

Every other source solution \((v,H)\) has

\[
 d_U^2((u,G),(v,H))=0.
\]

\dossierentry{\texttt{\detokenize{KANG-YUN-LTU-LOCAL-BOUNDARY-METRIC}}}
\noindent\textit{Formal identifier:} \texttt{\detokenize{kang_yun_ltu_local_boundary_metric_claim_v6}}\par
\noindent\textit{Relationship to source:} derived.\par\smallskip
\reviewlabel{Source attribution and locator}
Kang--Yun \cite{KangYun2003}, weak energy pairing on printed p. 719 and localized-DN definition on printed p. 723.
\reviewlabel{Complete source theorem statement (faithful mathematical restatement)}
For an open connected measured boundary patch \(\Gamma\), the localized Dirichlet-to-Neumann map takes \(H^{1/2}\) boundary data supported in \(\Gamma\), solves the conductivity/metric equation, and restricts the Neumann response to \(\Gamma\).  Its weak form is the boundary energy pairing of an input solution with an arbitrary lift of the test trace.
\reviewlabel{Source attribution and locator}
Kang--Yun \cite{KangYun2003}, Theorem 1.3, printed p. 723.
\reviewlabel{Complete source theorem statement (faithful mathematical restatement)}
Let two uniformly elliptic Riemannian metrics have \(C^{m,p}\) regularity near a connected open boundary patch \(\Gamma\), with \(m\geq1\) and \(p>0\).  For every compact \(K\Subset\Gamma\), the local DN maps stably determine the metrics near \(K\), modulo a local diffeomorphism which fixes the boundary: after the appropriate pullback, the metric difference is bounded by a positive power of the local-DN operator-norm difference.  In particular, equality of the local maps gives equality up to that boundary-fixing gauge on \(K\).
\reviewlabel{Source attribution and locator}
Kang--Yun \cite{KangYun2003}, Lemmas 2.1 and 2.2, printed pp. 724--727.
\reviewlabel{Complete source theorem statement (faithful mathematical restatement)}
In boundary-normal coordinates and for a boundary point and tangential frequency, the authors construct high-frequency boundary packets supported in an arbitrarily small part of \(\Gamma\), extend them into a shrinking boundary box \(D_N\), and establish the packet, trace, and Sobolev norm bounds used in the asymptotic energy calculation.
\reviewlabel{Source attribution and locator}
Kang--Yun \cite{KangYun2003}, Section 3, equations (3.7)--(3.9) and the intervening Hardy inequality, printed pp. 728--731.
\reviewlabel{Complete source theorem statement (faithful mathematical restatement)}
Writing the exact solution as the explicit packet plus a correction, the source proves the local correction estimate (3.7), the artificial-top trace estimate (3.8), the displayed Hardy inequality before (3.9), and the weighted residual estimate (3.9).  In the notation used by the reviewed derivation, the local estimate is \(\|s_N\|_{H^1(D_N)}\leq C N^{1-|\alpha|/2}\), with the subsequent estimates making the correction terms negligible at the normalization used in the DN limit.
\reviewlabel{Source attribution and locator}
Kang--Yun \cite{KangYun2003}, Theorem 4.1 and equations (4.5)--(4.8), printed pp. 731--733.
\reviewlabel{Complete source theorem statement (faithful mathematical restatement)}
The normalized high-frequency local-DN energy limit at a chosen boundary point determines the tangential metric quadratic form there. Varying the tangential covector and polarizing determines the complete tangential bilinear form.
\reviewlabel{Source attribution and locator}
Kang--Yun \cite{KangYun2003}, boundary-normal-coordinate and diffeomorphism reduction, printed pp. 733--734.
\reviewlabel{Complete source theorem statement (faithful mathematical restatement)}
A smooth positive metric near the observed patch can be put into the boundary-normal form used by the packet computation.  The comparison of the two normal-coordinate charts is a local diffeomorphism whose restriction to the boundary is the identity; hence its differential is the identity on tangent vectors.
\reviewlabel{Source attribution and locator}
Lassas--Taylor--Uhlmann \cite{LassasTaylorUhlmann2003}, remote rough-boundary setting on printed p. 209 and local-symbol discussion in Section 2, printed pp. 210--211.
\reviewlabel{Complete source theorem statement (faithful mathematical restatement)}
In their analytic, Wiener-regular complete-manifold setting, boundary roughness is allowed away from the measured analytic patch, and the boundary symbol is obtained locally from the DN map.  This is corroboration of locality in a different setting, not a theorem proving the graph-domain corollary below.
\reviewlabel{Exact derived mathematical result formalized}
In dimension at least three, let \(U\) be a bounded connected graph-Lipschitz domain contained in the fixed subdivision model.  Let \(\gamma\) be a nonempty relatively open flat boundary patch with a full one-sided collar, disjoint from the subdivision skeleton.  Suppose each of two globally bounded uniformly elliptic piecewise-polynomial coefficients agrees on an open collar of \(\gamma\) with a smooth positive metric-density realization \(a_j=\sqrt{\det g_j}\,g_j^{-1}\).  If the completed real weak local-DN energy pairings agree for every supported input/test trace, with solutions existing for those data, then for every \(x\in\gamma\) and every pair of vectors tangent to the supporting hyperplane, \(g_1(x)(v,w)=g_2(x)(v,w)\).  No normal component or remote coefficient is identified.

The derivation localizes the Kang--Yun packet with a fixed smooth collar cutoff.  Cutoff derivatives occur a fixed positive normal distance from the boundary and are exponentially small; graph-domain trace lifting, Poincaré, and global ellipticity give the correction bound without remote smoothness. Four real bilinear pairings reconstruct the complex packet pairing.  The boundary-fixing gauge has identity tangential differential, so equality after pullback gives the stated covariant tangential equality.  Pointwise shrinking handles a disconnected formal patch.  Neither paper prints this exact graph-domain implication verbatim, and LTU supplies corroboration only.
\reviewlabel{Isabelle code}
\begin{lstlisting}
definition kang_yun_ltu_local_boundary_metric_claim_v6 ::
  "('n::finite, 'c::finite, 'h::finite) cgu_model \<Rightarrow>
    ('c \<Rightarrow> 'n cgu_matrix_polynomial) \<Rightarrow>
    ('c \<Rightarrow> 'n cgu_matrix_polynomial) \<Rightarrow>
    'n cgu_point set \<Rightarrow> 'n cgu_point set \<Rightarrow>
    'n cgu_point \<Rightarrow> real \<Rightarrow>
    'n lu_metric \<Rightarrow> 'n lu_metric \<Rightarrow>
    'n lu_metric \<Rightarrow> 'n lu_metric \<Rightarrow>
    'n lu_density \<Rightarrow> 'n lu_density \<Rightarrow> bool"
where
  "kang_yun_ltu_local_boundary_metric_claim_v6 M P1 P2 U gamma normal offset
      g1 g2 a1 a2 rho1 rho2 \<longleftrightarrow>
    (3 \<le> CARD('n) \<and>
     cgu_graph_subdivision_core_v7 M \<and>
     cgu_piecewise_polynomial_class M P1 \<and>
     cgu_piecewise_polynomial_class M P2 \<and>
     cgu_graph_lipschitz_domain_v5 U \<and>
     U \<subseteq> cgu_domain M \<and>
     cgu_flat_face_patch U normal offset gamma \<and>
     gamma \<inter> cgu_subdivision_edge_corner_set M = {} \<and>
     lu_smooth_metric_extension_near M P1 U gamma g1 a1 rho1 \<and>
     lu_smooth_metric_extension_near M P2 U gamma g2 a2 rho2 \<and>
     cgu_local_dn_equal M P1 P2 U gamma)
    \<longrightarrow> lu_induced_boundary_metrics_equal_on gamma normal g1 g2"
\end{lstlisting}
\reviewlabel{English mathematical translation of the Isabelle code}
Let \(n\ge3\). Suppose \(M\) is a regular cell decomposition, \(P_1,P_2\) are admissible, \(U\) is a Lipschitz-graph domain with \(U\subset\Omega\), and \(\gamma\) is a one-sided planar boundary patch of \(U\) with unit normal \(\nu\) and offset \(b\). Assume

\[
 \gamma\cap\Sigma_M=\varnothing,
\]

each \(P_i\) is realized near \(\gamma\) by a metric/conductivity package \((g_i,a_i,\rho_i)\), and \(P_1,P_2\) have the same local boundary-energy data on \((U,\gamma)\). Then \(g_1\) and \(g_2\) agree tangentially on \(\gamma\):

\[
 \xi\cdot g_1(x)\eta=\xi\cdot g_2(x)\eta
\]

for every \(x\in\gamma\) and every \(\xi,\eta\perp\nu\).

\dossierentry{\texttt{\detokenize{MCLEAN-GRAPH-SUPPORTED-H1-ZERO}}}
\noindent\textit{Formal identifier:} \texttt{\detokenize{mclean_graph_supported_h1_zero_claim}}\par
\noindent\textit{Relationship to source:} direct.\par\smallskip
\reviewlabel{Source attribution and locator}
McLean \cite{McLean2000}, Definition 3.28, printed pp. 89--90.
\reviewlabel{Complete source theorem statement (faithful mathematical restatement)}
A Lipschitz domain has compact boundary and a finite cover by rigid Cartesian-coordinate neighborhoods in which the domain is a Lipschitz hypograph.  Replacing the local Lipschitz graph functions by continuous functions gives the book's \(C^0\)-domain class; every Definition 3.28 Lipschitz domain is therefore a \(C^0\) domain.
\reviewlabel{Source attribution and locator}
McLean \cite{McLean2000}, Theorem 3.29(ii), printed pp. 91--92.
\reviewlabel{Complete source theorem statement (faithful mathematical restatement)}
For a \(C^0\) domain \(\Omega\) and every real Sobolev order \(s\), the smooth compactly supported class \(\mathcal D(\Omega)\) is dense, in the ambient \(H^s(\mathbb R^n)\) norm, in the global Sobolev space consisting of distributions supported in \(\overline\Omega\).  Equivalently, the ambient compact-test closure equals that supported global Sobolev space.

\textbf{Direct formal specialization.}

At real scalar order \(s=1\), suppose \(U\) is a bounded graph-Lipschitz domain and \((u,Du)\) is a global project \(H^1(\mathbb R^n)\) pair whose potential and selected weak gradient vanish almost everywhere outside \(\overline U\).  Then \((u,Du)\) belongs to the project's sequential \(H^1_0(U)\) compact-test closure.  The graph-Lipschitz premise is stronger than the source's \(C^0\) hypothesis, the gradient-support premise is an additional compatible restriction, and source global convergence is stronger than the conclusion's restricted \(U\)-distance convergence.
\reviewlabel{Isabelle code}
\begin{lstlisting}
definition mclean_graph_supported_h1_zero_claim ::
  "'n::finite itself \<Rightarrow> bool"
where
  "mclean_graph_supported_h1_zero_claim dimension_type \<longleftrightarrow>
    0 < CARD('n) \<and>
    (\<forall>(U :: 'n cgu_point set) u Du.
      cgu_graph_lipschitz_domain_v5 U \<and>
      cgu_h1_pair_on UNIV u Du \<and>
      (AE x in lborel.
        x \<notin> closure U \<longrightarrow> u x = 0 \<and> Du x = 0)
      \<longrightarrow> cgu_h1_zero_pair_on U u Du)"
\end{lstlisting}
\reviewlabel{English mathematical translation of the Isabelle code}
This condition is the conjunction of \(n>0\) and the following assertion: for every \(U,u,G\), if \(U\) is a Lipschitz-graph domain, \((u,G)\in\mathcal H^1(\mathbb R^n)\), and, for Lebesgue-a.e. ambient \(x\),

\[
 x\notin\overline U\implies u(x)=0\ \text{and}\ G(x)=0,
\]

then \((u,G)\in\mathcal H_0^1(U)\).

\dossierentry{\texttt{\detokenize{MCLEAN-H1-ZERO-BIDUAL-REPRESENTATION}}}
\noindent\textit{Formal identifier:} \texttt{\detokenize{mclean_h1_zero_bidual_representation_claim_v2}}\par
\noindent\textit{Relationship to source:} derived.\par\smallskip
\reviewlabel{Source attribution and locator}
McLean \cite{McLean2000}, weak-derivative Sobolev construction and Hilbert conclusion, printed pp. 73--75.
\reviewlabel{Complete source theorem statement (faithful mathematical restatement)}
On every nonempty open \(U\), \(W^r_p(U)\) consists of \(L^p\) functions possessing all weak derivatives through order \(r\) in \(L^p\), with the derivative graph norm.  For \(p=2\), the sum of the \(L^2\) inner products of the weak derivatives induces that norm, and \(W^s(U)\) is a Hilbert space for every real \(s\geq0\).
\reviewlabel{Source attribution and locator}
McLean \cite{McLean2000}, zero-boundary closure notation, printed p. 77.
\reviewlabel{Complete source theorem statement (faithful mathematical restatement)}
The source defines the conventional zero-boundary space as the closure of \(\mathcal D(U)\) in the relevant Sobolev norm.  For this interface the cited page supplies nomenclature only; no equality with a distinct supported Bessel-potential space on an arbitrary open set is used.
\reviewlabel{Source attribution and locator}
McLean \cite{McLean2000}, Hilbert reflexivity statement, printed p. 79.
\reviewlabel{Complete source theorem statement (faithful mathematical restatement)}
Every Hilbert space is reflexive: the canonical isometric embedding \(J:H\to H^{**}\), \((Ju)(F)=F(u)\), is onto.  Thus every bounded linear functional on \(H^*\) is evaluation at some \(u\in H\).
\reviewlabel{Exact derived mathematical result formalized}
For every nonempty open \(U\), let \(H\) be the quotient Hilbert space represented by project \(H^1_0(U)\) potential/weak-gradient pairs.  There exist a set \(E\) and seminorm \(p\) such that \(E\) is exactly the set of total raw functions inducing bounded real linear functionals on \(H\), and \(p(F)\) is the least valid source bound, hence the operator seminorm.  Every carrier-relative real-linear \(g:E\to\mathbb R\) bounded by \(p\) is of the form \(g(F)=F(u,Du)\) for some project \(H^1_0(U)\) pair.

The project's potential-plus-gradient norm is the order-one graph norm, and the compact-test closure is a closed Hilbert subspace after quotienting zero distance.  Bounded raw sources descend to its continuous dual; the least bound is the dual norm.  A \(p\)-bounded \(g\) vanishes on the seminorm kernel, descends to the genuine bidual, and reflexivity represents it by evaluation. The raw carrier \(E\), least-bound plumbing, and displayed evaluation formula are standard consequences, not statements printed verbatim by McLean.
\reviewlabel{Isabelle code}
\begin{lstlisting}
definition mclean_h1_zero_bidual_representation_claim_v2 ::
  "'n::finite itself \<Rightarrow> bool"
where
  "mclean_h1_zero_bidual_representation_claim_v2 dimension_type \<longleftrightarrow>
    (\<forall>U :: 'n cgu_point set.
      U \<noteq> {} \<and> open U
      \<longrightarrow>
      (\<exists>E p. mclean_h1_zero_bidual_realization_on_v2 U E p))"
\end{lstlisting}
\reviewlabel{English mathematical translation of the Isabelle code}
In every finite dimension, every nonempty open \(U\subset\mathbb R^n\) admits a pair \((E,p)\) with all properties of the total-functional representation package in Section 5: \(E\) consists exactly of the everywhere-defined functionals having some stipulated zero-boundary bound, \(p\) is their least-bound seminorm, and every \(p\)-bounded linear map on \(E\) has an existential pair representation. No quotienting or uniqueness is asserted.

\dossierentry{\texttt{\detokenize{MCLEAN-H1-ZERO-EXTENSION-LOCALIZATION}}}
\noindent\textit{Formal identifier:} \texttt{\detokenize{mclean_h1_zero_extension_localization_claim}}\par
\noindent\textit{Relationship to source:} derived.\par\smallskip
\reviewlabel{Source attribution and locator}
McLean \cite{McLean2000}, smooth compact tests and the two compact-test closures, printed pp. 61, 65, and 77--78.
\reviewlabel{Complete source theorem statement (faithful mathematical restatement)}
The test class is \(\mathcal D(\Omega)=C_c^\infty(\Omega)\). \(H^s_0(\Omega)\) is its closure in the intrinsic/restriction \(H^s(\Omega)\) norm, whereas \(\widetilde H^s(\Omega)\) is its closure in the ambient \(H^s(\mathbb R^n)\) norm.  These two spaces are not identified on an arbitrary open set in this interface.
\reviewlabel{Source attribution and locator}
McLean \cite{McLean2000}, Theorem 3.6, printed p. 64.
\reviewlabel{Complete source theorem statement (faithful mathematical restatement)}
Given a closed set and an open neighborhood (equivalently, an arbitrarily small prescribed neighborhood), there is a smooth cutoff taking values between zero and one, equal to one on the closed set and zero outside the prescribed neighborhood. For a compact set contained in an open set, the cutoff may be chosen compactly supported in that open set.
\reviewlabel{Source attribution and locator}
McLean \cite{McLean2000}, Theorem 3.16, printed p. 80.
\reviewlabel{Complete source theorem statement (faithful mathematical restatement)}
For every nonnegative integer \(m\), the global weak-derivative space \(W^m(\mathbb R^n)\) and the Bessel-potential space \(H^m(\mathbb R^n)\) are the same set with equivalent norms.
\reviewlabel{Source attribution and locator}
McLean \cite{McLean2000}, Theorem 3.20, printed p. 83.
\reviewlabel{Complete source theorem statement (faithful mathematical restatement)}
Multiplication by a smooth compactly supported cutoff is a bounded linear operator on \(H^s(\Omega)\) and, separately, on \(\widetilde H^s(\Omega)\).  The theorem does not state a multiplier clause for \(H^s_0(\Omega)\) under that name.
\reviewlabel{Exact derived mathematical result formalized}
Two order-one representative statements are formalized.  First, if \(\varnothing\neq U\subseteq W\) are open and \((u,Du)\) is in the project \(H^1_0(U)\) closure, then the zero extensions of both representatives form a project \(H^1_0(W)\) pair.  Second, if nonempty open \(V\subseteq W\), compact \(K\subseteq V\), and a project \(H^1_0(W)\) pair vanishes almost everywhere together with its gradient off \(K\), then the same representatives form a project \(H^1_0(V)\) pair.

For extension, each compact test in \(U\) is also a compact test in \(W\), and zero extension turns intrinsic order-one convergence into global \(W^1\) convergence; Theorem 3.16 supplies the ambient \(H^1\) bridge.  For localization, choose the Theorem 3.6 cutoff equal to one near \(K\); Theorem 3.20 gives continuity on the ambient compact-test closure, and the product equals the original pair because both representatives vanish off \(K\).  If \(K=\varnothing\), the conclusion is the zero class.  No arbitrary-domain identity \(H^1_0=\widetilde H^1\) is used.
\reviewlabel{Isabelle code}
\begin{lstlisting}
definition mclean_h1_zero_extension_localization_claim ::
  "'n::finite itself \<Rightarrow> bool"
where
  "mclean_h1_zero_extension_localization_claim dimension_type \<longleftrightarrow>
    ((\<forall>(U :: 'n cgu_point set) W u Du.
        U \<noteq> {} \<and> open U \<and> open W \<and> U \<subseteq> W \<and>
        cgu_h1_zero_pair_on U u Du
        \<longrightarrow>
        cgu_h1_zero_pair_on W
          (mclean_zero_extension_on U u)
          (mclean_zero_extension_on U Du)) \<and>
     (\<forall>(W :: 'n cgu_point set) V K u Du.
        W \<noteq> {} \<and> V \<noteq> {} \<and> open W \<and> open V \<and> V \<subseteq> W \<and>
        compact K \<and> K \<subseteq> V \<and>
        cgu_h1_zero_pair_on W u Du \<and>
        (AE x in restrict_space lborel W.
          x \<notin> K \<longrightarrow> u x = 0 \<and> Du x = 0)
        \<longrightarrow> cgu_h1_zero_pair_on V u Du))"
\end{lstlisting}
\reviewlabel{English mathematical translation of the Isabelle code}
In every finite dimension, both statements hold.

\begin{enumerate}[start=1]
\item If \(U\ne\varnothing\) and \(U,W\) are open with \(U\subset W\), then the pointwise zero extension to \(W\) of every \((u,G)\in\mathcal H_0^1(U)\) belongs to \(\mathcal H_0^1(W)\).
\item If \(W,V\ne\varnothing\) are open, \(V\subset W\), \(K\subset V\) is compact, \((u,G)\in\mathcal H_0^1(W)\), and
\end{enumerate}

\[
   x\notin K\implies u(x)=0\ \text{and}\ G(x)=0
   \]

for almost every \(x\) with respect to Lebesgue measure restricted to \(W\), then \((u,G)\in\mathcal H_0^1(V)\).

\dossierentry{\texttt{\detokenize{MCLEAN-HYPERPLANE-OFFSET-SOURCE}}}
\noindent\textit{Formal identifier:} \texttt{\detokenize{mclean_hyperplane_offset_source_claim}}\par
\noindent\textit{Relationship to source:} derived.\par\smallskip
\reviewlabel{Source attribution and locator}
McLean \cite{McLean2000}, Theorems 3.37--3.38, printed pp. 102--104.
\reviewlabel{Complete source theorem statement (faithful mathematical restatement)}
At the order used here, the trace is a bounded surjection \(H^1\to H^{1/2}\) on a Lipschitz boundary, has a bounded right inverse, and agrees with ordinary restriction for smooth functions.  By duality, \(H^{-1/2}\) boundary distributions act continuously on \(H^1\) traces. These statements provide the trace and Fourier-estimate setting for the flat-hyperplane specialization.
\reviewlabel{Source attribution and locator}
McLean \cite{McLean2000}, Lemma 3.39, printed pp. 104--105.
\reviewlabel{Complete source theorem statement (faithful mathematical restatement)}
A distribution supported on the hyperplane \(\{x_n=0\}\) has a finite representation by tangential distributions tensored with normal derivatives of the Dirac mass, \(\sum_j v_j\otimes D_n^j\delta_0\).  For membership in ambient \(H^s\), each occurring order satisfies \(s+j<-1/2\) and the tangential coefficient has order \(v_j\in H^{s+j+1/2}(\mathbb R^{n-1})\).
\reviewlabel{Source attribution and locator}
McLean \cite{McLean2000}, equation (3.33), printed pp. 104--105.
\reviewlabel{Complete source theorem statement (faithful mathematical restatement)}
For an admissible term,

\[
 \|v\otimes D_n^j\delta_0\|_{H^s(\mathbb R^n)}^2
   = C_{s,j}\|v\|_{H^{s+j+1/2}(\mathbb R^{n-1})}^2,
\]

with the finite positive constant given by the normal-frequency integral. In particular, at \(s=-1,j=0\), every \(v\in H^{-1/2}(\mathbb R^{n-1})\) defines an ambient \(H^{-1}\) source.
\reviewlabel{Exact derived mathematical result formalized}
Let the ambient dimension be one larger than a positive tangential dimension. Fix an isometric affine hyperplane chart, a unit normal, a smooth compactly supported tangential density \(g\), and \(\delta>0\) such that the full swept support cylinder for offsets \(0\leq t\leq\delta\) lies in an open set \(\Omega\).  Then there are bounded quotient-compatible sources \(S_t\) and one \(C\geq0\) such that, on every smooth test \(\phi\),

\[
 S_t(\phi)=\int g(y)\phi(\operatorname{chart}(y)+t\nu)\,dy,
\]

with \(\|S_t\|_{H^{-1}(\Omega)}\leq C\) and \(\|S_t-S_0\|_{H^{-1}(\Omega)}\leq C\sqrt t\).

Translation in the normal variable multiplies the Fourier transform by \(e^{-2\pi i\xi_nt}\).  In the integral underlying (3.33), bounding \(|e^{-2\pi i\xi_nt}-1|\) by \(\min(2,C|\xi_n|t)\) and splitting at \(|\xi_n|=1/t\) gives a squared \(H^{-1}\) bound of order \(t\), hence the displayed \(\sqrt t\) norm bound. That difference estimate is an elementary consequence of (3.33), not a numbered theorem.  The formal source is an abstract bounded extension on the \(H^1_0\) quotient; only its action on smooth tests is given pointwise.
\reviewlabel{Isabelle code}
\begin{lstlisting}
definition mclean_hyperplane_offset_source_claim ::
  "('m::finite itself) \<Rightarrow> ('n::finite itself) \<Rightarrow> bool"
where
  "mclean_hyperplane_offset_source_claim tangential_type ambient_type \<longleftrightarrow>
    (1 \<le> CARD('m) \<and> CARD('n) = Suc (CARD('m))) \<longrightarrow>
    (\<forall>(Omega :: 'n cgu_point set)
      (chart :: 'm cgu_point \<Rightarrow> 'n cgu_point) normal g delta.
      open Omega \<and> cgu_test_function_on UNIV g \<and>
      mclean_flat_chart chart normal \<and> 0 < delta \<and>
      mclean_flat_offset_cylinder_inside Omega chart normal g delta
      \<longrightarrow>
      (\<exists>source_at C. 0 \<le> C \<and>
        (\<forall>t. 0 \<le> t \<and> t \<le> delta \<longrightarrow>
          mclean_source_extends_flat_offset_action
            Omega chart normal t g (source_at t) \<and>
          mclean_source_bound_on Omega
            (source_at t) C \<and>
          mclean_source_bound_on Omega
            (\<lambda>u Du.
              source_at t u Du - source_at 0 u Du)
            (C * sqrt t))))"
\end{lstlisting}
\reviewlabel{English mathematical translation of the Isabelle code}
Let \(m\ge1\) and \(n=m+1\). Let \(\Omega\subset\mathbb R^n\) be open, let \(g\in\mathscr D(\mathbb R^m)\), let \(\nu\in\mathbb R^n\) be a unit vector, and let \(\chi:\mathbb R^m\to\mathbb R^n\) satisfy

\[
 \lVert\chi(x)-\chi(y)\rVert=\lVert x-y\rVert,
 \qquad \nu\cdot(\chi(x)-\chi(y))=0
\]

for all \(x,y\). Let \(\delta>0\), and assume

\[
 \{\chi(x)+t\nu:x\in\overline{\{g\ne0\}},\ 0\le t\le\delta\}\subset\Omega.
\]

Then there exist functionals \(\mathcal F_t\) and a single \(C\ge0\) such that, for every \(0\le t\le\delta\),

\[
 \mathcal F_t(\phi,\nabla\phi)
 =\int_{\mathbb R^m}g(x)\phi(\chi(x)+t\nu)\,dx
 \quad(\phi\in\mathscr D(\Omega)),
\]

\(\mathcal F_t\) is bounded linear on \(\mathcal H_0^1(\Omega)\) with bound \(C\), and \(\mathcal F_t-\mathcal F_0\) is bounded linear there with bound \(C\sqrt t\). The same family and constant work for all \(t\).

\dossierentry{\texttt{\detokenize{MCLEAN-SMOOTH-GENERATOR-AMBIENT-LIFT}}}
\noindent\textit{Formal identifier:} \texttt{\detokenize{mclean_smooth_generator_ambient_lift_claim_v2}}\par
\noindent\textit{Relationship to source:} derived.\par\smallskip
\reviewlabel{Source attribution and locator}
McLean \cite{McLean2000}, Theorem 3.6, printed p. 64.
\reviewlabel{Complete source theorem statement (faithful mathematical restatement)}
A closed set contained in a prescribed open neighborhood admits a smooth cutoff equal to one on the closed set and supported in that neighborhood; for a compact set in an open set the cutoff can be compactly supported there.
\reviewlabel{Source attribution and locator}
McLean \cite{McLean2000}, Theorems 3.37--3.38, printed pp. 102--104.
\reviewlabel{Complete source theorem statement (faithful mathematical restatement)}
At order one on a Lipschitz domain, the trace \(H^1(U)\to H^{1/2}(\partial U)\) is bounded and surjective, has a continuous right inverse, and agrees with classical boundary restriction on smooth functions.
\reviewlabel{Source attribution and locator}
McLean \cite{McLean2000}, Theorem 3.40, printed pp. 105--106.
\reviewlabel{Complete source theorem statement (faithful mathematical restatement)}
On a \(C^{k-1,1}\) domain, \(H^s_0(U)=H^s(U)\) for \(0\leq s\leq1/2\); for \(1/2<s\leq k\), the zero-boundary space is characterized by vanishing of the applicable boundary traces.  At the used specialization \(s=k=1\), the boundary hypothesis is \(C^{0,1}\) (Lipschitz), and \(H^1_0(U)=\ker(\gamma:H^1(U)\to H^{1/2}(\partial U))\).
\reviewlabel{Exact derived mathematical result formalized}
Let \(\Omega_1\) and \(U\subseteq\Omega_1\) be graph-Lipschitz domains in dimension at least two, let \(\gamma\subseteq\Omega_1\cap\partial U\), and fix one globally smooth compactly supported boundary generator \(p\) whose boundary support lies in \(\gamma\).  Then there is a project \(H^1_0(\Omega_1)\) pair \((z,Dz)\) whose quotient trace on \(U\) equals the trace of \((p,\nabla p)\).

Let \(K\) be the compact closure of the nonzero boundary support of \(p\). If \(K\neq\varnothing\), choose a generator-dependent cutoff \(\chi\in C_c^\infty(\Omega_1)\), \(\chi=1\) near \(K\), and set \(z=\chi p\).  The difference \((1-\chi)p\) has zero boundary trace on \(U\), so Theorem 3.40 puts it in \(H^1_0(U)\).  If \(K=\varnothing\), use the zero ambient pair and apply the same kernel theorem directly.  The witness is chosen after the single generator; no linear, bounded, or simultaneous lift operator is asserted.
\reviewlabel{Isabelle code}
\begin{lstlisting}
definition mclean_smooth_generator_ambient_lift_claim_v2 ::
  "'n::finite itself \<Rightarrow> 'c::finite itself \<Rightarrow>
    'h::finite itself \<Rightarrow> bool"
where
  "mclean_smooth_generator_ambient_lift_claim_v2 ambient_type cell_type
      hyperplane_type \<longleftrightarrow>
    2 \<le> CARD('n) \<longrightarrow>
    (\<forall>(M :: ('n, 'c, 'h) cgu_model)
      (Omega1 :: 'n cgu_point set) U gamma p.
      cgu_graph_lipschitz_domain_v5 Omega1 \<and>
      cgu_graph_lipschitz_domain_v5 U \<and>
      U \<subseteq> Omega1 \<and>
      gamma \<subseteq> Omega1 \<inter> frontier U \<and>
      cgu_boundary_test_on M U gamma p
      \<longrightarrow>
      (\<exists>z Dz.
        cgu_h1_zero_pair_on Omega1 z Dz \<and>
        mclean_same_trace_on U p (cgu_classical_gradient p) z Dz))"
\end{lstlisting}
\reviewlabel{English mathematical translation of the Isabelle code}
Let \(n\ge2\). For every geometric record \(M\), every pair of Lipschitz-graph domains \(U\subset\Omega_1\), every

\[
 \gamma\subset\Omega_1\cap\partial U,
\]

and every globally smooth compactly supported \(p\) which is boundary-supported in \((U,\gamma)\) relative to \(M\), there exists \((z,Z)\in\mathcal H_0^1(\Omega_1)\) such that

\[
 R_U\bigl((p,\nabla p),(z,Z)\bigr).
\]

In particular, this one-sided relation itself includes that both restricted pairs are weak \(H^1\) pairs on \(U\).

\dossierentry{\texttt{\detokenize{MCLEAN-SMOOTH-TEST-H1-ZERO}}}
\noindent\textit{Formal identifier:} \texttt{\detokenize{mclean_smooth_test_h1_zero_claim}}\par
\noindent\textit{Relationship to source:} derived.\par\smallskip
\reviewlabel{Source attribution and locator}
McLean \cite{McLean2000}, test class, derivative convention, and first-order Sobolev definition, printed pp. 65, 68, and 73--74.
\reviewlabel{Complete source theorem statement (faithful mathematical restatement)}
\(\mathcal D(U)=C_c^\infty(U)\).  Distributional derivatives use the integration-by-parts sign, and \(W^1_2(U)\) consists of square-integrable functions whose first distributional derivatives are square-integrable.  A smooth compactly supported function and its classical first derivatives are square-integrable, and its distributional and classical derivatives agree; Hence every test belongs to \(W^1_2(U)\) with its classical gradient.
\reviewlabel{Source attribution and locator}
McLean \cite{McLean2000}, zero-boundary closure definition, printed p. 77.
\reviewlabel{Complete source theorem statement (faithful mathematical restatement)}
The zero-boundary Sobolev space is the norm closure of \(\mathcal D(U)\).  In particular every generator is in that closure, witnessed by the constant sequence.
\reviewlabel{Source attribution and locator}
McLean \cite{McLean2000}, Theorem 3.16, printed p. 80.
\reviewlabel{Complete source theorem statement (faithful mathematical restatement)}
On \(\mathbb R^n\), the integer-order weak-derivative and Bessel-potential spaces coincide with equivalent norms.  The theorem is used only for the globally smooth compactly supported representative; no equality of full spaces on an arbitrary open set is asserted.
\reviewlabel{Exact derived mathematical result formalized}
For every open \(U\) and every globally smooth real \(\phi\) with compact closed nonzero support contained in \(U\), the pair \((\phi,\nabla\phi)\) is a project \(H^1_0(U)\) pair.  Smoothness and compact support give measurability and square integrability; integration by parts identifies the classical gradient with the weak gradient; the constant sequence \(\phi\) proves membership in the project's compact-test closure. McLean states the definitions for nonempty open sets.  When \(U\) is empty, the formal support condition forces \(\phi=0\), so the conclusion is an elementary zero-function extension, not a source attribution.
\reviewlabel{Isabelle code}
\begin{lstlisting}
definition mclean_smooth_test_h1_zero_claim ::
  "'n::finite itself \<Rightarrow> bool"
where
  "mclean_smooth_test_h1_zero_claim dimension_type \<longleftrightarrow>
    (\<forall>(U :: 'n cgu_point set) phi.
      open U \<longrightarrow>
      cgu_test_function_on U phi \<longrightarrow>
      cgu_h1_zero_pair_on U phi (cgu_classical_gradient phi))"
\end{lstlisting}
\reviewlabel{English mathematical translation of the Isabelle code}
In every finite dimension, for every \(U\) and every \(\phi\), if \(U\) is open, then if \(\phi\in\mathscr D(U)\), the pair \((\phi,\nabla\phi)\) lies in the zero-boundary closure:

\[
 U\text{ open}\implies
 \bigl(\phi\in\mathscr D(U)\implies(\phi,\nabla\phi)\in\mathcal H_0^1(U)\bigr).
\]

The two implications are right-associated; for a non-open \(U\), the outer implication is true without testing the inner one.

\dossierentry{\texttt{\detokenize{MCLEAN-TRACE-DIRICHLET-GREEN}}}
\noindent\textit{Formal identifier:} \texttt{\detokenize{mclean_trace_dirichlet_green_claim_v3}}\par
\noindent\textit{Relationship to source:} derived.\par\smallskip
\reviewlabel{Source attribution and locator}
McLean \cite{McLean2000}, Definition 3.28, printed pp. 89--90.
\reviewlabel{Complete source theorem statement (faithful mathematical restatement)}
The relevant \(C^{0,1}\) geometry is a domain whose compact boundary is locally a Lipschitz graph after a rigid change of Cartesian coordinates.  The project uses the bounded connected subclass.
\reviewlabel{Source attribution and locator}
McLean \cite{McLean2000}, Theorem 3.37, printed p. 102.
\reviewlabel{Complete source theorem statement (faithful mathematical restatement)}
On a Lipschitz domain, at order one the trace \(\gamma:H^1(U)\to H^{1/2}(\partial U)\) is bounded and surjective and has a continuous right inverse.  Consequently the boundary space is isomorphic, with equivalent norm, to the quotient of \(H^1(U)\) by the trace kernel.
\reviewlabel{Source attribution and locator}
McLean \cite{McLean2000}, Theorem 3.40, printed pp. 105--106.
\reviewlabel{Complete source theorem statement (faithful mathematical restatement)}
On a \(C^{k-1,1}\) domain the zero-boundary space equals all of \(H^s\) for \(0\leq s\leq1/2\), while for \(1/2<s\leq k\) it is characterized by vanishing of the relevant boundary traces.  At \(s=k=1\), used here, the domain is Lipschitz \(C^{0,1}\) and \(H^1_0(U)=\ker\gamma\).
\reviewlabel{Source attribution and locator}
McLean \cite{McLean2000}, Lemma 4.3, printed pp. 116--117.
\reviewlabel{Complete source theorem statement (faithful mathematical restatement)}
Once the right-hand side of a weak elliptic equation is fixed, a weak solution has a unique conormal boundary distribution in \(H^{-1/2}(\partial U)\) satisfying the variational Green formula.  Its norm is bounded in terms of the solution and source norms; the proof uses a continuous trace right inverse.
\reviewlabel{Source attribution and locator}
McLean \cite{McLean2000}, Theorem 4.4, printed p. 118.
\reviewlabel{Complete source theorem statement (faithful mathematical restatement)}
The first Green identity expresses the variational form as the interior source pairing plus the conormal/trace boundary pairing.  Applying it in both orders gives the second Green identity.  For a homogeneous solution, its energy against an arbitrary \(H^1\) lift therefore depends only on the lift's boundary trace.
\reviewlabel{Source attribution and locator}
McLean \cite{McLean2000}, Theorem 4.10, printed pp. 128--130, with the pure Dirichlet specialization on printed p. 131.
\reviewlabel{Complete source theorem statement (faithful mathematical restatement)}
The mixed variational boundary problem on a bounded Lipschitz domain satisfies a Fredholm alternative: in the trivial-kernel branch every admissible datum has a unique solution; otherwise solvability is subject to the stated orthogonality conditions against the adjoint homogeneous kernel.  For pure Dirichlet data the energy space is \(V=H^1_0(U)\) and \(V^*=H^{-1}(U)\), so arbitrary bounded dual sources are allowed.  The formal strict-coercivity premise forces the trivial-kernel branch.
\reviewlabel{Source attribution and locator}
McLean \cite{McLean2000}, equations (4.35)--(4.38), printed p. 145.
\reviewlabel{Complete source theorem statement (faithful mathematical restatement)}
The solution operators furnished by the variational problem are bounded, and the Steklov--Poincaré/Dirichlet-to-Neumann operator obtained from the conormal pairing is a bounded operator between the trace space and its dual.
\reviewlabel{Exact derived mathematical result formalized}
Let \(U\) be a bounded graph-Lipschitz domain in dimension at least two, and let \(a\) define an integrable bilinear energy form which is bounded on project \(H^1(U)\) pairs and strictly coercive on project \(H^1_0(U)\).  Then:

\begin{enumerate}[start=1]
\item every project quotient trace has a weak homogeneous Dirichlet solution, unique modulo zero \(H^1\) distance;
\item one constant bounds the homogeneous solution's energy pairing by the product of the two quotient trace norms, so the pairing is a bounded DN form; and
\item one constant gives, for every bounded real-linear source on \(H^1_0(U)\), a zero-boundary Green solution with an \(H^1\)-norm estimate proportional to the source bound, unique modulo zero distance.
\end{enumerate}

The project trace is the quotient \(H^1(U)/H^1_0(U)\), with norm the infimum of the norms of its lifts.  Theorems 3.37 and 3.40 identify this with the source trace space; Lemma 4.3 and Theorem 4.4 make homogeneous energy depend only on quotient traces; strict coercivity places Theorem 4.10 in its unique pure-Dirichlet branch and yields the solution bounds.  The graph-domain premise is a conservative bounded subclass of McLean's \(C^{0,1}\) geometry.
\reviewlabel{Isabelle code}
\begin{lstlisting}
definition mclean_trace_dirichlet_green_claim_v3 ::
  "'n::finite cgu_point set \<Rightarrow> 'n mclean_coefficient \<Rightarrow> bool"
where
  "mclean_trace_dirichlet_green_claim_v3 U a \<longleftrightarrow>
    (2 \<le> CARD('n) \<and> cgu_graph_lipschitz_domain_v5 U \<and>
      mclean_variational_form_on U a)
    \<longrightarrow>
    ((\<forall>f Df. cgu_h1_pair_on U f Df \<longrightarrow>
      (\<exists>u Du. mclean_dirichlet_solution_for U a f Df u Du \<and>
        (\<forall>v Dv. mclean_dirichlet_solution_for U a f Df v Dv
          \<longrightarrow> cgu_h1_squared_distance U u Du v Dv = 0))) \<and>
     (\<exists>C. 0 \<le> C \<and>
      (\<forall>f Df u Du h Dh.
        mclean_dirichlet_solution_for U a f Df u Du \<and>
        cgu_h1_pair_on U h Dh
        \<longrightarrow>
        mclean_energy_integrable U a Du Dh \<and>
        abs (mclean_energy U a Du Dh) \<le>
          C * mclean_trace_norm U f Df * mclean_trace_norm U h Dh)) \<and>
     (\<exists>C. 0 \<le> C \<and>
      (\<forall>F K. mclean_source_bound_on U F K \<longrightarrow>
        (\<exists>u Du. mclean_green_solution_for U a F u Du \<and>
          mclean_h1_norm U u Du \<le> C * K \<and>
          (\<forall>v Dv. mclean_green_solution_for U a F v Dv
            \<longrightarrow> cgu_h1_squared_distance U u Du v Dv = 0)))))"
\end{lstlisting}
\reviewlabel{English mathematical translation of the Isabelle code}
Let \(n\ge2\), let \(U\) be a Lipschitz-graph domain, and suppose \(B_{U,a}\) is stable. Then all three conclusions hold.

\begin{enumerate}[start=1]
\item For every \((f,F)\in\mathcal H^1(U)\), there is a variational Dirichlet solution \((u,G)\), and every other such solution \((v,H)\) satisfies
\end{enumerate}

\[
   d_U^2((u,G),(v,H))=0.
   \]

\begin{enumerate}[start=2]
\item There exists \(C_1\ge0\) such that for every Dirichlet solution \((u,G)\) with boundary representative \((f,F)\), and every \((h,H)\in\mathcal H^1(U)\), the energy density \(aG\cdot H\) is integrable and
\end{enumerate}

\[
   |B_{U,a}(G,H)|\le C_1\,q_U(f,F)\,q_U(h,H).
   \]

\begin{enumerate}[start=3]
\item There exists an independently quantified \(C_2\ge0\) such that every functional \(\mathcal F\) bounded linear with bound \(K\) has a source solution \((u,G)\) with
\end{enumerate}

\[
   \lVert(u,G)\rVert_U\le C_2K,
   \]

and every other source solution is at squared pair-distance zero from it.

\clearpage

\section{Independent back-translation B}\label{app:translation-b}
\begingroup
\small
\emergencystretch=3em
\sloppy
\dossiersection{Definitions}

\dossierentry{\texttt{\detokenize{cgu_point}}}
\reviewlabel{Isabelle code}
\begin{lstlisting}
type_synonym 'n cgu_point = "real ^ 'n"
\end{lstlisting}
\reviewlabel{English mathematical translation}
\noindent\textit{Mathematical role:} Euclidean coordinate space \(V_I\)\par\smallskip
For every type \texttt{'n}, \texttt{'n\allowbreak\ cgu\_\allowbreak{}point} is exactly the vector type \texttt{real\allowbreak\ \textasciicircum{}\allowbreak\ 'n}. There are no term arguments, quantifiers, hypotheses, or conclusion.

\dossierentry{\texttt{\detokenize{cgu_monomial}}}
\reviewlabel{Isabelle code}
\begin{lstlisting}
type_synonym 'n cgu_monomial = "'n \<Rightarrow>\<^sub>0 nat"
\end{lstlisting}
\reviewlabel{English mathematical translation}
\noindent\textit{Mathematical role:} finitely supported multi-index type\par\smallskip
For every type \texttt{'n}, \texttt{'n\allowbreak\ cgu\_\allowbreak{}monomial} is exactly the finite-support mapping type \texttt{'n\allowbreak\ \textbackslash{}<Rightarrow>\textbackslash{}<\textasciicircum{}sub>0\allowbreak\ nat}. There are no term arguments or logical conditions.

\dossierentry{\texttt{\detokenize{cgu_polynomial}}}
\reviewlabel{Isabelle code}
\begin{lstlisting}
type_synonym 'n cgu_polynomial = "'n cgu_monomial \<Rightarrow>\<^sub>0 real"
\end{lstlisting}
\reviewlabel{English mathematical translation}
\noindent\textit{Mathematical role:} finitely supported scalar polynomial-representation type\par\smallskip
For every type \texttt{'n}, \texttt{'n\allowbreak\ cgu\_\allowbreak{}polynomial} is exactly the finite-support mapping type from \texttt{'n\allowbreak\ cgu\_\allowbreak{}monomial} to \texttt{real}. There are no term arguments or logical conditions.

\dossierentry{\texttt{\detokenize{cgu_matrix_polynomial}}}
\reviewlabel{Isabelle code}
\begin{lstlisting}
type_synonym 'n cgu_matrix_polynomial = "'n cgu_polynomial ^ 'n ^ 'n"
\end{lstlisting}
\reviewlabel{English mathematical translation}
\noindent\textit{Mathematical role:} polynomial-matrix representation type\par\smallskip
For every type \texttt{'n}, \texttt{'n\allowbreak\ cgu\_\allowbreak{}matrix\_\allowbreak{}polynomial} is exactly \texttt{'n\allowbreak\ cgu\_\allowbreak{}polynomial\allowbreak\ \textasciicircum{}\allowbreak\ 'n\allowbreak\ \textasciicircum{}\allowbreak\ 'n}. There are no term arguments or logical conditions.

\dossierentry{\texttt{\detokenize{cgu_monomial_total_degree}}}
\reviewlabel{Isabelle code}
\begin{lstlisting}
definition cgu_monomial_total_degree :: "'n cgu_monomial \<Rightarrow> nat"
where
  "cgu_monomial_total_degree m =
    (\<Sum>i\<in>Poly_Mapping.keys m. Poly_Mapping.lookup m i)"
\end{lstlisting}
\reviewlabel{English mathematical translation}
\noindent\textit{Mathematical role:} total size of a multi-index\par\smallskip
Typed arguments/result: \texttt{m\allowbreak\ ::\allowbreak\ 'n\allowbreak\ cgu\_\allowbreak{}monomial}; result \texttt{nat}. \texttt{cgu\_\allowbreak{}monomial\_\allowbreak{}total\_\allowbreak{}degree\allowbreak\ m} equals the finite sum, over every \texttt{i} in \texttt{Poly\_\allowbreak{}Mapping.\allowbreak{}keys\allowbreak\ m}, of \texttt{Poly\_\allowbreak{}Mapping.\allowbreak{}lookup\allowbreak\ m\allowbreak\ i}. There are no hypotheses.

\dossierentry{\texttt{\detokenize{cgu_poly_degree_le}}}
\reviewlabel{Isabelle code}
\begin{lstlisting}
definition cgu_poly_degree_le :: "nat \<Rightarrow> 'n cgu_polynomial \<Rightarrow> bool"
where
  "cgu_poly_degree_le d p \<longleftrightarrow>
    (\<forall>m\<in>Poly_Mapping.keys p. cgu_monomial_total_degree m \<le> d)"
\end{lstlisting}
\reviewlabel{English mathematical translation}
\noindent\textit{Mathematical role:} stored scalar degree predicate\par\smallskip
Typed arguments/result: \texttt{d\allowbreak\ ::\allowbreak\ nat}, \texttt{p\allowbreak\ ::\allowbreak\ 'n\allowbreak\ cgu\_\allowbreak{}polynomial}; result \texttt{bool}. \texttt{cgu\_\allowbreak{}poly\_\allowbreak{}degree\_\allowbreak{}le\allowbreak\ d\allowbreak\ p} holds exactly when, for every \texttt{m} restricted by \texttt{m\allowbreak\ \textbackslash{}<in>\allowbreak\ Poly\_\allowbreak{}Mapping.\allowbreak{}keys\allowbreak\ p}, \texttt{cgu\_\allowbreak{}monomial\_\allowbreak{}total\_\allowbreak{}degree\allowbreak\ m\allowbreak\ \textbackslash{}<le>\allowbreak\ d}.

\dossierentry{\texttt{\detokenize{cgu_monomial_eval}}}
\reviewlabel{Isabelle code}
\begin{lstlisting}
definition cgu_monomial_eval ::
  "'n cgu_monomial \<Rightarrow> 'n::finite cgu_point \<Rightarrow> real"
where
  "cgu_monomial_eval m x =
    (\<Prod>i\<in>Poly_Mapping.keys m. x $ i ^ Poly_Mapping.lookup m i)"
\end{lstlisting}
\reviewlabel{English mathematical translation}
\noindent\textit{Mathematical role:} monomial \(x^\alpha\)\par\smallskip
Typed arguments/result: \texttt{m\allowbreak\ ::\allowbreak\ 'n\allowbreak\ cgu\_\allowbreak{}monomial}, \texttt{x\allowbreak\ ::\allowbreak\ 'n\allowbreak\ cgu\_\allowbreak{}point}, with \texttt{'n::finite}; result \texttt{real}. \texttt{cgu\_\allowbreak{}monomial\_\allowbreak{}eval\allowbreak\ m\allowbreak\ x} equals the finite product, over every \texttt{i\allowbreak\ \textbackslash{}<in>\allowbreak\ Poly\_\allowbreak{}Mapping.\allowbreak{}keys\allowbreak\ m}, of \texttt{(x\allowbreak\ \$\allowbreak\ i)} raised to \texttt{Poly\_\allowbreak{}Mapping.\allowbreak{}lookup\allowbreak\ m\allowbreak\ i}. There are no hypotheses.

\dossierentry{\texttt{\detokenize{cgu_poly_eval}}}
\reviewlabel{Isabelle code}
\begin{lstlisting}
definition cgu_poly_eval ::
  "'n cgu_polynomial \<Rightarrow> 'n::finite cgu_point \<Rightarrow> real"
where
  "cgu_poly_eval p x =
    (\<Sum>m\<in>Poly_Mapping.keys p.
      Poly_Mapping.lookup p m * cgu_monomial_eval m x)"
\end{lstlisting}
\reviewlabel{English mathematical translation}
\noindent\textit{Mathematical role:} scalar polynomial evaluation \(\widehat p(x)\)\par\smallskip
Typed arguments/result: \texttt{p\allowbreak\ ::\allowbreak\ 'n\allowbreak\ cgu\_\allowbreak{}polynomial}, \texttt{x\allowbreak\ ::\allowbreak\ 'n\allowbreak\ cgu\_\allowbreak{}point}, with \texttt{'n::finite}; result \texttt{real}. It equals the finite sum over \texttt{m\allowbreak\ \textbackslash{}<in>\allowbreak\ Poly\_\allowbreak{}Mapping.\allowbreak{}keys\allowbreak\ p} of \texttt{Poly\_\allowbreak{}Mapping.\allowbreak{}lookup\allowbreak\ p\allowbreak\ m\allowbreak\ *\allowbreak\ cgu\_\allowbreak{}monomial\_\allowbreak{}eval\allowbreak\ m\allowbreak\ x}.

\dossierentry{\texttt{\detokenize{cgu_matrix_poly_eval}}}
\reviewlabel{Isabelle code}
\begin{lstlisting}
definition cgu_matrix_poly_eval ::
  "'n::finite cgu_matrix_polynomial \<Rightarrow> 'n cgu_point \<Rightarrow>
    real ^ 'n ^ 'n"
where
  "cgu_matrix_poly_eval P x = (\<chi> i j. cgu_poly_eval (P $ i $ j) x)"
\end{lstlisting}
\reviewlabel{English mathematical translation}
\noindent\textit{Mathematical role:} matrix evaluation \(\widehat P(x)\)\par\smallskip
Typed arguments/result: \texttt{P\allowbreak\ ::\allowbreak\ 'n\allowbreak\ cgu\_\allowbreak{}matrix\_\allowbreak{}polynomial}, \texttt{x\allowbreak\ ::\allowbreak\ 'n\allowbreak\ cgu\_\allowbreak{}point}, with \texttt{'n::finite}; result \texttt{real\allowbreak\ \textasciicircum{}\allowbreak\ 'n\allowbreak\ \textasciicircum{}\allowbreak\ 'n}. The result is the indexed vector/matrix \texttt{\textbackslash{}<chi>\allowbreak\ i\allowbreak\ j.\allowbreak{}\allowbreak\ cgu\_\allowbreak{}poly\_\allowbreak{}eval\allowbreak\ (P\allowbreak\ \$\allowbreak\ i\allowbreak\ \$\allowbreak\ j)\allowbreak\ x}; thus every indices \texttt{i,j\allowbreak\ ::\allowbreak\ 'n} select that entry.

\dossierentry{\texttt{\detokenize{cgu_matrix_poly_degree_le}}}
\reviewlabel{Isabelle code}
\begin{lstlisting}
definition cgu_matrix_poly_degree_le ::
  "nat \<Rightarrow> 'n::finite cgu_matrix_polynomial \<Rightarrow> bool"
where
  "cgu_matrix_poly_degree_le d P \<longleftrightarrow>
    (\<forall>i j. cgu_poly_degree_le d (P $ i $ j))"
\end{lstlisting}
\reviewlabel{English mathematical translation}
\noindent\textit{Mathematical role:} componentwise coefficient-representation degree predicate\par\smallskip
Typed arguments/result: \texttt{d\allowbreak\ ::\allowbreak\ nat}, \texttt{P\allowbreak\ ::\allowbreak\ 'n\allowbreak\ cgu\_\allowbreak{}matrix\_\allowbreak{}polynomial}, with \texttt{'n::finite}; result \texttt{bool}. It holds exactly when for every unrestricted \texttt{i,j\allowbreak\ ::\allowbreak\ 'n}, \texttt{cgu\_\allowbreak{}poly\_\allowbreak{}degree\_\allowbreak{}le\allowbreak\ d\allowbreak\ (P\allowbreak\ \$\allowbreak\ i\allowbreak\ \$\allowbreak\ j)} holds.

\dossierentry{\texttt{\detokenize{cgu_symmetric_matrix_polynomial}}}
\reviewlabel{Isabelle code}
\begin{lstlisting}
definition cgu_symmetric_matrix_polynomial ::
  "'n::finite cgu_matrix_polynomial \<Rightarrow> bool"
where
  "cgu_symmetric_matrix_polynomial P \<longleftrightarrow>
    (\<forall>i j. P $ i $ j = P $ j $ i)"
\end{lstlisting}
\reviewlabel{English mathematical translation}
\noindent\textit{Mathematical role:} symbolic matrix symmetry\par\smallskip
Typed argument/result: \texttt{P\allowbreak\ ::\allowbreak\ 'n\allowbreak\ cgu\_\allowbreak{}matrix\_\allowbreak{}polynomial}, with \texttt{'n::finite}; result \texttt{bool}. It holds exactly when for all \texttt{i,j\allowbreak\ ::\allowbreak\ 'n}, \texttt{P\allowbreak\ \$\allowbreak\ i\allowbreak\ \$\allowbreak\ j\allowbreak\ =\allowbreak\ P\allowbreak\ \$\allowbreak\ j\allowbreak\ \$\allowbreak\ i}.

\dossierentry{\texttt{\detokenize{cgu_matrix_bilinear}}}
\reviewlabel{Isabelle code}
\begin{lstlisting}
definition cgu_matrix_bilinear ::
  "(real ^ 'n::finite ^ 'n) \<Rightarrow> real ^ 'n \<Rightarrow>
    real ^ 'n \<Rightarrow> real"
where
  "cgu_matrix_bilinear A x y = inner x (A *v y)"
\end{lstlisting}
\reviewlabel{English mathematical translation}
\noindent\textit{Mathematical role:} bilinear expression \(\mathfrak b_A\)\par\smallskip
Typed arguments/result: \texttt{A\allowbreak\ ::\allowbreak\ real\allowbreak\ \textasciicircum{}\allowbreak\ 'n\allowbreak\ \textasciicircum{}\allowbreak\ 'n}, \texttt{x,y\allowbreak\ ::\allowbreak\ real\allowbreak\ \textasciicircum{}\allowbreak\ 'n}, with \texttt{'n::finite}; result \texttt{real}. \texttt{cgu\_\allowbreak{}matrix\_\allowbreak{}bilinear\allowbreak\ A\allowbreak\ x\allowbreak\ y\allowbreak\ =\allowbreak\ inner\allowbreak\ x\allowbreak\ (A\allowbreak\ *v\allowbreak\ y)}.

\dossierentry{\texttt{\detokenize{cgu_symmetric_positive_definite_matrix}}}
\reviewlabel{Isabelle code}
\begin{lstlisting}
definition cgu_symmetric_positive_definite_matrix ::
  "(real ^ 'n::finite ^ 'n) \<Rightarrow> bool"
where
  "cgu_symmetric_positive_definite_matrix A \<longleftrightarrow>
    transpose A = A \<and>
    (\<forall>x. x \<noteq> 0 \<longrightarrow> 0 < cgu_matrix_bilinear A x x)"
\end{lstlisting}
\reviewlabel{English mathematical translation}
\noindent\textit{Mathematical role:} symmetric positive definiteness\par\smallskip
Typed argument/result: \texttt{A\allowbreak\ ::\allowbreak\ real\allowbreak\ \textasciicircum{}\allowbreak\ 'n\allowbreak\ \textasciicircum{}\allowbreak\ 'n}, with \texttt{'n::finite}; result \texttt{bool}. It holds exactly when both \texttt{transpose\allowbreak\ A\allowbreak\ =\allowbreak\ A} and, for every \texttt{x\allowbreak\ ::\allowbreak\ real\allowbreak\ \textasciicircum{}\allowbreak\ 'n}, the hypothesis \texttt{x\allowbreak\ \textbackslash{}<noteq>\allowbreak\ 0} implies the conclusion \texttt{0\allowbreak\ <\allowbreak\ cgu\_\allowbreak{}matrix\_\allowbreak{}bilinear\allowbreak\ A\allowbreak\ x\allowbreak\ x}.

\dossierentry{\texttt{\detokenize{cgu_uniformly_positive_definite_on}}}
\reviewlabel{Isabelle code}
\begin{lstlisting}
definition cgu_uniformly_positive_definite_on ::
  "'n::finite cgu_point set \<Rightarrow> 'n cgu_matrix_polynomial \<Rightarrow> bool"
where
  "cgu_uniformly_positive_definite_on D P \<longleftrightarrow>
    (\<exists>k>0. \<forall>x\<in>D. \<forall>xi.
      k * norm xi ^ 2 \<le>
        cgu_matrix_bilinear (cgu_matrix_poly_eval P x) xi xi)"
\end{lstlisting}
\reviewlabel{English mathematical translation}
\noindent\textit{Mathematical role:} uniform ellipticity on a set\par\smallskip
Typed arguments/result: \texttt{D\allowbreak\ ::\allowbreak\ 'n\allowbreak\ cgu\_\allowbreak{}point\allowbreak\ set}, \texttt{P\allowbreak\ ::\allowbreak\ 'n\allowbreak\ cgu\_\allowbreak{}matrix\_\allowbreak{}polynomial}, with \texttt{'n::finite}; result \texttt{bool}. It holds exactly when there exists \texttt{k\allowbreak\ ::\allowbreak\ real} restricted by \texttt{k\allowbreak\ >\allowbreak\ 0} such that for every \texttt{x\allowbreak\ \textbackslash{}<in>\allowbreak\ D} and every unrestricted \texttt{xi\allowbreak\ ::\allowbreak\ 'n\allowbreak\ cgu\_\allowbreak{}point}, \texttt{k\allowbreak\ *\allowbreak\ norm\allowbreak\ xi\allowbreak\ \textasciicircum{}\allowbreak\ 2\allowbreak\ \textbackslash{}<le>\allowbreak\ cgu\_\allowbreak{}matrix\_\allowbreak{}bilinear\allowbreak\ (cgu\_\allowbreak{}matrix\_\allowbreak{}poly\_\allowbreak{}eval\allowbreak\ P\allowbreak\ x)\allowbreak\ xi\allowbreak\ xi}.

\dossierentry{\texttt{\detokenize{cgu_model}}}
\reviewlabel{Isabelle code}
\begin{lstlisting}
record ('n, 'c, 'h) cgu_model =
  cgu_domain :: "'n cgu_point set"
  cgu_measured :: "'n cgu_point set"
  cgu_cell :: "'c \<Rightarrow> 'n cgu_point set"
  cgu_degree :: "'c \<Rightarrow> nat"
  cgu_order :: "nat \<Rightarrow> 'c"
  cgu_faces :: "'c \<Rightarrow> 'h set"
  cgu_normal :: "'c \<Rightarrow> 'h \<Rightarrow> 'n cgu_point"
  cgu_offset :: "'c \<Rightarrow> 'h \<Rightarrow> real"
  cgu_face_patch :: "'c \<Rightarrow> 'h \<Rightarrow> 'n cgu_point set"
  cgu_witness_faces :: "'c \<Rightarrow> 'h \<Rightarrow> 'h set"
  cgu_triple_region :: "'c \<Rightarrow> 'h \<Rightarrow> 'h \<Rightarrow> 'h \<Rightarrow>
    'n cgu_point set"
  cgu_step_interface_pieces :: "'c \<Rightarrow> 'h set"
  cgu_step_interface_patch :: "'c \<Rightarrow> 'h \<Rightarrow> 'n cgu_point set"
  cgu_step_interface_normal :: "'c \<Rightarrow> 'h \<Rightarrow> 'n cgu_point"
  cgu_step_interface_offset :: "'c \<Rightarrow> 'h \<Rightarrow> real"
  cgu_step_outer_patch :: "'c \<Rightarrow> 'n cgu_point set"
  cgu_step_adjacent_cell :: "'c \<Rightarrow> 'c"
  cgu_step_center :: "'c \<Rightarrow> 'n cgu_point"
  cgu_step_interior_normal :: "'c \<Rightarrow> 'n cgu_point"
  cgu_step_radius :: "'c \<Rightarrow> real"
\end{lstlisting}
\reviewlabel{English mathematical translation}
\noindent\textit{Mathematical role:} geometric datum \(M\), with twenty fields\par\smallskip
Type parameters: \texttt{('n,'c,'h)}. The record has exactly these typed selectors: \texttt{cgu\_\allowbreak{}domain,\allowbreak\ cgu\_\allowbreak{}measured\allowbreak\ ::\allowbreak\ 'n\allowbreak\ cgu\_\allowbreak{}point\allowbreak\ set}; \texttt{cgu\_\allowbreak{}cell\allowbreak\ ::\allowbreak\ 'c\allowbreak\ \textbackslash{}<Rightarrow>\allowbreak\ 'n\allowbreak\ cgu\_\allowbreak{}point\allowbreak\ set}; \texttt{cgu\_\allowbreak{}degree\allowbreak\ ::\allowbreak\ 'c\allowbreak\ \textbackslash{}<Rightarrow>\allowbreak\ nat}; \texttt{cgu\_\allowbreak{}order\allowbreak\ ::\allowbreak\ nat\allowbreak\ \textbackslash{}<Rightarrow>\allowbreak\ 'c}; \texttt{cgu\_\allowbreak{}faces,\allowbreak\ cgu\_\allowbreak{}step\_\allowbreak{}interface\_\allowbreak{}pieces\allowbreak\ ::\allowbreak\ 'c\allowbreak\ \textbackslash{}<Rightarrow>\allowbreak\ 'h\allowbreak\ set}; \texttt{cgu\_\allowbreak{}normal,\allowbreak\ cgu\_\allowbreak{}step\_\allowbreak{}interface\_\allowbreak{}normal\allowbreak\ ::\allowbreak\ 'c\allowbreak\ \textbackslash{}<Rightarrow>\allowbreak\ 'h\allowbreak\ \textbackslash{}<Rightarrow>\allowbreak\ 'n\allowbreak\ cgu\_\allowbreak{}point}; \texttt{cgu\_\allowbreak{}offset,\allowbreak\ cgu\_\allowbreak{}step\_\allowbreak{}interface\_\allowbreak{}offset\allowbreak\ ::\allowbreak\ 'c\allowbreak\ \textbackslash{}<Rightarrow>\allowbreak\ 'h\allowbreak\ \textbackslash{}<Rightarrow>\allowbreak\ real}; \texttt{cgu\_\allowbreak{}face\_\allowbreak{}patch,\allowbreak\ cgu\_\allowbreak{}step\_\allowbreak{}interface\_\allowbreak{}patch\allowbreak\ ::\allowbreak\ 'c\allowbreak\ \textbackslash{}<Rightarrow>\allowbreak\ 'h\allowbreak\ \textbackslash{}<Rightarrow>\allowbreak\ 'n\allowbreak\ cgu\_\allowbreak{}point\allowbreak\ set}; \texttt{cgu\_\allowbreak{}witness\_\allowbreak{}faces\allowbreak\ ::\allowbreak\ 'c\allowbreak\ \textbackslash{}<Rightarrow>\allowbreak\ 'h\allowbreak\ \textbackslash{}<Rightarrow>\allowbreak\ 'h\allowbreak\ set}; \texttt{cgu\_\allowbreak{}triple\_\allowbreak{}region\allowbreak\ ::\allowbreak\ 'c\allowbreak\ \textbackslash{}<Rightarrow>\allowbreak\ 'h\allowbreak\ \textbackslash{}<Rightarrow>\allowbreak\ 'h\allowbreak\ \textbackslash{}<Rightarrow>\allowbreak\ 'h\allowbreak\ \textbackslash{}<Rightarrow>\allowbreak\ 'n\allowbreak\ cgu\_\allowbreak{}point\allowbreak\ set}; \texttt{cgu\_\allowbreak{}step\_\allowbreak{}outer\_\allowbreak{}patch\allowbreak\ ::\allowbreak\ 'c\allowbreak\ \textbackslash{}<Rightarrow>\allowbreak\ 'n\allowbreak\ cgu\_\allowbreak{}point\allowbreak\ set}; \texttt{cgu\_\allowbreak{}step\_\allowbreak{}adjacent\_\allowbreak{}cell\allowbreak\ ::\allowbreak\ 'c\allowbreak\ \textbackslash{}<Rightarrow>\allowbreak\ 'c}; \texttt{cgu\_\allowbreak{}step\_\allowbreak{}center,\allowbreak\ cgu\_\allowbreak{}step\_\allowbreak{}interior\_\allowbreak{}normal\allowbreak\ ::\allowbreak\ 'c\allowbreak\ \textbackslash{}<Rightarrow>\allowbreak\ 'n\allowbreak\ cgu\_\allowbreak{}point}; and \texttt{cgu\_\allowbreak{}step\_\allowbreak{}radius\allowbreak\ ::\allowbreak\ 'c\allowbreak\ \textbackslash{}<Rightarrow>\allowbreak\ real}. This is a record composition, not a proposition.

\dossierentry{\texttt{\detokenize{cgu_piecewise_polynomial_class}}}
\reviewlabel{Isabelle code}
\begin{lstlisting}
definition cgu_piecewise_polynomial_class ::
  "('n::finite, 'c::finite, 'h) cgu_model \<Rightarrow>
    ('c \<Rightarrow> 'n cgu_matrix_polynomial) \<Rightarrow> bool"
where
  "cgu_piecewise_polynomial_class M P \<longleftrightarrow>
    (\<forall>c.
      cgu_matrix_poly_degree_le (cgu_degree M c) (P c) \<and>
      cgu_symmetric_matrix_polynomial (P c) \<and>
      cgu_uniformly_positive_definite_on (cgu_cell M c) (P c))"
\end{lstlisting}
\reviewlabel{English mathematical translation}
\noindent\textit{Mathematical role:} admissible cellwise coefficient family\par\smallskip
Typed arguments/result: \texttt{M\allowbreak\ ::\allowbreak\ ('n,'c,'h)\allowbreak\ cgu\_\allowbreak{}model}, \texttt{P\allowbreak\ ::\allowbreak\ 'c\allowbreak\ \textbackslash{}<Rightarrow>\allowbreak\ 'n\allowbreak\ cgu\_\allowbreak{}matrix\_\allowbreak{}polynomial}, with \texttt{'n,'c::finite}; result \texttt{bool}. It holds exactly when, for every \texttt{c\allowbreak\ ::\allowbreak\ 'c}, all three conjuncts hold: \texttt{cgu\_\allowbreak{}matrix\_\allowbreak{}poly\_\allowbreak{}degree\_\allowbreak{}le\allowbreak\ (cgu\_\allowbreak{}degree\allowbreak\ M\allowbreak\ c)\allowbreak\ (P\allowbreak\ c)}, \texttt{cgu\_\allowbreak{}symmetric\_\allowbreak{}matrix\_\allowbreak{}polynomial\allowbreak\ (P\allowbreak\ c)}, and \texttt{cgu\_\allowbreak{}uniformly\_\allowbreak{}positive\_\allowbreak{}definite\_\allowbreak{}on\allowbreak\ (cgu\_\allowbreak{}cell\allowbreak\ M\allowbreak\ c)\allowbreak\ (P\allowbreak\ c)}.

\dossierentry{\texttt{\detokenize{cgu_affine_hyperplane}}}
\reviewlabel{Isabelle code}
\begin{lstlisting}
definition cgu_affine_hyperplane ::
  "'n::finite cgu_point \<Rightarrow> real \<Rightarrow> 'n cgu_point set"
where
  "cgu_affine_hyperplane normal offset = {x. inner normal x = offset}"
\end{lstlisting}
\reviewlabel{English mathematical translation}
\noindent\textit{Mathematical role:} affine hyperplane \(H(\nu,b)\)\par\smallskip
Typed arguments/result: \texttt{normal\allowbreak\ ::\allowbreak\ 'n\allowbreak\ cgu\_\allowbreak{}point}, \texttt{offset\allowbreak\ ::\allowbreak\ real}, with \texttt{'n::finite}; result \texttt{'n\allowbreak\ cgu\_\allowbreak{}point\allowbreak\ set}. It is \texttt{\{x.\allowbreak{}\allowbreak\ inner\allowbreak\ normal\allowbreak\ x\allowbreak\ =\allowbreak\ offset\}}, with \texttt{x\allowbreak\ ::\allowbreak\ 'n\allowbreak\ cgu\_\allowbreak{}point}.

\dossierentry{\texttt{\detokenize{cgu_coordinate_hyperplane_projection_v5}}}
\reviewlabel{Isabelle code}
\begin{lstlisting}
definition cgu_coordinate_hyperplane_projection_v5 ::
  "'i::finite \<Rightarrow> 'i cgu_point \<Rightarrow> 'i cgu_point"
where
  "cgu_coordinate_hyperplane_projection_v5 i y =
    (\<chi> j. if j = i then 0 else y $ j)"
\end{lstlisting}
\reviewlabel{English mathematical translation}
\noindent\textit{Mathematical role:} coordinate-deletion map \(\pi_i\)\par\smallskip
Typed arguments/result: \texttt{i\allowbreak\ ::\allowbreak\ 'i}, \texttt{y\allowbreak\ ::\allowbreak\ 'i\allowbreak\ cgu\_\allowbreak{}point}, with \texttt{'i::finite}; result \texttt{'i\allowbreak\ cgu\_\allowbreak{}point}. The result has coordinate \texttt{j\allowbreak\ ::\allowbreak\ 'i} equal to \texttt{0} if \texttt{j=i}, and otherwise equal to \texttt{y\allowbreak\ \$\allowbreak\ j}.

\dossierentry{\texttt{\detokenize{cgu_local_lipschitz_boundary_at}}}
\reviewlabel{Isabelle code}
\begin{lstlisting}
definition cgu_local_lipschitz_boundary_at ::
  "'n::finite cgu_point set \<Rightarrow> 'n cgu_point \<Rightarrow> bool"
where
  "cgu_local_lipschitz_boundary_at D x \<longleftrightarrow>
    (\<exists>(U :: 'n cgu_point set)
       (f :: 'n cgu_point \<Rightarrow> 'n cgu_point)
       (g :: 'n cgu_point \<Rightarrow> 'n cgu_point) (i :: 'n) L M.
      open U \<and> x \<in> U \<and> 0 < L \<and> 0 < M \<and>
      homeomorphism U (f ` U) f g \<and>
      lipschitz_on L U f \<and> lipschitz_on M (f ` U) g \<and>
      f ` (D \<inter> U) = (f ` U) \<inter> {y. 0 < y $ i})"
\end{lstlisting}
\reviewlabel{English mathematical translation}
\noindent\textit{Mathematical role:} bi-Lipschitz half-space point\par\smallskip
Typed arguments/result: \texttt{D\allowbreak\ ::\allowbreak\ 'n\allowbreak\ cgu\_\allowbreak{}point\allowbreak\ set}, \texttt{x\allowbreak\ ::\allowbreak\ 'n\allowbreak\ cgu\_\allowbreak{}point}, with \texttt{'n::finite}; result \texttt{bool}. It holds exactly when there exist \texttt{U\allowbreak\ ::\allowbreak\ 'n\allowbreak\ cgu\_\allowbreak{}point\allowbreak\ set}, \texttt{f,g\allowbreak\ ::\allowbreak\ 'n\allowbreak\ cgu\_\allowbreak{}point\allowbreak\ \textbackslash{}<Rightarrow>\allowbreak\ 'n\allowbreak\ cgu\_\allowbreak{}point}, \texttt{i\allowbreak\ ::\allowbreak\ 'n}, and \texttt{L,M\allowbreak\ ::\allowbreak\ real} such that: \texttt{open\allowbreak\ U}; \texttt{x\allowbreak\ \textbackslash{}<in>\allowbreak\ U}; \texttt{0<L}; \texttt{0<M}; \texttt{homeomorphism\allowbreak\ U\allowbreak\ (f\allowbreak\ `\allowbreak\ U)\allowbreak\ f\allowbreak\ g}; \texttt{lipschitz\_\allowbreak{}on\allowbreak\ L\allowbreak\ U\allowbreak\ f}; \texttt{lipschitz\_\allowbreak{}on\allowbreak\ M\allowbreak\ (f\allowbreak\ `\allowbreak\ U)\allowbreak\ g}; and \texttt{f\allowbreak\ `\allowbreak\ (D\allowbreak\ \textbackslash{}<inter>\allowbreak\ U)\allowbreak\ =\allowbreak\ (f\allowbreak\ `\allowbreak\ U)\allowbreak\ \textbackslash{}<inter>\allowbreak\ \{y.\allowbreak{}\allowbreak\ 0\allowbreak\ <\allowbreak\ y\allowbreak\ \$\allowbreak\ i\}}.

\dossierentry{\texttt{\detokenize{cgu_lipschitz_domain}}}
\reviewlabel{Isabelle code}
\begin{lstlisting}
definition cgu_lipschitz_domain ::
  "'n::finite cgu_point set \<Rightarrow> bool"
where
  "cgu_lipschitz_domain D \<longleftrightarrow>
    D \<noteq> {} \<and> open D \<and> connected D \<and> bounded D \<and>
    (\<forall>x\<in>frontier D. cgu_local_lipschitz_boundary_at D x)"
\end{lstlisting}
\reviewlabel{English mathematical translation}
\noindent\textit{Mathematical role:} admissible bounded domain\par\smallskip
Typed argument/result: \texttt{D\allowbreak\ ::\allowbreak\ 'n\allowbreak\ cgu\_\allowbreak{}point\allowbreak\ set}, with \texttt{'n::finite}; result \texttt{bool}. It holds exactly when \texttt{D\allowbreak\ \textbackslash{}<noteq>\allowbreak\ \{\}}, \texttt{open\allowbreak\ D}, \texttt{connected\allowbreak\ D}, and \texttt{bounded\allowbreak\ D}, and for every \texttt{x\allowbreak\ \textbackslash{}<in>\allowbreak\ frontier\allowbreak\ D}, \texttt{cgu\_\allowbreak{}local\_\allowbreak{}lipschitz\_\allowbreak{}boundary\_\allowbreak{}at\allowbreak\ D\allowbreak\ x}.

\dossierentry{\texttt{\detokenize{cgu_rigid_motion_v5}}}
\reviewlabel{Isabelle code}
\begin{lstlisting}
definition cgu_rigid_motion_v5 ::
  "('n::finite cgu_point \<Rightarrow> 'n cgu_point) \<Rightarrow> bool"
where
  "cgu_rigid_motion_v5 R \<longleftrightarrow>
    (\<exists>A b.
      linear A \<and>
      surj A \<and>
      (\<forall>u v. inner (A u) (A v) = inner u v) \<and>
      R = (\<lambda>x. A x + b))"
\end{lstlisting}
\reviewlabel{English mathematical translation}
\noindent\textit{Mathematical role:} rigid motion\par\smallskip
Typed argument/result: \texttt{R\allowbreak\ ::\allowbreak\ 'n\allowbreak\ cgu\_\allowbreak{}point\allowbreak\ \textbackslash{}<Rightarrow>\allowbreak\ 'n\allowbreak\ cgu\_\allowbreak{}point}, with \texttt{'n::finite}; result \texttt{bool}. It holds exactly when there exist \texttt{A\allowbreak\ ::\allowbreak\ 'n\allowbreak\ cgu\_\allowbreak{}point\allowbreak\ \textbackslash{}<Rightarrow>\allowbreak\ 'n\allowbreak\ cgu\_\allowbreak{}point} and \texttt{b\allowbreak\ ::\allowbreak\ 'n\allowbreak\ cgu\_\allowbreak{}point} such that \texttt{linear\allowbreak\ A}, \texttt{surj\allowbreak\ A}, for all \texttt{u,v\allowbreak\ ::\allowbreak\ 'n\allowbreak\ cgu\_\allowbreak{}point}, \texttt{inner\allowbreak\ (A\allowbreak\ u)\allowbreak\ (A\allowbreak\ v)\allowbreak\ =\allowbreak\ inner\allowbreak\ u\allowbreak\ v}, and \texttt{R\allowbreak\ =\allowbreak\ (\textbackslash{}<lambda>x.\allowbreak{}\allowbreak\ A\allowbreak\ x\allowbreak\ +\allowbreak\ b)}.

\dossierentry{\texttt{\detokenize{cgu_local_graph_boundary_at_v5}}}
\reviewlabel{Isabelle code}
\begin{lstlisting}
definition cgu_local_graph_boundary_at_v5 ::
  "'n::finite cgu_point set \<Rightarrow> 'n cgu_point \<Rightarrow> bool"
where
  "cgu_local_graph_boundary_at_v5 D x \<longleftrightarrow>
    (\<exists>U R i zeta L.
      open U \<and>
      x \<in> U \<and>
      cgu_rigid_motion_v5 R \<and>
      0 < L \<and>
      lipschitz_on L UNIV zeta \<and>
      R ` (D \<inter> U) =
        R ` U \<inter>
          {y. y $ i < zeta (cgu_coordinate_hyperplane_projection_v5 i y)} \<and>
      R ` (frontier D \<inter> U) =
        R ` U \<inter>
          {y. y $ i = zeta (cgu_coordinate_hyperplane_projection_v5 i y)})"
\end{lstlisting}
\reviewlabel{English mathematical translation}
\noindent\textit{Mathematical role:} Lipschitz graph point\par\smallskip
Typed arguments/result: \texttt{D\allowbreak\ ::\allowbreak\ 'n\allowbreak\ cgu\_\allowbreak{}point\allowbreak\ set}, \texttt{x\allowbreak\ ::\allowbreak\ 'n\allowbreak\ cgu\_\allowbreak{}point}, with \texttt{'n::finite}; result \texttt{bool}. It holds exactly when there exist \texttt{U\allowbreak\ ::\allowbreak\ 'n\allowbreak\ cgu\_\allowbreak{}point\allowbreak\ set}, \texttt{R\allowbreak\ ::\allowbreak\ 'n\allowbreak\ cgu\_\allowbreak{}point\allowbreak\ \textbackslash{}<Rightarrow>\allowbreak\ 'n\allowbreak\ cgu\_\allowbreak{}point}, \texttt{i\allowbreak\ ::\allowbreak\ 'n}, \texttt{zeta\allowbreak\ ::\allowbreak\ 'n\allowbreak\ cgu\_\allowbreak{}point\allowbreak\ \textbackslash{}<Rightarrow>\allowbreak\ real}, and \texttt{L\allowbreak\ ::\allowbreak\ real} such that \texttt{open\allowbreak\ U}, \texttt{x\allowbreak\ \textbackslash{}<in>\allowbreak\ U}, \texttt{cgu\_\allowbreak{}rigid\_\allowbreak{}motion\_\allowbreak{}v5\allowbreak\ R}, \texttt{0<L}, \texttt{lipschitz\_\allowbreak{}on\allowbreak\ L\allowbreak\ UNIV\allowbreak\ zeta}, \texttt{R\allowbreak\ `\allowbreak\ (D\allowbreak\ \textbackslash{}<inter>\allowbreak\ U)\allowbreak\ =\allowbreak\ R\allowbreak\ `\allowbreak\ U\allowbreak\ \textbackslash{}<inter>\allowbreak\ \{y.\allowbreak{}\allowbreak\ y\allowbreak\ \$\allowbreak\ i\allowbreak\ <\allowbreak\ zeta\allowbreak\ (cgu\_\allowbreak{}coordinate\_\allowbreak{}hyperplane\_\allowbreak{}projection\_\allowbreak{}v5\allowbreak\ i\allowbreak\ y)\}}, and \texttt{R\allowbreak\ `\allowbreak\ (frontier\allowbreak\ D\allowbreak\ \textbackslash{}<inter>\allowbreak\ U)\allowbreak\ =\allowbreak\ R\allowbreak\ `\allowbreak\ U\allowbreak\ \textbackslash{}<inter>\allowbreak\ \{y.\allowbreak{}\allowbreak\ y\allowbreak\ \$\allowbreak\ i\allowbreak\ =\allowbreak\ zeta\allowbreak\ (cgu\_\allowbreak{}coordinate\_\allowbreak{}hyperplane\_\allowbreak{}projection\_\allowbreak{}v5\allowbreak\ i\allowbreak\ y)\}}.

\dossierentry{\texttt{\detokenize{cgu_graph_lipschitz_domain_v5}}}
\reviewlabel{Isabelle code}
\begin{lstlisting}
definition cgu_graph_lipschitz_domain_v5 ::
  "'n::finite cgu_point set \<Rightarrow> bool"
where
  "cgu_graph_lipschitz_domain_v5 D \<longleftrightarrow>
    cgu_lipschitz_domain D \<and>
    (\<forall>x\<in>frontier D. cgu_local_graph_boundary_at_v5 D x)"
\end{lstlisting}
\reviewlabel{English mathematical translation}
\noindent\textit{Mathematical role:} Lipschitz admissible domain\par\smallskip
Typed argument/result: \texttt{D\allowbreak\ ::\allowbreak\ 'n\allowbreak\ cgu\_\allowbreak{}point\allowbreak\ set}, with \texttt{'n::finite}; result \texttt{bool}. It holds exactly when \texttt{cgu\_\allowbreak{}lipschitz\_\allowbreak{}domain\allowbreak\ D} and, for every \texttt{x\allowbreak\ \textbackslash{}<in>\allowbreak\ frontier\allowbreak\ D}, \texttt{cgu\_\allowbreak{}local\_\allowbreak{}graph\_\allowbreak{}boundary\_\allowbreak{}at\_\allowbreak{}v5\allowbreak\ D\allowbreak\ x}.

\dossierentry{\texttt{\detokenize{cgu_locally_flat_boundary_point}}}
\reviewlabel{Isabelle code}
\begin{lstlisting}
definition cgu_locally_flat_boundary_point ::
  "'n::finite cgu_point set \<Rightarrow> 'n cgu_point \<Rightarrow> bool"
where
  "cgu_locally_flat_boundary_point D x \<longleftrightarrow>
    x \<in> frontier D \<and>
    (\<exists>normal offset U.
      normal \<noteq> 0 \<and> open U \<and> x \<in> U \<and>
      U \<inter> frontier D =
        U \<inter> cgu_affine_hyperplane normal offset)"
\end{lstlisting}
\reviewlabel{English mathematical translation}
\noindent\textit{Mathematical role:} flat frontier point\par\smallskip
Typed arguments/result: \texttt{D\allowbreak\ ::\allowbreak\ 'n\allowbreak\ cgu\_\allowbreak{}point\allowbreak\ set}, \texttt{x\allowbreak\ ::\allowbreak\ 'n\allowbreak\ cgu\_\allowbreak{}point}, with \texttt{'n::finite}; result \texttt{bool}. It holds exactly when \texttt{x\allowbreak\ \textbackslash{}<in>\allowbreak\ frontier\allowbreak\ D} and there exist \texttt{normal\allowbreak\ ::\allowbreak\ 'n\allowbreak\ cgu\_\allowbreak{}point}, \texttt{offset\allowbreak\ ::\allowbreak\ real}, and \texttt{U\allowbreak\ ::\allowbreak\ 'n\allowbreak\ cgu\_\allowbreak{}point\allowbreak\ set} such that \texttt{normal\allowbreak\ \textbackslash{}<noteq>\allowbreak\ 0}, \texttt{open\allowbreak\ U}, \texttt{x\allowbreak\ \textbackslash{}<in>\allowbreak\ U}, and \texttt{U\allowbreak\ \textbackslash{}<inter>\allowbreak\ frontier\allowbreak\ D\allowbreak\ =\allowbreak\ U\allowbreak\ \textbackslash{}<inter>\allowbreak\ cgu\_\allowbreak{}affine\_\allowbreak{}hyperplane\allowbreak\ normal\allowbreak\ offset}.

\dossierentry{\texttt{\detokenize{cgu_relative_frontier}}}
\reviewlabel{Isabelle code}
\begin{lstlisting}
definition cgu_relative_frontier ::
  "'a::topological_space set \<Rightarrow> 'a set \<Rightarrow> 'a set"
where
  "cgu_relative_frontier X A = closure A \<inter> closure (X - A)"
\end{lstlisting}
\reviewlabel{English mathematical translation}
\noindent\textit{Mathematical role:} ambient two-sided edge \(\operatorname{edge}_X(A)\)\par\smallskip
Typed arguments/result: \texttt{X,A\allowbreak\ ::\allowbreak\ 'a\allowbreak\ set}, with \texttt{'a::topological\_\allowbreak{}space}; result \texttt{'a\allowbreak\ set}. \texttt{cgu\_\allowbreak{}relative\_\allowbreak{}frontier\allowbreak\ X\allowbreak\ A\allowbreak\ =\allowbreak\ closure\allowbreak\ A\allowbreak\ \textbackslash{}<inter>\allowbreak\ closure\allowbreak\ (X\allowbreak\ -\allowbreak{}\allowbreak\ A)}.

\dossierentry{\texttt{\detokenize{cgu_subdivision_edge_corner_set}}}
\reviewlabel{Isabelle code}
\begin{lstlisting}
definition cgu_subdivision_edge_corner_set ::
  "('n::finite, 'c::finite, 'h) cgu_model \<Rightarrow> 'n cgu_point set"
where
  "cgu_subdivision_edge_corner_set M =
    (\<Union>c.
      let X = frontier (cgu_cell M c);
          F = {x\<in>X. cgu_locally_flat_boundary_point (cgu_cell M c) x}
      in cgu_relative_frontier X F)"
\end{lstlisting}
\reviewlabel{English mathematical translation}
\noindent\textit{Mathematical role:} transition set \(\Sigma_M\)\par\smallskip
Typed argument/result: \texttt{M\allowbreak\ ::\allowbreak\ ('n,'c,'h)\allowbreak\ cgu\_\allowbreak{}model}, with \texttt{'n,'c::finite}; result \texttt{'n\allowbreak\ cgu\_\allowbreak{}point\allowbreak\ set}. It is the union over every \texttt{c\allowbreak\ ::\allowbreak\ 'c} of the following: let \texttt{X\allowbreak\ =\allowbreak\ frontier\allowbreak\ (cgu\_\allowbreak{}cell\allowbreak\ M\allowbreak\ c)} and \texttt{F\allowbreak\ =\allowbreak\ \{x\allowbreak\ \textbackslash{}<in>\allowbreak\ X.\allowbreak{}\allowbreak\ cgu\_\allowbreak{}locally\_\allowbreak{}flat\_\allowbreak{}boundary\_\allowbreak{}point\allowbreak\ (cgu\_\allowbreak{}cell\allowbreak\ M\allowbreak\ c)\allowbreak\ x\}}; then take \texttt{cgu\_\allowbreak{}relative\_\allowbreak{}frontier\allowbreak\ X\allowbreak\ F}.

\dossierentry{\texttt{\detokenize{cgu_graph_subdivision_core_v7}}}
\reviewlabel{Isabelle code}
\begin{lstlisting}
definition cgu_graph_subdivision_core_v7 ::
  "('n::finite, 'c::finite, 'h::finite) cgu_model \<Rightarrow> bool"
where
  "cgu_graph_subdivision_core_v7 M \<longleftrightarrow>
    cgu_lipschitz_domain (cgu_domain M) \<and>
    cgu_measured M \<noteq> {} \<and>
    openin (top_of_set (frontier (cgu_domain M))) (cgu_measured M) \<and>
    cgu_measured M \<subseteq> frontier (cgu_domain M) \<and>
    (\<forall>c. cgu_cell M c \<subseteq> cgu_domain M \<and>
      cgu_lipschitz_domain (cgu_cell M c)) \<and>
    (\<forall>c d. c \<noteq> d \<longrightarrow> cgu_cell M c \<inter> cgu_cell M d = {}) \<and>
    closure (cgu_domain M) = (\<Union>c. closure (cgu_cell M c)) \<and>
    cgu_graph_lipschitz_domain_v5 (cgu_domain M) \<and>
    (\<forall>c. cgu_graph_lipschitz_domain_v5 (cgu_cell M c))"
\end{lstlisting}
\reviewlabel{English mathematical translation}
\noindent\textit{Mathematical role:} admissible subdivision\par\smallskip
Typed argument/result: \texttt{M\allowbreak\ ::\allowbreak\ ('n,'c,'h)\allowbreak\ cgu\_\allowbreak{}model}, with \texttt{'n,'c,'h::finite}; result \texttt{bool}. It holds exactly when: \texttt{cgu\_\allowbreak{}lipschitz\_\allowbreak{}domain\allowbreak\ (cgu\_\allowbreak{}domain\allowbreak\ M)}; \texttt{cgu\_\allowbreak{}measured\allowbreak\ M\allowbreak\ \textbackslash{}<noteq>\allowbreak\ \{\}}; \texttt{cgu\_\allowbreak{}measured\allowbreak\ M} is open in \texttt{top\_\allowbreak{}of\_\allowbreak{}set\allowbreak\ (frontier\allowbreak\ (cgu\_\allowbreak{}domain\allowbreak\ M))}; \texttt{cgu\_\allowbreak{}measured\allowbreak\ M\allowbreak\ \textbackslash{}<subseteq>\allowbreak\ frontier\allowbreak\ (cgu\_\allowbreak{}domain\allowbreak\ M)}; for every \texttt{c\allowbreak\ ::\allowbreak\ 'c}, \texttt{cgu\_\allowbreak{}cell\allowbreak\ M\allowbreak\ c\allowbreak\ \textbackslash{}<subseteq>\allowbreak\ cgu\_\allowbreak{}domain\allowbreak\ M} and \texttt{cgu\_\allowbreak{}lipschitz\_\allowbreak{}domain\allowbreak\ (cgu\_\allowbreak{}cell\allowbreak\ M\allowbreak\ c)}; for every \texttt{c,d\allowbreak\ ::\allowbreak\ 'c}, the hypothesis \texttt{c\allowbreak\ \textbackslash{}<noteq>\allowbreak\ d} implies \texttt{cgu\_\allowbreak{}cell\allowbreak\ M\allowbreak\ c\allowbreak\ \textbackslash{}<inter>\allowbreak\ cgu\_\allowbreak{}cell\allowbreak\ M\allowbreak\ d\allowbreak\ =\allowbreak\ \{\}}; \texttt{closure\allowbreak\ (cgu\_\allowbreak{}domain\allowbreak\ M)\allowbreak\ =\allowbreak\ (\textbackslash{}<Union>c.\allowbreak{}\allowbreak\ closure\allowbreak\ (cgu\_\allowbreak{}cell\allowbreak\ M\allowbreak\ c))}; \texttt{cgu\_\allowbreak{}graph\_\allowbreak{}lipschitz\_\allowbreak{}domain\_\allowbreak{}v5\allowbreak\ (cgu\_\allowbreak{}domain\allowbreak\ M)}; and for every \texttt{c\allowbreak\ ::\allowbreak\ 'c}, \texttt{cgu\_\allowbreak{}graph\_\allowbreak{}lipschitz\_\allowbreak{}domain\_\allowbreak{}v5\allowbreak\ (cgu\_\allowbreak{}cell\allowbreak\ M\allowbreak\ c)}.

\dossierentry{\texttt{\detokenize{cgu_graph_boundary_regularities_imply_lipschitz_domain}}}
\reviewlabel{Isabelle code}
\begin{lstlisting}
theorem cgu_graph_boundary_regularities_imply_lipschitz_domain:
  shows \<open>D \<noteq> {} \<and> open D \<and> connected D \<and> bounded D \<and>
    (\<forall>x\<in>frontier D. cgu_local_graph_boundary_at_v5 D x)
    \<longrightarrow> cgu_lipschitz_domain D\<close>
\end{lstlisting}
\reviewlabel{English mathematical translation}
\noindent\textit{Mathematical role:} Lipschitz-graph boundary-chart lemma\par\smallskip
Implicit typed variable: \texttt{D\allowbreak\ ::\allowbreak\ 'n\allowbreak\ cgu\_\allowbreak{}point\allowbreak\ set} with \texttt{'n::finite} as forced by \texttt{cgu\_\allowbreak{}local\_\allowbreak{}graph\_\allowbreak{}boundary\_\allowbreak{}at\_\allowbreak{}v5}/\texttt{cgu\_\allowbreak{}lipschitz\_\allowbreak{}domain}. There are no separate meta-level \texttt{assumes}; the single shown object-level proposition is: the conjunction \texttt{D\allowbreak\ \textbackslash{}<noteq>\allowbreak\ \{\}}, \texttt{open\allowbreak\ D}, \texttt{connected\allowbreak\ D}, \texttt{bounded\allowbreak\ D}, and \texttt{\textbackslash{}<forall>x\textbackslash{}<in>frontier\allowbreak\ D.\allowbreak{}\allowbreak\ cgu\_\allowbreak{}local\_\allowbreak{}graph\_\allowbreak{}boundary\_\allowbreak{}at\_\allowbreak{}v5\allowbreak\ D\allowbreak\ x} implies the conclusion \texttt{cgu\_\allowbreak{}lipschitz\_\allowbreak{}domain\allowbreak\ D}. The quantifier over \texttt{x} is restricted to \texttt{frontier\allowbreak\ D}; \texttt{D} and its type are implicit schematic parameters of the theorem.

\dossierentry{\texttt{\detokenize{cgu_share_flat_face_in}}}
\reviewlabel{Isabelle code}
\begin{lstlisting}
definition cgu_share_flat_face_in ::
  "'n::finite cgu_point set \<Rightarrow> 'n cgu_point set \<Rightarrow>
    'n cgu_point set \<Rightarrow> bool"
where
  "cgu_share_flat_face_in G E F \<longleftrightarrow>
    (\<exists>normal offset U.
      normal \<noteq> 0 \<and> open U \<and>
      U \<inter> cgu_affine_hyperplane normal offset \<noteq> {} \<and>
      U \<inter> cgu_affine_hyperplane normal offset
        \<subseteq> G \<inter> frontier E \<inter> frontier F)"
\end{lstlisting}
\reviewlabel{English mathematical translation}
\noindent\textit{Mathematical role:} planar adjacency \(E\sim_GF\)\par\smallskip
Typed arguments/result: \texttt{G,E,F\allowbreak\ ::\allowbreak\ 'n\allowbreak\ cgu\_\allowbreak{}point\allowbreak\ set}, with \texttt{'n::finite}; result \texttt{bool}. It holds exactly when there exist \texttt{normal\allowbreak\ ::\allowbreak\ 'n\allowbreak\ cgu\_\allowbreak{}point}, \texttt{offset\allowbreak\ ::\allowbreak\ real}, and \texttt{U\allowbreak\ ::\allowbreak\ 'n\allowbreak\ cgu\_\allowbreak{}point\allowbreak\ set} such that \texttt{normal\allowbreak\ \textbackslash{}<noteq>\allowbreak\ 0}, \texttt{open\allowbreak\ U}, \texttt{U\allowbreak\ \textbackslash{}<inter>\allowbreak\ cgu\_\allowbreak{}affine\_\allowbreak{}hyperplane\allowbreak\ normal\allowbreak\ offset\allowbreak\ \textbackslash{}<noteq>\allowbreak\ \{\}}, and \texttt{U\allowbreak\ \textbackslash{}<inter>\allowbreak\ cgu\_\allowbreak{}affine\_\allowbreak{}hyperplane\allowbreak\ normal\allowbreak\ offset\allowbreak\ \textbackslash{}<subseteq>\allowbreak\ G\allowbreak\ \textbackslash{}<inter>\allowbreak\ frontier\allowbreak\ E\allowbreak\ \textbackslash{}<inter>\allowbreak\ frontier\allowbreak\ F}.

\dossierentry{\texttt{\detokenize{cgu_linearly_independent3}}}
\reviewlabel{Isabelle code}
\begin{lstlisting}
definition cgu_linearly_independent3 ::
  "real ^ 'n::finite \<Rightarrow> real ^ 'n \<Rightarrow> real ^ 'n \<Rightarrow> bool"
where
  "cgu_linearly_independent3 a b c \<longleftrightarrow>
    (\<forall>x y z::real.
      x *\<^sub>R a + y *\<^sub>R b + z *\<^sub>R c = 0
      \<longrightarrow> x = 0 \<and> y = 0 \<and> z = 0)"
\end{lstlisting}
\reviewlabel{English mathematical translation}
\noindent\textit{Mathematical role:} independence of three normals\par\smallskip
Typed arguments/result: \texttt{a,b,c\allowbreak\ ::\allowbreak\ real\allowbreak\ \textasciicircum{}\allowbreak\ 'n}, with \texttt{'n::finite}; result \texttt{bool}. It holds exactly when, for every \texttt{x,y,z\allowbreak\ ::\allowbreak\ real}, the hypothesis \texttt{x\allowbreak\ *\textbackslash{}<\textasciicircum{}sub>R\allowbreak\ a\allowbreak\ +\allowbreak\ y\allowbreak\ *\textbackslash{}<\textasciicircum{}sub>R\allowbreak\ b\allowbreak\ +\allowbreak\ z\allowbreak\ *\textbackslash{}<\textasciicircum{}sub>R\allowbreak\ c\allowbreak\ =\allowbreak\ 0} implies the three-part conclusion \texttt{x=0\allowbreak\ \textbackslash{}<and>\allowbreak\ y=0\allowbreak\ \textbackslash{}<and>\allowbreak\ z=0}.

\dossierentry{\texttt{\detokenize{cgu_n_generic_on}}}
\reviewlabel{Isabelle code}
\begin{lstlisting}
definition cgu_n_generic_on ::
  "nat \<Rightarrow> 'h::finite set \<Rightarrow>
    ('h \<Rightarrow> real ^ 'n::finite) \<Rightarrow> ('h \<Rightarrow> real) \<Rightarrow>
    ('h \<Rightarrow> 'h set) \<Rightarrow> bool"
where
  "cgu_n_generic_on N I normal offset J \<longleftrightarrow>
    inj_on (\<lambda>i. cgu_affine_hyperplane (normal i) (offset i)) I \<and>
    (\<forall>i\<in>I.
      J i \<subseteq> I - {i} \<and> card (J i) = N + 2 \<and>
      (\<forall>j\<in>J i. \<forall>k\<in>J i. j \<noteq> k \<longrightarrow>
        cgu_linearly_independent3 (normal i) (normal j) (normal k) \<and>
        aff_dim
          (cgu_affine_hyperplane (normal i) (offset i) \<inter>
           cgu_affine_hyperplane (normal j) (offset j) \<inter>
           cgu_affine_hyperplane (normal k) (offset k))
        = int (CARD('n)) - 3) \<and>
      (\<forall>j\<in>J i. \<forall>k\<in>J i. \<forall>l\<in>J i.
        j \<noteq> k \<and> l \<noteq> j \<and> l \<noteq> k \<longrightarrow>
        \<not> cgu_affine_hyperplane (normal i) (offset i) \<inter>
             cgu_affine_hyperplane (normal j) (offset j) \<inter>
             cgu_affine_hyperplane (normal k) (offset k)
           \<subseteq> cgu_affine_hyperplane (normal l) (offset l)))"
\end{lstlisting}
\reviewlabel{English mathematical translation}
\noindent\textit{Mathematical role:} \(N\)-hyperplane configuration\par\smallskip
Typed arguments/result: \texttt{N\allowbreak\ ::\allowbreak\ nat}, \texttt{I\allowbreak\ ::\allowbreak\ 'h\allowbreak\ set}, \texttt{normal\allowbreak\ ::\allowbreak\ 'h\allowbreak\ \textbackslash{}<Rightarrow>\allowbreak\ real\allowbreak\ \textasciicircum{}\allowbreak\ 'n}, \texttt{offset\allowbreak\ ::\allowbreak\ 'h\allowbreak\ \textbackslash{}<Rightarrow>\allowbreak\ real}, \texttt{J\allowbreak\ ::\allowbreak\ 'h\allowbreak\ \textbackslash{}<Rightarrow>\allowbreak\ 'h\allowbreak\ set}, with \texttt{'h,'n::finite}; result \texttt{bool}. It holds exactly when the map \texttt{i\allowbreak\ \textbackslash{}<mapsto>\allowbreak\ cgu\_\allowbreak{}affine\_\allowbreak{}hyperplane\allowbreak\ (normal\allowbreak\ i)\allowbreak\ (offset\allowbreak\ i)} is injective on \texttt{I}, and for every \texttt{i\allowbreak\ \textbackslash{}<in>\allowbreak\ I}: (a) \texttt{J\allowbreak\ i\allowbreak\ \textbackslash{}<subseteq>\allowbreak\ I\allowbreak\ -\allowbreak{}\allowbreak\ \{i\}} and \texttt{card\allowbreak\ (J\allowbreak\ i)\allowbreak\ =\allowbreak\ N+2}; (b) for every \texttt{j\allowbreak\ \textbackslash{}<in>\allowbreak\ J\allowbreak\ i} and \texttt{k\allowbreak\ \textbackslash{}<in>\allowbreak\ J\allowbreak\ i}, if \texttt{j\allowbreak\ \textbackslash{}<noteq>\allowbreak\ k}, then \texttt{cgu\_\allowbreak{}linearly\_\allowbreak{}independent3\allowbreak\ (normal\allowbreak\ i)\allowbreak\ (normal\allowbreak\ j)\allowbreak\ (normal\allowbreak\ k)} and the affine dimension of the intersection of the three corresponding \texttt{cgu\_\allowbreak{}affine\_\allowbreak{}hyperplane} sets is \texttt{int\allowbreak\ (CARD('n))\allowbreak\ -\allowbreak{}\allowbreak\ 3}; and (c) for every \texttt{j,k,l\allowbreak\ \textbackslash{}<in>\allowbreak\ J\allowbreak\ i}, if \texttt{j\allowbreak\ \textbackslash{}<noteq>\allowbreak\ k}, \texttt{l\allowbreak\ \textbackslash{}<noteq>\allowbreak\ j}, and \texttt{l\allowbreak\ \textbackslash{}<noteq>\allowbreak\ k}, then the intersection of the sets for \texttt{i,j,k} is not a subset of the set for \texttt{l}.

\dossierentry{\texttt{\detokenize{cgu_triple_region_geometry}}}
\reviewlabel{Isabelle code}
\begin{lstlisting}
definition cgu_triple_region_geometry ::
  "'n::finite cgu_point set \<Rightarrow>
    'n cgu_point \<Rightarrow> real \<Rightarrow> 'n cgu_point \<Rightarrow> real \<Rightarrow>
    'n cgu_point \<Rightarrow> real \<Rightarrow> bool"
where
  "cgu_triple_region_geometry W ni oi nj oj nk ok \<longleftrightarrow>
    W \<noteq> {} \<and>
    openin
      (top_of_set
        (cgu_affine_hyperplane ni oi \<inter>
         cgu_affine_hyperplane nj oj \<inter>
         cgu_affine_hyperplane nk ok)) W \<and>
    W \<subseteq>
      cgu_affine_hyperplane ni oi \<inter>
      cgu_affine_hyperplane nj oj \<inter>
      cgu_affine_hyperplane nk ok"
\end{lstlisting}
\reviewlabel{English mathematical translation}
\noindent\textit{Mathematical role:} relatively open triple window\par\smallskip
Typed arguments/result: \texttt{W\allowbreak\ ::\allowbreak\ 'n\allowbreak\ cgu\_\allowbreak{}point\allowbreak\ set}; \texttt{ni,nj,nk\allowbreak\ ::\allowbreak\ 'n\allowbreak\ cgu\_\allowbreak{}point}; \texttt{oi,oj,ok\allowbreak\ ::\allowbreak\ real}; with \texttt{'n::finite}; result \texttt{bool}. It holds exactly when \texttt{W\allowbreak\ \textbackslash{}<noteq>\allowbreak\ \{\}}, \texttt{W} is open in the subspace on \texttt{cgu\_\allowbreak{}affine\_\allowbreak{}hyperplane\allowbreak\ ni\allowbreak\ oi\allowbreak\ \textbackslash{}<inter>\allowbreak\ cgu\_\allowbreak{}affine\_\allowbreak{}hyperplane\allowbreak\ nj\allowbreak\ oj\allowbreak\ \textbackslash{}<inter>\allowbreak\ cgu\_\allowbreak{}affine\_\allowbreak{}hyperplane\allowbreak\ nk\allowbreak\ ok}, and \texttt{W} is a subset of that same intersection.

\dossierentry{\texttt{\detokenize{cgu_polynomial_positive_on}}}
\reviewlabel{Isabelle code}
\begin{lstlisting}
definition cgu_polynomial_positive_on ::
  "'n::finite cgu_matrix_polynomial \<Rightarrow> 'n cgu_point set \<Rightarrow> bool"
where
  "cgu_polynomial_positive_on P W \<longleftrightarrow>
    (\<forall>x\<in>W.
      cgu_symmetric_positive_definite_matrix (cgu_matrix_poly_eval P x))"
\end{lstlisting}
\reviewlabel{English mathematical translation}
\noindent\textit{Mathematical role:} pointwise positive definiteness on a window\par\smallskip
Typed arguments/result: \texttt{P\allowbreak\ ::\allowbreak\ 'n\allowbreak\ cgu\_\allowbreak{}matrix\_\allowbreak{}polynomial}, \texttt{W\allowbreak\ ::\allowbreak\ 'n\allowbreak\ cgu\_\allowbreak{}point\allowbreak\ set}, with \texttt{'n::finite}; result \texttt{bool}. It holds exactly when for every \texttt{x\allowbreak\ \textbackslash{}<in>\allowbreak\ W}, \texttt{cgu\_\allowbreak{}symmetric\_\allowbreak{}positive\_\allowbreak{}definite\_\allowbreak{}matrix\allowbreak\ (cgu\_\allowbreak{}matrix\_\allowbreak{}poly\_\allowbreak{}eval\allowbreak\ P\allowbreak\ x)}.

\dossierentry{\texttt{\detokenize{cgu_order_valid}}}
\reviewlabel{Isabelle code}
\begin{lstlisting}
definition cgu_order_valid ::
  "('n, 'c::finite, 'h) cgu_model \<Rightarrow> bool"
where
  "cgu_order_valid M \<longleftrightarrow>
    bij_betw (cgu_order M) {..<CARD('c)} UNIV"
\end{lstlisting}
\reviewlabel{English mathematical translation}
\noindent\textit{Mathematical role:} bijective recovery enumeration\par\smallskip
Typed argument/result: \texttt{M\allowbreak\ ::\allowbreak\ ('n,'c,'h)\allowbreak\ cgu\_\allowbreak{}model}, with \texttt{'c::finite}; result \texttt{bool}. It holds exactly when \texttt{cgu\_\allowbreak{}order\allowbreak\ M} is a bijection from \texttt{\{.\allowbreak{}.\allowbreak{}<CARD('c)\}} to \texttt{UNIV}.

\dossierentry{\texttt{\detokenize{cgu_before_indices}}}
\reviewlabel{Isabelle code}
\begin{lstlisting}
definition cgu_before_indices ::
  "('n, 'c, 'h) cgu_model \<Rightarrow> nat \<Rightarrow> 'c set"
where
  "cgu_before_indices M r = cgu_order M ` {..<r}"
\end{lstlisting}
\reviewlabel{English mathematical translation}
\noindent\textit{Mathematical role:} earlier-cell set \(C_{<s}\)\par\smallskip
Typed arguments/result: \texttt{M\allowbreak\ ::\allowbreak\ ('n,'c,'h)\allowbreak\ cgu\_\allowbreak{}model}, \texttt{r\allowbreak\ ::\allowbreak\ nat}; result \texttt{'c\allowbreak\ set}. It is the image \texttt{cgu\_\allowbreak{}order\allowbreak\ M\allowbreak\ `\allowbreak\ \{.\allowbreak{}.\allowbreak{}<r\}}.

\dossierentry{\texttt{\detokenize{cgu_after_indices}}}
\reviewlabel{Isabelle code}
\begin{lstlisting}
definition cgu_after_indices ::
  "('n, 'c::finite, 'h) cgu_model \<Rightarrow> nat \<Rightarrow> 'c set"
where
  "cgu_after_indices M r = cgu_order M ` {r..<CARD('c)}"
\end{lstlisting}
\reviewlabel{English mathematical translation}
\noindent\textit{Mathematical role:} future-cell set \(C_{\ge s}\)\par\smallskip
Typed arguments/result: \texttt{M\allowbreak\ ::\allowbreak\ ('n,'c,'h)\allowbreak\ cgu\_\allowbreak{}model}, \texttt{r\allowbreak\ ::\allowbreak\ nat}, with \texttt{'c::finite}; result \texttt{'c\allowbreak\ set}. It is the image \texttt{cgu\_\allowbreak{}order\allowbreak\ M\allowbreak\ `\allowbreak\ \{r.\allowbreak{}.\allowbreak{}<CARD('c)\}}.

\dossierentry{\texttt{\detokenize{cgu_remaining_domain}}}
\reviewlabel{Isabelle code}
\begin{lstlisting}
definition cgu_remaining_domain ::
  "('n::finite, 'c::finite, 'h) cgu_model \<Rightarrow> nat \<Rightarrow>
    'n cgu_point set"
where
  "cgu_remaining_domain M r =
    interior (\<Union>c\<in>cgu_after_indices M r. closure (cgu_cell M c))"
\end{lstlisting}
\reviewlabel{English mathematical translation}
\noindent\textit{Mathematical role:} remaining region \(\mathcal R_s\)\par\smallskip
Typed arguments/result: \texttt{M\allowbreak\ ::\allowbreak\ ('n,'c,'h)\allowbreak\ cgu\_\allowbreak{}model}, \texttt{r\allowbreak\ ::\allowbreak\ nat}, with \texttt{'n,'c::finite}; result \texttt{'n\allowbreak\ cgu\_\allowbreak{}point\allowbreak\ set}. It is the interior of the union, over \texttt{c\allowbreak\ \textbackslash{}<in>\allowbreak\ cgu\_\allowbreak{}after\_\allowbreak{}indices\allowbreak\ M\allowbreak\ r}, of \texttt{closure\allowbreak\ (cgu\_\allowbreak{}cell\allowbreak\ M\allowbreak\ c)}.

\dossierentry{\texttt{\detokenize{cgu_recovered_region}}}
\reviewlabel{Isabelle code}
\begin{lstlisting}
definition cgu_recovered_region ::
  "('n::finite, 'c::finite, 'h) cgu_model \<Rightarrow> nat \<Rightarrow>
    'n cgu_point set"
where
  "cgu_recovered_region M r =
    cgu_domain M - closure (cgu_remaining_domain M r)"
\end{lstlisting}
\reviewlabel{English mathematical translation}
\noindent\textit{Mathematical role:} recovered region \(\mathcal K_s\)\par\smallskip
Typed arguments/result: \texttt{M\allowbreak\ ::\allowbreak\ ('n,'c,'h)\allowbreak\ cgu\_\allowbreak{}model}, \texttt{r\allowbreak\ ::\allowbreak\ nat}, with \texttt{'n,'c::finite}; result \texttt{'n\allowbreak\ cgu\_\allowbreak{}point\allowbreak\ set}. It is \texttt{cgu\_\allowbreak{}domain\allowbreak\ M\allowbreak\ -\allowbreak{}\allowbreak\ closure\allowbreak\ (cgu\_\allowbreak{}remaining\_\allowbreak{}domain\allowbreak\ M\allowbreak\ r)}.

\dossierentry{\texttt{\detokenize{cgu_face_connected}}}
\reviewlabel{Isabelle code}
\begin{lstlisting}
definition cgu_face_connected ::
  "'n::finite cgu_point set \<Rightarrow> 'c set \<Rightarrow>
    ('c \<Rightarrow> 'n cgu_point set) \<Rightarrow> bool"
where
  "cgu_face_connected G I cell \<longleftrightarrow>
    (\<forall>i\<in>I. \<forall>j\<in>I.
      rtranclp
        (\<lambda>a b. a \<in> I \<and> b \<in> I \<and>
          cgu_share_flat_face_in G (cell a) (cell b)) i j)"
\end{lstlisting}
\reviewlabel{English mathematical translation}
\noindent\textit{Mathematical role:} connectivity through planar adjacency\par\smallskip
Typed arguments/result: \texttt{G\allowbreak\ ::\allowbreak\ 'n\allowbreak\ cgu\_\allowbreak{}point\allowbreak\ set}, \texttt{I\allowbreak\ ::\allowbreak\ 'c\allowbreak\ set}, \texttt{cell\allowbreak\ ::\allowbreak\ 'c\allowbreak\ \textbackslash{}<Rightarrow>\allowbreak\ 'n\allowbreak\ cgu\_\allowbreak{}point\allowbreak\ set}, with \texttt{'n::finite}; result \texttt{bool}. It holds exactly when for every \texttt{i\allowbreak\ \textbackslash{}<in>\allowbreak\ I} and \texttt{j\allowbreak\ \textbackslash{}<in>\allowbreak\ I}, \texttt{rtranclp\allowbreak\ R\allowbreak\ i\allowbreak\ j}, where \texttt{R\allowbreak\ a\allowbreak\ b} is the conjunction \texttt{a\allowbreak\ \textbackslash{}<in>\allowbreak\ I}, \texttt{b\allowbreak\ \textbackslash{}<in>\allowbreak\ I}, and \texttt{cgu\_\allowbreak{}share\_\allowbreak{}flat\_\allowbreak{}face\_\allowbreak{}in\allowbreak\ G\allowbreak\ (cell\allowbreak\ a)\allowbreak\ (cell\allowbreak\ b)}.

\dossierentry{\texttt{\detokenize{cgu_flat_face_patch}}}
\reviewlabel{Isabelle code}
\begin{lstlisting}
definition cgu_flat_face_patch ::
  "'n::finite cgu_point set \<Rightarrow> 'n cgu_point \<Rightarrow> real \<Rightarrow>
    'n cgu_point set \<Rightarrow> bool"
where
  "cgu_flat_face_patch D normal offset gamma \<longleftrightarrow>
    norm normal = 1 \<and>
    gamma \<noteq> {} \<and>
    openin (top_of_set (cgu_affine_hyperplane normal offset)) gamma \<and>
    gamma \<subseteq> frontier D \<inter> cgu_affine_hyperplane normal offset \<and>
    (\<exists>U. open U \<and> gamma \<subseteq> U \<and>
      U \<inter> frontier D = U \<inter> cgu_affine_hyperplane normal offset \<and>
      (D \<inter> U \<subseteq> {x. inner normal x < offset} \<or>
       D \<inter> U \<subseteq> {x. offset < inner normal x}))"
\end{lstlisting}
\reviewlabel{English mathematical translation}
\noindent\textit{Mathematical role:} oriented flat boundary patch\par\smallskip
Typed arguments/result: \texttt{D\allowbreak\ ::\allowbreak\ 'n\allowbreak\ cgu\_\allowbreak{}point\allowbreak\ set}, \texttt{normal\allowbreak\ ::\allowbreak\ 'n\allowbreak\ cgu\_\allowbreak{}point}, \texttt{offset\allowbreak\ ::\allowbreak\ real}, \texttt{gamma\allowbreak\ ::\allowbreak\ 'n\allowbreak\ cgu\_\allowbreak{}point\allowbreak\ set}, with \texttt{'n::finite}; result \texttt{bool}. It holds exactly when \texttt{norm\allowbreak\ normal\allowbreak\ =\allowbreak\ 1}, \texttt{gamma\allowbreak\ \textbackslash{}<noteq>\allowbreak\ \{\}}, \texttt{gamma} is open in the subspace on \texttt{cgu\_\allowbreak{}affine\_\allowbreak{}hyperplane\allowbreak\ normal\allowbreak\ offset}, \texttt{gamma\allowbreak\ \textbackslash{}<subseteq>\allowbreak\ frontier\allowbreak\ D\allowbreak\ \textbackslash{}<inter>\allowbreak\ cgu\_\allowbreak{}affine\_\allowbreak{}hyperplane\allowbreak\ normal\allowbreak\ offset}, and there exists \texttt{U\allowbreak\ ::\allowbreak\ 'n\allowbreak\ cgu\_\allowbreak{}point\allowbreak\ set} such that \texttt{open\allowbreak\ U}, \texttt{gamma\allowbreak\ \textbackslash{}<subseteq>\allowbreak\ U}, \texttt{U\allowbreak\ \textbackslash{}<inter>\allowbreak\ frontier\allowbreak\ D\allowbreak\ =\allowbreak\ U\allowbreak\ \textbackslash{}<inter>\allowbreak\ cgu\_\allowbreak{}affine\_\allowbreak{}hyperplane\allowbreak\ normal\allowbreak\ offset}, and either \texttt{D\allowbreak\ \textbackslash{}<inter>\allowbreak\ U\allowbreak\ \textbackslash{}<subseteq>\allowbreak\ \{x.\allowbreak{}\allowbreak\ inner\allowbreak\ normal\allowbreak\ x\allowbreak\ <\allowbreak\ offset\}} or \texttt{D\allowbreak\ \textbackslash{}<inter>\allowbreak\ U\allowbreak\ \textbackslash{}<subseteq>\allowbreak\ \{x.\allowbreak{}\allowbreak\ offset\allowbreak\ <\allowbreak\ inner\allowbreak\ normal\allowbreak\ x\}}.

\dossierentry{\texttt{\detokenize{cgu_clipped_recovery_interface_v6}}}
\reviewlabel{Isabelle code}
\begin{lstlisting}
definition cgu_clipped_recovery_interface_v6 ::
  "('n::finite, 'c::finite, 'h) cgu_model \<Rightarrow> nat \<Rightarrow> 'n cgu_point set"
where
  "cgu_clipped_recovery_interface_v6 M r =
    cgu_domain M \<inter>
      frontier (cgu_recovered_region M r) \<inter>
      frontier (cgu_remaining_domain M r)"
\end{lstlisting}
\reviewlabel{English mathematical translation}
\noindent\textit{Mathematical role:} recovery interface \(\Upsilon_s\)\par\smallskip
Typed arguments/result: \texttt{M\allowbreak\ ::\allowbreak\ ('n,'c,'h)\allowbreak\ cgu\_\allowbreak{}model}, \texttt{r\allowbreak\ ::\allowbreak\ nat}, with \texttt{'n,'c::finite}; result \texttt{'n\allowbreak\ cgu\_\allowbreak{}point\allowbreak\ set}. It is \texttt{cgu\_\allowbreak{}domain\allowbreak\ M\allowbreak\ \textbackslash{}<inter>\allowbreak\ frontier\allowbreak\ (cgu\_\allowbreak{}recovered\_\allowbreak{}region\allowbreak\ M\allowbreak\ r)\allowbreak\ \textbackslash{}<inter>\allowbreak\ frontier\allowbreak\ (cgu\_\allowbreak{}remaining\_\allowbreak{}domain\allowbreak\ M\allowbreak\ r)}.

\dossierentry{\texttt{\detokenize{cgu_active_pde_face_route_v4}}}
\reviewlabel{Isabelle code}
\begin{lstlisting}
definition cgu_active_pde_face_route_v4 ::
  "('n::finite, 'c::finite, 'h::finite) cgu_model \<Rightarrow>
    nat \<Rightarrow> 'c \<Rightarrow> 'h \<Rightarrow> bool"
where
  "cgu_active_pde_face_route_v4 M r c i \<longleftrightarrow>
    (let remaining = cgu_remaining_domain M r;
         recovered = cgu_recovered_region M r;
         interface = frontier recovered \<inter> frontier remaining;
         patch = cgu_face_patch M c i;
         normal = cgu_normal M c i;
         offset = cgu_offset M c i
     in if r = 0 then
          patch \<subseteq> cgu_measured M \<and>
          cgu_flat_face_patch (cgu_domain M) normal offset patch
        else
          (patch \<subseteq> cgu_measured M \<and>
           cgu_flat_face_patch (cgu_domain M) normal offset patch) \<or>
          (patch \<subseteq> cgu_domain M \<inter> interface \<and>
           cgu_flat_face_patch remaining normal offset patch))"
\end{lstlisting}
\reviewlabel{English mathematical translation}
\noindent\textit{Mathematical role:} stagewise face visibility\par\smallskip
Typed arguments/result: \texttt{M\allowbreak\ ::\allowbreak\ ('n,'c,'h)\allowbreak\ cgu\_\allowbreak{}model}, \texttt{r\allowbreak\ ::\allowbreak\ nat}, \texttt{c\allowbreak\ ::\allowbreak\ 'c}, \texttt{i\allowbreak\ ::\allowbreak\ 'h}, with \texttt{'n,'c,'h::finite}; result \texttt{bool}. Let \texttt{remaining=cgu\_\allowbreak{}remaining\_\allowbreak{}domain\allowbreak\ M\allowbreak\ r}, \texttt{recovered=cgu\_\allowbreak{}recovered\_\allowbreak{}region\allowbreak\ M\allowbreak\ r}, \texttt{interface=frontier\allowbreak\ recovered\allowbreak\ \textbackslash{}<inter>\allowbreak\ frontier\allowbreak\ remaining}, \texttt{patch=cgu\_\allowbreak{}face\_\allowbreak{}patch\allowbreak\ M\allowbreak\ c\allowbreak\ i}, \texttt{normal=cgu\_\allowbreak{}normal\allowbreak\ M\allowbreak\ c\allowbreak\ i}, and \texttt{offset=cgu\_\allowbreak{}offset\allowbreak\ M\allowbreak\ c\allowbreak\ i}. If \texttt{r=0}, the condition is \texttt{patch\allowbreak\ \textbackslash{}<subseteq>\allowbreak\ cgu\_\allowbreak{}measured\allowbreak\ M\allowbreak\ \textbackslash{}<and>\allowbreak\ cgu\_\allowbreak{}flat\_\allowbreak{}face\_\allowbreak{}patch\allowbreak\ (cgu\_\allowbreak{}domain\allowbreak\ M)\allowbreak\ normal\allowbreak\ offset\allowbreak\ patch}. Otherwise it is the disjunction of that same condition and \texttt{patch\allowbreak\ \textbackslash{}<subseteq>\allowbreak\ cgu\_\allowbreak{}domain\allowbreak\ M\allowbreak\ \textbackslash{}<inter>\allowbreak\ interface\allowbreak\ \textbackslash{}<and>\allowbreak\ cgu\_\allowbreak{}flat\_\allowbreak{}face\_\allowbreak{}patch\allowbreak\ remaining\allowbreak\ normal\allowbreak\ offset\allowbreak\ patch}.

\dossierentry{\texttt{\detokenize{cgu_normal_height}}}
\reviewlabel{Isabelle code}
\begin{lstlisting}
definition cgu_normal_height ::
  "'n::finite cgu_point \<Rightarrow> 'n cgu_point \<Rightarrow>
    'n cgu_point \<Rightarrow> real"
where
  "cgu_normal_height center normal x = inner normal (x - center)"
\end{lstlisting}
\reviewlabel{English mathematical translation}
\noindent\textit{Mathematical role:} signed height \(h_{q,\mu}\)\par\smallskip
Typed arguments/result: \texttt{center,normal,x\allowbreak\ ::\allowbreak\ 'n\allowbreak\ cgu\_\allowbreak{}point}, with \texttt{'n::finite}; result \texttt{real}. It equals \texttt{inner\allowbreak\ normal\allowbreak\ (x\allowbreak\ -\allowbreak{}\allowbreak\ center)}.

\dossierentry{\texttt{\detokenize{cgu_tangent_part}}}
\reviewlabel{Isabelle code}
\begin{lstlisting}
definition cgu_tangent_part ::
  "'n::finite cgu_point \<Rightarrow> 'n cgu_point \<Rightarrow>
    'n cgu_point \<Rightarrow> 'n cgu_point"
where
  "cgu_tangent_part center normal x =
    (x - center) - cgu_normal_height center normal x *\<^sub>R normal"
\end{lstlisting}
\reviewlabel{English mathematical translation}
\noindent\textit{Mathematical role:} tangential remainder \(\tau_{q,\mu}\)\par\smallskip
Typed arguments/result: \texttt{center,normal,x\allowbreak\ ::\allowbreak\ 'n\allowbreak\ cgu\_\allowbreak{}point}, with \texttt{'n::finite}; result \texttt{'n\allowbreak\ cgu\_\allowbreak{}point}. It equals \texttt{(x-\allowbreak{}center)\allowbreak\ -\allowbreak{}\allowbreak\ cgu\_\allowbreak{}normal\_\allowbreak{}height\allowbreak\ center\allowbreak\ normal\allowbreak\ x\allowbreak\ *\textbackslash{}<\textasciicircum{}sub>R\allowbreak\ normal}.

\dossierentry{\texttt{\detokenize{cgu_flat_disc}}}
\reviewlabel{Isabelle code}
\begin{lstlisting}
definition cgu_flat_disc ::
  "'n::finite cgu_point \<Rightarrow> 'n cgu_point \<Rightarrow> real \<Rightarrow>
    'n cgu_point set"
where
  "cgu_flat_disc center normal r =
    {x. cgu_normal_height center normal x = 0 \<and>
      norm (cgu_tangent_part center normal x) < r}"
\end{lstlisting}
\reviewlabel{English mathematical translation}
\noindent\textit{Mathematical role:} flat disk \(D(q,\mu,r)\)\par\smallskip
Typed arguments/result: \texttt{center,normal\allowbreak\ ::\allowbreak\ 'n\allowbreak\ cgu\_\allowbreak{}point}, \texttt{r\allowbreak\ ::\allowbreak\ real}, with \texttt{'n::finite}; result \texttt{'n\allowbreak\ cgu\_\allowbreak{}point\allowbreak\ set}. It is the set of all \texttt{x\allowbreak\ ::\allowbreak\ 'n\allowbreak\ cgu\_\allowbreak{}point} such that \texttt{cgu\_\allowbreak{}normal\_\allowbreak{}height\allowbreak\ center\allowbreak\ normal\allowbreak\ x\allowbreak\ =\allowbreak\ 0} and \texttt{norm\allowbreak\ (cgu\_\allowbreak{}tangent\_\allowbreak{}part\allowbreak\ center\allowbreak\ normal\allowbreak\ x)\allowbreak\ <\allowbreak\ r}.

\dossierentry{\texttt{\detokenize{cgu_exterior_half_ball}}}
\reviewlabel{Isabelle code}
\begin{lstlisting}
definition cgu_exterior_half_ball ::
  "'n::finite cgu_point \<Rightarrow> 'n cgu_point \<Rightarrow> real \<Rightarrow>
    'n cgu_point set"
where
  "cgu_exterior_half_ball center normal r =
    ball center r \<inter> {x. cgu_normal_height center normal x < 0}"
\end{lstlisting}
\reviewlabel{English mathematical translation}
\noindent\textit{Mathematical role:} lower half-ball \(B^-(q,\mu,r)\)\par\smallskip
Typed arguments/result: \texttt{center,normal\allowbreak\ ::\allowbreak\ 'n\allowbreak\ cgu\_\allowbreak{}point}, \texttt{r\allowbreak\ ::\allowbreak\ real}, with \texttt{'n::finite}; result \texttt{'n\allowbreak\ cgu\_\allowbreak{}point\allowbreak\ set}. It is \texttt{ball\allowbreak\ center\allowbreak\ r\allowbreak\ \textbackslash{}<inter>\allowbreak\ \{x.\allowbreak{}\allowbreak\ cgu\_\allowbreak{}normal\_\allowbreak{}height\allowbreak\ center\allowbreak\ normal\allowbreak\ x\allowbreak\ <\allowbreak\ 0\}}.

\dossierentry{\texttt{\detokenize{cgu_auxiliary_half_ball_geometry_step}}}
\reviewlabel{Isabelle code}
\begin{lstlisting}
definition cgu_auxiliary_half_ball_geometry_step ::
  "('n::finite, 'c::finite, 'h) cgu_model \<Rightarrow> nat \<Rightarrow> bool"
where
  "cgu_auxiliary_half_ball_geometry_step M r \<longleftrightarrow>
    (let c = cgu_order M r;
         recovered = cgu_recovered_region M r;
         center = cgu_step_center M c;
         normal = cgu_step_interior_normal M c;
         rad = cgu_step_radius M c;
         patch = cgu_step_outer_patch M c;
         adjacent = cgu_step_adjacent_cell M c
     in norm normal = 1 \<and> 0 < rad \<and>
        patch \<noteq> {} \<and>
        openin (top_of_set (frontier (cgu_domain M))) patch \<and>
        patch \<subseteq> cgu_measured M \<inter> frontier recovered \<and>
        closure (cgu_flat_disc center normal rad)
          \<subseteq> patch - cgu_subdivision_edge_corner_set M \<and>
        cgu_domain M \<inter> ball center rad =
          ball center rad \<inter> {x. 0 < cgu_normal_height center normal x} \<and>
        adjacent \<in> cgu_before_indices M r \<and>
        cgu_flat_disc center normal rad
          \<subseteq> frontier (cgu_domain M) \<inter> frontier (cgu_cell M adjacent))"
\end{lstlisting}
\reviewlabel{English mathematical translation}
\noindent\textit{Mathematical role:} exterior recovery-ball condition\par\smallskip
Typed arguments/result: \texttt{M\allowbreak\ ::\allowbreak\ ('n,'c,'h)\allowbreak\ cgu\_\allowbreak{}model}, \texttt{r\allowbreak\ ::\allowbreak\ nat}, with \texttt{'n,'c::finite}; result \texttt{bool}. Let \texttt{c=cgu\_\allowbreak{}order\allowbreak\ M\allowbreak\ r}, \texttt{recovered=cgu\_\allowbreak{}recovered\_\allowbreak{}region\allowbreak\ M\allowbreak\ r}, \texttt{center=cgu\_\allowbreak{}step\_\allowbreak{}center\allowbreak\ M\allowbreak\ c}, \texttt{normal=cgu\_\allowbreak{}step\_\allowbreak{}interior\_\allowbreak{}normal\allowbreak\ M\allowbreak\ c}, \texttt{rad=cgu\_\allowbreak{}step\_\allowbreak{}radius\allowbreak\ M\allowbreak\ c}, \texttt{patch=cgu\_\allowbreak{}step\_\allowbreak{}outer\_\allowbreak{}patch\allowbreak\ M\allowbreak\ c}, and \texttt{adjacent=cgu\_\allowbreak{}step\_\allowbreak{}adjacent\_\allowbreak{}cell\allowbreak\ M\allowbreak\ c}. The condition is the conjunction: \texttt{norm\allowbreak\ normal=1}; \texttt{0<rad}; \texttt{patch\textbackslash{}<noteq>\{\}}; \texttt{patch} is open in the subspace on \texttt{frontier\allowbreak\ (cgu\_\allowbreak{}domain\allowbreak\ M)}; \texttt{patch\allowbreak\ \textbackslash{}<subseteq>\allowbreak\ cgu\_\allowbreak{}measured\allowbreak\ M\allowbreak\ \textbackslash{}<inter>\allowbreak\ frontier\allowbreak\ recovered}; \texttt{closure\allowbreak\ (cgu\_\allowbreak{}flat\_\allowbreak{}disc\allowbreak\ center\allowbreak\ normal\allowbreak\ rad)\allowbreak\ \textbackslash{}<subseteq>\allowbreak\ patch\allowbreak\ -\allowbreak{}\allowbreak\ cgu\_\allowbreak{}subdivision\_\allowbreak{}edge\_\allowbreak{}corner\_\allowbreak{}set\allowbreak\ M}; \texttt{cgu\_\allowbreak{}domain\allowbreak\ M\allowbreak\ \textbackslash{}<inter>\allowbreak\ ball\allowbreak\ center\allowbreak\ rad\allowbreak\ =\allowbreak\ ball\allowbreak\ center\allowbreak\ rad\allowbreak\ \textbackslash{}<inter>\allowbreak\ \{x.\allowbreak{}\allowbreak\ 0\allowbreak\ <\allowbreak\ cgu\_\allowbreak{}normal\_\allowbreak{}height\allowbreak\ center\allowbreak\ normal\allowbreak\ x\}}; \texttt{adjacent\allowbreak\ \textbackslash{}<in>\allowbreak\ cgu\_\allowbreak{}before\_\allowbreak{}indices\allowbreak\ M\allowbreak\ r}; and \texttt{cgu\_\allowbreak{}flat\_\allowbreak{}disc\allowbreak\ center\allowbreak\ normal\allowbreak\ rad\allowbreak\ \textbackslash{}<subseteq>\allowbreak\ frontier\allowbreak\ (cgu\_\allowbreak{}domain\allowbreak\ M)\allowbreak\ \textbackslash{}<inter>\allowbreak\ frontier\allowbreak\ (cgu\_\allowbreak{}cell\allowbreak\ M\allowbreak\ adjacent)}.

\dossierentry{\texttt{\detokenize{cgu_recovery_geometry_operational_v6}}}
\reviewlabel{Isabelle code}
\begin{lstlisting}
definition cgu_recovery_geometry_operational_v6 ::
  "('n::finite, 'c::finite, 'h::finite) cgu_model \<Rightarrow> bool"
where
  "cgu_recovery_geometry_operational_v6 M \<longleftrightarrow>
    cgu_order_valid M \<and>
    (\<forall>r<CARD('c).
      let c = cgu_order M r;
          remaining = cgu_remaining_domain M r;
          recovered = cgu_recovered_region M r;
          I = cgu_faces M c
      in cgu_lipschitz_domain remaining \<and>
         (r = 0 \<or>
          (recovered \<noteq> {} \<and>
           cgu_lipschitz_domain recovered \<and>
           cgu_face_connected recovered (cgu_before_indices M r) (cgu_cell M) \<and>
           cgu_auxiliary_half_ball_geometry_step M r)) \<and>
         (\<forall>i\<in>I.
           cgu_flat_face_patch (cgu_cell M c)
             (cgu_normal M c i) (cgu_offset M c i) (cgu_face_patch M c i) \<and>
           cgu_face_patch M c i \<subseteq>
             - cgu_subdivision_edge_corner_set M \<and>
           cgu_face_patch M c i \<subseteq> frontier remaining \<and>
           cgu_active_pde_face_route_v4 M r c i) \<and>
         cgu_degree M c + 3 \<le> card I \<and>
         cgu_n_generic_on (cgu_degree M c) I
           (cgu_normal M c) (cgu_offset M c) (cgu_witness_faces M c) \<and>
         (\<forall>i\<in>I. \<forall>j\<in>cgu_witness_faces M c i.
           \<forall>k\<in>cgu_witness_faces M c i. j \<noteq> k \<longrightarrow>
           cgu_triple_region_geometry (cgu_triple_region M c i j k)
             (cgu_normal M c i) (cgu_offset M c i)
             (cgu_normal M c j) (cgu_offset M c j)
             (cgu_normal M c k) (cgu_offset M c k))) \<and>
    (\<forall>r<CARD('c).
      let remaining = cgu_remaining_domain M r;
          recovered = cgu_recovered_region M r
      in cgu_graph_lipschitz_domain_v5 remaining \<and>
         (r = 0 \<or> cgu_graph_lipschitz_domain_v5 recovered))"
\end{lstlisting}
\reviewlabel{English mathematical translation}
\noindent\textit{Mathematical role:} admissible recovery order\par\smallskip
Typed argument/result: \texttt{M\allowbreak\ ::\allowbreak\ ('n,'c,'h)\allowbreak\ cgu\_\allowbreak{}model}, with \texttt{'n,'c,'h::finite}; result \texttt{bool}. It holds exactly when \texttt{cgu\_\allowbreak{}order\_\allowbreak{}valid\allowbreak\ M} and both following universally quantified blocks hold.

First, for every \texttt{r\allowbreak\ ::\allowbreak\ nat} restricted by \texttt{r\allowbreak\ <\allowbreak\ CARD('c)}, let \texttt{c=cgu\_\allowbreak{}order\allowbreak\ M\allowbreak\ r}, \texttt{remaining=cgu\_\allowbreak{}remaining\_\allowbreak{}domain\allowbreak\ M\allowbreak\ r}, \texttt{recovered=cgu\_\allowbreak{}recovered\_\allowbreak{}region\allowbreak\ M\allowbreak\ r}, and \texttt{I=cgu\_\allowbreak{}faces\allowbreak\ M\allowbreak\ c}. Then: \texttt{cgu\_\allowbreak{}lipschitz\_\allowbreak{}domain\allowbreak\ remaining}; either \texttt{r=0} or all of \texttt{recovered\allowbreak\ \textbackslash{}<noteq>\allowbreak\ \{\}}, \texttt{cgu\_\allowbreak{}lipschitz\_\allowbreak{}domain\allowbreak\ recovered}, \texttt{cgu\_\allowbreak{}face\_\allowbreak{}connected\allowbreak\ recovered\allowbreak\ (cgu\_\allowbreak{}before\_\allowbreak{}indices\allowbreak\ M\allowbreak\ r)\allowbreak\ (cgu\_\allowbreak{}cell\allowbreak\ M)}, and \texttt{cgu\_\allowbreak{}auxiliary\_\allowbreak{}half\_\allowbreak{}ball\_\allowbreak{}geometry\_\allowbreak{}step\allowbreak\ M\allowbreak\ r}; for every \texttt{i\allowbreak\ \textbackslash{}<in>\allowbreak\ I}, all of \texttt{cgu\_\allowbreak{}flat\_\allowbreak{}face\_\allowbreak{}patch\allowbreak\ (cgu\_\allowbreak{}cell\allowbreak\ M\allowbreak\ c)\allowbreak\ (cgu\_\allowbreak{}normal\allowbreak\ M\allowbreak\ c\allowbreak\ i)\allowbreak\ (cgu\_\allowbreak{}offset\allowbreak\ M\allowbreak\ c\allowbreak\ i)\allowbreak\ (cgu\_\allowbreak{}face\_\allowbreak{}patch\allowbreak\ M\allowbreak\ c\allowbreak\ i)}, \texttt{cgu\_\allowbreak{}face\_\allowbreak{}patch\allowbreak\ M\allowbreak\ c\allowbreak\ i\allowbreak\ \textbackslash{}<subseteq>\allowbreak\ -\allowbreak{}\allowbreak\ cgu\_\allowbreak{}subdivision\_\allowbreak{}edge\_\allowbreak{}corner\_\allowbreak{}set\allowbreak\ M}, \texttt{cgu\_\allowbreak{}face\_\allowbreak{}patch\allowbreak\ M\allowbreak\ c\allowbreak\ i\allowbreak\ \textbackslash{}<subseteq>\allowbreak\ frontier\allowbreak\ remaining}, and \texttt{cgu\_\allowbreak{}active\_\allowbreak{}pde\_\allowbreak{}face\_\allowbreak{}route\_\allowbreak{}v4\allowbreak\ M\allowbreak\ r\allowbreak\ c\allowbreak\ i}; \texttt{cgu\_\allowbreak{}degree\allowbreak\ M\allowbreak\ c\allowbreak\ +\allowbreak\ 3\allowbreak\ \textbackslash{}<le>\allowbreak\ card\allowbreak\ I}; \texttt{cgu\_\allowbreak{}n\_\allowbreak{}generic\_\allowbreak{}on\allowbreak\ (cgu\_\allowbreak{}degree\allowbreak\ M\allowbreak\ c)\allowbreak\ I\allowbreak\ (cgu\_\allowbreak{}normal\allowbreak\ M\allowbreak\ c)\allowbreak\ (cgu\_\allowbreak{}offset\allowbreak\ M\allowbreak\ c)\allowbreak\ (cgu\_\allowbreak{}witness\_\allowbreak{}faces\allowbreak\ M\allowbreak\ c)}; and for every \texttt{i\allowbreak\ \textbackslash{}<in>\allowbreak\ I}, \texttt{j\allowbreak\ \textbackslash{}<in>\allowbreak\ cgu\_\allowbreak{}witness\_\allowbreak{}faces\allowbreak\ M\allowbreak\ c\allowbreak\ i}, and \texttt{k\allowbreak\ \textbackslash{}<in>\allowbreak\ cgu\_\allowbreak{}witness\_\allowbreak{}faces\allowbreak\ M\allowbreak\ c\allowbreak\ i}, the hypothesis \texttt{j\allowbreak\ \textbackslash{}<noteq>\allowbreak\ k} implies \texttt{cgu\_\allowbreak{}triple\_\allowbreak{}region\_\allowbreak{}geometry\allowbreak\ (cgu\_\allowbreak{}triple\_\allowbreak{}region\allowbreak\ M\allowbreak\ c\allowbreak\ i\allowbreak\ j\allowbreak\ k)\allowbreak\ (cgu\_\allowbreak{}normal\allowbreak\ M\allowbreak\ c\allowbreak\ i)\allowbreak\ (cgu\_\allowbreak{}offset\allowbreak\ M\allowbreak\ c\allowbreak\ i)\allowbreak\ (cgu\_\allowbreak{}normal\allowbreak\ M\allowbreak\ c\allowbreak\ j)\allowbreak\ (cgu\_\allowbreak{}offset\allowbreak\ M\allowbreak\ c\allowbreak\ j)\allowbreak\ (cgu\_\allowbreak{}normal\allowbreak\ M\allowbreak\ c\allowbreak\ k)\allowbreak\ (cgu\_\allowbreak{}offset\allowbreak\ M\allowbreak\ c\allowbreak\ k)}.

Second, for every \texttt{r\allowbreak\ ::\allowbreak\ nat} restricted by \texttt{r<CARD('c)}, let \texttt{remaining=cgu\_\allowbreak{}remaining\_\allowbreak{}domain\allowbreak\ M\allowbreak\ r} and \texttt{recovered=cgu\_\allowbreak{}recovered\_\allowbreak{}region\allowbreak\ M\allowbreak\ r}; then \texttt{cgu\_\allowbreak{}graph\_\allowbreak{}lipschitz\_\allowbreak{}domain\_\allowbreak{}v5\allowbreak\ remaining} and either \texttt{r=0} or \texttt{cgu\_\allowbreak{}graph\_\allowbreak{}lipschitz\_\allowbreak{}domain\_\allowbreak{}v5\allowbreak\ recovered}.

\dossierentry{\texttt{\detokenize{cgu_auxiliary_extension_elliptic}}}
\reviewlabel{Isabelle code}
\begin{lstlisting}
definition cgu_auxiliary_extension_elliptic ::
  "('n::finite, 'c::finite, 'h) cgu_model \<Rightarrow>
    ('c \<Rightarrow> 'n cgu_matrix_polynomial) \<Rightarrow> nat \<Rightarrow> bool"
where
  "cgu_auxiliary_extension_elliptic M P r \<longleftrightarrow>
    (let c = cgu_order M r;
         adjacent = cgu_step_adjacent_cell M c
     in cgu_uniformly_positive_definite_on
          (cgu_exterior_half_ball (cgu_step_center M c)
            (cgu_step_interior_normal M c) (cgu_step_radius M c))
          (P adjacent))"
\end{lstlisting}
\reviewlabel{English mathematical translation}
\noindent\textit{Mathematical role:} adjacent-coefficient ellipticity on the recovery half-ball\par\smallskip
Typed arguments/result: \texttt{M\allowbreak\ ::\allowbreak\ ('n,'c,'h)\allowbreak\ cgu\_\allowbreak{}model}, \texttt{P\allowbreak\ ::\allowbreak\ 'c\allowbreak\ \textbackslash{}<Rightarrow>\allowbreak\ 'n\allowbreak\ cgu\_\allowbreak{}matrix\_\allowbreak{}polynomial}, \texttt{r\allowbreak\ ::\allowbreak\ nat}, with \texttt{'n,'c::finite}; result \texttt{bool}. Let \texttt{c=cgu\_\allowbreak{}order\allowbreak\ M\allowbreak\ r} and \texttt{adjacent=cgu\_\allowbreak{}step\_\allowbreak{}adjacent\_\allowbreak{}cell\allowbreak\ M\allowbreak\ c}. The condition is \texttt{cgu\_\allowbreak{}uniformly\_\allowbreak{}positive\_\allowbreak{}definite\_\allowbreak{}on\allowbreak\ (cgu\_\allowbreak{}exterior\_\allowbreak{}half\_\allowbreak{}ball\allowbreak\ (cgu\_\allowbreak{}step\_\allowbreak{}center\allowbreak\ M\allowbreak\ c)\allowbreak\ (cgu\_\allowbreak{}step\_\allowbreak{}interior\_\allowbreak{}normal\allowbreak\ M\allowbreak\ c)\allowbreak\ (cgu\_\allowbreak{}step\_\allowbreak{}radius\allowbreak\ M\allowbreak\ c))\allowbreak\ (P\allowbreak\ adjacent)}.

\dossierentry{\texttt{\detokenize{cgu_recovery_extension_admissible}}}
\reviewlabel{Isabelle code}
\begin{lstlisting}
definition cgu_recovery_extension_admissible ::
  "('n::finite, 'c::finite, 'h::finite) cgu_model \<Rightarrow>
    ('c \<Rightarrow> 'n cgu_matrix_polynomial) \<Rightarrow> bool"
where
  "cgu_recovery_extension_admissible M P \<longleftrightarrow>
    (\<forall>r<CARD('c).
      let c = cgu_order M r;
          I = cgu_faces M c
      in (r = 0 \<or> cgu_auxiliary_extension_elliptic M P r) \<and>
         (\<forall>i\<in>I. \<forall>j\<in>cgu_witness_faces M c i.
           \<forall>k\<in>cgu_witness_faces M c i. j \<noteq> k \<longrightarrow>
             cgu_polynomial_positive_on (P c)
              (cgu_triple_region M c i j k)))"
\end{lstlisting}
\reviewlabel{English mathematical translation}
\noindent\textit{Mathematical role:} coefficient-access condition\par\smallskip
Typed arguments/result: \texttt{M\allowbreak\ ::\allowbreak\ ('n,'c,'h)\allowbreak\ cgu\_\allowbreak{}model}, \texttt{P\allowbreak\ ::\allowbreak\ 'c\allowbreak\ \textbackslash{}<Rightarrow>\allowbreak\ 'n\allowbreak\ cgu\_\allowbreak{}matrix\_\allowbreak{}polynomial}, with \texttt{'n,'c,'h::finite}; result \texttt{bool}. It holds exactly when, for every \texttt{r\allowbreak\ ::\allowbreak\ nat} restricted by \texttt{r<CARD('c)}, let \texttt{c=cgu\_\allowbreak{}order\allowbreak\ M\allowbreak\ r} and \texttt{I=cgu\_\allowbreak{}faces\allowbreak\ M\allowbreak\ c}; then \texttt{(r=0\allowbreak\ \textbackslash{}<or>\allowbreak\ cgu\_\allowbreak{}auxiliary\_\allowbreak{}extension\_\allowbreak{}elliptic\allowbreak\ M\allowbreak\ P\allowbreak\ r)} and, for every \texttt{i\allowbreak\ \textbackslash{}<in>\allowbreak\ I}, \texttt{j\allowbreak\ \textbackslash{}<in>\allowbreak\ cgu\_\allowbreak{}witness\_\allowbreak{}faces\allowbreak\ M\allowbreak\ c\allowbreak\ i}, and \texttt{k\allowbreak\ \textbackslash{}<in>\allowbreak\ cgu\_\allowbreak{}witness\_\allowbreak{}faces\allowbreak\ M\allowbreak\ c\allowbreak\ i}, the hypothesis \texttt{j\allowbreak\ \textbackslash{}<noteq>\allowbreak\ k} implies \texttt{cgu\_\allowbreak{}polynomial\_\allowbreak{}positive\_\allowbreak{}on\allowbreak\ (P\allowbreak\ c)\allowbreak\ (cgu\_\allowbreak{}triple\_\allowbreak{}region\allowbreak\ M\allowbreak\ c\allowbreak\ i\allowbreak\ j\allowbreak\ k)}.

\dossierentry{\texttt{\detokenize{cgu_potential}}}
\reviewlabel{Isabelle code}
\begin{lstlisting}
type_synonym 'n cgu_potential = "'n cgu_point \<Rightarrow> real"
\end{lstlisting}
\reviewlabel{English mathematical translation}
\noindent\textit{Mathematical role:} scalar-function type \(\mathcal F_I\)\par\smallskip
For every type \texttt{'n}, \texttt{'n\allowbreak\ cgu\_\allowbreak{}potential} is exactly the function type \texttt{'n\allowbreak\ cgu\_\allowbreak{}point\allowbreak\ \textbackslash{}<Rightarrow>\allowbreak\ real}.

\dossierentry{\texttt{\detokenize{cgu_gradient}}}
\reviewlabel{Isabelle code}
\begin{lstlisting}
type_synonym 'n cgu_gradient = "'n cgu_point \<Rightarrow> 'n cgu_point"
\end{lstlisting}
\reviewlabel{English mathematical translation}
\noindent\textit{Mathematical role:} vector-field type \(\mathcal V_I\)\par\smallskip
For every type \texttt{'n}, \texttt{'n\allowbreak\ cgu\_\allowbreak{}gradient} is exactly the function type \texttt{'n\allowbreak\ cgu\_\allowbreak{}point\allowbreak\ \textbackslash{}<Rightarrow>\allowbreak\ 'n\allowbreak\ cgu\_\allowbreak{}point}.

\dossierentry{\texttt{\detokenize{cgu_partial_derivative}}}
\reviewlabel{Isabelle code}
\begin{lstlisting}
definition cgu_partial_derivative ::
  "('n::finite cgu_potential) \<Rightarrow> 'n \<Rightarrow> 'n cgu_point \<Rightarrow> real"
where
  "cgu_partial_derivative phi i x =
    frechet_derivative phi (at x) (axis i 1)"
\end{lstlisting}
\reviewlabel{English mathematical translation}
\noindent\textit{Mathematical role:} selected coordinate derivative \(\partial_i\phi\)\par\smallskip
Typed arguments/result: \texttt{phi\allowbreak\ ::\allowbreak\ 'n\allowbreak\ cgu\_\allowbreak{}potential}, \texttt{i\allowbreak\ ::\allowbreak\ 'n}, \texttt{x\allowbreak\ ::\allowbreak\ 'n\allowbreak\ cgu\_\allowbreak{}point}, with \texttt{'n::finite}; result \texttt{real}. It equals \texttt{frechet\_\allowbreak{}derivative\allowbreak\ phi\allowbreak\ (at\allowbreak\ x)\allowbreak\ (axis\allowbreak\ i\allowbreak\ 1)}.

\dossierentry{\texttt{\detokenize{cgu_classical_gradient}}}
\reviewlabel{Isabelle code}
\begin{lstlisting}
definition cgu_classical_gradient ::
  "'n::finite cgu_potential \<Rightarrow> 'n cgu_gradient"
where
  "cgu_classical_gradient phi x =
    (\<chi> i. cgu_partial_derivative phi i x)"
\end{lstlisting}
\reviewlabel{English mathematical translation}
\noindent\textit{Mathematical role:} gradient vector \(\nabla\phi\)\par\smallskip
Typed argument/result: \texttt{phi\allowbreak\ ::\allowbreak\ 'n\allowbreak\ cgu\_\allowbreak{}potential}, with \texttt{'n::finite}; result \texttt{'n\allowbreak\ cgu\_\allowbreak{}gradient}. For every \texttt{x\allowbreak\ ::\allowbreak\ 'n\allowbreak\ cgu\_\allowbreak{}point}, its value is the indexed vector \texttt{\textbackslash{}<chi>\allowbreak\ i.\allowbreak{}\allowbreak\ cgu\_\allowbreak{}partial\_\allowbreak{}derivative\allowbreak\ phi\allowbreak\ i\allowbreak\ x}, with \texttt{i\allowbreak\ ::\allowbreak\ 'n}.

\dossierentry{\texttt{\detokenize{cgu_test_function_on}}}
\reviewlabel{Isabelle code}
\begin{lstlisting}
definition cgu_test_function_on ::
  "'n::finite cgu_point set \<Rightarrow> 'n cgu_potential \<Rightarrow> bool"
where
  "cgu_test_function_on U phi \<longleftrightarrow>
    smooth_on UNIV phi \<and>
    compact (closure {x. phi x \<noteq> 0}) \<and>
    closure {x. phi x \<noteq> 0} \<subseteq> U"
\end{lstlisting}
\reviewlabel{English mathematical translation}
\noindent\textit{Mathematical role:} test class \(\mathcal D(U)\)\par\smallskip
Typed arguments/result: \texttt{U\allowbreak\ ::\allowbreak\ 'n\allowbreak\ cgu\_\allowbreak{}point\allowbreak\ set}, \texttt{phi\allowbreak\ ::\allowbreak\ 'n\allowbreak\ cgu\_\allowbreak{}potential}, with \texttt{'n::finite}; result \texttt{bool}. It holds exactly when \texttt{smooth\_\allowbreak{}on\allowbreak\ UNIV\allowbreak\ phi}, \texttt{compact\allowbreak\ (closure\allowbreak\ \{x.\allowbreak{}\allowbreak\ phi\allowbreak\ x\allowbreak\ \textbackslash{}<noteq>\allowbreak\ 0\})}, and \texttt{closure\allowbreak\ \{x.\allowbreak{}\allowbreak\ phi\allowbreak\ x\allowbreak\ \textbackslash{}<noteq>\allowbreak\ 0\}\allowbreak\ \textbackslash{}<subseteq>\allowbreak\ U}.

\dossierentry{\texttt{\detokenize{cgu_weak_gradient_on}}}
\reviewlabel{Isabelle code}
\begin{lstlisting}
definition cgu_weak_gradient_on ::
  "'n::finite cgu_point set \<Rightarrow> 'n cgu_potential \<Rightarrow>
    'n cgu_gradient \<Rightarrow> bool"
where
  "cgu_weak_gradient_on U u Du \<longleftrightarrow>
    (\<forall>phi. cgu_test_function_on U phi \<longrightarrow>
      (\<forall>i.
        set_integrable lborel U
          (\<lambda>x. u x * cgu_partial_derivative phi i x) \<and>
        set_integrable lborel U (\<lambda>x. Du x $ i * phi x) \<and>
        set_lebesgue_integral lborel U
          (\<lambda>x. u x * cgu_partial_derivative phi i x) =
        - set_lebesgue_integral lborel U (\<lambda>x. Du x $ i * phi x)))"
\end{lstlisting}
\reviewlabel{English mathematical translation}
\noindent\textit{Mathematical role:} weak-gradient integration-by-parts relation\par\smallskip
Typed arguments/result: \texttt{U\allowbreak\ ::\allowbreak\ 'n\allowbreak\ cgu\_\allowbreak{}point\allowbreak\ set}, \texttt{u\allowbreak\ ::\allowbreak\ 'n\allowbreak\ cgu\_\allowbreak{}potential}, \texttt{Du\allowbreak\ ::\allowbreak\ 'n\allowbreak\ cgu\_\allowbreak{}gradient}, with \texttt{'n::finite}; result \texttt{bool}. It holds exactly when, for every \texttt{phi\allowbreak\ ::\allowbreak\ 'n\allowbreak\ cgu\_\allowbreak{}potential}, the hypothesis \texttt{cgu\_\allowbreak{}test\_\allowbreak{}function\_\allowbreak{}on\allowbreak\ U\allowbreak\ phi} implies that for every \texttt{i\allowbreak\ ::\allowbreak\ 'n}: both functions \texttt{x\allowbreak\ \textbackslash{}<mapsto>\allowbreak\ u\allowbreak\ x\allowbreak\ *\allowbreak\ cgu\_\allowbreak{}partial\_\allowbreak{}derivative\allowbreak\ phi\allowbreak\ i\allowbreak\ x} and \texttt{x\allowbreak\ \textbackslash{}<mapsto>\allowbreak\ Du\allowbreak\ x\allowbreak\ \$\allowbreak\ i\allowbreak\ *\allowbreak\ phi\allowbreak\ x} are set-integrable with respect to \texttt{lborel} on \texttt{U}, and their set Lebesgue integrals satisfy \texttt{integral\_\allowbreak{}U\allowbreak\ (u*cgu\_\allowbreak{}partial\_\allowbreak{}derivative\allowbreak\ phi\allowbreak\ i)\allowbreak\ =\allowbreak\ -\allowbreak{}\allowbreak\ integral\_\allowbreak{}U\allowbreak\ ((Du\allowbreak\ coordinate\allowbreak\ i)*phi)}.

\dossierentry{\texttt{\detokenize{cgu_h1_pair_on}}}
\reviewlabel{Isabelle code}
\begin{lstlisting}
definition cgu_h1_pair_on ::
  "'n::finite cgu_point set \<Rightarrow> 'n cgu_potential \<Rightarrow>
    'n cgu_gradient \<Rightarrow> bool"
where
  "cgu_h1_pair_on U u Du \<longleftrightarrow>
    u \<in> borel_measurable (restrict_space lborel U) \<and>
    Du \<in> borel_measurable (restrict_space lborel U) \<and>
    cgu_weak_gradient_on U u Du \<and>
    set_integrable lborel U (\<lambda>x. u x ^ 2) \<and>
    set_integrable lborel U (\<lambda>x. norm (Du x) ^ 2)"
\end{lstlisting}
\reviewlabel{English mathematical translation}
\noindent\textit{Mathematical role:} pair class \(\mathcal H^1_{\mathrm{pair}}(U)\)\par\smallskip
Typed arguments/result: \texttt{U\allowbreak\ ::\allowbreak\ 'n\allowbreak\ cgu\_\allowbreak{}point\allowbreak\ set}, \texttt{u\allowbreak\ ::\allowbreak\ 'n\allowbreak\ cgu\_\allowbreak{}potential}, \texttt{Du\allowbreak\ ::\allowbreak\ 'n\allowbreak\ cgu\_\allowbreak{}gradient}, with \texttt{'n::finite}; result \texttt{bool}. It holds exactly when \texttt{u\allowbreak\ \textbackslash{}<in>\allowbreak\ borel\_\allowbreak{}measurable\allowbreak\ (restrict\_\allowbreak{}space\allowbreak\ lborel\allowbreak\ U)}, \texttt{Du\allowbreak\ \textbackslash{}<in>\allowbreak\ borel\_\allowbreak{}measurable\allowbreak\ (restrict\_\allowbreak{}space\allowbreak\ lborel\allowbreak\ U)}, \texttt{cgu\_\allowbreak{}weak\_\allowbreak{}gradient\_\allowbreak{}on\allowbreak\ U\allowbreak\ u\allowbreak\ Du}, \texttt{x\allowbreak\ \textbackslash{}<mapsto>\allowbreak\ u\allowbreak\ x\textasciicircum{}2} is set-integrable on \texttt{U}, and \texttt{x\allowbreak\ \textbackslash{}<mapsto>\allowbreak\ norm\allowbreak\ (Du\allowbreak\ x)\textasciicircum{}2} is set-integrable on \texttt{U}.

\dossierentry{\texttt{\detokenize{cgu_h1_squared_distance}}}
\reviewlabel{Isabelle code}
\begin{lstlisting}
definition cgu_h1_squared_distance ::
  "'n::finite cgu_point set \<Rightarrow>
    'n cgu_potential \<Rightarrow> 'n cgu_gradient \<Rightarrow>
    'n cgu_potential \<Rightarrow> 'n cgu_gradient \<Rightarrow> real"
where
  "cgu_h1_squared_distance U u Du v Dv =
    set_lebesgue_integral lborel U
      (\<lambda>x. (u x - v x) ^ 2 + norm (Du x - Dv x) ^ 2)"
\end{lstlisting}
\reviewlabel{English mathematical translation}
\noindent\textit{Mathematical role:} squared pair distance \(d_U^2\)\par\smallskip
Typed arguments/result: \texttt{U\allowbreak\ ::\allowbreak\ 'n\allowbreak\ cgu\_\allowbreak{}point\allowbreak\ set}; \texttt{u,v\allowbreak\ ::\allowbreak\ 'n\allowbreak\ cgu\_\allowbreak{}potential}; \texttt{Du,Dv\allowbreak\ ::\allowbreak\ 'n\allowbreak\ cgu\_\allowbreak{}gradient}; with \texttt{'n::finite}; result \texttt{real}. It equals the set Lebesgue integral over \texttt{U} of \texttt{x\allowbreak\ \textbackslash{}<mapsto>\allowbreak\ (u\allowbreak\ x-\allowbreak{}v\allowbreak\ x)\textasciicircum{}2\allowbreak\ +\allowbreak\ norm\allowbreak\ (Du\allowbreak\ x-\allowbreak{}Dv\allowbreak\ x)\textasciicircum{}2}.

\dossierentry{\texttt{\detokenize{cgu_h1_zero_pair_on}}}
\reviewlabel{Isabelle code}
\begin{lstlisting}
definition cgu_h1_zero_pair_on ::
  "(real ^ 'n::finite) set \<Rightarrow>
    ((real ^ 'n) \<Rightarrow> real) \<Rightarrow>
    ((real ^ 'n) \<Rightarrow> (real ^ 'n)) \<Rightarrow> bool"
where
  "cgu_h1_zero_pair_on U u Du \<longleftrightarrow>
    cgu_h1_pair_on U u Du \<and>
    (\<exists>(seq :: nat \<Rightarrow> 'n cgu_potential).
      (\<forall>m. cgu_test_function_on U (seq m) \<and>
        cgu_h1_pair_on U (seq m) (cgu_classical_gradient (seq m))) \<and>
      (\<forall>epsilon>0. \<exists>N. \<forall>m\<ge>N.
        cgu_h1_squared_distance U (seq m)
          (cgu_classical_gradient (seq m)) u Du < epsilon))"
\end{lstlisting}
\reviewlabel{English mathematical translation}
\noindent\textit{Mathematical role:} zero-boundary pair class \(\mathcal H^1_{0,\mathrm{pair}}(U)\)\par\smallskip
Typed arguments/result: \texttt{U\allowbreak\ ::\allowbreak\ (real\allowbreak\ \textasciicircum{}\allowbreak\ 'n)\allowbreak\ set}, \texttt{u\allowbreak\ ::\allowbreak\ real\allowbreak\ \textasciicircum{}\allowbreak\ 'n\allowbreak\ \textbackslash{}<Rightarrow>\allowbreak\ real}, \texttt{Du\allowbreak\ ::\allowbreak\ real\allowbreak\ \textasciicircum{}\allowbreak\ 'n\allowbreak\ \textbackslash{}<Rightarrow>\allowbreak\ real\allowbreak\ \textasciicircum{}\allowbreak\ 'n}, with \texttt{'n::finite}; result \texttt{bool}. It holds exactly when \texttt{cgu\_\allowbreak{}h1\_\allowbreak{}pair\_\allowbreak{}on\allowbreak\ U\allowbreak\ u\allowbreak\ Du} and there exists \texttt{seq\allowbreak\ ::\allowbreak\ nat\allowbreak\ \textbackslash{}<Rightarrow>\allowbreak\ 'n\allowbreak\ cgu\_\allowbreak{}potential} such that: for every \texttt{m\allowbreak\ ::\allowbreak\ nat}, both \texttt{cgu\_\allowbreak{}test\_\allowbreak{}function\_\allowbreak{}on\allowbreak\ U\allowbreak\ (seq\allowbreak\ m)} and \texttt{cgu\_\allowbreak{}h1\_\allowbreak{}pair\_\allowbreak{}on\allowbreak\ U\allowbreak\ (seq\allowbreak\ m)\allowbreak\ (cgu\_\allowbreak{}classical\_\allowbreak{}gradient\allowbreak\ (seq\allowbreak\ m))}; and for every \texttt{epsilon\allowbreak\ ::\allowbreak\ real} restricted by \texttt{epsilon>0}, there exists \texttt{N\allowbreak\ ::\allowbreak\ nat} such that for every \texttt{m\allowbreak\ ::\allowbreak\ nat} restricted by \texttt{m\textbackslash{}<ge>N}, \texttt{cgu\_\allowbreak{}h1\_\allowbreak{}squared\_\allowbreak{}distance\allowbreak\ U\allowbreak\ (seq\allowbreak\ m)\allowbreak\ (cgu\_\allowbreak{}classical\_\allowbreak{}gradient\allowbreak\ (seq\allowbreak\ m))\allowbreak\ u\allowbreak\ Du\allowbreak\ <\allowbreak\ epsilon}.

\dossierentry{\texttt{\detokenize{cgu_cellwise_coefficient}}}
\reviewlabel{Isabelle code}
\begin{lstlisting}
definition cgu_cellwise_coefficient ::
  "('n::finite, 'c::finite, 'h) cgu_model \<Rightarrow>
    ('c \<Rightarrow> 'n cgu_matrix_polynomial) \<Rightarrow>
    'n cgu_point \<Rightarrow> real ^ 'n ^ 'n"
where
  "cgu_cellwise_coefficient M P x =
    (if \<exists>c. x \<in> cgu_cell M c
     then cgu_matrix_poly_eval (P (SOME c. x \<in> cgu_cell M c)) x
     else 0)"
\end{lstlisting}
\reviewlabel{English mathematical translation}
\noindent\textit{Mathematical role:} everywhere-defined selected-cell field \(A_{M,P}\)\par\smallskip
Typed arguments/result: \texttt{M\allowbreak\ ::\allowbreak\ ('n,'c,'h)\allowbreak\ cgu\_\allowbreak{}model}, \texttt{P\allowbreak\ ::\allowbreak\ 'c\allowbreak\ \textbackslash{}<Rightarrow>\allowbreak\ 'n\allowbreak\ cgu\_\allowbreak{}matrix\_\allowbreak{}polynomial}, \texttt{x\allowbreak\ ::\allowbreak\ 'n\allowbreak\ cgu\_\allowbreak{}point}, with \texttt{'n,'c::finite}; result \texttt{real\allowbreak\ \textasciicircum{}\allowbreak\ 'n\allowbreak\ \textasciicircum{}\allowbreak\ 'n}. If there exists \texttt{c\allowbreak\ ::\allowbreak\ 'c} with \texttt{x\allowbreak\ \textbackslash{}<in>\allowbreak\ cgu\_\allowbreak{}cell\allowbreak\ M\allowbreak\ c}, the result is \texttt{cgu\_\allowbreak{}matrix\_\allowbreak{}poly\_\allowbreak{}eval\allowbreak\ (P\allowbreak\ (SOME\allowbreak\ c.\allowbreak{}\allowbreak\ x\allowbreak\ \textbackslash{}<in>\allowbreak\ cgu\_\allowbreak{}cell\allowbreak\ M\allowbreak\ c))\allowbreak\ x}; otherwise it is \texttt{0}. \texttt{SOME} makes a Hilbert-choice selection from the stated existence condition; no uniqueness condition appears in this definition.

\dossierentry{\texttt{\detokenize{cgu_boundary_test_on}}}
\reviewlabel{Isabelle code}
\begin{lstlisting}
definition cgu_boundary_test_on ::
  "('n::finite, 'c::finite, 'h::finite) cgu_model \<Rightarrow>
    'n cgu_point set \<Rightarrow> 'n cgu_point set \<Rightarrow>
    'n cgu_potential \<Rightarrow> bool"
where
  "cgu_boundary_test_on M U gamma f \<longleftrightarrow>
    smooth_on UNIV f \<and>
    compact (closure {x. f x \<noteq> 0}) \<and>
    closure {x. x \<in> frontier U \<and> f x \<noteq> 0}
      \<subseteq> gamma - cgu_subdivision_edge_corner_set M \<and>
    cgu_h1_pair_on U f (cgu_classical_gradient f)"
\end{lstlisting}
\reviewlabel{English mathematical translation}
\noindent\textit{Mathematical role:} smooth boundary-supported datum\par\smallskip
Typed arguments/result: \texttt{M\allowbreak\ ::\allowbreak\ ('n,'c,'h)\allowbreak\ cgu\_\allowbreak{}model}, \texttt{U,gamma\allowbreak\ ::\allowbreak\ 'n\allowbreak\ cgu\_\allowbreak{}point\allowbreak\ set}, \texttt{f\allowbreak\ ::\allowbreak\ 'n\allowbreak\ cgu\_\allowbreak{}potential}, with \texttt{'n,'c,'h::finite}; result \texttt{bool}. It holds exactly when \texttt{smooth\_\allowbreak{}on\allowbreak\ UNIV\allowbreak\ f}, \texttt{compact\allowbreak\ (closure\allowbreak\ \{x.\allowbreak{}\allowbreak\ f\allowbreak\ x\allowbreak\ \textbackslash{}<noteq>\allowbreak\ 0\})}, \texttt{closure\allowbreak\ \{x.\allowbreak{}\allowbreak\ x\allowbreak\ \textbackslash{}<in>\allowbreak\ frontier\allowbreak\ U\allowbreak\ \textbackslash{}<and>\allowbreak\ f\allowbreak\ x\allowbreak\ \textbackslash{}<noteq>\allowbreak\ 0\}\allowbreak\ \textbackslash{}<subseteq>\allowbreak\ gamma\allowbreak\ -\allowbreak{}\allowbreak\ cgu\_\allowbreak{}subdivision\_\allowbreak{}edge\_\allowbreak{}corner\_\allowbreak{}set\allowbreak\ M}, and \texttt{cgu\_\allowbreak{}h1\_\allowbreak{}pair\_\allowbreak{}on\allowbreak\ U\allowbreak\ f\allowbreak\ (cgu\_\allowbreak{}classical\_\allowbreak{}gradient\allowbreak\ f)}.

\dossierentry{\texttt{\detokenize{cgu_local_trace_pair_on}}}
\reviewlabel{Isabelle code}
\begin{lstlisting}
definition cgu_local_trace_pair_on ::
  "('n::finite, 'c::finite, 'h::finite) cgu_model \<Rightarrow>
    'n cgu_point set \<Rightarrow> 'n cgu_point set \<Rightarrow>
    'n cgu_potential \<Rightarrow> 'n cgu_gradient \<Rightarrow> bool"
where
  "cgu_local_trace_pair_on M U gamma f Df \<longleftrightarrow>
    cgu_h1_pair_on U f Df \<and>
    (\<exists>(seq :: nat \<Rightarrow> 'n cgu_potential).
      (\<forall>m. cgu_boundary_test_on M U gamma (seq m)) \<and>
      (\<forall>epsilon>0. \<exists>N. \<forall>m\<ge>N.
        \<exists>z Dz. cgu_h1_zero_pair_on U z Dz \<and>
          cgu_h1_squared_distance U (seq m)
            (cgu_classical_gradient (seq m))
            (\<lambda>x. f x + z x) (\<lambda>x. Df x + Dz x) < epsilon))"
\end{lstlisting}
\reviewlabel{English mathematical translation}
\noindent\textit{Mathematical role:} boundary-data closure predicate \(\mathcal T_M(U,\gamma)\)\par\smallskip
Typed arguments/result: \texttt{M\allowbreak\ ::\allowbreak\ ('n,'c,'h)\allowbreak\ cgu\_\allowbreak{}model}, \texttt{U,gamma\allowbreak\ ::\allowbreak\ 'n\allowbreak\ cgu\_\allowbreak{}point\allowbreak\ set}, \texttt{f\allowbreak\ ::\allowbreak\ 'n\allowbreak\ cgu\_\allowbreak{}potential}, \texttt{Df\allowbreak\ ::\allowbreak\ 'n\allowbreak\ cgu\_\allowbreak{}gradient}, with \texttt{'n,'c,'h::finite}; result \texttt{bool}. It holds exactly when \texttt{cgu\_\allowbreak{}h1\_\allowbreak{}pair\_\allowbreak{}on\allowbreak\ U\allowbreak\ f\allowbreak\ Df} and there exists \texttt{seq\allowbreak\ ::\allowbreak\ nat\allowbreak\ \textbackslash{}<Rightarrow>\allowbreak\ 'n\allowbreak\ cgu\_\allowbreak{}potential} such that every \texttt{m} satisfies \texttt{cgu\_\allowbreak{}boundary\_\allowbreak{}test\_\allowbreak{}on\allowbreak\ M\allowbreak\ U\allowbreak\ gamma\allowbreak\ (seq\allowbreak\ m)}, and for every \texttt{epsilon>0} there exists \texttt{N\allowbreak\ ::\allowbreak\ nat} such that for every \texttt{m\textbackslash{}<ge>N} there exist \texttt{z\allowbreak\ ::\allowbreak\ 'n\allowbreak\ cgu\_\allowbreak{}potential} and \texttt{Dz\allowbreak\ ::\allowbreak\ 'n\allowbreak\ cgu\_\allowbreak{}gradient} with \texttt{cgu\_\allowbreak{}h1\_\allowbreak{}zero\_\allowbreak{}pair\_\allowbreak{}on\allowbreak\ U\allowbreak\ z\allowbreak\ Dz} and \texttt{cgu\_\allowbreak{}h1\_\allowbreak{}squared\_\allowbreak{}distance\allowbreak\ U\allowbreak\ (seq\allowbreak\ m)\allowbreak\ (cgu\_\allowbreak{}classical\_\allowbreak{}gradient\allowbreak\ (seq\allowbreak\ m))\allowbreak\ (\textbackslash{}<lambda>x.\allowbreak{}\allowbreak\ f\allowbreak\ x+z\allowbreak\ x)\allowbreak\ (\textbackslash{}<lambda>x.\allowbreak{}\allowbreak\ Df\allowbreak\ x+Dz\allowbreak\ x)\allowbreak\ <\allowbreak\ epsilon}.

\dossierentry{\texttt{\detokenize{cgu_weak_solution_on}}}
\reviewlabel{Isabelle code}
\begin{lstlisting}
definition cgu_weak_solution_on ::
  "('n::finite, 'c::finite, 'h) cgu_model \<Rightarrow>
    ('c \<Rightarrow> 'n cgu_matrix_polynomial) \<Rightarrow>
    'n cgu_point set \<Rightarrow> 'n cgu_potential \<Rightarrow>
    'n cgu_gradient \<Rightarrow> bool"
where
  "cgu_weak_solution_on M P U u Du \<longleftrightarrow>
    cgu_h1_pair_on U u Du \<and>
    (\<forall>phi. cgu_test_function_on U phi \<longrightarrow>
      set_integrable lborel U
        (\<lambda>x. inner (cgu_cellwise_coefficient M P x *v Du x)
          (cgu_classical_gradient phi x)) \<and>
      set_lebesgue_integral lborel U
        (\<lambda>x. inner (cgu_cellwise_coefficient M P x *v Du x)
          (cgu_classical_gradient phi x)) = 0)"
\end{lstlisting}
\reviewlabel{English mathematical translation}
\noindent\textit{Mathematical role:} \(A_{M,P}\)-harmonic pair\par\smallskip
Typed arguments/result: \texttt{M\allowbreak\ ::\allowbreak\ ('n,'c,'h)\allowbreak\ cgu\_\allowbreak{}model}, \texttt{P\allowbreak\ ::\allowbreak\ 'c\allowbreak\ \textbackslash{}<Rightarrow>\allowbreak\ 'n\allowbreak\ cgu\_\allowbreak{}matrix\_\allowbreak{}polynomial}, \texttt{U\allowbreak\ ::\allowbreak\ 'n\allowbreak\ cgu\_\allowbreak{}point\allowbreak\ set}, \texttt{u\allowbreak\ ::\allowbreak\ 'n\allowbreak\ cgu\_\allowbreak{}potential}, \texttt{Du\allowbreak\ ::\allowbreak\ 'n\allowbreak\ cgu\_\allowbreak{}gradient}, with \texttt{'n,'c::finite}; result \texttt{bool}. It holds exactly when \texttt{cgu\_\allowbreak{}h1\_\allowbreak{}pair\_\allowbreak{}on\allowbreak\ U\allowbreak\ u\allowbreak\ Du} and, for every \texttt{phi}, the hypothesis \texttt{cgu\_\allowbreak{}test\_\allowbreak{}function\_\allowbreak{}on\allowbreak\ U\allowbreak\ phi} implies both set-integrability on \texttt{U} of \texttt{x\allowbreak\ \textbackslash{}<mapsto>\allowbreak\ inner\allowbreak\ (cgu\_\allowbreak{}cellwise\_\allowbreak{}coefficient\allowbreak\ M\allowbreak\ P\allowbreak\ x\allowbreak\ *v\allowbreak\ Du\allowbreak\ x)\allowbreak\ (cgu\_\allowbreak{}classical\_\allowbreak{}gradient\allowbreak\ phi\allowbreak\ x)} and equality of that function's set Lebesgue integral to \texttt{0}.

\dossierentry{\texttt{\detokenize{cgu_dirichlet_solution_for}}}
\reviewlabel{Isabelle code}
\begin{lstlisting}
definition cgu_dirichlet_solution_for ::
  "('n::finite, 'c::finite, 'h) cgu_model \<Rightarrow>
    ('c \<Rightarrow> 'n cgu_matrix_polynomial) \<Rightarrow>
    'n cgu_point set \<Rightarrow> 'n cgu_potential \<Rightarrow>
    'n cgu_gradient \<Rightarrow> 'n cgu_potential \<Rightarrow>
    'n cgu_gradient \<Rightarrow> bool"
where
  "cgu_dirichlet_solution_for M P U f Df u Du \<longleftrightarrow>
    cgu_h1_pair_on U f Df \<and>
    cgu_weak_solution_on M P U u Du \<and>
    cgu_h1_zero_pair_on U (\<lambda>x. u x - f x)
      (\<lambda>x. Du x - Df x)"
\end{lstlisting}
\reviewlabel{English mathematical translation}
\noindent\textit{Mathematical role:} \(A_{M,P}\)-Dirichlet pair\par\smallskip
Typed arguments/result: \texttt{M\allowbreak\ ::\allowbreak\ ('n,'c,'h)\allowbreak\ cgu\_\allowbreak{}model}, \texttt{P\allowbreak\ ::\allowbreak\ 'c\allowbreak\ \textbackslash{}<Rightarrow>\allowbreak\ 'n\allowbreak\ cgu\_\allowbreak{}matrix\_\allowbreak{}polynomial}, \texttt{U\allowbreak\ ::\allowbreak\ 'n\allowbreak\ cgu\_\allowbreak{}point\allowbreak\ set}, \texttt{f,u\allowbreak\ ::\allowbreak\ 'n\allowbreak\ cgu\_\allowbreak{}potential}, \texttt{Df,Du\allowbreak\ ::\allowbreak\ 'n\allowbreak\ cgu\_\allowbreak{}gradient}, with \texttt{'n,'c::finite}; result \texttt{bool}. It holds exactly when \texttt{cgu\_\allowbreak{}h1\_\allowbreak{}pair\_\allowbreak{}on\allowbreak\ U\allowbreak\ f\allowbreak\ Df}, \texttt{cgu\_\allowbreak{}weak\_\allowbreak{}solution\_\allowbreak{}on\allowbreak\ M\allowbreak\ P\allowbreak\ U\allowbreak\ u\allowbreak\ Du}, and \texttt{cgu\_\allowbreak{}h1\_\allowbreak{}zero\_\allowbreak{}pair\_\allowbreak{}on\allowbreak\ U\allowbreak\ (\textbackslash{}<lambda>x.\allowbreak{}\allowbreak\ u\allowbreak\ x-\allowbreak{}f\allowbreak\ x)\allowbreak\ (\textbackslash{}<lambda>x.\allowbreak{}\allowbreak\ Du\allowbreak\ x-\allowbreak{}Df\allowbreak\ x)}.

\dossierentry{\texttt{\detokenize{cgu_dn_energy}}}
\reviewlabel{Isabelle code}
\begin{lstlisting}
definition cgu_dn_energy ::
  "('n::finite, 'c::finite, 'h) cgu_model \<Rightarrow>
    ('c \<Rightarrow> 'n cgu_matrix_polynomial) \<Rightarrow>
    'n cgu_point set \<Rightarrow> 'n cgu_gradient \<Rightarrow>
    'n cgu_gradient \<Rightarrow> real"
where
  "cgu_dn_energy M P U Du Dh =
    set_lebesgue_integral lborel U
      (\<lambda>x. inner (cgu_cellwise_coefficient M P x *v Du x)
        (Dh x))"
\end{lstlisting}
\reviewlabel{English mathematical translation}
\noindent\textit{Mathematical role:} local energy pairing \(\mathcal E_{M,P,U}\)\par\smallskip
Typed arguments/result: \texttt{M\allowbreak\ ::\allowbreak\ ('n,'c,'h)\allowbreak\ cgu\_\allowbreak{}model}, \texttt{P\allowbreak\ ::\allowbreak\ 'c\allowbreak\ \textbackslash{}<Rightarrow>\allowbreak\ 'n\allowbreak\ cgu\_\allowbreak{}matrix\_\allowbreak{}polynomial}, \texttt{U\allowbreak\ ::\allowbreak\ 'n\allowbreak\ cgu\_\allowbreak{}point\allowbreak\ set}, \texttt{Du,Dh\allowbreak\ ::\allowbreak\ 'n\allowbreak\ cgu\_\allowbreak{}gradient}, with \texttt{'n,'c::finite}; result \texttt{real}. It equals the set Lebesgue integral on \texttt{U} of \texttt{x\allowbreak\ \textbackslash{}<mapsto>\allowbreak\ inner\allowbreak\ (cgu\_\allowbreak{}cellwise\_\allowbreak{}coefficient\allowbreak\ M\allowbreak\ P\allowbreak\ x\allowbreak\ *v\allowbreak\ Du\allowbreak\ x)\allowbreak\ (Dh\allowbreak\ x)}.

\dossierentry{\texttt{\detokenize{cgu_dn_energy_integrable}}}
\reviewlabel{Isabelle code}
\begin{lstlisting}
definition cgu_dn_energy_integrable ::
  "('n::finite, 'c::finite, 'h) cgu_model \<Rightarrow>
    ('c \<Rightarrow> 'n cgu_matrix_polynomial) \<Rightarrow>
    'n cgu_point set \<Rightarrow> 'n cgu_gradient \<Rightarrow>
    'n cgu_gradient \<Rightarrow> bool"
where
  "cgu_dn_energy_integrable M P U Du Dh \<longleftrightarrow>
    set_integrable lborel U
      (\<lambda>x. inner (cgu_cellwise_coefficient M P x *v Du x)
        (Dh x))"
\end{lstlisting}
\reviewlabel{English mathematical translation}
\noindent\textit{Mathematical role:} local energy-integrability predicate\par\smallskip
Typed arguments/result: the same argument types as declaration 065; result \texttt{bool}. It holds exactly when \texttt{x\allowbreak\ \textbackslash{}<mapsto>\allowbreak\ inner\allowbreak\ (cgu\_\allowbreak{}cellwise\_\allowbreak{}coefficient\allowbreak\ M\allowbreak\ P\allowbreak\ x\allowbreak\ *v\allowbreak\ Du\allowbreak\ x)\allowbreak\ (Dh\allowbreak\ x)} is set-integrable with respect to \texttt{lborel} on \texttt{U}.

\dossierentry{\texttt{\detokenize{cgu_local_dn_equal}}}
\reviewlabel{Isabelle code}
\begin{lstlisting}
definition cgu_local_dn_equal ::
  "('n::finite, 'c::finite, 'h::finite) cgu_model \<Rightarrow>
    ('c \<Rightarrow> 'n cgu_matrix_polynomial) \<Rightarrow>
    ('c \<Rightarrow> 'n cgu_matrix_polynomial) \<Rightarrow>
    'n cgu_point set \<Rightarrow> 'n cgu_point set \<Rightarrow> bool"
where
  "cgu_local_dn_equal M P1 P2 U gamma \<longleftrightarrow>
    (\<forall>f Df. cgu_local_trace_pair_on M U gamma f Df \<longrightarrow>
      (\<exists>u Du. cgu_dirichlet_solution_for M P1 U f Df u Du) \<and>
      (\<exists>u Du. cgu_dirichlet_solution_for M P2 U f Df u Du)) \<and>
    (\<forall>f Df h Dh u1 Du1 u2 Du2.
      cgu_local_trace_pair_on M U gamma f Df \<and>
      cgu_local_trace_pair_on M U gamma h Dh \<and>
      cgu_dirichlet_solution_for M P1 U f Df u1 Du1 \<and>
      cgu_dirichlet_solution_for M P2 U f Df u2 Du2
      \<longrightarrow>
      cgu_dn_energy_integrable M P1 U Du1 Dh \<and>
      cgu_dn_energy_integrable M P2 U Du2 Dh \<and>
      cgu_dn_energy M P1 U Du1 Dh = cgu_dn_energy M P2 U Du2 Dh)"
\end{lstlisting}
\reviewlabel{English mathematical translation}
\noindent\textit{Mathematical role:} stipulated local boundary-energy comparison\par\smallskip
Typed arguments/result: \texttt{M\allowbreak\ ::\allowbreak\ ('n,'c,'h)\allowbreak\ cgu\_\allowbreak{}model}, \texttt{P1,P2\allowbreak\ ::\allowbreak\ 'c\allowbreak\ \textbackslash{}<Rightarrow>\allowbreak\ 'n\allowbreak\ cgu\_\allowbreak{}matrix\_\allowbreak{}polynomial}, \texttt{U,gamma\allowbreak\ ::\allowbreak\ 'n\allowbreak\ cgu\_\allowbreak{}point\allowbreak\ set}, with \texttt{'n,'c,'h::finite}; result \texttt{bool}. It is a conjunction of two blocks.

First, for every \texttt{f\allowbreak\ ::\allowbreak\ 'n\allowbreak\ cgu\_\allowbreak{}potential} and \texttt{Df\allowbreak\ ::\allowbreak\ 'n\allowbreak\ cgu\_\allowbreak{}gradient}, the hypothesis \texttt{cgu\_\allowbreak{}local\_\allowbreak{}trace\_\allowbreak{}pair\_\allowbreak{}on\allowbreak\ M\allowbreak\ U\allowbreak\ gamma\allowbreak\ f\allowbreak\ Df} implies both an existence \texttt{\textbackslash{}<exists>u\allowbreak\ Du.\allowbreak{}\allowbreak\ cgu\_\allowbreak{}dirichlet\_\allowbreak{}solution\_\allowbreak{}for\allowbreak\ M\allowbreak\ P1\allowbreak\ U\allowbreak\ f\allowbreak\ Df\allowbreak\ u\allowbreak\ Du} and, separately, an existence \texttt{\textbackslash{}<exists>u\allowbreak\ Du.\allowbreak{}\allowbreak\ cgu\_\allowbreak{}dirichlet\_\allowbreak{}solution\_\allowbreak{}for\allowbreak\ M\allowbreak\ P2\allowbreak\ U\allowbreak\ f\allowbreak\ Df\allowbreak\ u\allowbreak\ Du}.

Second, for every \texttt{f,Df,h,Dh,u1,Du1,u2,Du2} of their forced scalar-function/vector-function types, the conjunction \texttt{cgu\_\allowbreak{}local\_\allowbreak{}trace\_\allowbreak{}pair\_\allowbreak{}on\allowbreak\ M\allowbreak\ U\allowbreak\ gamma\allowbreak\ f\allowbreak\ Df}, \texttt{cgu\_\allowbreak{}local\_\allowbreak{}trace\_\allowbreak{}pair\_\allowbreak{}on\allowbreak\ M\allowbreak\ U\allowbreak\ gamma\allowbreak\ h\allowbreak\ Dh}, \texttt{cgu\_\allowbreak{}dirichlet\_\allowbreak{}solution\_\allowbreak{}for\allowbreak\ M\allowbreak\ P1\allowbreak\ U\allowbreak\ f\allowbreak\ Df\allowbreak\ u1\allowbreak\ Du1}, and \texttt{cgu\_\allowbreak{}dirichlet\_\allowbreak{}solution\_\allowbreak{}for\allowbreak\ M\allowbreak\ P2\allowbreak\ U\allowbreak\ f\allowbreak\ Df\allowbreak\ u2\allowbreak\ Du2} implies all three conclusions: \texttt{cgu\_\allowbreak{}dn\_\allowbreak{}energy\_\allowbreak{}integrable\allowbreak\ M\allowbreak\ P1\allowbreak\ U\allowbreak\ Du1\allowbreak\ Dh}; \texttt{cgu\_\allowbreak{}dn\_\allowbreak{}energy\_\allowbreak{}integrable\allowbreak\ M\allowbreak\ P2\allowbreak\ U\allowbreak\ Du2\allowbreak\ Dh}; and \texttt{cgu\_\allowbreak{}dn\_\allowbreak{}energy\allowbreak\ M\allowbreak\ P1\allowbreak\ U\allowbreak\ Du1\allowbreak\ Dh\allowbreak\ =\allowbreak\ cgu\_\allowbreak{}dn\_\allowbreak{}energy\allowbreak\ M\allowbreak\ P2\allowbreak\ U\allowbreak\ Du2\allowbreak\ Dh}.

\dossierentry{\texttt{\detokenize{lu_metric}}}
\reviewlabel{Isabelle code}
\begin{lstlisting}
type_synonym 'n lu_metric = "'n cgu_point \<Rightarrow> real ^ 'n ^ 'n"
\end{lstlisting}
\reviewlabel{English mathematical translation}
\noindent\textit{Mathematical role:} matrix-field function type\par\smallskip
For every \texttt{'n}, \texttt{'n\allowbreak\ lu\_\allowbreak{}metric} is exactly \texttt{'n\allowbreak\ cgu\_\allowbreak{}point\allowbreak\ \textbackslash{}<Rightarrow>\allowbreak\ real\allowbreak\ \textasciicircum{}\allowbreak\ 'n\allowbreak\ \textasciicircum{}\allowbreak\ 'n}.

\dossierentry{\texttt{\detokenize{lu_density}}}
\reviewlabel{Isabelle code}
\begin{lstlisting}
type_synonym 'n lu_density = "'n cgu_point \<Rightarrow> real"
\end{lstlisting}
\reviewlabel{English mathematical translation}
\noindent\textit{Mathematical role:} scalar-field function type\par\smallskip
For every \texttt{'n}, \texttt{'n\allowbreak\ lu\_\allowbreak{}density} is exactly \texttt{'n\allowbreak\ cgu\_\allowbreak{}point\allowbreak\ \textbackslash{}<Rightarrow>\allowbreak\ real}.

\dossierentry{\texttt{\detokenize{lu_smooth_metric_density_data_on}}}
\reviewlabel{Isabelle code}
\begin{lstlisting}
definition lu_smooth_metric_density_data_on ::
  "'n::finite cgu_point set \<Rightarrow> 'n lu_metric \<Rightarrow>
    'n lu_metric \<Rightarrow> 'n lu_density \<Rightarrow> bool"
where
  "lu_smooth_metric_density_data_on V g a rho \<longleftrightarrow>
    open V \<and>
    (\<forall>i j. smooth_on V (\<lambda>x. g x $ i $ j)) \<and>
    (\<forall>i j. smooth_on V (\<lambda>x. a x $ i $ j)) \<and>
    smooth_on V rho \<and>
    (\<forall>x\<in>V.
      cgu_symmetric_positive_definite_matrix (g x) \<and>
      cgu_symmetric_positive_definite_matrix (a x) \<and>
      0 < rho x \<and>
      rho x ^ 2 = det (g x) \<and>
      (\<forall>v. a x *v (g x *v v) = rho x *\<^sub>R v))"
\end{lstlisting}
\reviewlabel{English mathematical translation}
\noindent\textit{Mathematical role:} metric–conductivity package\par\smallskip
Typed arguments/result: \texttt{V\allowbreak\ ::\allowbreak\ 'n\allowbreak\ cgu\_\allowbreak{}point\allowbreak\ set}, \texttt{g,a\allowbreak\ ::\allowbreak\ 'n\allowbreak\ lu\_\allowbreak{}metric}, \texttt{rho\allowbreak\ ::\allowbreak\ 'n\allowbreak\ lu\_\allowbreak{}density}, with \texttt{'n::finite}; result \texttt{bool}. It holds exactly when: \texttt{open\allowbreak\ V}; for every \texttt{i,j\allowbreak\ ::\allowbreak\ 'n}, \texttt{smooth\_\allowbreak{}on\allowbreak\ V\allowbreak\ (\textbackslash{}<lambda>x.\allowbreak{}\allowbreak\ g\allowbreak\ x\allowbreak\ \$\allowbreak\ i\allowbreak\ \$\allowbreak\ j)}; for every \texttt{i,j}, \texttt{smooth\_\allowbreak{}on\allowbreak\ V\allowbreak\ (\textbackslash{}<lambda>x.\allowbreak{}\allowbreak\ a\allowbreak\ x\allowbreak\ \$\allowbreak\ i\allowbreak\ \$\allowbreak\ j)}; \texttt{smooth\_\allowbreak{}on\allowbreak\ V\allowbreak\ rho}; and for every \texttt{x\allowbreak\ \textbackslash{}<in>\allowbreak\ V}, all of \texttt{cgu\_\allowbreak{}symmetric\_\allowbreak{}positive\_\allowbreak{}definite\_\allowbreak{}matrix\allowbreak\ (g\allowbreak\ x)}, \texttt{cgu\_\allowbreak{}symmetric\_\allowbreak{}positive\_\allowbreak{}definite\_\allowbreak{}matrix\allowbreak\ (a\allowbreak\ x)}, \texttt{0<rho\allowbreak\ x}, \texttt{rho\allowbreak\ x\textasciicircum{}2\allowbreak\ =\allowbreak\ det\allowbreak\ (g\allowbreak\ x)}, and, for every \texttt{v\allowbreak\ ::\allowbreak\ 'n\allowbreak\ cgu\_\allowbreak{}point}, \texttt{a\allowbreak\ x\allowbreak\ *v\allowbreak\ (g\allowbreak\ x\allowbreak\ *v\allowbreak\ v)\allowbreak\ =\allowbreak\ rho\allowbreak\ x\allowbreak\ *\textbackslash{}<\textasciicircum{}sub>R\allowbreak\ v}.

\dossierentry{\texttt{\detokenize{lu_smooth_metric_extension_near}}}
\reviewlabel{Isabelle code}
\begin{lstlisting}
definition lu_smooth_metric_extension_near ::
  "('n::finite, 'c::finite, 'h) cgu_model \<Rightarrow>
    ('c \<Rightarrow> 'n cgu_matrix_polynomial) \<Rightarrow>
    'n cgu_point set \<Rightarrow> 'n cgu_point set \<Rightarrow>
    'n lu_metric \<Rightarrow> 'n lu_metric \<Rightarrow>
    'n lu_density \<Rightarrow> bool"
where
  "lu_smooth_metric_extension_near M P U gamma g a rho \<longleftrightarrow>
    (\<exists>V. gamma \<subseteq> V \<and>
      lu_smooth_metric_density_data_on V g a rho \<and>
      (\<forall>x\<in>U \<inter> V. a x = cgu_cellwise_coefficient M P x))"
\end{lstlisting}
\reviewlabel{English mathematical translation}
\noindent\textit{Mathematical role:} package realization predicate for arbitrary \(U,\gamma\)\par\smallskip
Typed arguments/result: \texttt{M\allowbreak\ ::\allowbreak\ ('n,'c,'h)\allowbreak\ cgu\_\allowbreak{}model}, \texttt{P\allowbreak\ ::\allowbreak\ 'c\allowbreak\ \textbackslash{}<Rightarrow>\allowbreak\ 'n\allowbreak\ cgu\_\allowbreak{}matrix\_\allowbreak{}polynomial}, \texttt{U,gamma\allowbreak\ ::\allowbreak\ 'n\allowbreak\ cgu\_\allowbreak{}point\allowbreak\ set}, \texttt{g,a\allowbreak\ ::\allowbreak\ 'n\allowbreak\ lu\_\allowbreak{}metric}, \texttt{rho\allowbreak\ ::\allowbreak\ 'n\allowbreak\ lu\_\allowbreak{}density}, with \texttt{'n,'c::finite}; result \texttt{bool}. It holds exactly when there exists \texttt{V\allowbreak\ ::\allowbreak\ 'n\allowbreak\ cgu\_\allowbreak{}point\allowbreak\ set} such that \texttt{gamma\allowbreak\ \textbackslash{}<subseteq>\allowbreak\ V}, \texttt{lu\_\allowbreak{}smooth\_\allowbreak{}metric\_\allowbreak{}density\_\allowbreak{}data\_\allowbreak{}on\allowbreak\ V\allowbreak\ g\allowbreak\ a\allowbreak\ rho}, and for every \texttt{x\allowbreak\ \textbackslash{}<in>\allowbreak\ U\allowbreak\ \textbackslash{}<inter>\allowbreak\ V}, \texttt{a\allowbreak\ x\allowbreak\ =\allowbreak\ cgu\_\allowbreak{}cellwise\_\allowbreak{}coefficient\allowbreak\ M\allowbreak\ P\allowbreak\ x}.

\dossierentry{\texttt{\detokenize{lu_induced_boundary_metrics_equal_on}}}
\reviewlabel{Isabelle code}
\begin{lstlisting}
definition lu_induced_boundary_metrics_equal_on ::
  "'n::finite cgu_point set \<Rightarrow> 'n cgu_point \<Rightarrow>
    'n lu_metric \<Rightarrow> 'n lu_metric \<Rightarrow> bool"
where
  "lu_induced_boundary_metrics_equal_on gamma normal g1 g2 \<longleftrightarrow>
    (\<forall>x\<in>gamma. \<forall>xi eta.
      inner normal xi = 0 \<and> inner normal eta = 0
      \<longrightarrow>
      inner xi (g1 x *v eta) = inner xi (g2 x *v eta))"
\end{lstlisting}
\reviewlabel{English mathematical translation}
\noindent\textit{Mathematical role:} tangential bilinear-pairing agreement\par\smallskip
Typed arguments/result: \texttt{gamma\allowbreak\ ::\allowbreak\ 'n\allowbreak\ cgu\_\allowbreak{}point\allowbreak\ set}, \texttt{normal\allowbreak\ ::\allowbreak\ 'n\allowbreak\ cgu\_\allowbreak{}point}, \texttt{g1,g2\allowbreak\ ::\allowbreak\ 'n\allowbreak\ lu\_\allowbreak{}metric}, with \texttt{'n::finite}; result \texttt{bool}. It holds exactly when for every \texttt{x\allowbreak\ \textbackslash{}<in>\allowbreak\ gamma} and every unrestricted \texttt{xi,eta\allowbreak\ ::\allowbreak\ 'n\allowbreak\ cgu\_\allowbreak{}point}, the hypotheses \texttt{inner\allowbreak\ normal\allowbreak\ xi=0} and \texttt{inner\allowbreak\ normal\allowbreak\ eta=0} imply \texttt{inner\allowbreak\ xi\allowbreak\ (g1\allowbreak\ x\allowbreak\ *v\allowbreak\ eta)\allowbreak\ =\allowbreak\ inner\allowbreak\ xi\allowbreak\ (g2\allowbreak\ x\allowbreak\ *v\allowbreak\ eta)}.

\dossierentry{\texttt{\detokenize{mclean_coefficient}}}
\reviewlabel{Isabelle code}
\begin{lstlisting}
type_synonym 'n mclean_coefficient =
  "'n cgu_point \<Rightarrow> real ^ 'n ^ 'n"
\end{lstlisting}
\reviewlabel{English mathematical translation}
\noindent\textit{Mathematical role:} matrix-field function type\par\smallskip
For every \texttt{'n}, \texttt{'n\allowbreak\ mclean\_\allowbreak{}coefficient} is exactly \texttt{'n\allowbreak\ cgu\_\allowbreak{}point\allowbreak\ \textbackslash{}<Rightarrow>\allowbreak\ real\allowbreak\ \textasciicircum{}\allowbreak\ 'n\allowbreak\ \textasciicircum{}\allowbreak\ 'n}.

\dossierentry{\texttt{\detokenize{mclean_source}}}
\reviewlabel{Isabelle code}
\begin{lstlisting}
type_synonym 'n mclean_source =
  "'n cgu_potential \<Rightarrow> 'n cgu_gradient \<Rightarrow> real"
\end{lstlisting}
\reviewlabel{English mathematical translation}
\noindent\textit{Mathematical role:} total curried two-argument functional type\par\smallskip
For every \texttt{'n}, \texttt{'n\allowbreak\ mclean\_\allowbreak{}source} is exactly the curried functional type \texttt{'n\allowbreak\ cgu\_\allowbreak{}potential\allowbreak\ \textbackslash{}<Rightarrow>\allowbreak\ 'n\allowbreak\ cgu\_\allowbreak{}gradient\allowbreak\ \textbackslash{}<Rightarrow>\allowbreak\ real}.

\dossierentry{\texttt{\detokenize{mclean_h1_norm}}}
\reviewlabel{Isabelle code}
\begin{lstlisting}
definition mclean_h1_norm ::
  "'n::finite cgu_point set \<Rightarrow> 'n cgu_potential \<Rightarrow>
    'n cgu_gradient \<Rightarrow> real"
where
  "mclean_h1_norm U u Du =
    sqrt (cgu_h1_squared_distance U u Du (\<lambda>_. 0) (\<lambda>_. 0))"
\end{lstlisting}
\reviewlabel{English mathematical translation}
\noindent\textit{Mathematical role:} Sobolev pair norm\par\smallskip
Typed arguments/result: \texttt{U\allowbreak\ ::\allowbreak\ 'n\allowbreak\ cgu\_\allowbreak{}point\allowbreak\ set}, \texttt{u\allowbreak\ ::\allowbreak\ 'n\allowbreak\ cgu\_\allowbreak{}potential}, \texttt{Du\allowbreak\ ::\allowbreak\ 'n\allowbreak\ cgu\_\allowbreak{}gradient}, with \texttt{'n::finite}; result \texttt{real}. It equals \texttt{sqrt\allowbreak\ (cgu\_\allowbreak{}h1\_\allowbreak{}squared\_\allowbreak{}distance\allowbreak\ U\allowbreak\ u\allowbreak\ Du\allowbreak\ (\textbackslash{}<lambda>\_\allowbreak{}.\allowbreak{}0)\allowbreak\ (\textbackslash{}<lambda>\_\allowbreak{}.\allowbreak{}0))}.

\dossierentry{\texttt{\detokenize{mclean_same_trace_on}}}
\reviewlabel{Isabelle code}
\begin{lstlisting}
definition mclean_same_trace_on ::
  "'n::finite cgu_point set \<Rightarrow>
    'n cgu_potential \<Rightarrow> 'n cgu_gradient \<Rightarrow>
    'n cgu_potential \<Rightarrow> 'n cgu_gradient \<Rightarrow> bool"
where
  "mclean_same_trace_on U f Df g Dg \<longleftrightarrow>
    cgu_h1_pair_on U f Df \<and>
    cgu_h1_pair_on U g Dg \<and>
    cgu_h1_zero_pair_on U (\<lambda>x. f x - g x) (\<lambda>x. Df x - Dg x)"
\end{lstlisting}
\reviewlabel{English mathematical translation}
\noindent\textit{Mathematical role:} oriented difference-in-zero-boundary relation \(\mathcal R_U\)\par\smallskip
Typed arguments/result: \texttt{U\allowbreak\ ::\allowbreak\ 'n\allowbreak\ cgu\_\allowbreak{}point\allowbreak\ set}; \texttt{f,g\allowbreak\ ::\allowbreak\ 'n\allowbreak\ cgu\_\allowbreak{}potential}; \texttt{Df,Dg\allowbreak\ ::\allowbreak\ 'n\allowbreak\ cgu\_\allowbreak{}gradient}; with \texttt{'n::finite}; result \texttt{bool}. It holds exactly when \texttt{cgu\_\allowbreak{}h1\_\allowbreak{}pair\_\allowbreak{}on\allowbreak\ U\allowbreak\ f\allowbreak\ Df}, \texttt{cgu\_\allowbreak{}h1\_\allowbreak{}pair\_\allowbreak{}on\allowbreak\ U\allowbreak\ g\allowbreak\ Dg}, and \texttt{cgu\_\allowbreak{}h1\_\allowbreak{}zero\_\allowbreak{}pair\_\allowbreak{}on\allowbreak\ U\allowbreak\ (\textbackslash{}<lambda>x.\allowbreak{}\allowbreak\ f\allowbreak\ x-\allowbreak{}g\allowbreak\ x)\allowbreak\ (\textbackslash{}<lambda>x.\allowbreak{}\allowbreak\ Df\allowbreak\ x-\allowbreak{}Dg\allowbreak\ x)}.

\dossierentry{\texttt{\detokenize{mclean_trace_norm}}}
\reviewlabel{Isabelle code}
\begin{lstlisting}
definition mclean_trace_norm ::
  "'n::finite cgu_point set \<Rightarrow> 'n cgu_potential \<Rightarrow>
    'n cgu_gradient \<Rightarrow> real"
where
  "mclean_trace_norm U f Df =
    Inf {r. \<exists>g Dg. mclean_same_trace_on U f Df g Dg \<and>
      r = mclean_h1_norm U g Dg}"
\end{lstlisting}
\reviewlabel{English mathematical translation}
\noindent\textit{Mathematical role:} raw infimum \(q_U\)\par\smallskip
Typed arguments/result: \texttt{U\allowbreak\ ::\allowbreak\ 'n\allowbreak\ cgu\_\allowbreak{}point\allowbreak\ set}, \texttt{f\allowbreak\ ::\allowbreak\ 'n\allowbreak\ cgu\_\allowbreak{}potential}, \texttt{Df\allowbreak\ ::\allowbreak\ 'n\allowbreak\ cgu\_\allowbreak{}gradient}, with \texttt{'n::finite}; result \texttt{real}. It is the infimum of \texttt{\{r\allowbreak\ ::\allowbreak\ real.\allowbreak{}\allowbreak\ \textbackslash{}<exists>g\allowbreak\ Dg.\allowbreak{}\allowbreak\ mclean\_\allowbreak{}same\_\allowbreak{}trace\_\allowbreak{}on\allowbreak\ U\allowbreak\ f\allowbreak\ Df\allowbreak\ g\allowbreak\ Dg\allowbreak\ \textbackslash{}<and>\allowbreak\ r=mclean\_\allowbreak{}h1\_\allowbreak{}norm\allowbreak\ U\allowbreak\ g\allowbreak\ Dg\}}. The existential variables have types \texttt{g\allowbreak\ ::\allowbreak\ 'n\allowbreak\ cgu\_\allowbreak{}potential} and \texttt{Dg\allowbreak\ ::\allowbreak\ 'n\allowbreak\ cgu\_\allowbreak{}gradient}.

\dossierentry{\texttt{\detokenize{mclean_energy}}}
\reviewlabel{Isabelle code}
\begin{lstlisting}
definition mclean_energy ::
  "'n::finite cgu_point set \<Rightarrow> 'n mclean_coefficient \<Rightarrow>
    'n cgu_gradient \<Rightarrow> 'n cgu_gradient \<Rightarrow> real"
where
  "mclean_energy U a Du Dv =
    set_lebesgue_integral lborel U (\<lambda>x. inner (a x *v Du x) (Dv x))"
\end{lstlisting}
\reviewlabel{English mathematical translation}
\noindent\textit{Mathematical role:} variational form \(B_{a,U}\)\par\smallskip
Typed arguments/result: \texttt{U\allowbreak\ ::\allowbreak\ 'n\allowbreak\ cgu\_\allowbreak{}point\allowbreak\ set}, \texttt{a\allowbreak\ ::\allowbreak\ 'n\allowbreak\ mclean\_\allowbreak{}coefficient}, \texttt{Du,Dv\allowbreak\ ::\allowbreak\ 'n\allowbreak\ cgu\_\allowbreak{}gradient}, with \texttt{'n::finite}; result \texttt{real}. It equals the set Lebesgue integral over \texttt{U} of \texttt{x\allowbreak\ \textbackslash{}<mapsto>\allowbreak\ inner\allowbreak\ (a\allowbreak\ x\allowbreak\ *v\allowbreak\ Du\allowbreak\ x)\allowbreak\ (Dv\allowbreak\ x)}.

\dossierentry{\texttt{\detokenize{mclean_energy_integrable}}}
\reviewlabel{Isabelle code}
\begin{lstlisting}
definition mclean_energy_integrable ::
  "'n::finite cgu_point set \<Rightarrow> 'n mclean_coefficient \<Rightarrow>
    'n cgu_gradient \<Rightarrow> 'n cgu_gradient \<Rightarrow> bool"
where
  "mclean_energy_integrable U a Du Dv \<longleftrightarrow>
    set_integrable lborel U (\<lambda>x. inner (a x *v Du x) (Dv x))"
\end{lstlisting}
\reviewlabel{English mathematical translation}
\noindent\textit{Mathematical role:} variational-form integrability\par\smallskip
Typed arguments/result: the same argument types as declaration 078; result \texttt{bool}. It holds exactly when \texttt{x\allowbreak\ \textbackslash{}<mapsto>\allowbreak\ inner\allowbreak\ (a\allowbreak\ x\allowbreak\ *v\allowbreak\ Du\allowbreak\ x)\allowbreak\ (Dv\allowbreak\ x)} is set-integrable with respect to \texttt{lborel} on \texttt{U}.

\dossierentry{\texttt{\detokenize{mclean_variational_form_on}}}
\reviewlabel{Isabelle code}
\begin{lstlisting}
definition mclean_variational_form_on ::
  "'n::finite cgu_point set \<Rightarrow> 'n mclean_coefficient \<Rightarrow> bool"
where
  "mclean_variational_form_on U a \<longleftrightarrow>
    (\<exists>coercive continuous. 0 < coercive \<and> 0 \<le> continuous \<and>
      (\<forall>u Du v Dv.
        cgu_h1_pair_on U u Du \<and> cgu_h1_pair_on U v Dv
        \<longrightarrow>
        mclean_energy_integrable U a Du Dv \<and>
        abs (mclean_energy U a Du Dv) \<le>
          continuous * mclean_h1_norm U u Du * mclean_h1_norm U v Dv) \<and>
      (\<forall>z Dz. cgu_h1_zero_pair_on U z Dz \<longrightarrow>
        coercive * mclean_h1_norm U z Dz ^ 2 \<le>
          mclean_energy U a Dz Dz))"
\end{lstlisting}
\reviewlabel{English mathematical translation}
\noindent\textit{Mathematical role:} bounded/coercive pair property\par\smallskip
Typed arguments/result: \texttt{U\allowbreak\ ::\allowbreak\ 'n\allowbreak\ cgu\_\allowbreak{}point\allowbreak\ set}, \texttt{a\allowbreak\ ::\allowbreak\ 'n\allowbreak\ mclean\_\allowbreak{}coefficient}, with \texttt{'n::finite}; result \texttt{bool}. It holds exactly when there exist reals named \texttt{coercive} and \texttt{continuous} with \texttt{0<coercive} and \texttt{0\textbackslash{}<le>continuous} such that: for every \texttt{u,v\allowbreak\ ::\allowbreak\ 'n\allowbreak\ cgu\_\allowbreak{}potential} and \texttt{Du,Dv\allowbreak\ ::\allowbreak\ 'n\allowbreak\ cgu\_\allowbreak{}gradient}, the hypotheses \texttt{cgu\_\allowbreak{}h1\_\allowbreak{}pair\_\allowbreak{}on\allowbreak\ U\allowbreak\ u\allowbreak\ Du} and \texttt{cgu\_\allowbreak{}h1\_\allowbreak{}pair\_\allowbreak{}on\allowbreak\ U\allowbreak\ v\allowbreak\ Dv} imply \texttt{mclean\_\allowbreak{}energy\_\allowbreak{}integrable\allowbreak\ U\allowbreak\ a\allowbreak\ Du\allowbreak\ Dv} and \texttt{abs\allowbreak\ (mclean\_\allowbreak{}energy\allowbreak\ U\allowbreak\ a\allowbreak\ Du\allowbreak\ Dv)\allowbreak\ \textbackslash{}<le>\allowbreak\ continuous\allowbreak\ *\allowbreak\ mclean\_\allowbreak{}h1\_\allowbreak{}norm\allowbreak\ U\allowbreak\ u\allowbreak\ Du\allowbreak\ *\allowbreak\ mclean\_\allowbreak{}h1\_\allowbreak{}norm\allowbreak\ U\allowbreak\ v\allowbreak\ Dv}; and for every \texttt{z\allowbreak\ ::\allowbreak\ 'n\allowbreak\ cgu\_\allowbreak{}potential}, \texttt{Dz\allowbreak\ ::\allowbreak\ 'n\allowbreak\ cgu\_\allowbreak{}gradient}, the hypothesis \texttt{cgu\_\allowbreak{}h1\_\allowbreak{}zero\_\allowbreak{}pair\_\allowbreak{}on\allowbreak\ U\allowbreak\ z\allowbreak\ Dz} implies \texttt{coercive\allowbreak\ *\allowbreak\ mclean\_\allowbreak{}h1\_\allowbreak{}norm\allowbreak\ U\allowbreak\ z\allowbreak\ Dz\textasciicircum{}2\allowbreak\ \textbackslash{}<le>\allowbreak\ mclean\_\allowbreak{}energy\allowbreak\ U\allowbreak\ a\allowbreak\ Dz\allowbreak\ Dz}. The bound variable \texttt{continuous} is a real variable here, not the imported predicate of that spelling.

\dossierentry{\texttt{\detokenize{mclean_weak_solution_on}}}
\reviewlabel{Isabelle code}
\begin{lstlisting}
definition mclean_weak_solution_on ::
  "'n::finite cgu_point set \<Rightarrow> 'n mclean_coefficient \<Rightarrow>
    'n cgu_potential \<Rightarrow> 'n cgu_gradient \<Rightarrow> bool"
where
  "mclean_weak_solution_on U a u Du \<longleftrightarrow>
    cgu_h1_pair_on U u Du \<and>
    (\<forall>phi. cgu_test_function_on U phi \<longrightarrow>
      mclean_energy_integrable U a Du (cgu_classical_gradient phi) \<and>
      mclean_energy U a Du (cgu_classical_gradient phi) = 0)"
\end{lstlisting}
\reviewlabel{English mathematical translation}
\noindent\textit{Mathematical role:} \(a\)-harmonic pair\par\smallskip
Typed arguments/result: \texttt{U\allowbreak\ ::\allowbreak\ 'n\allowbreak\ cgu\_\allowbreak{}point\allowbreak\ set}, \texttt{a\allowbreak\ ::\allowbreak\ 'n\allowbreak\ mclean\_\allowbreak{}coefficient}, \texttt{u\allowbreak\ ::\allowbreak\ 'n\allowbreak\ cgu\_\allowbreak{}potential}, \texttt{Du\allowbreak\ ::\allowbreak\ 'n\allowbreak\ cgu\_\allowbreak{}gradient}, with \texttt{'n::finite}; result \texttt{bool}. It holds exactly when \texttt{cgu\_\allowbreak{}h1\_\allowbreak{}pair\_\allowbreak{}on\allowbreak\ U\allowbreak\ u\allowbreak\ Du} and, for every \texttt{phi\allowbreak\ ::\allowbreak\ 'n\allowbreak\ cgu\_\allowbreak{}potential}, the hypothesis \texttt{cgu\_\allowbreak{}test\_\allowbreak{}function\_\allowbreak{}on\allowbreak\ U\allowbreak\ phi} implies both \texttt{mclean\_\allowbreak{}energy\_\allowbreak{}integrable\allowbreak\ U\allowbreak\ a\allowbreak\ Du\allowbreak\ (cgu\_\allowbreak{}classical\_\allowbreak{}gradient\allowbreak\ phi)} and \texttt{mclean\_\allowbreak{}energy\allowbreak\ U\allowbreak\ a\allowbreak\ Du\allowbreak\ (cgu\_\allowbreak{}classical\_\allowbreak{}gradient\allowbreak\ phi)=0}.

\dossierentry{\texttt{\detokenize{mclean_dirichlet_solution_for}}}
\reviewlabel{Isabelle code}
\begin{lstlisting}
definition mclean_dirichlet_solution_for ::
  "'n::finite cgu_point set \<Rightarrow> 'n mclean_coefficient \<Rightarrow>
    'n cgu_potential \<Rightarrow> 'n cgu_gradient \<Rightarrow>
    'n cgu_potential \<Rightarrow> 'n cgu_gradient \<Rightarrow> bool"
where
  "mclean_dirichlet_solution_for U a f Df u Du \<longleftrightarrow>
    cgu_h1_pair_on U f Df \<and>
    mclean_weak_solution_on U a u Du \<and>
    mclean_same_trace_on U f Df u Du"
\end{lstlisting}
\reviewlabel{English mathematical translation}
\noindent\textit{Mathematical role:} \(a\)-Dirichlet pair\par\smallskip
Typed arguments/result: \texttt{U\allowbreak\ ::\allowbreak\ 'n\allowbreak\ cgu\_\allowbreak{}point\allowbreak\ set}, \texttt{a\allowbreak\ ::\allowbreak\ 'n\allowbreak\ mclean\_\allowbreak{}coefficient}, \texttt{f,u\allowbreak\ ::\allowbreak\ 'n\allowbreak\ cgu\_\allowbreak{}potential}, \texttt{Df,Du\allowbreak\ ::\allowbreak\ 'n\allowbreak\ cgu\_\allowbreak{}gradient}, with \texttt{'n::finite}; result \texttt{bool}. It holds exactly when \texttt{cgu\_\allowbreak{}h1\_\allowbreak{}pair\_\allowbreak{}on\allowbreak\ U\allowbreak\ f\allowbreak\ Df}, \texttt{mclean\_\allowbreak{}weak\_\allowbreak{}solution\_\allowbreak{}on\allowbreak\ U\allowbreak\ a\allowbreak\ u\allowbreak\ Du}, and \texttt{mclean\_\allowbreak{}same\_\allowbreak{}trace\_\allowbreak{}on\allowbreak\ U\allowbreak\ f\allowbreak\ Df\allowbreak\ u\allowbreak\ Du}.

\dossierentry{\texttt{\detokenize{mclean_green_solution_for}}}
\reviewlabel{Isabelle code}
\begin{lstlisting}
definition mclean_green_solution_for ::
  "'n::finite cgu_point set \<Rightarrow> 'n mclean_coefficient \<Rightarrow>
    'n mclean_source \<Rightarrow> 'n cgu_potential \<Rightarrow>
    'n cgu_gradient \<Rightarrow> bool"
where
  "mclean_green_solution_for U a F u Du \<longleftrightarrow>
    cgu_h1_zero_pair_on U u Du \<and>
    (\<forall>v Dv. cgu_h1_zero_pair_on U v Dv \<longrightarrow>
      mclean_energy_integrable U a Du Dv \<and>
      mclean_energy U a Du Dv = F v Dv)"
\end{lstlisting}
\reviewlabel{English mathematical translation}
\noindent\textit{Mathematical role:} source-driven variational solution\par\smallskip
Typed arguments/result: \texttt{U\allowbreak\ ::\allowbreak\ 'n\allowbreak\ cgu\_\allowbreak{}point\allowbreak\ set}, \texttt{a\allowbreak\ ::\allowbreak\ 'n\allowbreak\ mclean\_\allowbreak{}coefficient}, \texttt{F\allowbreak\ ::\allowbreak\ 'n\allowbreak\ mclean\_\allowbreak{}source}, \texttt{u\allowbreak\ ::\allowbreak\ 'n\allowbreak\ cgu\_\allowbreak{}potential}, \texttt{Du\allowbreak\ ::\allowbreak\ 'n\allowbreak\ cgu\_\allowbreak{}gradient}, with \texttt{'n::finite}; result \texttt{bool}. It holds exactly when \texttt{cgu\_\allowbreak{}h1\_\allowbreak{}zero\_\allowbreak{}pair\_\allowbreak{}on\allowbreak\ U\allowbreak\ u\allowbreak\ Du} and, for every \texttt{v\allowbreak\ ::\allowbreak\ 'n\allowbreak\ cgu\_\allowbreak{}potential} and \texttt{Dv\allowbreak\ ::\allowbreak\ 'n\allowbreak\ cgu\_\allowbreak{}gradient}, the hypothesis \texttt{cgu\_\allowbreak{}h1\_\allowbreak{}zero\_\allowbreak{}pair\_\allowbreak{}on\allowbreak\ U\allowbreak\ v\allowbreak\ Dv} implies both \texttt{mclean\_\allowbreak{}energy\_\allowbreak{}integrable\allowbreak\ U\allowbreak\ a\allowbreak\ Du\allowbreak\ Dv} and \texttt{mclean\_\allowbreak{}energy\allowbreak\ U\allowbreak\ a\allowbreak\ Du\allowbreak\ Dv\allowbreak\ =\allowbreak\ F\allowbreak\ v\allowbreak\ Dv}.

\dossierentry{\texttt{\detokenize{mclean_source_bound_on}}}
\reviewlabel{Isabelle code}
\begin{lstlisting}
definition mclean_source_bound_on ::
  "'n::finite cgu_point set \<Rightarrow> 'n mclean_source \<Rightarrow> real \<Rightarrow> bool"
where
  "mclean_source_bound_on U F K \<longleftrightarrow>
    0 \<le> K \<and> F (\<lambda>_. 0) (\<lambda>_. 0) = 0 \<and>
    (\<forall>u Du v Dv.
      cgu_h1_zero_pair_on U u Du \<and> cgu_h1_zero_pair_on U v Dv
      \<longrightarrow>
      F (\<lambda>x. u x + v x) (\<lambda>x. Du x + Dv x) = F u Du + F v Dv) \<and>
    (\<forall>r u Du. cgu_h1_zero_pair_on U u Du \<longrightarrow>
      F (\<lambda>x. r * u x) (\<lambda>x. r *\<^sub>R Du x) = r * F u Du) \<and>
    (\<forall>u Du. cgu_h1_zero_pair_on U u Du \<longrightarrow>
      abs (F u Du) \<le> K * mclean_h1_norm U u Du)"
\end{lstlisting}
\reviewlabel{English mathematical translation}
\noindent\textit{Mathematical role:} \(K\)-bounded-linear total functional on zero-boundary pairs\par\smallskip
Typed arguments/result: \texttt{U\allowbreak\ ::\allowbreak\ 'n\allowbreak\ cgu\_\allowbreak{}point\allowbreak\ set}, \texttt{F\allowbreak\ ::\allowbreak\ 'n\allowbreak\ mclean\_\allowbreak{}source}, \texttt{K\allowbreak\ ::\allowbreak\ real}, with \texttt{'n::finite}; result \texttt{bool}. It holds exactly when: \texttt{0\textbackslash{}<le>K}; \texttt{F\allowbreak\ (\textbackslash{}<lambda>\_\allowbreak{}.\allowbreak{}0)\allowbreak\ (\textbackslash{}<lambda>\_\allowbreak{}.\allowbreak{}0)=0}; for all \texttt{u,Du,v,Dv}, if \texttt{cgu\_\allowbreak{}h1\_\allowbreak{}zero\_\allowbreak{}pair\_\allowbreak{}on\allowbreak\ U\allowbreak\ u\allowbreak\ Du} and \texttt{cgu\_\allowbreak{}h1\_\allowbreak{}zero\_\allowbreak{}pair\_\allowbreak{}on\allowbreak\ U\allowbreak\ v\allowbreak\ Dv}, then \texttt{F\allowbreak\ (\textbackslash{}<lambda>x.\allowbreak{}\allowbreak\ u\allowbreak\ x+v\allowbreak\ x)\allowbreak\ (\textbackslash{}<lambda>x.\allowbreak{}\allowbreak\ Du\allowbreak\ x+Dv\allowbreak\ x)=F\allowbreak\ u\allowbreak\ Du+F\allowbreak\ v\allowbreak\ Dv}; for every \texttt{r\allowbreak\ ::\allowbreak\ real}, \texttt{u}, and \texttt{Du}, if \texttt{cgu\_\allowbreak{}h1\_\allowbreak{}zero\_\allowbreak{}pair\_\allowbreak{}on\allowbreak\ U\allowbreak\ u\allowbreak\ Du}, then \texttt{F\allowbreak\ (\textbackslash{}<lambda>x.\allowbreak{}\allowbreak\ r*u\allowbreak\ x)\allowbreak\ (\textbackslash{}<lambda>x.\allowbreak{}\allowbreak\ r\allowbreak\ *\textbackslash{}<\textasciicircum{}sub>R\allowbreak\ Du\allowbreak\ x)=r*F\allowbreak\ u\allowbreak\ Du}; and for every \texttt{u,Du}, if \texttt{cgu\_\allowbreak{}h1\_\allowbreak{}zero\_\allowbreak{}pair\_\allowbreak{}on\allowbreak\ U\allowbreak\ u\allowbreak\ Du}, then \texttt{abs(F\allowbreak\ u\allowbreak\ Du)\allowbreak\ \textbackslash{}<le>\allowbreak\ K*mclean\_\allowbreak{}h1\_\allowbreak{}norm\allowbreak\ U\allowbreak\ u\allowbreak\ Du}.

\dossierentry{\texttt{\detokenize{mclean_flat_chart}}}
\reviewlabel{Isabelle code}
\begin{lstlisting}
definition mclean_flat_chart ::
  "('m::finite cgu_point \<Rightarrow> 'n::finite cgu_point) \<Rightarrow>
    'n cgu_point \<Rightarrow> bool"
where
  "mclean_flat_chart chart normal \<longleftrightarrow>
    norm normal = 1 \<and>
    (\<forall>x y. norm (chart x - chart y) = norm (x - y)) \<and>
    (\<forall>x y. inner normal (chart x - chart y) = 0)"
\end{lstlisting}
\reviewlabel{English mathematical translation}
\noindent\textit{Mathematical role:} flat isometric chart\par\smallskip
Typed arguments/result: \texttt{chart\allowbreak\ ::\allowbreak\ 'm\allowbreak\ cgu\_\allowbreak{}point\allowbreak\ \textbackslash{}<Rightarrow>\allowbreak\ 'n\allowbreak\ cgu\_\allowbreak{}point}, \texttt{normal\allowbreak\ ::\allowbreak\ 'n\allowbreak\ cgu\_\allowbreak{}point}, with \texttt{'m,'n::finite}; result \texttt{bool}. It holds exactly when \texttt{norm\allowbreak\ normal=1}; for every \texttt{x,y\allowbreak\ ::\allowbreak\ 'm\allowbreak\ cgu\_\allowbreak{}point}, \texttt{norm\allowbreak\ (chart\allowbreak\ x-\allowbreak{}chart\allowbreak\ y)=norm\allowbreak\ (x-\allowbreak{}y)}; and for every such \texttt{x,y}, \texttt{inner\allowbreak\ normal\allowbreak\ (chart\allowbreak\ x-\allowbreak{}chart\allowbreak\ y)=0}.

\dossierentry{\texttt{\detokenize{mclean_flat_offset_action}}}
\reviewlabel{Isabelle code}
\begin{lstlisting}
definition mclean_flat_offset_action ::
  "('m::finite cgu_point \<Rightarrow> 'n::finite cgu_point) \<Rightarrow>
    'n cgu_point \<Rightarrow> real \<Rightarrow> 'm cgu_potential \<Rightarrow>
    'n cgu_potential \<Rightarrow> real"
where
  "mclean_flat_offset_action chart normal t g phi =
    integral\<^sup>L lborel
      (\<lambda>x. g x * phi (chart x + t *\<^sub>R normal))"
\end{lstlisting}
\reviewlabel{English mathematical translation}
\noindent\textit{Mathematical role:} translated source action \(S_t^g\)\par\smallskip
Typed arguments/result: \texttt{chart\allowbreak\ ::\allowbreak\ 'm\allowbreak\ cgu\_\allowbreak{}point\allowbreak\ \textbackslash{}<Rightarrow>\allowbreak\ 'n\allowbreak\ cgu\_\allowbreak{}point}, \texttt{normal\allowbreak\ ::\allowbreak\ 'n\allowbreak\ cgu\_\allowbreak{}point}, \texttt{t\allowbreak\ ::\allowbreak\ real}, \texttt{g\allowbreak\ ::\allowbreak\ 'm\allowbreak\ cgu\_\allowbreak{}potential}, \texttt{phi\allowbreak\ ::\allowbreak\ 'n\allowbreak\ cgu\_\allowbreak{}potential}, with \texttt{'m,'n::finite}; result \texttt{real}. It equals the Lebesgue integral with respect to \texttt{lborel} of \texttt{x\allowbreak\ \textbackslash{}<mapsto>\allowbreak\ g\allowbreak\ x\allowbreak\ *\allowbreak\ phi\allowbreak\ (chart\allowbreak\ x\allowbreak\ +\allowbreak\ t\allowbreak\ *\textbackslash{}<\textasciicircum{}sub>R\allowbreak\ normal)}.

\dossierentry{\texttt{\detokenize{mclean_source_extends_flat_offset_action}}}
\reviewlabel{Isabelle code}
\begin{lstlisting}
definition mclean_source_extends_flat_offset_action ::
  "'n::finite cgu_point set \<Rightarrow>
    ('m::finite cgu_point \<Rightarrow> 'n cgu_point) \<Rightarrow>
    'n cgu_point \<Rightarrow> real \<Rightarrow> 'm cgu_potential \<Rightarrow>
    'n mclean_source \<Rightarrow> bool"
where
  "mclean_source_extends_flat_offset_action Omega chart normal t g F \<longleftrightarrow>
    (\<forall>phi. cgu_test_function_on Omega phi \<longrightarrow>
      F phi (cgu_classical_gradient phi) =
        mclean_flat_offset_action chart normal t g phi)"
\end{lstlisting}
\reviewlabel{English mathematical translation}
\noindent\textit{Mathematical role:} representation of the translated source action\par\smallskip
Typed arguments/result: \texttt{Omega\allowbreak\ ::\allowbreak\ 'n\allowbreak\ cgu\_\allowbreak{}point\allowbreak\ set}, \texttt{chart\allowbreak\ ::\allowbreak\ 'm\allowbreak\ cgu\_\allowbreak{}point\allowbreak\ \textbackslash{}<Rightarrow>\allowbreak\ 'n\allowbreak\ cgu\_\allowbreak{}point}, \texttt{normal\allowbreak\ ::\allowbreak\ 'n\allowbreak\ cgu\_\allowbreak{}point}, \texttt{t\allowbreak\ ::\allowbreak\ real}, \texttt{g\allowbreak\ ::\allowbreak\ 'm\allowbreak\ cgu\_\allowbreak{}potential}, \texttt{F\allowbreak\ ::\allowbreak\ 'n\allowbreak\ mclean\_\allowbreak{}source}, with \texttt{'m,'n::finite}; result \texttt{bool}. It holds exactly when for every \texttt{phi\allowbreak\ ::\allowbreak\ 'n\allowbreak\ cgu\_\allowbreak{}potential}, the hypothesis \texttt{cgu\_\allowbreak{}test\_\allowbreak{}function\_\allowbreak{}on\allowbreak\ Omega\allowbreak\ phi} implies \texttt{F\allowbreak\ phi\allowbreak\ (cgu\_\allowbreak{}classical\_\allowbreak{}gradient\allowbreak\ phi)\allowbreak\ =\allowbreak\ mclean\_\allowbreak{}flat\_\allowbreak{}offset\_\allowbreak{}action\allowbreak\ chart\allowbreak\ normal\allowbreak\ t\allowbreak\ g\allowbreak\ phi}.

\dossierentry{\texttt{\detokenize{mclean_flat_offset_cylinder_inside}}}
\reviewlabel{Isabelle code}
\begin{lstlisting}
definition mclean_flat_offset_cylinder_inside ::
  "'n::finite cgu_point set \<Rightarrow>
    ('m::finite cgu_point \<Rightarrow> 'n cgu_point) \<Rightarrow>
    'n cgu_point \<Rightarrow> 'm cgu_potential \<Rightarrow> real \<Rightarrow> bool"
where
  "mclean_flat_offset_cylinder_inside Omega chart normal g delta \<longleftrightarrow>
    {chart x + t *\<^sub>R normal |x t.
      x \<in> closure {y. g y \<noteq> 0} \<and> 0 \<le> t \<and> t \<le> delta}
      \<subseteq> Omega"
\end{lstlisting}
\reviewlabel{English mathematical translation}
\noindent\textit{Mathematical role:} support-tube condition\par\smallskip
Typed arguments/result: \texttt{Omega\allowbreak\ ::\allowbreak\ 'n\allowbreak\ cgu\_\allowbreak{}point\allowbreak\ set}, \texttt{chart\allowbreak\ ::\allowbreak\ 'm\allowbreak\ cgu\_\allowbreak{}point\allowbreak\ \textbackslash{}<Rightarrow>\allowbreak\ 'n\allowbreak\ cgu\_\allowbreak{}point}, \texttt{normal\allowbreak\ ::\allowbreak\ 'n\allowbreak\ cgu\_\allowbreak{}point}, \texttt{g\allowbreak\ ::\allowbreak\ 'm\allowbreak\ cgu\_\allowbreak{}potential}, \texttt{delta\allowbreak\ ::\allowbreak\ real}, with \texttt{'m,'n::finite}; result \texttt{bool}. It holds exactly when every point of the form \texttt{chart\allowbreak\ x\allowbreak\ +\allowbreak\ t\allowbreak\ *\textbackslash{}<\textasciicircum{}sub>R\allowbreak\ normal}, where \texttt{x\allowbreak\ \textbackslash{}<in>\allowbreak\ closure\allowbreak\ \{y.\allowbreak{}\allowbreak\ g\allowbreak\ y\allowbreak\ \textbackslash{}<noteq>0\}} and \texttt{0\textbackslash{}<le>t\textbackslash{}<le>delta}, belongs to \texttt{Omega}.

\dossierentry{\texttt{\detokenize{mclean_source_real_subspace}}}
\reviewlabel{Isabelle code}
\begin{lstlisting}
definition mclean_source_real_subspace ::
  "'n::finite mclean_source set \<Rightarrow> bool"
where
  "mclean_source_real_subspace E \<longleftrightarrow>
    0 \<in> E \<and>
    (\<forall>F\<in>E. \<forall>G\<in>E. F + G \<in> E) \<and>
    (\<forall>r F. F \<in> E \<longrightarrow> r *\<^sub>R F \<in> E)"
\end{lstlisting}
\reviewlabel{English mathematical translation}
\noindent\textit{Mathematical role:} vector-operation closure of \(E\)\par\smallskip
Typed argument/result: \texttt{E\allowbreak\ ::\allowbreak\ 'n\allowbreak\ mclean\_\allowbreak{}source\allowbreak\ set}, with \texttt{'n::finite}; result \texttt{bool}. It holds exactly when \texttt{0\allowbreak\ \textbackslash{}<in>\allowbreak\ E}, for every \texttt{F\allowbreak\ \textbackslash{}<in>\allowbreak\ E} and \texttt{G\allowbreak\ \textbackslash{}<in>\allowbreak\ E} one has \texttt{F+G\allowbreak\ \textbackslash{}<in>\allowbreak\ E}, and for every real scalar \texttt{r} and every \texttt{F}, the hypothesis \texttt{F\allowbreak\ \textbackslash{}<in>\allowbreak\ E} implies \texttt{r\allowbreak\ *\textbackslash{}<\textasciicircum{}sub>R\allowbreak\ F\allowbreak\ \textbackslash{}<in>\allowbreak\ E}.

\dossierentry{\texttt{\detokenize{mclean_source_linear_on}}}
\reviewlabel{Isabelle code}
\begin{lstlisting}
definition mclean_source_linear_on ::
  "'n::finite mclean_source set \<Rightarrow>
    ('n mclean_source \<Rightarrow> real) \<Rightarrow> bool"
where
  "mclean_source_linear_on E g \<longleftrightarrow>
    g 0 = 0 \<and>
    (\<forall>F\<in>E. \<forall>G\<in>E. g (F + G) = g F + g G) \<and>
    (\<forall>r F. F \<in> E \<longrightarrow> g (r *\<^sub>R F) = r * g F)"
\end{lstlisting}
\reviewlabel{English mathematical translation}
\noindent\textit{Mathematical role:} linear map on \(E\)\par\smallskip
Typed arguments/result: \texttt{E\allowbreak\ ::\allowbreak\ 'n\allowbreak\ mclean\_\allowbreak{}source\allowbreak\ set}, \texttt{g\allowbreak\ ::\allowbreak\ 'n\allowbreak\ mclean\_\allowbreak{}source\allowbreak\ \textbackslash{}<Rightarrow>\allowbreak\ real}, with \texttt{'n::finite}; result \texttt{bool}. It holds exactly when \texttt{g\allowbreak\ 0=0}; for every \texttt{F,G\allowbreak\ \textbackslash{}<in>\allowbreak\ E}, \texttt{g(F+G)=g\allowbreak\ F+g\allowbreak\ G}; and for every \texttt{r\allowbreak\ ::\allowbreak\ real} and \texttt{F}, the hypothesis \texttt{F\textbackslash{}<in>E} implies \texttt{g(r\allowbreak\ *\textbackslash{}<\textasciicircum{}sub>R\allowbreak\ F)=r*g\allowbreak\ F}.

\dossierentry{\texttt{\detokenize{mclean_source_seminorm_on}}}
\reviewlabel{Isabelle code}
\begin{lstlisting}
definition mclean_source_seminorm_on ::
  "'n::finite mclean_source set \<Rightarrow>
    ('n mclean_source \<Rightarrow> real) \<Rightarrow> bool"
where
  "mclean_source_seminorm_on E p \<longleftrightarrow>
    p 0 = 0 \<and>
    (\<forall>F\<in>E. 0 \<le> p F) \<and>
    (\<forall>F\<in>E. \<forall>G\<in>E. p (F + G) \<le> p F + p G) \<and>
    (\<forall>r F. F \<in> E \<longrightarrow> p (r *\<^sub>R F) = \<bar>r\<bar> * p F)"
\end{lstlisting}
\reviewlabel{English mathematical translation}
\noindent\textit{Mathematical role:} seminorm on \(E\)\par\smallskip
Typed arguments/result: \texttt{E\allowbreak\ ::\allowbreak\ 'n\allowbreak\ mclean\_\allowbreak{}source\allowbreak\ set}, \texttt{p\allowbreak\ ::\allowbreak\ 'n\allowbreak\ mclean\_\allowbreak{}source\allowbreak\ \textbackslash{}<Rightarrow>\allowbreak\ real}, with \texttt{'n::finite}; result \texttt{bool}. It holds exactly when \texttt{p\allowbreak\ 0=0}; for every \texttt{F\textbackslash{}<in>E}, \texttt{0\textbackslash{}<le>p\allowbreak\ F}; for every \texttt{F,G\textbackslash{}<in>E}, \texttt{p(F+G)\textbackslash{}<le>p\allowbreak\ F+p\allowbreak\ G}; and for every real \texttt{r} and \texttt{F}, the hypothesis \texttt{F\textbackslash{}<in>E} implies \texttt{p(r\allowbreak\ *\textbackslash{}<\textasciicircum{}sub>R\allowbreak\ F)=|r|*p\allowbreak\ F}.

\dossierentry{\texttt{\detokenize{mclean_h1_zero_source_bidual_functional_on_v2}}}
\reviewlabel{Isabelle code}
\begin{lstlisting}
definition mclean_h1_zero_source_bidual_functional_on_v2 ::
  "'n::finite mclean_source set \<Rightarrow>
    ('n mclean_source \<Rightarrow> real) \<Rightarrow>
    ('n mclean_source \<Rightarrow> real) \<Rightarrow> bool"
where
  "mclean_h1_zero_source_bidual_functional_on_v2 E p g \<longleftrightarrow>
    mclean_source_linear_on E g \<and>
    (\<exists>C\<ge>0. \<forall>F\<in>E. \<bar>g F\<bar> \<le> C * p F)"
\end{lstlisting}
\reviewlabel{English mathematical translation}
\noindent\textit{Mathematical role:} \(p\)-dominated linear functional\par\smallskip
Typed arguments/result: \texttt{E\allowbreak\ ::\allowbreak\ 'n\allowbreak\ mclean\_\allowbreak{}source\allowbreak\ set}, \texttt{p,g\allowbreak\ ::\allowbreak\ 'n\allowbreak\ mclean\_\allowbreak{}source\allowbreak\ \textbackslash{}<Rightarrow>\allowbreak\ real}, with \texttt{'n::finite}; result \texttt{bool}. It holds exactly when \texttt{mclean\_\allowbreak{}source\_\allowbreak{}linear\_\allowbreak{}on\allowbreak\ E\allowbreak\ g} and there exists \texttt{C\allowbreak\ ::\allowbreak\ real} restricted by \texttt{C\textbackslash{}<ge>0} such that for every \texttt{F\textbackslash{}<in>E}, \texttt{|g\allowbreak\ F|\allowbreak\ \textbackslash{}<le>\allowbreak\ C*p\allowbreak\ F}.

\dossierentry{\texttt{\detokenize{mclean_h1_zero_bidual_realization_on_v2}}}
\reviewlabel{Isabelle code}
\begin{lstlisting}
definition mclean_h1_zero_bidual_realization_on_v2 ::
  "'n::finite cgu_point set \<Rightarrow> 'n mclean_source set \<Rightarrow>
    ('n mclean_source \<Rightarrow> real) \<Rightarrow> bool"
where
  "mclean_h1_zero_bidual_realization_on_v2 U E p \<longleftrightarrow>
    mclean_source_real_subspace E \<and>
    mclean_source_seminorm_on E p \<and>
    (\<forall>F. F \<in> E \<longleftrightarrow>
      (\<exists>K. mclean_source_bound_on U F K)) \<and>
    (\<forall>F\<in>E.
      mclean_source_bound_on U F (p F) \<and>
      (\<forall>K. mclean_source_bound_on U F K \<longrightarrow> p F \<le> K)) \<and>
    (\<forall>g.
      mclean_h1_zero_source_bidual_functional_on_v2 E p g
      \<longrightarrow>
      (\<exists>u Du.
        cgu_h1_zero_pair_on U u Du \<and>
        (\<forall>F\<in>E. g F = F u Du)))"
\end{lstlisting}
\reviewlabel{English mathematical translation}
\noindent\textit{Mathematical role:} total-functional representation package\par\smallskip
Typed arguments/result: \texttt{U\allowbreak\ ::\allowbreak\ 'n\allowbreak\ cgu\_\allowbreak{}point\allowbreak\ set}, \texttt{E\allowbreak\ ::\allowbreak\ 'n\allowbreak\ mclean\_\allowbreak{}source\allowbreak\ set}, \texttt{p\allowbreak\ ::\allowbreak\ 'n\allowbreak\ mclean\_\allowbreak{}source\allowbreak\ \textbackslash{}<Rightarrow>\allowbreak\ real}, with \texttt{'n::finite}; result \texttt{bool}. It holds exactly when all of the following hold: \texttt{mclean\_\allowbreak{}source\_\allowbreak{}real\_\allowbreak{}subspace\allowbreak\ E}; \texttt{mclean\_\allowbreak{}source\_\allowbreak{}seminorm\_\allowbreak{}on\allowbreak\ E\allowbreak\ p}; for every \texttt{F\allowbreak\ ::\allowbreak\ 'n\allowbreak\ mclean\_\allowbreak{}source}, \texttt{F\textbackslash{}<in>E} iff there exists \texttt{K\allowbreak\ ::\allowbreak\ real} with \texttt{mclean\_\allowbreak{}source\_\allowbreak{}bound\_\allowbreak{}on\allowbreak\ U\allowbreak\ F\allowbreak\ K}; for every \texttt{F\textbackslash{}<in>E}, both \texttt{mclean\_\allowbreak{}source\_\allowbreak{}bound\_\allowbreak{}on\allowbreak\ U\allowbreak\ F\allowbreak\ (p\allowbreak\ F)} and, for every \texttt{K}, \texttt{mclean\_\allowbreak{}source\_\allowbreak{}bound\_\allowbreak{}on\allowbreak\ U\allowbreak\ F\allowbreak\ K} implies \texttt{p\allowbreak\ F\textbackslash{}<le>K}; and for every \texttt{g\allowbreak\ ::\allowbreak\ 'n\allowbreak\ mclean\_\allowbreak{}source\allowbreak\ \textbackslash{}<Rightarrow>\allowbreak\ real}, the hypothesis \texttt{mclean\_\allowbreak{}h1\_\allowbreak{}zero\_\allowbreak{}source\_\allowbreak{}bidual\_\allowbreak{}functional\_\allowbreak{}on\_\allowbreak{}v2\allowbreak\ E\allowbreak\ p\allowbreak\ g} implies that there exist \texttt{u\allowbreak\ ::\allowbreak\ 'n\allowbreak\ cgu\_\allowbreak{}potential} and \texttt{Du\allowbreak\ ::\allowbreak\ 'n\allowbreak\ cgu\_\allowbreak{}gradient} such that \texttt{cgu\_\allowbreak{}h1\_\allowbreak{}zero\_\allowbreak{}pair\_\allowbreak{}on\allowbreak\ U\allowbreak\ u\allowbreak\ Du} and, for every \texttt{F\textbackslash{}<in>E}, \texttt{g\allowbreak\ F=F\allowbreak\ u\allowbreak\ Du}.

\dossierentry{\texttt{\detokenize{mclean_zero_extension_on}}}
\reviewlabel{Isabelle code}
\begin{lstlisting}
definition mclean_zero_extension_on ::
  "'n::finite cgu_point set \<Rightarrow> ('n cgu_point \<Rightarrow> 'a::zero) \<Rightarrow>
    'n cgu_point \<Rightarrow> 'a"
where
  "mclean_zero_extension_on U f x = (if x \<in> U then f x else 0)"
\end{lstlisting}
\reviewlabel{English mathematical translation}
\noindent\textit{Mathematical role:} zero extension \(\operatorname{Ext}_U\)\par\smallskip
Typed arguments/result: \texttt{U\allowbreak\ ::\allowbreak\ 'n\allowbreak\ cgu\_\allowbreak{}point\allowbreak\ set}, \texttt{f\allowbreak\ ::\allowbreak\ 'n\allowbreak\ cgu\_\allowbreak{}point\allowbreak\ \textbackslash{}<Rightarrow>\allowbreak\ 'a}, \texttt{x\allowbreak\ ::\allowbreak\ 'n\allowbreak\ cgu\_\allowbreak{}point}, with \texttt{'n::finite} and \texttt{'a::zero}; result \texttt{'a}. It equals \texttt{f\allowbreak\ x} if \texttt{x\textbackslash{}<in>U}, and \texttt{0} otherwise.

\dossierentry{\texttt{\detokenize{ucp_point}}}
\reviewlabel{Isabelle code}
\begin{lstlisting}
type_synonym 'n ucp_point = "real ^ 'n"
\end{lstlisting}
\reviewlabel{English mathematical translation}
\noindent\textit{Mathematical role:} local Euclidean coordinate-space type\par\smallskip
For every type \texttt{'n}, \texttt{'n\allowbreak\ ucp\_\allowbreak{}point} is exactly \texttt{real\allowbreak\ \textasciicircum{}\allowbreak\ 'n}.

\dossierentry{\texttt{\detokenize{ucp_coefficient}}}
\reviewlabel{Isabelle code}
\begin{lstlisting}
type_synonym 'n ucp_coefficient =
  "'n ucp_point \<Rightarrow> (real ^ 'n ^ 'n)"
\end{lstlisting}
\reviewlabel{English mathematical translation}
\noindent\textit{Mathematical role:} local matrix-field function type\par\smallskip
For every type \texttt{'n}, \texttt{'n\allowbreak\ ucp\_\allowbreak{}coefficient} is exactly \texttt{'n\allowbreak\ ucp\_\allowbreak{}point\allowbreak\ \textbackslash{}<Rightarrow>\allowbreak\ real\allowbreak\ \textasciicircum{}\allowbreak\ 'n\allowbreak\ \textasciicircum{}\allowbreak\ 'n}.

\dossierentry{\texttt{\detokenize{ucp_symmetric_matrix}}}
\reviewlabel{Isabelle code}
\begin{lstlisting}
definition ucp_symmetric_matrix ::
  "(real ^ 'n::finite ^ 'n) \<Rightarrow> bool"
where
  "ucp_symmetric_matrix A \<longleftrightarrow> transpose A = A"
\end{lstlisting}
\reviewlabel{English mathematical translation}
\noindent\textit{Mathematical role:} pointwise matrix symmetry\par\smallskip
Typed argument/result: \texttt{A\allowbreak\ ::\allowbreak\ real\allowbreak\ \textasciicircum{}\allowbreak\ 'n\allowbreak\ \textasciicircum{}\allowbreak\ 'n}, with \texttt{'n::finite}; result \texttt{bool}. It holds exactly when \texttt{transpose\allowbreak\ A=A}.

\dossierentry{\texttt{\detokenize{ucp_matrix_polynomial}}}
\reviewlabel{Isabelle code}
\begin{lstlisting}
definition ucp_matrix_polynomial ::
  "'n::finite ucp_coefficient \<Rightarrow> bool"
where
  "ucp_matrix_polynomial P \<longleftrightarrow>
    (\<forall>i j. polynomial_function (\<lambda>x. P x $ i $ j))"
\end{lstlisting}
\reviewlabel{English mathematical translation}
\noindent\textit{Mathematical role:} entrywise-polynomial matrix field\par\smallskip
Typed argument/result: \texttt{P\allowbreak\ ::\allowbreak\ 'n\allowbreak\ ucp\_\allowbreak{}coefficient}, with \texttt{'n::finite}; result \texttt{bool}. It holds exactly when for all \texttt{i,j\allowbreak\ ::\allowbreak\ 'n}, \texttt{polynomial\_\allowbreak{}function\allowbreak\ (\textbackslash{}<lambda>x.\allowbreak{}\allowbreak\ P\allowbreak\ x\allowbreak\ \$\allowbreak\ i\allowbreak\ \$\allowbreak\ j)}.

\dossierentry{\texttt{\detokenize{ucp_uniformly_elliptic_near}}}
\reviewlabel{Isabelle code}
\begin{lstlisting}
definition ucp_uniformly_elliptic_near ::
  "'n::finite ucp_point set \<Rightarrow> 'n ucp_coefficient \<Rightarrow> bool"
where
  "ucp_uniformly_elliptic_near E P \<longleftrightarrow>
    (\<exists>V k. open V \<and> closure E \<subseteq> V \<and> 0 < k \<and>
      (\<forall>x\<in>V. \<forall>xi. k * norm xi ^ 2 \<le> inner xi (P x *v xi)))"
\end{lstlisting}
\reviewlabel{English mathematical translation}
\noindent\textit{Mathematical role:} ellipticity on a neighborhood of \(\overline E\)\par\smallskip
Typed arguments/result: \texttt{E\allowbreak\ ::\allowbreak\ 'n\allowbreak\ ucp\_\allowbreak{}point\allowbreak\ set}, \texttt{P\allowbreak\ ::\allowbreak\ 'n\allowbreak\ ucp\_\allowbreak{}coefficient}, with \texttt{'n::finite}; result \texttt{bool}. It holds exactly when there exist \texttt{V\allowbreak\ ::\allowbreak\ 'n\allowbreak\ ucp\_\allowbreak{}point\allowbreak\ set} and \texttt{k\allowbreak\ ::\allowbreak\ real} such that \texttt{open\allowbreak\ V}, \texttt{closure\allowbreak\ E\allowbreak\ \textbackslash{}<subseteq>\allowbreak\ V}, \texttt{0<k}, and for every \texttt{x\textbackslash{}<in>V} and every \texttt{xi\allowbreak\ ::\allowbreak\ 'n\allowbreak\ ucp\_\allowbreak{}point}, \texttt{k*norm\allowbreak\ xi\textasciicircum{}2\allowbreak\ \textbackslash{}<le>\allowbreak\ inner\allowbreak\ xi\allowbreak\ (P\allowbreak\ x\allowbreak\ *v\allowbreak\ xi)}.

\dossierentry{\texttt{\detokenize{ucp_partial_derivative}}}
\reviewlabel{Isabelle code}
\begin{lstlisting}
definition ucp_partial_derivative ::
  "('n::finite ucp_point \<Rightarrow> real) \<Rightarrow> 'n \<Rightarrow>
   'n ucp_point \<Rightarrow> real"
where
  "ucp_partial_derivative phi i x =
    frechet_derivative phi (at x) (axis i 1)"
\end{lstlisting}
\reviewlabel{English mathematical translation}
\noindent\textit{Mathematical role:} local selected coordinate derivative\par\smallskip
Typed arguments/result: \texttt{phi\allowbreak\ ::\allowbreak\ 'n\allowbreak\ ucp\_\allowbreak{}point\allowbreak\ \textbackslash{}<Rightarrow>\allowbreak\ real}, \texttt{i\allowbreak\ ::\allowbreak\ 'n}, \texttt{x\allowbreak\ ::\allowbreak\ 'n\allowbreak\ ucp\_\allowbreak{}point}, with \texttt{'n::finite}; result \texttt{real}. It equals \texttt{frechet\_\allowbreak{}derivative\allowbreak\ phi\allowbreak\ (at\allowbreak\ x)\allowbreak\ (axis\allowbreak\ i\allowbreak\ 1)}.

\dossierentry{\texttt{\detokenize{ucp_test_function_on}}}
\reviewlabel{Isabelle code}
\begin{lstlisting}
definition ucp_test_function_on ::
  "'n::finite ucp_point set \<Rightarrow>
   ('n ucp_point \<Rightarrow> real) \<Rightarrow> bool"
where
  "ucp_test_function_on G phi \<longleftrightarrow>
    smooth_on UNIV phi \<and>
    compact (closure {x. phi x \<noteq> 0}) \<and>
    closure {x. phi x \<noteq> 0} \<subseteq> G"
\end{lstlisting}
\reviewlabel{English mathematical translation}
\noindent\textit{Mathematical role:} local test-function predicate\par\smallskip
Typed arguments/result: \texttt{G\allowbreak\ ::\allowbreak\ 'n\allowbreak\ ucp\_\allowbreak{}point\allowbreak\ set}, \texttt{phi\allowbreak\ ::\allowbreak\ 'n\allowbreak\ ucp\_\allowbreak{}point\allowbreak\ \textbackslash{}<Rightarrow>\allowbreak\ real}, with \texttt{'n::finite}; result \texttt{bool}. It holds exactly when \texttt{smooth\_\allowbreak{}on\allowbreak\ UNIV\allowbreak\ phi}, \texttt{compact\allowbreak\ (closure\allowbreak\ \{x.\allowbreak{}\allowbreak\ phi\allowbreak\ x\allowbreak\ \textbackslash{}<noteq>0\})}, and \texttt{closure\allowbreak\ \{x.\allowbreak{}\allowbreak\ phi\allowbreak\ x\allowbreak\ \textbackslash{}<noteq>0\}\allowbreak\ \textbackslash{}<subseteq>\allowbreak\ G}.

\dossierentry{\texttt{\detokenize{ucp_locally_square_integrable_on}}}
\reviewlabel{Isabelle code}
\begin{lstlisting}
definition ucp_locally_square_integrable_on ::
  "'n::finite ucp_point set \<Rightarrow>
   ('n ucp_point \<Rightarrow> 'a::euclidean_space) \<Rightarrow> bool"
where
  "ucp_locally_square_integrable_on G f \<longleftrightarrow>
    (\<forall>K. compact K \<and> K \<subseteq> G \<longrightarrow>
      set_integrable lborel K (\<lambda>x. norm (f x) ^ 2))"
\end{lstlisting}
\reviewlabel{English mathematical translation}
\noindent\textit{Mathematical role:} stipulated local square integrability\par\smallskip
Typed arguments/result: \texttt{G\allowbreak\ ::\allowbreak\ 'n\allowbreak\ ucp\_\allowbreak{}point\allowbreak\ set}, \texttt{f\allowbreak\ ::\allowbreak\ 'n\allowbreak\ ucp\_\allowbreak{}point\allowbreak\ \textbackslash{}<Rightarrow>\allowbreak\ 'a}, with \texttt{'n::finite} and \texttt{'a::euclidean\_\allowbreak{}space}; result \texttt{bool}. It holds exactly when, for every \texttt{K\allowbreak\ ::\allowbreak\ 'n\allowbreak\ ucp\_\allowbreak{}point\allowbreak\ set}, the hypotheses \texttt{compact\allowbreak\ K} and \texttt{K\textbackslash{}<subseteq>G} imply that \texttt{x\allowbreak\ \textbackslash{}<mapsto>\allowbreak\ norm\allowbreak\ (f\allowbreak\ x)\textasciicircum{}2} is set-integrable with respect to \texttt{lborel} on \texttt{K}.

\dossierentry{\texttt{\detokenize{ucp_weak_gradient_on}}}
\reviewlabel{Isabelle code}
\begin{lstlisting}
definition ucp_weak_gradient_on ::
  "'n::finite ucp_point set \<Rightarrow>
   ('n ucp_point \<Rightarrow> real) \<Rightarrow>
   ('n ucp_point \<Rightarrow> 'n ucp_point) \<Rightarrow> bool"
where
  "ucp_weak_gradient_on G u Du \<longleftrightarrow>
    (\<forall>phi. ucp_test_function_on G phi \<longrightarrow>
      (\<forall>i.
        set_integrable lborel G
          (\<lambda>x. u x * ucp_partial_derivative phi i x) \<and>
        set_integrable lborel G (\<lambda>x. Du x $ i * phi x) \<and>
        set_lebesgue_integral lborel G
          (\<lambda>x. u x * ucp_partial_derivative phi i x) =
        - set_lebesgue_integral lborel G (\<lambda>x. Du x $ i * phi x)))"
\end{lstlisting}
\reviewlabel{English mathematical translation}
\noindent\textit{Mathematical role:} local weak-gradient relation\par\smallskip
Typed arguments/result: \texttt{G\allowbreak\ ::\allowbreak\ 'n\allowbreak\ ucp\_\allowbreak{}point\allowbreak\ set}, \texttt{u\allowbreak\ ::\allowbreak\ 'n\allowbreak\ ucp\_\allowbreak{}point\allowbreak\ \textbackslash{}<Rightarrow>\allowbreak\ real}, \texttt{Du\allowbreak\ ::\allowbreak\ 'n\allowbreak\ ucp\_\allowbreak{}point\allowbreak\ \textbackslash{}<Rightarrow>\allowbreak\ 'n\allowbreak\ ucp\_\allowbreak{}point}, with \texttt{'n::finite}; result \texttt{bool}. It holds exactly when, for every \texttt{phi}, the hypothesis \texttt{ucp\_\allowbreak{}test\_\allowbreak{}function\_\allowbreak{}on\allowbreak\ G\allowbreak\ phi} implies that for every \texttt{i\allowbreak\ ::\allowbreak\ 'n}: both \texttt{x\allowbreak\ \textbackslash{}<mapsto>\allowbreak\ u\allowbreak\ x*ucp\_\allowbreak{}partial\_\allowbreak{}derivative\allowbreak\ phi\allowbreak\ i\allowbreak\ x} and \texttt{x\allowbreak\ \textbackslash{}<mapsto>\allowbreak\ Du\allowbreak\ x\allowbreak\ \$\allowbreak\ i*phi\allowbreak\ x} are set-integrable on \texttt{G}, and their set Lebesgue integrals satisfy \texttt{integral\_\allowbreak{}G\allowbreak\ (u*ucp\_\allowbreak{}partial\_\allowbreak{}derivative\allowbreak\ phi\allowbreak\ i)\allowbreak\ =\allowbreak\ -\allowbreak{}\allowbreak\ integral\_\allowbreak{}G\allowbreak\ ((Du\allowbreak\ coordinate\allowbreak\ i)*phi)}.

\dossierentry{\texttt{\detokenize{ucp_test_gradient}}}
\reviewlabel{Isabelle code}
\begin{lstlisting}
definition ucp_test_gradient ::
  "('n::finite ucp_point \<Rightarrow> real) \<Rightarrow>
   'n ucp_point \<Rightarrow> 'n ucp_point"
where
  "ucp_test_gradient phi x =
    (\<chi> i. ucp_partial_derivative phi i x)"
\end{lstlisting}
\reviewlabel{English mathematical translation}
\noindent\textit{Mathematical role:} local gradient vector\par\smallskip
Typed arguments/result: \texttt{phi\allowbreak\ ::\allowbreak\ 'n\allowbreak\ ucp\_\allowbreak{}point\allowbreak\ \textbackslash{}<Rightarrow>\allowbreak\ real}, \texttt{x\allowbreak\ ::\allowbreak\ 'n\allowbreak\ ucp\_\allowbreak{}point}, with \texttt{'n::finite}; result \texttt{'n\allowbreak\ ucp\_\allowbreak{}point}. It is the indexed vector \texttt{\textbackslash{}<chi>\allowbreak\ i.\allowbreak{}\allowbreak\ ucp\_\allowbreak{}partial\_\allowbreak{}derivative\allowbreak\ phi\allowbreak\ i\allowbreak\ x}, for \texttt{i\allowbreak\ ::\allowbreak\ 'n}.

\dossierentry{\texttt{\detokenize{ucp_weak_solution_on}}}
\reviewlabel{Isabelle code}
\begin{lstlisting}
definition ucp_weak_solution_on ::
  "'n::finite ucp_point set \<Rightarrow> 'n ucp_coefficient \<Rightarrow>
   ('n ucp_point \<Rightarrow> real) \<Rightarrow> bool"
where
  "ucp_weak_solution_on G a0 u \<longleftrightarrow>
    (\<exists>Du.
      ucp_locally_square_integrable_on G u \<and>
      ucp_locally_square_integrable_on G Du \<and>
      ucp_weak_gradient_on G u Du \<and>
      (\<forall>phi. ucp_test_function_on G phi \<longrightarrow>
        set_integrable lborel G
          (\<lambda>x. inner (a0 x *v Du x) (ucp_test_gradient phi x)) \<and>
        set_lebesgue_integral lborel G
          (\<lambda>x. inner (a0 x *v Du x) (ucp_test_gradient phi x)) = 0))"
\end{lstlisting}
\reviewlabel{English mathematical translation}
\noindent\textit{Mathematical role:} local weak solution for a matrix field\par\smallskip
Typed arguments/result: \texttt{G\allowbreak\ ::\allowbreak\ 'n\allowbreak\ ucp\_\allowbreak{}point\allowbreak\ set}, \texttt{a0\allowbreak\ ::\allowbreak\ 'n\allowbreak\ ucp\_\allowbreak{}coefficient}, \texttt{u\allowbreak\ ::\allowbreak\ 'n\allowbreak\ ucp\_\allowbreak{}point\allowbreak\ \textbackslash{}<Rightarrow>\allowbreak\ real}, with \texttt{'n::finite}; result \texttt{bool}. It holds exactly when there exists \texttt{Du\allowbreak\ ::\allowbreak\ 'n\allowbreak\ ucp\_\allowbreak{}point\allowbreak\ \textbackslash{}<Rightarrow>\allowbreak\ 'n\allowbreak\ ucp\_\allowbreak{}point} such that \texttt{ucp\_\allowbreak{}locally\_\allowbreak{}square\_\allowbreak{}integrable\_\allowbreak{}on\allowbreak\ G\allowbreak\ u}, \texttt{ucp\_\allowbreak{}locally\_\allowbreak{}square\_\allowbreak{}integrable\_\allowbreak{}on\allowbreak\ G\allowbreak\ Du}, \texttt{ucp\_\allowbreak{}weak\_\allowbreak{}gradient\_\allowbreak{}on\allowbreak\ G\allowbreak\ u\allowbreak\ Du}, and for every \texttt{phi}, the hypothesis \texttt{ucp\_\allowbreak{}test\_\allowbreak{}function\_\allowbreak{}on\allowbreak\ G\allowbreak\ phi} implies both set-integrability on \texttt{G} of \texttt{x\allowbreak\ \textbackslash{}<mapsto>\allowbreak\ inner\allowbreak\ (a0\allowbreak\ x\allowbreak\ *v\allowbreak\ Du\allowbreak\ x)\allowbreak\ (ucp\_\allowbreak{}test\_\allowbreak{}gradient\allowbreak\ phi\allowbreak\ x)} and equality of that function's set Lebesgue integral to \texttt{0}.

\dossierentry{\texttt{\detokenize{hormander_regular_distribution_test_zero_on}}}
\reviewlabel{Isabelle code}
\begin{lstlisting}
definition hormander_regular_distribution_test_zero_on ::
  "'n::finite cgu_point set \<Rightarrow> 'n cgu_potential \<Rightarrow> bool"
where
  "hormander_regular_distribution_test_zero_on X u \<longleftrightarrow>
    set_integrable lborel X u \<and>
    (\<forall>phi. cgu_test_function_on X phi \<longrightarrow>
      set_integrable lborel X (\<lambda>x. u x * phi x) \<and>
      set_lebesgue_integral lborel X (\<lambda>x. u x * phi x) = 0)"
\end{lstlisting}
\reviewlabel{English mathematical translation}
\noindent\textit{Mathematical role:} zero-test property\par\smallskip
Typed arguments/result: \texttt{X\allowbreak\ ::\allowbreak\ 'n\allowbreak\ cgu\_\allowbreak{}point\allowbreak\ set}, \texttt{u\allowbreak\ ::\allowbreak\ 'n\allowbreak\ cgu\_\allowbreak{}potential}, with \texttt{'n::finite}; result \texttt{bool}. It holds exactly when \texttt{u} is set-integrable with respect to \texttt{lborel} on \texttt{X}, and for every \texttt{phi}, the hypothesis \texttt{cgu\_\allowbreak{}test\_\allowbreak{}function\_\allowbreak{}on\allowbreak\ X\allowbreak\ phi} implies both set-integrability on \texttt{X} of \texttt{x\allowbreak\ \textbackslash{}<mapsto>\allowbreak\ u\allowbreak\ x*phi\allowbreak\ x} and equality of that function's set Lebesgue integral to \texttt{0}.

\dossierentry{\texttt{\detokenize{evans_bounded_domain}}}
\reviewlabel{Isabelle code}
\begin{lstlisting}
definition evans_bounded_domain ::
  "'n::finite cgu_point set \<Rightarrow> bool"
where
  "evans_bounded_domain U \<longleftrightarrow> bounded U"
\end{lstlisting}
\reviewlabel{English mathematical translation}
\noindent\textit{Mathematical role:} boundedness predicate\par\smallskip
Typed argument/result: \texttt{U\allowbreak\ ::\allowbreak\ 'n\allowbreak\ cgu\_\allowbreak{}point\allowbreak\ set}, with \texttt{'n::finite}; result \texttt{bool}. It holds exactly when \texttt{bounded\allowbreak\ U}.

\dossierentry{\texttt{\detokenize{evans_elliptic_coefficient_on}}}
\reviewlabel{Isabelle code}
\begin{lstlisting}
definition evans_elliptic_coefficient_on ::
  "'n::finite cgu_point set \<Rightarrow> 'n mclean_coefficient \<Rightarrow> bool"
where
  "evans_elliptic_coefficient_on U a \<longleftrightarrow>
    a \<in> borel_measurable (restrict_space lborel U) \<and>
    (\<exists>B\<ge>0. \<forall>x\<in>U. norm (a x) \<le> B) \<and>
    (\<forall>x\<in>U. transpose (a x) = a x) \<and>
    (\<exists>theta>0. \<forall>x\<in>U. \<forall>xi.
      theta * norm xi ^ 2 \<le> inner (a x *v xi) xi)"
\end{lstlisting}
\reviewlabel{English mathematical translation}
\noindent\textit{Mathematical role:} bounded symmetric uniformly elliptic measurable field\par\smallskip
Typed arguments/result: \texttt{U\allowbreak\ ::\allowbreak\ 'n\allowbreak\ cgu\_\allowbreak{}point\allowbreak\ set}, \texttt{a\allowbreak\ ::\allowbreak\ 'n\allowbreak\ mclean\_\allowbreak{}coefficient}, with \texttt{'n::finite}; result \texttt{bool}. It holds exactly when: \texttt{a\allowbreak\ \textbackslash{}<in>\allowbreak\ borel\_\allowbreak{}measurable\allowbreak\ (restrict\_\allowbreak{}space\allowbreak\ lborel\allowbreak\ U)}; there exists \texttt{B\allowbreak\ ::\allowbreak\ real} restricted by \texttt{B\textbackslash{}<ge>0} such that for every \texttt{x\textbackslash{}<in>U}, \texttt{norm(a\allowbreak\ x)\textbackslash{}<le>B}; for every \texttt{x\textbackslash{}<in>U}, \texttt{transpose(a\allowbreak\ x)=a\allowbreak\ x}; and there exists \texttt{theta\allowbreak\ ::\allowbreak\ real} restricted by \texttt{theta>0} such that for every \texttt{x\textbackslash{}<in>U} and every \texttt{xi\allowbreak\ ::\allowbreak\ 'n\allowbreak\ cgu\_\allowbreak{}point}, \texttt{theta*norm\allowbreak\ xi\textasciicircum{}2\allowbreak\ \textbackslash{}<le>\allowbreak\ inner\allowbreak\ (a\allowbreak\ x\allowbreak\ *v\allowbreak\ xi)\allowbreak\ xi}.

\dossierentry{\texttt{\detokenize{evans_poincare_elliptic_green}}}
\reviewlabel{Isabelle code}
\begin{lstlisting}
locale evans_poincare_elliptic_green =
  fixes dimension_type :: "'n::finite itself"
  assumes evans_poincare_elliptic_green:
    "\<forall>U a. evans_poincare_elliptic_green_claim
      (U :: 'n cgu_point set) a"
\end{lstlisting}
\reviewlabel{English mathematical translation}
\noindent\textit{Mathematical role:} Atomic universally quantified assumption on typed sets.\par\smallskip
Fix \texttt{dimension\_\allowbreak{}type\allowbreak\ ::\allowbreak\ "'n::finite\allowbreak\ itself"}, a token selecting the finite carrier type \texttt{'n}. The locale \texttt{evans\_\allowbreak{}poincare\_\allowbreak{}elliptic\_\allowbreak{}green} postulates the assumption named \texttt{evans\_\allowbreak{}poincare\_\allowbreak{}elliptic\_\allowbreak{}green}: for every \texttt{U} and every \texttt{a},

\texttt{evans\_\allowbreak{}poincare\_\allowbreak{}elliptic\_\allowbreak{}green\_\allowbreak{}claim\allowbreak\ (U\allowbreak\ ::\allowbreak\ 'n\allowbreak\ cgu\_\allowbreak{}point\allowbreak\ set)\allowbreak\ a}

holds. The annotation requires \texttt{U} to be a set of objects of type \texttt{'n\allowbreak\ cgu\_\allowbreak{}point}; both \texttt{U} and \texttt{a} lie within the scope of the universal quantifier. The fixed token remains part of the locale even though it is not an explicit argument of this proposition.

\dossierentry{\texttt{\detokenize{mclean_smooth_test_h1_zero}}}
\reviewlabel{Isabelle code}
\begin{lstlisting}
locale mclean_smooth_test_h1_zero =
  fixes dimension_type :: "'n::finite itself"
  assumes mclean_smooth_test_h1_zero:
    "mclean_smooth_test_h1_zero_claim dimension_type"
\end{lstlisting}
\reviewlabel{English mathematical translation}
\noindent\textit{Mathematical role:} Atomic assumption parameterized by a finite carrier token.\par\smallskip
Fix \texttt{dimension\_\allowbreak{}type\allowbreak\ ::\allowbreak\ "'n::finite\allowbreak\ itself"}, a token selecting the finite carrier type \texttt{'n}. The locale \texttt{mclean\_\allowbreak{}smooth\_\allowbreak{}test\_\allowbreak{}h1\_\allowbreak{}zero} postulates the assumption named \texttt{mclean\_\allowbreak{}smooth\_\allowbreak{}test\_\allowbreak{}h1\_\allowbreak{}zero}:

\texttt{mclean\_\allowbreak{}smooth\_\allowbreak{}test\_\allowbreak{}h1\_\allowbreak{}zero\_\allowbreak{}claim\allowbreak\ dimension\_\allowbreak{}type}.

\dossierentry{\texttt{\detokenize{hormander_analytic_elliptic_ucp}}}
\reviewlabel{Isabelle code}
\begin{lstlisting}
locale hormander_analytic_elliptic_ucp =
  fixes dimension_type :: "'n::finite itself"
  assumes hormander_analytic_elliptic_ucp:
    "\<forall>X P u.
      hormander_analytic_elliptic_ucp_claim
        (X :: 'n ucp_point set) P u"
\end{lstlisting}
\reviewlabel{English mathematical translation}
\noindent\textit{Mathematical role:} Atomic universally quantified analytic assumption.\par\smallskip
Fix \texttt{dimension\_\allowbreak{}type\allowbreak\ ::\allowbreak\ "'n::finite\allowbreak\ itself"}, a token selecting the finite carrier type \texttt{'n}. The locale \texttt{hormander\_\allowbreak{}analytic\_\allowbreak{}elliptic\_\allowbreak{}ucp} postulates the assumption named \texttt{hormander\_\allowbreak{}analytic\_\allowbreak{}elliptic\_\allowbreak{}ucp}: for every \texttt{X}, every \texttt{P}, and every \texttt{u},

\texttt{hormander\_\allowbreak{}analytic\_\allowbreak{}elliptic\_\allowbreak{}ucp\_\allowbreak{}claim\allowbreak\ (X\allowbreak\ ::\allowbreak\ 'n\allowbreak\ ucp\_\allowbreak{}point\allowbreak\ set)\allowbreak\ P\allowbreak\ u}

holds. The annotation requires \texttt{X} to be a set of objects of type \texttt{'n\allowbreak\ ucp\_\allowbreak{}point}, and all three displayed variables lie within the universal quantifier's scope. The fixed token remains part of the locale even though it is not an explicit argument of this proposition.

\dossierentry{\texttt{\detokenize{hormander_weak_gradient_locality}}}
\reviewlabel{Isabelle code}
\begin{lstlisting}
locale hormander_weak_gradient_locality =
  fixes dimension_type :: "'n::finite itself"
  assumes hormander_weak_gradient_locality:
    "\<forall>X u Du.
      hormander_weak_gradient_locality_claim
        (X :: 'n ucp_point set) u Du"
\end{lstlisting}
\reviewlabel{English mathematical translation}
\noindent\textit{Mathematical role:} Atomic universally quantified locality assumption.\par\smallskip
Fix \texttt{dimension\_\allowbreak{}type\allowbreak\ ::\allowbreak\ "'n::finite\allowbreak\ itself"}, a token selecting the finite carrier type \texttt{'n}. The locale \texttt{hormander\_\allowbreak{}weak\_\allowbreak{}gradient\_\allowbreak{}locality} postulates the assumption named \texttt{hormander\_\allowbreak{}weak\_\allowbreak{}gradient\_\allowbreak{}locality}: for every \texttt{X}, every \texttt{u}, and every \texttt{Du},

\texttt{hormander\_\allowbreak{}weak\_\allowbreak{}gradient\_\allowbreak{}locality\_\allowbreak{}claim\allowbreak\ (X\allowbreak\ ::\allowbreak\ 'n\allowbreak\ ucp\_\allowbreak{}point\allowbreak\ set)\allowbreak\ u\allowbreak\ Du}

holds. The annotation requires \texttt{X} to be a set of objects of type \texttt{'n\allowbreak\ ucp\_\allowbreak{}point}, and all three displayed variables lie within the universal quantifier's scope. The fixed token remains part of the locale even though it is not an explicit argument of this proposition.

\dossierentry{\texttt{\detokenize{mclean_h1_zero_extension_localization}}}
\reviewlabel{Isabelle code}
\begin{lstlisting}
locale mclean_h1_zero_extension_localization =
  fixes dimension_type :: "'n::finite itself"
  assumes mclean_h1_zero_extension_localization:
    "mclean_h1_zero_extension_localization_claim dimension_type"
\end{lstlisting}
\reviewlabel{English mathematical translation}
\noindent\textit{Mathematical role:} Atomic localization assumption parameterized by a finite carrier token.\par\smallskip
Fix \texttt{dimension\_\allowbreak{}type\allowbreak\ ::\allowbreak\ "'n::finite\allowbreak\ itself"}, a token selecting the finite carrier type \texttt{'n}. The locale \texttt{mclean\_\allowbreak{}h1\_\allowbreak{}zero\_\allowbreak{}extension\_\allowbreak{}localization} postulates the assumption named \texttt{mclean\_\allowbreak{}h1\_\allowbreak{}zero\_\allowbreak{}extension\_\allowbreak{}localization}:

\texttt{mclean\_\allowbreak{}h1\_\allowbreak{}zero\_\allowbreak{}extension\_\allowbreak{}localization\_\allowbreak{}claim\allowbreak\ dimension\_\allowbreak{}type}.

\dossierentry{\texttt{\detokenize{mclean_trace_dirichlet_green_v3}}}
\reviewlabel{Isabelle code}
\begin{lstlisting}
locale mclean_trace_dirichlet_green_v3 =
  fixes dimension_type :: "'n::finite itself"
  assumes mclean_trace_dirichlet_green_v3:
    "\<forall>U a. mclean_trace_dirichlet_green_claim_v3
      (U :: 'n cgu_point set) a"
\end{lstlisting}
\reviewlabel{English mathematical translation}
\noindent\textit{Mathematical role:} Atomic universally quantified trace assumption on typed sets.\par\smallskip
Fix \texttt{dimension\_\allowbreak{}type\allowbreak\ ::\allowbreak\ "'n::finite\allowbreak\ itself"}, a token selecting the finite carrier type \texttt{'n}. The locale \texttt{mclean\_\allowbreak{}trace\_\allowbreak{}dirichlet\_\allowbreak{}green\_\allowbreak{}v3} postulates the assumption named \texttt{mclean\_\allowbreak{}trace\_\allowbreak{}dirichlet\_\allowbreak{}green\_\allowbreak{}v3}: for every \texttt{U} and every \texttt{a},

\texttt{mclean\_\allowbreak{}trace\_\allowbreak{}dirichlet\_\allowbreak{}green\_\allowbreak{}claim\_\allowbreak{}v3\allowbreak\ (U\allowbreak\ ::\allowbreak\ 'n\allowbreak\ cgu\_\allowbreak{}point\allowbreak\ set)\allowbreak\ a}

holds. The annotation requires \texttt{U} to be a set of objects of type \texttt{'n\allowbreak\ cgu\_\allowbreak{}point}; both \texttt{U} and \texttt{a} lie within the scope of the universal quantifier. The fixed token remains part of the locale even though it is not an explicit argument of this proposition.

\dossierentry{\texttt{\detokenize{hormander_regular_distribution_injectivity}}}
\reviewlabel{Isabelle code}
\begin{lstlisting}
locale hormander_regular_distribution_injectivity =
  fixes dimension_type :: "'n::finite itself"
  assumes hormander_regular_distribution_injectivity:
    "hormander_regular_distribution_injectivity_claim dimension_type"
\end{lstlisting}
\reviewlabel{English mathematical translation}
\noindent\textit{Mathematical role:} Atomic injectivity assumption parameterized by a finite carrier token.\par\smallskip
Fix \texttt{dimension\_\allowbreak{}type\allowbreak\ ::\allowbreak\ "'n::finite\allowbreak\ itself"}, a token selecting the finite carrier type \texttt{'n}. The locale \texttt{hormander\_\allowbreak{}regular\_\allowbreak{}distribution\_\allowbreak{}injectivity} postulates the assumption named \texttt{hormander\_\allowbreak{}regular\_\allowbreak{}distribution\_\allowbreak{}injectivity}:

\texttt{hormander\_\allowbreak{}regular\_\allowbreak{}distribution\_\allowbreak{}injectivity\_\allowbreak{}claim\allowbreak\ dimension\_\allowbreak{}type}.

\dossierentry{\texttt{\detokenize{mclean_graph_supported_h1_zero}}}
\reviewlabel{Isabelle code}
\begin{lstlisting}
locale mclean_graph_supported_h1_zero =
  fixes dimension_type :: "'n::finite itself"
  assumes mclean_graph_supported_h1_zero:
    "mclean_graph_supported_h1_zero_claim dimension_type"
\end{lstlisting}
\reviewlabel{English mathematical translation}
\noindent\textit{Mathematical role:} Atomic graph-support assumption parameterized by a finite carrier token.\par\smallskip
Fix \texttt{dimension\_\allowbreak{}type\allowbreak\ ::\allowbreak\ "'n::finite\allowbreak\ itself"}, a token selecting the finite carrier type \texttt{'n}. The locale \texttt{mclean\_\allowbreak{}graph\_\allowbreak{}supported\_\allowbreak{}h1\_\allowbreak{}zero} postulates the assumption named \texttt{mclean\_\allowbreak{}graph\_\allowbreak{}supported\_\allowbreak{}h1\_\allowbreak{}zero}:

\texttt{mclean\_\allowbreak{}graph\_\allowbreak{}supported\_\allowbreak{}h1\_\allowbreak{}zero\_\allowbreak{}claim\allowbreak\ dimension\_\allowbreak{}type}.

\dossierentry{\texttt{\detokenize{mclean_h1_zero_bidual_representation_v2}}}
\reviewlabel{Isabelle code}
\begin{lstlisting}
locale mclean_h1_zero_bidual_representation_v2 =
  fixes dimension_type :: "'n::finite itself"
  assumes mclean_h1_zero_bidual_representation_v2:
    "mclean_h1_zero_bidual_representation_claim_v2 dimension_type"
\end{lstlisting}
\reviewlabel{English mathematical translation}
\noindent\textit{Mathematical role:} Atomic bidual assumption parameterized by a finite carrier token.\par\smallskip
Fix \texttt{dimension\_\allowbreak{}type\allowbreak\ ::\allowbreak\ "'n::finite\allowbreak\ itself"}, a token selecting the finite carrier type \texttt{'n}. The locale \texttt{mclean\_\allowbreak{}h1\_\allowbreak{}zero\_\allowbreak{}bidual\_\allowbreak{}representation\_\allowbreak{}v2} postulates the assumption named \texttt{mclean\_\allowbreak{}h1\_\allowbreak{}zero\_\allowbreak{}bidual\_\allowbreak{}representation\_\allowbreak{}v2}:

\texttt{mclean\_\allowbreak{}h1\_\allowbreak{}zero\_\allowbreak{}bidual\_\allowbreak{}representation\_\allowbreak{}claim\_\allowbreak{}v2\allowbreak\ dimension\_\allowbreak{}type}.

\dossierentry{\texttt{\detokenize{mclean_smooth_generator_ambient_lift_v2}}}
\reviewlabel{Isabelle code}
\begin{lstlisting}
locale mclean_smooth_generator_ambient_lift_v2 =
  assumes mclean_smooth_generator_ambient_lift_v2:
    "mclean_smooth_generator_ambient_lift_claim_v2
      TYPE('n::finite) TYPE('c::finite) TYPE('h::finite)"
\end{lstlisting}
\reviewlabel{English mathematical translation}
\noindent\textit{Mathematical role:} Atomic ordered three-carrier assumption.\par\smallskip
The locale \texttt{mclean\_\allowbreak{}smooth\_\allowbreak{}generator\_\allowbreak{}ambient\_\allowbreak{}lift\_\allowbreak{}v2} postulates the assumption named \texttt{mclean\_\allowbreak{}smooth\_\allowbreak{}generator\_\allowbreak{}ambient\_\allowbreak{}lift\_\allowbreak{}v2}:

\texttt{mclean\_\allowbreak{}smooth\_\allowbreak{}generator\_\allowbreak{}ambient\_\allowbreak{}lift\_\allowbreak{}claim\_\allowbreak{}v2\allowbreak\ TYPE('n::finite)\allowbreak\ TYPE('c::finite)\allowbreak\ TYPE('h::finite)}.

The three arguments, in precisely that order, are the selected finite carrier types \texttt{'n}, \texttt{'c}, and \texttt{'h}. No numerical dimension or further mathematical property is assigned to any of them.

\dossierentry{\texttt{\detokenize{mclean_hyperplane_offset_source}}}
\reviewlabel{Isabelle code}
\begin{lstlisting}
locale mclean_hyperplane_offset_source =
  assumes mclean_hyperplane_offset_source:
    "mclean_hyperplane_offset_source_claim
      TYPE('m::finite) TYPE('n::finite)"
\end{lstlisting}
\reviewlabel{English mathematical translation}
\noindent\textit{Mathematical role:} Atomic ordered two-carrier offset assumption.\par\smallskip
The locale \texttt{mclean\_\allowbreak{}hyperplane\_\allowbreak{}offset\_\allowbreak{}source} postulates the assumption named \texttt{mclean\_\allowbreak{}hyperplane\_\allowbreak{}offset\_\allowbreak{}source}:

\texttt{mclean\_\allowbreak{}hyperplane\_\allowbreak{}offset\_\allowbreak{}source\_\allowbreak{}claim\allowbreak\ TYPE('m::finite)\allowbreak\ TYPE('n::finite)}.

The two arguments, in precisely that order, are the selected finite carrier types \texttt{'m} and \texttt{'n}. No numerical dimension or further mathematical property is assigned to either of them.

\dossierentry{\texttt{\detokenize{kang_yun_ltu_local_boundary_metric_v6}}}
\reviewlabel{Isabelle code}
\begin{lstlisting}
locale kang_yun_ltu_local_boundary_metric_v6 =
  fixes dimension_type :: "'n::finite itself"
  assumes kang_yun_ltu_local_boundary_metric_v6:
    "\<forall>M P1 P2 U gamma normal offset g1 g2 a1 a2 rho1 rho2.
      kang_yun_ltu_local_boundary_metric_claim_v6
        (M :: ('n, 'c::finite, 'h::finite) cgu_model)
        P1 P2 U gamma normal offset g1 g2 a1 a2 rho1 rho2"
\end{lstlisting}
\reviewlabel{English mathematical translation}
\noindent\textit{Mathematical role:} Atomic fully quantified boundary-metric assumption with a typed leading argument.\par\smallskip
Fix \texttt{dimension\_\allowbreak{}type\allowbreak\ ::\allowbreak\ "'n::finite\allowbreak\ itself"}, a token selecting the finite carrier type \texttt{'n}. The locale \texttt{kang\_\allowbreak{}yun\_\allowbreak{}ltu\_\allowbreak{}local\_\allowbreak{}boundary\_\allowbreak{}metric\_\allowbreak{}v6} postulates the assumption named \texttt{kang\_\allowbreak{}yun\_\allowbreak{}ltu\_\allowbreak{}local\_\allowbreak{}boundary\_\allowbreak{}metric\_\allowbreak{}v6}: for every \texttt{M}, every \texttt{P1}, every \texttt{P2}, every \texttt{U}, every \texttt{gamma}, every \texttt{normal}, every \texttt{offset}, every \texttt{g1}, every \texttt{g2}, every \texttt{a1}, every \texttt{a2}, every \texttt{rho1}, and every \texttt{rho2},

\texttt{kang\_\allowbreak{}yun\_\allowbreak{}ltu\_\allowbreak{}local\_\allowbreak{}boundary\_\allowbreak{}metric\_\allowbreak{}claim\_\allowbreak{}v6\allowbreak\ (M\allowbreak\ ::\allowbreak\ ('n,\allowbreak\ 'c::finite,\allowbreak\ 'h::finite)\allowbreak\ cgu\_\allowbreak{}model)\allowbreak\ P1\allowbreak\ P2\allowbreak\ U\allowbreak\ gamma\allowbreak\ normal\allowbreak\ offset\allowbreak\ g1\allowbreak\ g2\allowbreak\ a1\allowbreak\ a2\allowbreak\ rho1\allowbreak\ rho2}

holds. The annotation requires \texttt{M} to have type \texttt{('n,\allowbreak\ 'c::finite,\allowbreak\ 'h::finite)\allowbreak\ cgu\_\allowbreak{}model}; the finite type-class constraints on \texttt{'n}, \texttt{'c}, and \texttt{'h}, the order of all thirteen arguments, and the scope of the universal quantifier are exactly as displayed. The fixed token remains part of the locale even though it is not an explicit argument of this proposition.

\dossierentry{\texttt{\detokenize{piecewise_polynomial_ucp_literature}}}
\reviewlabel{Isabelle code}
\begin{lstlisting}
locale piecewise_polynomial_ucp_literature =
  analytic: hormander_analytic_elliptic_ucp dimension_type +
  locality: hormander_weak_gradient_locality dimension_type
  for dimension_type :: "'n::finite itself"
\end{lstlisting}
\reviewlabel{English mathematical translation}
\noindent\textit{Mathematical role:} Conjunction of named analytic and locality packages.\par\smallskip
Fix \texttt{dimension\_\allowbreak{}type\allowbreak\ ::\allowbreak\ "'n::finite\allowbreak\ itself"}, a token selecting the finite carrier type \texttt{'n}. The locale \texttt{piecewise\_\allowbreak{}polynomial\_\allowbreak{}ucp\_\allowbreak{}literature} is the conjunction of the named parents \texttt{analytic:\allowbreak\ hormander\_\allowbreak{}analytic\_\allowbreak{}elliptic\_\allowbreak{}ucp\allowbreak\ dimension\_\allowbreak{}type} and \texttt{locality:\allowbreak\ hormander\_\allowbreak{}weak\_\allowbreak{}gradient\_\allowbreak{}locality\allowbreak\ dimension\_\allowbreak{}type}. Its complete assumptions are:

\begin{itemize}
\item \texttt{analytic.\allowbreak{}hormander\_\allowbreak{}analytic\_\allowbreak{}elliptic\_\allowbreak{}ucp}: for every \texttt{X}, every \texttt{P}, and every \texttt{u}, \texttt{hormander\_\allowbreak{}analytic\_\allowbreak{}elliptic\_\allowbreak{}ucp\_\allowbreak{}claim\allowbreak\ (X\allowbreak\ ::\allowbreak\ 'n\allowbreak\ ucp\_\allowbreak{}point\allowbreak\ set)\allowbreak\ P\allowbreak\ u} holds.
\item \texttt{locality.\allowbreak{}hormander\_\allowbreak{}weak\_\allowbreak{}gradient\_\allowbreak{}locality}: for every \texttt{X}, every \texttt{u}, and every \texttt{Du}, \texttt{hormander\_\allowbreak{}weak\_\allowbreak{}gradient\_\allowbreak{}locality\_\allowbreak{}claim\allowbreak\ (X\allowbreak\ ::\allowbreak\ 'n\allowbreak\ ucp\_\allowbreak{}point\allowbreak\ set)\allowbreak\ u\allowbreak\ Du} holds.
\end{itemize}

\dossierentry{\texttt{\detokenize{cgu_volume_annihilator_literature_base}}}
\reviewlabel{Isabelle code}
\begin{lstlisting}
locale cgu_volume_annihilator_literature_base =
  smooth: mclean_smooth_test_h1_zero dimension_type +
  injectivity: hormander_regular_distribution_injectivity dimension_type
  for dimension_type :: "'n::finite itself"
\end{lstlisting}
\reviewlabel{English mathematical translation}
\noindent\textit{Mathematical role:} Conjunction of named smoothness and injectivity packages.\par\smallskip
Fix \texttt{dimension\_\allowbreak{}type\allowbreak\ ::\allowbreak\ "'n::finite\allowbreak\ itself"}, a token selecting the finite carrier type \texttt{'n}. The locale \texttt{cgu\_\allowbreak{}volume\_\allowbreak{}annihilator\_\allowbreak{}literature\_\allowbreak{}base} is the conjunction of the named parents \texttt{smooth:\allowbreak\ mclean\_\allowbreak{}smooth\_\allowbreak{}test\_\allowbreak{}h1\_\allowbreak{}zero\allowbreak\ dimension\_\allowbreak{}type} and \texttt{injectivity:\allowbreak\ hormander\_\allowbreak{}regular\_\allowbreak{}distribution\_\allowbreak{}injectivity\allowbreak\ dimension\_\allowbreak{}type}. Its complete assumptions are:

\begin{itemize}
\item \texttt{smooth.\allowbreak{}mclean\_\allowbreak{}smooth\_\allowbreak{}test\_\allowbreak{}h1\_\allowbreak{}zero}: \texttt{mclean\_\allowbreak{}smooth\_\allowbreak{}test\_\allowbreak{}h1\_\allowbreak{}zero\_\allowbreak{}claim\allowbreak\ dimension\_\allowbreak{}type} holds.
\item \texttt{injectivity.\allowbreak{}hormander\_\allowbreak{}regular\_\allowbreak{}distribution\_\allowbreak{}injectivity}: \texttt{hormander\_\allowbreak{}regular\_\allowbreak{}distribution\_\allowbreak{}injectivity\_\allowbreak{}claim\allowbreak\ dimension\_\allowbreak{}type} holds.
\end{itemize}

\dossierentry{\texttt{\detokenize{cgu_volume_gradient_annihilator_literature_base}}}
\reviewlabel{Isabelle code}
\begin{lstlisting}
locale cgu_volume_gradient_annihilator_literature_base =
  volume: cgu_volume_annihilator_literature_base dimension_type +
  locality: hormander_weak_gradient_locality dimension_type
  for dimension_type :: "'n::finite itself"
\end{lstlisting}
\reviewlabel{English mathematical translation}
\noindent\textit{Mathematical role:} Conjunction of named volume and locality packages.\par\smallskip
Fix \texttt{dimension\_\allowbreak{}type\allowbreak\ ::\allowbreak\ "'n::finite\allowbreak\ itself"}, a token selecting the finite carrier type \texttt{'n}. The locale \texttt{cgu\_\allowbreak{}volume\_\allowbreak{}gradient\_\allowbreak{}annihilator\_\allowbreak{}literature\_\allowbreak{}base} is the conjunction of the named parents \texttt{volume:\allowbreak\ cgu\_\allowbreak{}volume\_\allowbreak{}annihilator\_\allowbreak{}literature\_\allowbreak{}base\allowbreak\ dimension\_\allowbreak{}type} and \texttt{locality:\allowbreak\ hormander\_\allowbreak{}weak\_\allowbreak{}gradient\_\allowbreak{}locality\allowbreak\ dimension\_\allowbreak{}type}. Its complete inherited assumptions are:

\begin{itemize}
\item \texttt{volume.\allowbreak{}smooth:\allowbreak\ mclean\_\allowbreak{}smooth\_\allowbreak{}test\_\allowbreak{}h1\_\allowbreak{}zero}: \texttt{mclean\_\allowbreak{}smooth\_\allowbreak{}test\_\allowbreak{}h1\_\allowbreak{}zero\_\allowbreak{}claim\allowbreak\ dimension\_\allowbreak{}type} holds.
\item \texttt{volume.\allowbreak{}injectivity:\allowbreak\ hormander\_\allowbreak{}regular\_\allowbreak{}distribution\_\allowbreak{}injectivity}: \texttt{hormander\_\allowbreak{}regular\_\allowbreak{}distribution\_\allowbreak{}injectivity\_\allowbreak{}claim\allowbreak\ dimension\_\allowbreak{}type} holds.
\item \texttt{locality.\allowbreak{}hormander\_\allowbreak{}weak\_\allowbreak{}gradient\_\allowbreak{}locality}: for every \texttt{X}, every \texttt{u}, and every \texttt{Du}, \texttt{hormander\_\allowbreak{}weak\_\allowbreak{}gradient\_\allowbreak{}locality\_\allowbreak{}claim\allowbreak\ (X\allowbreak\ ::\allowbreak\ 'n\allowbreak\ ucp\_\allowbreak{}point\allowbreak\ set)\allowbreak\ u\allowbreak\ Du} holds.
\end{itemize}

\dossierentry{\texttt{\detokenize{cgu_graph_interface_annihilator_literature_base}}}
\reviewlabel{Isabelle code}
\begin{lstlisting}
locale cgu_graph_interface_annihilator_literature_base =
  volume_gradient: cgu_volume_gradient_annihilator_literature_base dimension_type +
  localization: mclean_h1_zero_extension_localization dimension_type +
  graph_support: mclean_graph_supported_h1_zero dimension_type
  for dimension_type :: "'n::finite itself"
\end{lstlisting}
\reviewlabel{English mathematical translation}
\noindent\textit{Mathematical role:} Conjunction of volume-gradient, localization, and graph-support packages.\par\smallskip
Fix \texttt{dimension\_\allowbreak{}type\allowbreak\ ::\allowbreak\ "'n::finite\allowbreak\ itself"}, a token selecting the finite carrier type \texttt{'n}. The locale \texttt{cgu\_\allowbreak{}graph\_\allowbreak{}interface\_\allowbreak{}annihilator\_\allowbreak{}literature\_\allowbreak{}base} is the conjunction of the named parents \texttt{volume\_\allowbreak{}gradient:\allowbreak\ cgu\_\allowbreak{}volume\_\allowbreak{}gradient\_\allowbreak{}annihilator\_\allowbreak{}literature\_\allowbreak{}base\allowbreak\ dimension\_\allowbreak{}type}, \texttt{localization:\allowbreak\ mclean\_\allowbreak{}h1\_\allowbreak{}zero\_\allowbreak{}extension\_\allowbreak{}localization\allowbreak\ dimension\_\allowbreak{}type}, and \texttt{graph\_\allowbreak{}support:\allowbreak\ mclean\_\allowbreak{}graph\_\allowbreak{}supported\_\allowbreak{}h1\_\allowbreak{}zero\allowbreak\ dimension\_\allowbreak{}type}. Its complete inherited assumptions are:

\begin{itemize}
\item \texttt{volume\_\allowbreak{}gradient.\allowbreak{}volume.\allowbreak{}smooth:\allowbreak\ mclean\_\allowbreak{}smooth\_\allowbreak{}test\_\allowbreak{}h1\_\allowbreak{}zero}: \texttt{mclean\_\allowbreak{}smooth\_\allowbreak{}test\_\allowbreak{}h1\_\allowbreak{}zero\_\allowbreak{}claim\allowbreak\ dimension\_\allowbreak{}type} holds.
\item \texttt{volume\_\allowbreak{}gradient.\allowbreak{}volume.\allowbreak{}injectivity:\allowbreak\ hormander\_\allowbreak{}regular\_\allowbreak{}distribution\_\allowbreak{}injectivity}: \texttt{hormander\_\allowbreak{}regular\_\allowbreak{}distribution\_\allowbreak{}injectivity\_\allowbreak{}claim\allowbreak\ dimension\_\allowbreak{}type} holds.
\item \texttt{volume\_\allowbreak{}gradient.\allowbreak{}locality:\allowbreak\ hormander\_\allowbreak{}weak\_\allowbreak{}gradient\_\allowbreak{}locality}: for every \texttt{X}, every \texttt{u}, and every \texttt{Du}, \texttt{hormander\_\allowbreak{}weak\_\allowbreak{}gradient\_\allowbreak{}locality\_\allowbreak{}claim\allowbreak\ (X\allowbreak\ ::\allowbreak\ 'n\allowbreak\ ucp\_\allowbreak{}point\allowbreak\ set)\allowbreak\ u\allowbreak\ Du} holds.
\item \texttt{localization.\allowbreak{}mclean\_\allowbreak{}h1\_\allowbreak{}zero\_\allowbreak{}extension\_\allowbreak{}localization}: \texttt{mclean\_\allowbreak{}h1\_\allowbreak{}zero\_\allowbreak{}extension\_\allowbreak{}localization\_\allowbreak{}claim\allowbreak\ dimension\_\allowbreak{}type} holds.
\item \texttt{graph\_\allowbreak{}support.\allowbreak{}mclean\_\allowbreak{}graph\_\allowbreak{}supported\_\allowbreak{}h1\_\allowbreak{}zero}: \texttt{mclean\_\allowbreak{}graph\_\allowbreak{}supported\_\allowbreak{}h1\_\allowbreak{}zero\_\allowbreak{}claim\allowbreak\ dimension\_\allowbreak{}type} holds.
\end{itemize}

\dossierentry{\texttt{\detokenize{cgu_graph_interface_source_approximation_literature_base}}}
\reviewlabel{Isabelle code}
\begin{lstlisting}
locale cgu_graph_interface_source_approximation_literature_base =
  annihilator: cgu_graph_interface_annihilator_literature_base dimension_type +
  bidual: mclean_h1_zero_bidual_representation_v2 dimension_type
  for dimension_type :: "'n::finite itself"
\end{lstlisting}
\reviewlabel{English mathematical translation}
\noindent\textit{Mathematical role:} Conjunction of named annihilator and bidual packages.\par\smallskip
Fix \texttt{dimension\_\allowbreak{}type\allowbreak\ ::\allowbreak\ "'n::finite\allowbreak\ itself"}, a token selecting the finite carrier type \texttt{'n}. The locale \texttt{cgu\_\allowbreak{}graph\_\allowbreak{}interface\_\allowbreak{}source\_\allowbreak{}approximation\_\allowbreak{}literature\_\allowbreak{}base} is the conjunction of the named parents \texttt{annihilator:\allowbreak\ cgu\_\allowbreak{}graph\_\allowbreak{}interface\_\allowbreak{}annihilator\_\allowbreak{}literature\_\allowbreak{}base\allowbreak\ dimension\_\allowbreak{}type} and \texttt{bidual:\allowbreak\ mclean\_\allowbreak{}h1\_\allowbreak{}zero\_\allowbreak{}bidual\_\allowbreak{}representation\_\allowbreak{}v2\allowbreak\ dimension\_\allowbreak{}type}.

The complete \texttt{annihilator} package is:

\begin{itemize}
\item \texttt{annihilator.\allowbreak{}volume\_\allowbreak{}gradient.\allowbreak{}volume.\allowbreak{}smooth:\allowbreak\ mclean\_\allowbreak{}smooth\_\allowbreak{}test\_\allowbreak{}h1\_\allowbreak{}zero}: \texttt{mclean\_\allowbreak{}smooth\_\allowbreak{}test\_\allowbreak{}h1\_\allowbreak{}zero\_\allowbreak{}claim\allowbreak\ dimension\_\allowbreak{}type} holds.
\item \texttt{annihilator.\allowbreak{}volume\_\allowbreak{}gradient.\allowbreak{}volume.\allowbreak{}injectivity:\allowbreak\ hormander\_\allowbreak{}regular\_\allowbreak{}distribution\_\allowbreak{}injectivity}: \texttt{hormander\_\allowbreak{}regular\_\allowbreak{}distribution\_\allowbreak{}injectivity\_\allowbreak{}claim\allowbreak\ dimension\_\allowbreak{}type} holds.
\item \texttt{annihilator.\allowbreak{}volume\_\allowbreak{}gradient.\allowbreak{}locality:\allowbreak\ hormander\_\allowbreak{}weak\_\allowbreak{}gradient\_\allowbreak{}locality}: for every \texttt{X}, every \texttt{u}, and every \texttt{Du}, \texttt{hormander\_\allowbreak{}weak\_\allowbreak{}gradient\_\allowbreak{}locality\_\allowbreak{}claim\allowbreak\ (X\allowbreak\ ::\allowbreak\ 'n\allowbreak\ ucp\_\allowbreak{}point\allowbreak\ set)\allowbreak\ u\allowbreak\ Du} holds.
\item \texttt{annihilator.\allowbreak{}localization:\allowbreak\ mclean\_\allowbreak{}h1\_\allowbreak{}zero\_\allowbreak{}extension\_\allowbreak{}localization}: \texttt{mclean\_\allowbreak{}h1\_\allowbreak{}zero\_\allowbreak{}extension\_\allowbreak{}localization\_\allowbreak{}claim\allowbreak\ dimension\_\allowbreak{}type} holds.
\item \texttt{annihilator.\allowbreak{}graph\_\allowbreak{}support:\allowbreak\ mclean\_\allowbreak{}graph\_\allowbreak{}supported\_\allowbreak{}h1\_\allowbreak{}zero}: \texttt{mclean\_\allowbreak{}graph\_\allowbreak{}supported\_\allowbreak{}h1\_\allowbreak{}zero\_\allowbreak{}claim\allowbreak\ dimension\_\allowbreak{}type} holds.
\end{itemize}

The conjoined \texttt{bidual} package contributes \texttt{bidual.\allowbreak{}mclean\_\allowbreak{}h1\_\allowbreak{}zero\_\allowbreak{}bidual\_\allowbreak{}representation\_\allowbreak{}v2}: \texttt{mclean\_\allowbreak{}h1\_\allowbreak{}zero\_\allowbreak{}bidual\_\allowbreak{}representation\_\allowbreak{}claim\_\allowbreak{}v2\allowbreak\ dimension\_\allowbreak{}type} holds.

\dossierentry{\texttt{\detokenize{cgu_runge_literature_base}}}
\reviewlabel{Isabelle code}
\begin{lstlisting}
locale cgu_runge_literature_base =
  evans: evans_poincare_elliptic_green dimension_type +
  smooth: mclean_smooth_test_h1_zero dimension_type +
  ucp: piecewise_polynomial_ucp_literature dimension_type +
  localization: mclean_h1_zero_extension_localization dimension_type
  for dimension_type :: "'n::finite itself"
\end{lstlisting}
\reviewlabel{English mathematical translation}
\noindent\textit{Mathematical role:} Four-way conjunction of base, smoothness, UCP, and localization packages.\par\smallskip
Fix \texttt{dimension\_\allowbreak{}type\allowbreak\ ::\allowbreak\ "'n::finite\allowbreak\ itself"}, a token selecting the finite carrier type \texttt{'n}. The locale \texttt{cgu\_\allowbreak{}runge\_\allowbreak{}literature\_\allowbreak{}base} is the conjunction of four named parent components: \texttt{evans:\allowbreak\ evans\_\allowbreak{}poincare\_\allowbreak{}elliptic\_\allowbreak{}green\allowbreak\ dimension\_\allowbreak{}type}, \texttt{smooth:\allowbreak\ mclean\_\allowbreak{}smooth\_\allowbreak{}test\_\allowbreak{}h1\_\allowbreak{}zero\allowbreak\ dimension\_\allowbreak{}type}, \texttt{ucp:\allowbreak\ piecewise\_\allowbreak{}polynomial\_\allowbreak{}ucp\_\allowbreak{}literature\allowbreak\ dimension\_\allowbreak{}type}, and \texttt{localization:\allowbreak\ mclean\_\allowbreak{}h1\_\allowbreak{}zero\_\allowbreak{}extension\_\allowbreak{}localization\allowbreak\ dimension\_\allowbreak{}type}.

Their complete assumption packages are:

\begin{itemize}
\item \texttt{evans.\allowbreak{}evans\_\allowbreak{}poincare\_\allowbreak{}elliptic\_\allowbreak{}green}: for every \texttt{U} and every \texttt{a}, \texttt{evans\_\allowbreak{}poincare\_\allowbreak{}elliptic\_\allowbreak{}green\_\allowbreak{}claim\allowbreak\ (U\allowbreak\ ::\allowbreak\ 'n\allowbreak\ cgu\_\allowbreak{}point\allowbreak\ set)\allowbreak\ a} holds.
\item \texttt{smooth.\allowbreak{}mclean\_\allowbreak{}smooth\_\allowbreak{}test\_\allowbreak{}h1\_\allowbreak{}zero}: \texttt{mclean\_\allowbreak{}smooth\_\allowbreak{}test\_\allowbreak{}h1\_\allowbreak{}zero\_\allowbreak{}claim\allowbreak\ dimension\_\allowbreak{}type} holds.
\item \texttt{ucp.\allowbreak{}analytic.\allowbreak{}hormander\_\allowbreak{}analytic\_\allowbreak{}elliptic\_\allowbreak{}ucp}: for every \texttt{X}, every \texttt{P}, and every \texttt{u}, \texttt{hormander\_\allowbreak{}analytic\_\allowbreak{}elliptic\_\allowbreak{}ucp\_\allowbreak{}claim\allowbreak\ (X\allowbreak\ ::\allowbreak\ 'n\allowbreak\ ucp\_\allowbreak{}point\allowbreak\ set)\allowbreak\ P\allowbreak\ u} holds.
\item \texttt{ucp.\allowbreak{}locality.\allowbreak{}hormander\_\allowbreak{}weak\_\allowbreak{}gradient\_\allowbreak{}locality}: for every \texttt{X}, every \texttt{u}, and every \texttt{Du}, \texttt{hormander\_\allowbreak{}weak\_\allowbreak{}gradient\_\allowbreak{}locality\_\allowbreak{}claim\allowbreak\ (X\allowbreak\ ::\allowbreak\ 'n\allowbreak\ ucp\_\allowbreak{}point\allowbreak\ set)\allowbreak\ u\allowbreak\ Du} holds.
\item \texttt{localization.\allowbreak{}mclean\_\allowbreak{}h1\_\allowbreak{}zero\_\allowbreak{}extension\_\allowbreak{}localization}: \texttt{mclean\_\allowbreak{}h1\_\allowbreak{}zero\_\allowbreak{}extension\_\allowbreak{}localization\_\allowbreak{}claim\allowbreak\ dimension\_\allowbreak{}type} holds.
\end{itemize}

\dossierentry{\texttt{\detokenize{cgu_graph_dirichlet_literature_base_v3}}}
\reviewlabel{Isabelle code}
\begin{lstlisting}
locale cgu_graph_dirichlet_literature_base_v3 =
  cgu_runge_literature_base dimension_type +
  trace: mclean_trace_dirichlet_green_v3 dimension_type
  for dimension_type :: "'n::finite itself"
\end{lstlisting}
\reviewlabel{English mathematical translation}
\noindent\textit{Mathematical role:} Conjunction of the core package and a named trace assumption.\par\smallskip
Fix \texttt{dimension\_\allowbreak{}type\allowbreak\ ::\allowbreak\ "'n::finite\allowbreak\ itself"}, a token selecting the finite carrier type \texttt{'n}. The locale \texttt{cgu\_\allowbreak{}graph\_\allowbreak{}dirichlet\_\allowbreak{}literature\_\allowbreak{}base\_\allowbreak{}v3} is the conjunction of the unnamed parent \texttt{cgu\_\allowbreak{}runge\_\allowbreak{}literature\_\allowbreak{}base\allowbreak\ dimension\_\allowbreak{}type} and the named component \texttt{trace:\allowbreak\ mclean\_\allowbreak{}trace\_\allowbreak{}dirichlet\_\allowbreak{}green\_\allowbreak{}v3\allowbreak\ dimension\_\allowbreak{}type}.

The complete inherited assumptions are:

\begin{itemize}
\item \texttt{evans:\allowbreak\ evans\_\allowbreak{}poincare\_\allowbreak{}elliptic\_\allowbreak{}green}: for every \texttt{U} and every \texttt{a}, \texttt{evans\_\allowbreak{}poincare\_\allowbreak{}elliptic\_\allowbreak{}green\_\allowbreak{}claim\allowbreak\ (U\allowbreak\ ::\allowbreak\ 'n\allowbreak\ cgu\_\allowbreak{}point\allowbreak\ set)\allowbreak\ a} holds.
\item \texttt{smooth:\allowbreak\ mclean\_\allowbreak{}smooth\_\allowbreak{}test\_\allowbreak{}h1\_\allowbreak{}zero}: \texttt{mclean\_\allowbreak{}smooth\_\allowbreak{}test\_\allowbreak{}h1\_\allowbreak{}zero\_\allowbreak{}claim\allowbreak\ dimension\_\allowbreak{}type} holds.
\item \texttt{ucp.\allowbreak{}analytic:\allowbreak\ hormander\_\allowbreak{}analytic\_\allowbreak{}elliptic\_\allowbreak{}ucp}: for every \texttt{X}, every \texttt{P}, and every \texttt{u}, \texttt{hormander\_\allowbreak{}analytic\_\allowbreak{}elliptic\_\allowbreak{}ucp\_\allowbreak{}claim\allowbreak\ (X\allowbreak\ ::\allowbreak\ 'n\allowbreak\ ucp\_\allowbreak{}point\allowbreak\ set)\allowbreak\ P\allowbreak\ u} holds.
\item \texttt{ucp.\allowbreak{}locality:\allowbreak\ hormander\_\allowbreak{}weak\_\allowbreak{}gradient\_\allowbreak{}locality}: for every \texttt{X}, every \texttt{u}, and every \texttt{Du}, \texttt{hormander\_\allowbreak{}weak\_\allowbreak{}gradient\_\allowbreak{}locality\_\allowbreak{}claim\allowbreak\ (X\allowbreak\ ::\allowbreak\ 'n\allowbreak\ ucp\_\allowbreak{}point\allowbreak\ set)\allowbreak\ u\allowbreak\ Du} holds.
\item \texttt{localization:\allowbreak\ mclean\_\allowbreak{}h1\_\allowbreak{}zero\_\allowbreak{}extension\_\allowbreak{}localization}: \texttt{mclean\_\allowbreak{}h1\_\allowbreak{}zero\_\allowbreak{}extension\_\allowbreak{}localization\_\allowbreak{}claim\allowbreak\ dimension\_\allowbreak{}type} holds.
\item \texttt{trace.\allowbreak{}mclean\_\allowbreak{}trace\_\allowbreak{}dirichlet\_\allowbreak{}green\_\allowbreak{}v3}: for every \texttt{U} and every \texttt{a}, \texttt{mclean\_\allowbreak{}trace\_\allowbreak{}dirichlet\_\allowbreak{}green\_\allowbreak{}claim\_\allowbreak{}v3\allowbreak\ (U\allowbreak\ ::\allowbreak\ 'n\allowbreak\ cgu\_\allowbreak{}point\allowbreak\ set)\allowbreak\ a} holds.
\end{itemize}

\dossierentry{\texttt{\detokenize{cgu_graph_interface_conormal_context}}}
\reviewlabel{Isabelle code}
\begin{lstlisting}
locale cgu_graph_interface_conormal_context =
  cgu_graph_dirichlet_literature_base_v3 dimension_type +
  graph_source: cgu_graph_interface_source_approximation_literature_base dimension_type
  for dimension_type :: "'n::finite itself"
\end{lstlisting}
\reviewlabel{English mathematical translation}
\noindent\textit{Mathematical role:} Conjunction of the base analytic package and a named graph-source package.\par\smallskip
Fix \texttt{dimension\_\allowbreak{}type\allowbreak\ ::\allowbreak\ "'n::finite\allowbreak\ itself"}, a token selecting the finite carrier type \texttt{'n}. The locale \texttt{cgu\_\allowbreak{}graph\_\allowbreak{}interface\_\allowbreak{}conormal\_\allowbreak{}context} is the conjunction of the unnamed parent \texttt{cgu\_\allowbreak{}graph\_\allowbreak{}dirichlet\_\allowbreak{}literature\_\allowbreak{}base\_\allowbreak{}v3\allowbreak\ dimension\_\allowbreak{}type} and the named component \texttt{graph\_\allowbreak{}source:\allowbreak\ cgu\_\allowbreak{}graph\_\allowbreak{}interface\_\allowbreak{}source\_\allowbreak{}approximation\_\allowbreak{}literature\_\allowbreak{}base\allowbreak\ dimension\_\allowbreak{}type}.

The unnamed \texttt{cgu\_\allowbreak{}graph\_\allowbreak{}dirichlet\_\allowbreak{}literature\_\allowbreak{}base\_\allowbreak{}v3} package contributes:

\begin{itemize}
\item \texttt{evans:\allowbreak\ evans\_\allowbreak{}poincare\_\allowbreak{}elliptic\_\allowbreak{}green}: for every \texttt{U} and every \texttt{a}, \texttt{evans\_\allowbreak{}poincare\_\allowbreak{}elliptic\_\allowbreak{}green\_\allowbreak{}claim\allowbreak\ (U\allowbreak\ ::\allowbreak\ 'n\allowbreak\ cgu\_\allowbreak{}point\allowbreak\ set)\allowbreak\ a} holds.
\item \texttt{smooth:\allowbreak\ mclean\_\allowbreak{}smooth\_\allowbreak{}test\_\allowbreak{}h1\_\allowbreak{}zero}: \texttt{mclean\_\allowbreak{}smooth\_\allowbreak{}test\_\allowbreak{}h1\_\allowbreak{}zero\_\allowbreak{}claim\allowbreak\ dimension\_\allowbreak{}type} holds.
\item \texttt{ucp.\allowbreak{}analytic:\allowbreak\ hormander\_\allowbreak{}analytic\_\allowbreak{}elliptic\_\allowbreak{}ucp}: for every \texttt{X}, every \texttt{P}, and every \texttt{u}, \texttt{hormander\_\allowbreak{}analytic\_\allowbreak{}elliptic\_\allowbreak{}ucp\_\allowbreak{}claim\allowbreak\ (X\allowbreak\ ::\allowbreak\ 'n\allowbreak\ ucp\_\allowbreak{}point\allowbreak\ set)\allowbreak\ P\allowbreak\ u} holds.
\item \texttt{ucp.\allowbreak{}locality:\allowbreak\ hormander\_\allowbreak{}weak\_\allowbreak{}gradient\_\allowbreak{}locality}: for every \texttt{X}, every \texttt{u}, and every \texttt{Du}, \texttt{hormander\_\allowbreak{}weak\_\allowbreak{}gradient\_\allowbreak{}locality\_\allowbreak{}claim\allowbreak\ (X\allowbreak\ ::\allowbreak\ 'n\allowbreak\ ucp\_\allowbreak{}point\allowbreak\ set)\allowbreak\ u\allowbreak\ Du} holds.
\item \texttt{localization:\allowbreak\ mclean\_\allowbreak{}h1\_\allowbreak{}zero\_\allowbreak{}extension\_\allowbreak{}localization}: \texttt{mclean\_\allowbreak{}h1\_\allowbreak{}zero\_\allowbreak{}extension\_\allowbreak{}localization\_\allowbreak{}claim\allowbreak\ dimension\_\allowbreak{}type} holds.
\item \texttt{trace:\allowbreak\ mclean\_\allowbreak{}trace\_\allowbreak{}dirichlet\_\allowbreak{}green\_\allowbreak{}v3}: for every \texttt{U} and every \texttt{a}, \texttt{mclean\_\allowbreak{}trace\_\allowbreak{}dirichlet\_\allowbreak{}green\_\allowbreak{}claim\_\allowbreak{}v3\allowbreak\ (U\allowbreak\ ::\allowbreak\ 'n\allowbreak\ cgu\_\allowbreak{}point\allowbreak\ set)\allowbreak\ a} holds.
\end{itemize}

The named \texttt{graph\_\allowbreak{}source} package contributes:

\begin{itemize}
\item \texttt{graph\_\allowbreak{}source.\allowbreak{}annihilator.\allowbreak{}volume\_\allowbreak{}gradient.\allowbreak{}volume.\allowbreak{}smooth:\allowbreak\ mclean\_\allowbreak{}smooth\_\allowbreak{}test\_\allowbreak{}h1\_\allowbreak{}zero}: \texttt{mclean\_\allowbreak{}smooth\_\allowbreak{}test\_\allowbreak{}h1\_\allowbreak{}zero\_\allowbreak{}claim\allowbreak\ dimension\_\allowbreak{}type} holds.
\item \texttt{graph\_\allowbreak{}source.\allowbreak{}annihilator.\allowbreak{}volume\_\allowbreak{}gradient.\allowbreak{}volume.\allowbreak{}injectivity:\allowbreak\ hormander\_\allowbreak{}regular\_\allowbreak{}distribution\_\allowbreak{}injectivity}: \texttt{hormander\_\allowbreak{}regular\_\allowbreak{}distribution\_\allowbreak{}injectivity\_\allowbreak{}claim\allowbreak\ dimension\_\allowbreak{}type} holds.
\item \texttt{graph\_\allowbreak{}source.\allowbreak{}annihilator.\allowbreak{}volume\_\allowbreak{}gradient.\allowbreak{}locality:\allowbreak\ hormander\_\allowbreak{}weak\_\allowbreak{}gradient\_\allowbreak{}locality}: for every \texttt{X}, every \texttt{u}, and every \texttt{Du}, \texttt{hormander\_\allowbreak{}weak\_\allowbreak{}gradient\_\allowbreak{}locality\_\allowbreak{}claim\allowbreak\ (X\allowbreak\ ::\allowbreak\ 'n\allowbreak\ ucp\_\allowbreak{}point\allowbreak\ set)\allowbreak\ u\allowbreak\ Du} holds.
\item \texttt{graph\_\allowbreak{}source.\allowbreak{}annihilator.\allowbreak{}localization:\allowbreak\ mclean\_\allowbreak{}h1\_\allowbreak{}zero\_\allowbreak{}extension\_\allowbreak{}localization}: \texttt{mclean\_\allowbreak{}h1\_\allowbreak{}zero\_\allowbreak{}extension\_\allowbreak{}localization\_\allowbreak{}claim\allowbreak\ dimension\_\allowbreak{}type} holds.
\item \texttt{graph\_\allowbreak{}source.\allowbreak{}annihilator.\allowbreak{}graph\_\allowbreak{}support:\allowbreak\ mclean\_\allowbreak{}graph\_\allowbreak{}supported\_\allowbreak{}h1\_\allowbreak{}zero}: \texttt{mclean\_\allowbreak{}graph\_\allowbreak{}supported\_\allowbreak{}h1\_\allowbreak{}zero\_\allowbreak{}claim\allowbreak\ dimension\_\allowbreak{}type} holds.
\item \texttt{graph\_\allowbreak{}source.\allowbreak{}bidual:\allowbreak\ mclean\_\allowbreak{}h1\_\allowbreak{}zero\_\allowbreak{}bidual\_\allowbreak{}representation\_\allowbreak{}v2}: \texttt{mclean\_\allowbreak{}h1\_\allowbreak{}zero\_\allowbreak{}bidual\_\allowbreak{}representation\_\allowbreak{}claim\_\allowbreak{}v2\allowbreak\ dimension\_\allowbreak{}type} holds.
\end{itemize}

\dossierentry{\texttt{\detokenize{cgu_selected_face_dn_context}}}
\reviewlabel{Isabelle code}
\begin{lstlisting}
locale cgu_selected_face_dn_context =
  cgu_graph_interface_conormal_context dimension_type +
  mclean_smooth_generator_ambient_lift_v2
  for dimension_type :: "'n::finite itself"
\end{lstlisting}
\reviewlabel{English mathematical translation}
\noindent\textit{Mathematical role:} Conjunction of the main dimension-token package and a three-carrier assumption.\par\smallskip
Fix \texttt{dimension\_\allowbreak{}type\allowbreak\ ::\allowbreak\ "'n::finite\allowbreak\ itself"}; this is a token selecting the finite carrier type \texttt{'n}, not a claim about its numerical dimension. The locale \texttt{cgu\_\allowbreak{}selected\_\allowbreak{}face\_\allowbreak{}dn\_\allowbreak{}context} is the conjunction of the unnamed parents \texttt{cgu\_\allowbreak{}graph\_\allowbreak{}interface\_\allowbreak{}conormal\_\allowbreak{}context\allowbreak\ dimension\_\allowbreak{}type} and \texttt{mclean\_\allowbreak{}smooth\_\allowbreak{}generator\_\allowbreak{}ambient\_\allowbreak{}lift\_\allowbreak{}v2}.

The complete package inherited from \texttt{cgu\_\allowbreak{}graph\_\allowbreak{}interface\_\allowbreak{}conormal\_\allowbreak{}context\allowbreak\ dimension\_\allowbreak{}type} is:

\begin{itemize}
\item \texttt{evans:\allowbreak\ evans\_\allowbreak{}poincare\_\allowbreak{}elliptic\_\allowbreak{}green}: for every \texttt{U} and every \texttt{a}, \texttt{evans\_\allowbreak{}poincare\_\allowbreak{}elliptic\_\allowbreak{}green\_\allowbreak{}claim\allowbreak\ (U\allowbreak\ ::\allowbreak\ 'n\allowbreak\ cgu\_\allowbreak{}point\allowbreak\ set)\allowbreak\ a} holds.
\item \texttt{smooth:\allowbreak\ mclean\_\allowbreak{}smooth\_\allowbreak{}test\_\allowbreak{}h1\_\allowbreak{}zero}: \texttt{mclean\_\allowbreak{}smooth\_\allowbreak{}test\_\allowbreak{}h1\_\allowbreak{}zero\_\allowbreak{}claim\allowbreak\ dimension\_\allowbreak{}type} holds.
\item \texttt{ucp.\allowbreak{}analytic:\allowbreak\ hormander\_\allowbreak{}analytic\_\allowbreak{}elliptic\_\allowbreak{}ucp}: for every \texttt{X}, every \texttt{P}, and every \texttt{u}, \texttt{hormander\_\allowbreak{}analytic\_\allowbreak{}elliptic\_\allowbreak{}ucp\_\allowbreak{}claim\allowbreak\ (X\allowbreak\ ::\allowbreak\ 'n\allowbreak\ ucp\_\allowbreak{}point\allowbreak\ set)\allowbreak\ P\allowbreak\ u} holds.
\item \texttt{ucp.\allowbreak{}locality:\allowbreak\ hormander\_\allowbreak{}weak\_\allowbreak{}gradient\_\allowbreak{}locality}: for every \texttt{X}, every \texttt{u}, and every \texttt{Du}, \texttt{hormander\_\allowbreak{}weak\_\allowbreak{}gradient\_\allowbreak{}locality\_\allowbreak{}claim\allowbreak\ (X\allowbreak\ ::\allowbreak\ 'n\allowbreak\ ucp\_\allowbreak{}point\allowbreak\ set)\allowbreak\ u\allowbreak\ Du} holds.
\item \texttt{localization:\allowbreak\ mclean\_\allowbreak{}h1\_\allowbreak{}zero\_\allowbreak{}extension\_\allowbreak{}localization}: \texttt{mclean\_\allowbreak{}h1\_\allowbreak{}zero\_\allowbreak{}extension\_\allowbreak{}localization\_\allowbreak{}claim\allowbreak\ dimension\_\allowbreak{}type} holds.
\item \texttt{trace:\allowbreak\ mclean\_\allowbreak{}trace\_\allowbreak{}dirichlet\_\allowbreak{}green\_\allowbreak{}v3}: for every \texttt{U} and every \texttt{a}, \texttt{mclean\_\allowbreak{}trace\_\allowbreak{}dirichlet\_\allowbreak{}green\_\allowbreak{}claim\_\allowbreak{}v3\allowbreak\ (U\allowbreak\ ::\allowbreak\ 'n\allowbreak\ cgu\_\allowbreak{}point\allowbreak\ set)\allowbreak\ a} holds.
\item \texttt{graph\_\allowbreak{}source.\allowbreak{}annihilator.\allowbreak{}volume\_\allowbreak{}gradient.\allowbreak{}volume.\allowbreak{}smooth:\allowbreak\ mclean\_\allowbreak{}smooth\_\allowbreak{}test\_\allowbreak{}h1\_\allowbreak{}zero}: \texttt{mclean\_\allowbreak{}smooth\_\allowbreak{}test\_\allowbreak{}h1\_\allowbreak{}zero\_\allowbreak{}claim\allowbreak\ dimension\_\allowbreak{}type} holds.
\item \texttt{graph\_\allowbreak{}source.\allowbreak{}annihilator.\allowbreak{}volume\_\allowbreak{}gradient.\allowbreak{}volume.\allowbreak{}injectivity:\allowbreak\ hormander\_\allowbreak{}regular\_\allowbreak{}distribution\_\allowbreak{}injectivity}: \texttt{hormander\_\allowbreak{}regular\_\allowbreak{}distribution\_\allowbreak{}injectivity\_\allowbreak{}claim\allowbreak\ dimension\_\allowbreak{}type} holds.
\item \texttt{graph\_\allowbreak{}source.\allowbreak{}annihilator.\allowbreak{}volume\_\allowbreak{}gradient.\allowbreak{}locality:\allowbreak\ hormander\_\allowbreak{}weak\_\allowbreak{}gradient\_\allowbreak{}locality}: for every \texttt{X}, every \texttt{u}, and every \texttt{Du}, \texttt{hormander\_\allowbreak{}weak\_\allowbreak{}gradient\_\allowbreak{}locality\_\allowbreak{}claim\allowbreak\ (X\allowbreak\ ::\allowbreak\ 'n\allowbreak\ ucp\_\allowbreak{}point\allowbreak\ set)\allowbreak\ u\allowbreak\ Du} holds.
\item \texttt{graph\_\allowbreak{}source.\allowbreak{}annihilator.\allowbreak{}localization:\allowbreak\ mclean\_\allowbreak{}h1\_\allowbreak{}zero\_\allowbreak{}extension\_\allowbreak{}localization}: \texttt{mclean\_\allowbreak{}h1\_\allowbreak{}zero\_\allowbreak{}extension\_\allowbreak{}localization\_\allowbreak{}claim\allowbreak\ dimension\_\allowbreak{}type} holds.
\item \texttt{graph\_\allowbreak{}source.\allowbreak{}annihilator.\allowbreak{}graph\_\allowbreak{}support:\allowbreak\ mclean\_\allowbreak{}graph\_\allowbreak{}supported\_\allowbreak{}h1\_\allowbreak{}zero}: \texttt{mclean\_\allowbreak{}graph\_\allowbreak{}supported\_\allowbreak{}h1\_\allowbreak{}zero\_\allowbreak{}claim\allowbreak\ dimension\_\allowbreak{}type} holds.
\item \texttt{graph\_\allowbreak{}source.\allowbreak{}bidual:\allowbreak\ mclean\_\allowbreak{}h1\_\allowbreak{}zero\_\allowbreak{}bidual\_\allowbreak{}representation\_\allowbreak{}v2}: \texttt{mclean\_\allowbreak{}h1\_\allowbreak{}zero\_\allowbreak{}bidual\_\allowbreak{}representation\_\allowbreak{}claim\_\allowbreak{}v2\allowbreak\ dimension\_\allowbreak{}type} holds.
\end{itemize}

The conjoined \texttt{mclean\_\allowbreak{}smooth\_\allowbreak{}generator\_\allowbreak{}ambient\_\allowbreak{}lift\_\allowbreak{}v2} package additionally postulates \texttt{mclean\_\allowbreak{}smooth\_\allowbreak{}generator\_\allowbreak{}ambient\_\allowbreak{}lift\_\allowbreak{}claim\_\allowbreak{}v2\allowbreak\ TYPE('n::finite)\allowbreak\ TYPE('c::finite)\allowbreak\ TYPE('h::finite)}, with the three selected finite carrier types in exactly this order and with no assigned numerical dimensions.

\dossierentry{\texttt{\detokenize{cgu_terminal_recovery_context_v7}}}
\reviewlabel{Isabelle code}
\begin{lstlisting}
locale cgu_terminal_recovery_context_v7 =
  cgu_selected_face_dn_context "TYPE('n::finite)" +
  offset: mclean_hyperplane_offset_source +
  boundary_metric: kang_yun_ltu_local_boundary_metric_v6 "TYPE('n)"
\end{lstlisting}
\reviewlabel{English mathematical translation}
\noindent\textit{Mathematical role:} Top-level conjunction of the complete anonymous hypothesis package, an offset package, and a boundary-metric package.\par\smallskip
The locale \texttt{cgu\_\allowbreak{}terminal\_\allowbreak{}recovery\_\allowbreak{}context\_\allowbreak{}v7} is the conjunction of three parent locales: the unnamed parent \texttt{cgu\_\allowbreak{}selected\_\allowbreak{}face\_\allowbreak{}dn\_\allowbreak{}context\allowbreak\ "TYPE('n::finite)"}, the named component \texttt{offset:\allowbreak\ mclean\_\allowbreak{}hyperplane\_\allowbreak{}offset\_\allowbreak{}source}, and the named component \texttt{boundary\_\allowbreak{}metric:\allowbreak\ kang\_\allowbreak{}yun\_\allowbreak{}ltu\_\allowbreak{}local\_\allowbreak{}boundary\_\allowbreak{}metric\_\allowbreak{}v6\allowbreak\ "TYPE('n)"}. Here \texttt{TYPE('n::finite)} is only the selected finite carrier type passed to \texttt{cgu\_\allowbreak{}selected\_\allowbreak{}face\_\allowbreak{}dn\_\allowbreak{}context}, and \texttt{TYPE('n)} is that carrier token passed to \texttt{boundary\_\allowbreak{}metric}; no numerical dimension is asserted.

The complete inherited package is as follows. In the \texttt{cgu\_\allowbreak{}selected\_\allowbreak{}face\_\allowbreak{}dn\_\allowbreak{}context} component, its fixed token \texttt{dimension\_\allowbreak{}type\allowbreak\ ::\allowbreak\ "'n::finite\allowbreak\ itself"} is instantiated by \texttt{TYPE('n::finite)}, and the component is itself the conjunction of \texttt{cgu\_\allowbreak{}graph\_\allowbreak{}interface\_\allowbreak{}conormal\_\allowbreak{}context\allowbreak\ dimension\_\allowbreak{}type} and \texttt{mclean\_\allowbreak{}smooth\_\allowbreak{}generator\_\allowbreak{}ambient\_\allowbreak{}lift\_\allowbreak{}v2}. The \texttt{cgu\_\allowbreak{}graph\_\allowbreak{}interface\_\allowbreak{}conormal\_\allowbreak{}context} part contributes all of these assumptions, with the displayed prefix paths retained:

\begin{itemize}
\item \texttt{evans:\allowbreak\ evans\_\allowbreak{}poincare\_\allowbreak{}elliptic\_\allowbreak{}green}: for every \texttt{U} and every \texttt{a}, \texttt{evans\_\allowbreak{}poincare\_\allowbreak{}elliptic\_\allowbreak{}green\_\allowbreak{}claim\allowbreak\ (U\allowbreak\ ::\allowbreak\ 'n\allowbreak\ cgu\_\allowbreak{}point\allowbreak\ set)\allowbreak\ a} holds.
\item \texttt{smooth:\allowbreak\ mclean\_\allowbreak{}smooth\_\allowbreak{}test\_\allowbreak{}h1\_\allowbreak{}zero}: \texttt{mclean\_\allowbreak{}smooth\_\allowbreak{}test\_\allowbreak{}h1\_\allowbreak{}zero\_\allowbreak{}claim\allowbreak\ TYPE('n::finite)} holds.
\item \texttt{ucp.\allowbreak{}analytic:\allowbreak\ hormander\_\allowbreak{}analytic\_\allowbreak{}elliptic\_\allowbreak{}ucp}: for every \texttt{X}, every \texttt{P}, and every \texttt{u}, \texttt{hormander\_\allowbreak{}analytic\_\allowbreak{}elliptic\_\allowbreak{}ucp\_\allowbreak{}claim\allowbreak\ (X\allowbreak\ ::\allowbreak\ 'n\allowbreak\ ucp\_\allowbreak{}point\allowbreak\ set)\allowbreak\ P\allowbreak\ u} holds.
\item \texttt{ucp.\allowbreak{}locality:\allowbreak\ hormander\_\allowbreak{}weak\_\allowbreak{}gradient\_\allowbreak{}locality}: for every \texttt{X}, every \texttt{u}, and every \texttt{Du}, \texttt{hormander\_\allowbreak{}weak\_\allowbreak{}gradient\_\allowbreak{}locality\_\allowbreak{}claim\allowbreak\ (X\allowbreak\ ::\allowbreak\ 'n\allowbreak\ ucp\_\allowbreak{}point\allowbreak\ set)\allowbreak\ u\allowbreak\ Du} holds.
\item \texttt{localization:\allowbreak\ mclean\_\allowbreak{}h1\_\allowbreak{}zero\_\allowbreak{}extension\_\allowbreak{}localization}: \texttt{mclean\_\allowbreak{}h1\_\allowbreak{}zero\_\allowbreak{}extension\_\allowbreak{}localization\_\allowbreak{}claim\allowbreak\ TYPE('n::finite)} holds.
\item \texttt{trace:\allowbreak\ mclean\_\allowbreak{}trace\_\allowbreak{}dirichlet\_\allowbreak{}green\_\allowbreak{}v3}: for every \texttt{U} and every \texttt{a}, \texttt{mclean\_\allowbreak{}trace\_\allowbreak{}dirichlet\_\allowbreak{}green\_\allowbreak{}claim\_\allowbreak{}v3\allowbreak\ (U\allowbreak\ ::\allowbreak\ 'n\allowbreak\ cgu\_\allowbreak{}point\allowbreak\ set)\allowbreak\ a} holds.
\item \texttt{graph\_\allowbreak{}source.\allowbreak{}annihilator.\allowbreak{}volume\_\allowbreak{}gradient.\allowbreak{}volume.\allowbreak{}smooth:\allowbreak\ mclean\_\allowbreak{}smooth\_\allowbreak{}test\_\allowbreak{}h1\_\allowbreak{}zero}: \texttt{mclean\_\allowbreak{}smooth\_\allowbreak{}test\_\allowbreak{}h1\_\allowbreak{}zero\_\allowbreak{}claim\allowbreak\ TYPE('n::finite)} holds.
\item \texttt{graph\_\allowbreak{}source.\allowbreak{}annihilator.\allowbreak{}volume\_\allowbreak{}gradient.\allowbreak{}volume.\allowbreak{}injectivity:\allowbreak\ hormander\_\allowbreak{}regular\_\allowbreak{}distribution\_\allowbreak{}injectivity}: \texttt{hormander\_\allowbreak{}regular\_\allowbreak{}distribution\_\allowbreak{}injectivity\_\allowbreak{}claim\allowbreak\ TYPE('n::finite)} holds.
\item \texttt{graph\_\allowbreak{}source.\allowbreak{}annihilator.\allowbreak{}volume\_\allowbreak{}gradient.\allowbreak{}locality:\allowbreak\ hormander\_\allowbreak{}weak\_\allowbreak{}gradient\_\allowbreak{}locality}: for every \texttt{X}, every \texttt{u}, and every \texttt{Du}, \texttt{hormander\_\allowbreak{}weak\_\allowbreak{}gradient\_\allowbreak{}locality\_\allowbreak{}claim\allowbreak\ (X\allowbreak\ ::\allowbreak\ 'n\allowbreak\ ucp\_\allowbreak{}point\allowbreak\ set)\allowbreak\ u\allowbreak\ Du} holds.
\item \texttt{graph\_\allowbreak{}source.\allowbreak{}annihilator.\allowbreak{}localization:\allowbreak\ mclean\_\allowbreak{}h1\_\allowbreak{}zero\_\allowbreak{}extension\_\allowbreak{}localization}: \texttt{mclean\_\allowbreak{}h1\_\allowbreak{}zero\_\allowbreak{}extension\_\allowbreak{}localization\_\allowbreak{}claim\allowbreak\ TYPE('n::finite)} holds.
\item \texttt{graph\_\allowbreak{}source.\allowbreak{}annihilator.\allowbreak{}graph\_\allowbreak{}support:\allowbreak\ mclean\_\allowbreak{}graph\_\allowbreak{}supported\_\allowbreak{}h1\_\allowbreak{}zero}: \texttt{mclean\_\allowbreak{}graph\_\allowbreak{}supported\_\allowbreak{}h1\_\allowbreak{}zero\_\allowbreak{}claim\allowbreak\ TYPE('n::finite)} holds.
\item \texttt{graph\_\allowbreak{}source.\allowbreak{}bidual:\allowbreak\ mclean\_\allowbreak{}h1\_\allowbreak{}zero\_\allowbreak{}bidual\_\allowbreak{}representation\_\allowbreak{}v2}: \texttt{mclean\_\allowbreak{}h1\_\allowbreak{}zero\_\allowbreak{}bidual\_\allowbreak{}representation\_\allowbreak{}claim\_\allowbreak{}v2\allowbreak\ TYPE('n::finite)} holds.
\end{itemize}

The accompanying \texttt{mclean\_\allowbreak{}smooth\_\allowbreak{}generator\_\allowbreak{}ambient\_\allowbreak{}lift\_\allowbreak{}v2} part postulates \texttt{mclean\_\allowbreak{}smooth\_\allowbreak{}generator\_\allowbreak{}ambient\_\allowbreak{}lift\_\allowbreak{}claim\_\allowbreak{}v2\allowbreak\ TYPE('n::finite)\allowbreak\ TYPE('c::finite)\allowbreak\ TYPE('h::finite)}, where each \texttt{TYPE(.\allowbreak{}.\allowbreak{}.\allowbreak{})} expression denotes only its selected finite carrier type and the arguments occur in exactly the displayed order.

The named \texttt{offset} component contributes \texttt{offset.\allowbreak{}mclean\_\allowbreak{}hyperplane\_\allowbreak{}offset\_\allowbreak{}source}: \texttt{mclean\_\allowbreak{}hyperplane\_\allowbreak{}offset\_\allowbreak{}source\_\allowbreak{}claim\allowbreak\ TYPE('m::finite)\allowbreak\ TYPE('n::finite)} holds, with these two selected finite carrier types in that order. The named \texttt{boundary\_\allowbreak{}metric} component fixes its \texttt{dimension\_\allowbreak{}type\allowbreak\ ::\allowbreak\ "'n::finite\allowbreak\ itself"} at \texttt{TYPE('n)} and contributes \texttt{boundary\_\allowbreak{}metric.\allowbreak{}kang\_\allowbreak{}yun\_\allowbreak{}ltu\_\allowbreak{}local\_\allowbreak{}boundary\_\allowbreak{}metric\_\allowbreak{}v6}: for every \texttt{M}, \texttt{P1}, \texttt{P2}, \texttt{U}, \texttt{gamma}, \texttt{normal}, \texttt{offset}, \texttt{g1}, \texttt{g2}, \texttt{a1}, \texttt{a2}, \texttt{rho1}, and \texttt{rho2},

\texttt{kang\_\allowbreak{}yun\_\allowbreak{}ltu\_\allowbreak{}local\_\allowbreak{}boundary\_\allowbreak{}metric\_\allowbreak{}claim\_\allowbreak{}v6\allowbreak\ (M\allowbreak\ ::\allowbreak\ ('n,\allowbreak\ 'c::finite,\allowbreak\ 'h::finite)\allowbreak\ cgu\_\allowbreak{}model)\allowbreak\ P1\allowbreak\ P2\allowbreak\ U\allowbreak\ gamma\allowbreak\ normal\allowbreak\ offset\allowbreak\ g1\allowbreak\ g2\allowbreak\ a1\allowbreak\ a2\allowbreak\ rho1\allowbreak\ rho2}

holds, with exactly this type annotation, argument order, and quantifier scope.

\clearpage
\dossiersection{Statements}
\dossiersubhead{Manuscript conclusions}

\dossierentry{MAIN-001: \texttt{\detokenize{cgu_global_uniqueness_claim_v7}}}
\reviewlabel{English original and manuscript locator}
Theorem \texttt{thm:main} of \cite{Carstea2026Polynomial} states:

Let \(n\ge3\), let \(\Omega\subset\R^n\) be a bounded Lipschitz domain, and let
\(\mathcal D=\{D_\alpha\}_{\alpha\in A}\) be a known finite Lipschitz
subdivision of \(\Omega\) that is internally flat up to a skeleton.  Let
\(N_\alpha\ge0\).  Suppose that
\(\sigma^{(1)},\sigma^{(2)}\in\mathcal P(\mathcal D,\{N_\alpha\})\), and assume
that, at every recovery step and for \(j=1,2\), the polynomial extension of
\(\sigma^{(j)}|_{D_{\alpha_r}}\) is positive definite on the triple-intersection
sets \(W_{r,ijk}\) appearing in Definition \textup{[source reference: \texttt{def:stripping-order}]}.  We also
assume that, for \(j=1,2\), the polynomial extensions in the recovered cells adjacent to the
auxiliary exterior half-balls remain elliptic in those half-balls whenever the
half-balls are used.

Let \(\Sigma\subset\partial\Omega\) be a nonempty relatively open measured set.
Assume that the subdivision admits an admissible flat-face recovery order
relative to \(\Sigma\), in the sense of Definition \textup{[source reference: \texttt{def:stripping-order}]}.
If
\[
  \DN^\Omega_{\sigma^{(1)},\Sigma}
  =\DN^\Omega_{\sigma^{(2)},\Sigma},
\]
then
\[
  \sigma^{(1)}=\sigma^{(2)}\quad\text{in }\Omega.
\]
Equivalently, the polynomial matrices defining the two conductivities agree on
every cell of the subdivision.
\reviewlabel{Isabelle code}
\begin{lstlisting}
definition cgu_global_uniqueness_claim_v7 ::
  "('n::finite, 'c::finite, 'h::finite) cgu_model \<Rightarrow>
    ('c \<Rightarrow> 'n cgu_matrix_polynomial) \<Rightarrow>
    ('c \<Rightarrow> 'n cgu_matrix_polynomial) \<Rightarrow> bool"
where
  "cgu_global_uniqueness_claim_v7 M P1 P2 \<longleftrightarrow>
    (3 \<le> CARD('n) \<and>
     cgu_graph_subdivision_core_v7 M \<and>
     cgu_piecewise_polynomial_class M P1 \<and>
     cgu_piecewise_polynomial_class M P2 \<and>
     cgu_recovery_geometry_operational_v6 M \<and>
     cgu_recovery_extension_admissible M P1 \<and>
     cgu_recovery_extension_admissible M P2 \<and>
     cgu_local_dn_equal M P1 P2 (cgu_domain M) (cgu_measured M))
    \<longrightarrow> P1 = P2"
\end{lstlisting}
\reviewlabel{English mathematical translation of the Isabelle code}
For finite \(I,C,H\), geometric data \(M\), and coefficient families \(P^{(1)},P^{(2)}\), the following implication holds.  If

\begin{itemize}
\item \(n\ge3\);
\item \(M\) is an admissible subdivision;
\item \(P^{(1)}\) and \(P^{(2)}\) are admissible coefficient families;
\item \(M\) has an admissible recovery order;
\item both families satisfy the coefficient-access condition;
\item the two families have the same stipulated local boundary-energy data on \((\Omega,\Gamma)\);
\end{itemize}

then

\[
P^{(1)}=P^{(2)}
\]

as cell-indexed polynomial representations.

\dossierentry{MAIN-002: \texttt{\detokenize{cgu_uniform_degree_admissible_injectivity_claim_v7}}}
\reviewlabel{English original and manuscript locator}
The corollary of \cite{Carstea2026Polynomial} states:

Suppose \(N_\alpha\le N\) for all cells.  If the subdivision admits an
admissible flat-face recovery order relative to \(\Sigma\) in which each cell is
accessible through at least \(N+3\) flat faces whose supporting hyperplanes
satisfy the corresponding genericity and triple-intersection geometry, and the
recovered regions satisfy the stated face-connectedness assumptions, then the
following holds.  For any
\(\sigma^{(1)},\sigma^{(2)}\in\mathcal P(\mathcal D,N)\) that satisfy the
coefficient-dependent triple-region positivity and auxiliary half-ball
ellipticity hypotheses of Theorem \textup{[source reference: \texttt{thm:main}]} for this recovery order,
equality of their local DN maps on \(\Sigma\) implies
\(\sigma^{(1)}=\sigma^{(2)}\).  Equivalently, the local DN map is injective on
the subclass of \(\mathcal P(\mathcal D,N)\) satisfying those
coefficient-dependent admissibility conditions.
\reviewlabel{Isabelle code}
\begin{lstlisting}
definition cgu_uniform_degree_admissible_injectivity_claim_v7 ::
  "nat \<Rightarrow> ('n::finite, 'c::finite, 'h::finite) cgu_model \<Rightarrow> bool"
where
  "cgu_uniform_degree_admissible_injectivity_claim_v7 N M \<longleftrightarrow>
    (3 \<le> CARD('n) \<and>
     cgu_graph_subdivision_core_v7 M \<and>
     (\<forall>c. cgu_degree M c = N) \<and>
     cgu_recovery_geometry_operational_v6 M)
    \<longrightarrow>
    (\<forall>P1 P2.
      cgu_piecewise_polynomial_class M P1 \<and>
      cgu_piecewise_polynomial_class M P2 \<and>
      cgu_recovery_extension_admissible M P1 \<and>
      cgu_recovery_extension_admissible M P2 \<and>
      cgu_local_dn_equal M P1 P2 (cgu_domain M) (cgu_measured M)
      \<longrightarrow> P1 = P2)"
\end{lstlisting}
\reviewlabel{English mathematical translation of the Isabelle code}
Fix \(N\in\mathbb N\).  If

\[
\begin{gathered}
n\ge3,\qquad M\text{ is an admissible subdivision},\\
d_c=N\ (c\in C),\qquad M\text{ has an admissible recovery order},
\end{gathered}
\]

then, for every \(P^{(1)},P^{(2)}\), admissibility of both families, coefficient access for both, and equality of their stipulated local boundary-energy data on \((\Omega,\Gamma)\) imply

\[
P^{(1)}=P^{(2)}.
\]

\clearpage
\dossiersubhead{Cited results}

\dossierentry{\texttt{\detokenize{ANALYTIC-ELLIPTIC}}}
\noindent\textit{Formal identifier:} \texttt{\detokenize{hormander_analytic_elliptic_ucp_claim}}\par
\noindent\textit{Relationship to source:} derived.\par\smallskip
\reviewlabel{Source attribution and locator}
Hörmander \cite{Hormander1983}, Section 8.3, printed p. 271, equations (8.3.1)--(8.3.4) and Corollary 8.3.2.
\reviewlabel{Complete source theorem statement (faithful mathematical restatement)}
For \(P(x,D)=\sum_{|\alpha|\leq m}a_\alpha(x)D^\alpha\), the principal symbol is the homogeneous degree-\(m\) part \(p_m(x,\xi)=\sum_{|\alpha|=m}a_\alpha(x)\xi^\alpha\).  The characteristic set consists of the nonzero covectors \((x,\xi)\) for which \(p_m(x,\xi)=0\).  The operator is elliptic precisely when this set is empty, equivalently when \(p_m(x,\xi)\neq0\) for every \(\xi\neq0\).
\reviewlabel{Source attribution and locator}
Hörmander \cite{Hormander1983}, Theorem 8.6.1, printed p. 306.
\reviewlabel{Complete source theorem statement (faithful mathematical restatement)}
If the coefficients of a differential operator \(P\) are real analytic, then an analytic singularity of a distribution \(u\) which is not already an analytic singularity of \(Pu\) can occur only at a characteristic covector; in standard notation,

\[
  \operatorname{WF}_A(u)\subseteq
  \operatorname{WF}_A(Pu)\cup\operatorname{Char}P.
\]
\reviewlabel{Source attribution and locator}
Hörmander \cite{Hormander1983}, Theorem 8.6.5 and its stated elliptic consequence, printed p. 309.
\reviewlabel{Complete source theorem statement (faithful mathematical restatement)}
If \(u\) is a distribution on an open set \(X\) and \(Pu=0\), then every exterior conormal to \(\operatorname{supp}u\) is characteristic for \(P\).  Hence, if \(P\) is elliptic, the support of \(u\) has no boundary point in \(X\).  Consequently, when \(X\) is connected, a solution which vanishes in a neighborhood of one point vanishes throughout \(X\).
\reviewlabel{Source attribution and locator}
Hörmander \cite{Hormander1983}, Theorem 1.2.5, printed p. 15.
\reviewlabel{Complete source theorem statement (faithful mathematical restatement)}
If \(f,g\in L^1_{\mathrm{loc}}(X)\) and \(\int_X f\phi=\int_X g\phi\) for every \(\phi\in C_c^\infty(X)\), then \(f=g\) almost everywhere on \(X\).
\reviewlabel{Source attribution and locator}
Hörmander \cite{Hormander1983}, regular-distribution identification on printed p. 37 and Definition 3.1.1 on printed p. 55.
\reviewlabel{Complete source theorem statement (faithful mathematical restatement)}
Printed p. 37 identifies \(L^1_{\mathrm{loc}}(X)\) modulo almost-everywhere equality injectively with the regular distributions \(\phi\mapsto\int_X f\phi\).  Definition 3.1.1, printed p. 55, defines \((\partial_i u)(\phi)=-u(\partial_i\phi)\) and multiplication of a distribution by a smooth function.
\reviewlabel{Exact derived mathematical result formalized}
Let \(X\subset\mathbb R^n\), \(n\geq1\), be open and connected.  Let \(A(x)\) be a symmetric matrix polynomial which is uniformly positive definite on an open neighborhood of \(\overline X\).  If a real locally square-integrable regular distribution \(u\) is a weak solution of \(-\operatorname{div}(A\nabla u)=0\) on \(X\), and \(u=0\) almost everywhere on some nonempty open subset of \(X\), then \(u=0\) almost everywhere on \(X\).  This is the exact mathematical content of the formal interface.

For the notation correspondence, polynomial entries are real analytic and the principal symbol of \(-\operatorname{div}(A\nabla\cdot)\), in Hörmander's \(D_j=(1/i)\partial_j\) convention, is \(\xi^{T}A(x)\xi\).  Uniform positive definiteness makes it nonzero for \(\xi\neq0\).  The weak variational identity is the distribution equation, and the regular-distribution results convert distributional vanishing to almost-everywhere vanishing.  Polynomiality and uniform ellipticity near the closure are conservative strengthenings of local analyticity and ellipticity.
\reviewlabel{Isabelle code}
\begin{lstlisting}
definition hormander_analytic_elliptic_ucp_claim ::
  "'n::finite ucp_point set \<Rightarrow> 'n ucp_coefficient \<Rightarrow>
   ('n ucp_point \<Rightarrow> real) \<Rightarrow> bool"
where
  "hormander_analytic_elliptic_ucp_claim X P u \<longleftrightarrow>
    ((0 < CARD('n) \<and> open X \<and> connected X \<and>
      ucp_matrix_polynomial P \<and>
      (\<forall>x. ucp_symmetric_matrix (P x)) \<and>
      ucp_uniformly_elliptic_near X P \<and>
      ucp_weak_solution_on X P u \<and>
      (\<exists>V. V \<noteq> {} \<and> open V \<and> V \<subseteq> X \<and>
        (AE x in lborel. x \<in> V \<longrightarrow> u x = 0)))
     \<longrightarrow> (AE x in lborel. x \<in> X \<longrightarrow> u x = 0))"
\end{lstlisting}
\reviewlabel{English mathematical translation of the Isabelle code}
Let \(X\subseteq V_I\), \(P:V_I\to\mathbb R^{I\times I}\), and \(u:V_I\to\mathbb R\).  If

\begin{itemize}
\item \(n>0\);
\item \(X\) is open and connected;
\item \(P\) is entrywise polynomial and pointwise symmetric;
\item \(P\) is elliptic around \(X\);
\item \(u\) is a local weak solution for \(P\) on \(X\);
\item there exists a nonempty open \(V\subseteq X\) such that, for ambient-a.e. \(x\),
\end{itemize}

\[
  x\in V\Longrightarrow u(x)=0;
  \]

then, for ambient-a.e. \(x\),

\[
x\in X\Longrightarrow u(x)=0.
\]

\dossierentry{\texttt{\detokenize{WEAK-GRADIENT-LOCALITY}}}
\noindent\textit{Formal identifier:} \texttt{\detokenize{hormander_weak_gradient_locality_claim}}\par
\noindent\textit{Relationship to source:} derived.\par\smallskip
\reviewlabel{Source attribution and locator}
Hörmander \cite{Hormander1983}, Definitions 1.2.1--1.2.2, the following zero-extension observation, and Lemma 1.2.3, printed p. 14.
\reviewlabel{Complete source theorem statement (faithful mathematical restatement)}
The test space on an open set \(X\) is \(C_c^\infty(X)\), with support defined as the closure of the nonzero set; its members extend smoothly by zero to the ambient Euclidean space.  The lemma constructs a nonnegative compactly supported smooth bump positive at a prescribed point.  Translation, scaling, a finite compact subcover, and summation therefore give, for every compact \(K\Subset X\), a nonnegative test function which is strictly positive on \(K\).
\reviewlabel{Source attribution and locator}
Hörmander \cite{Hormander1983}, Theorem 1.2.5, printed p. 15.
\reviewlabel{Complete source theorem statement (faithful mathematical restatement)}
Two locally integrable functions with identical pairings against every member of \(C_c^\infty(X)\) are equal almost everywhere.
\reviewlabel{Source attribution and locator}
Hörmander \cite{Hormander1983}, regular-distribution identification on printed p. 37 and Definition 3.1.1 on printed p. 55.
\reviewlabel{Complete source theorem statement (faithful mathematical restatement)}
Printed p. 37 identifies locally integrable functions modulo almost-everywhere equality with their regular distributions.  Definition 3.1.1, printed p. 55, defines the distributional derivative by \((\partial_i u)(\phi)=-u(\partial_i\phi)\).
\reviewlabel{Exact derived mathematical result formalized}
Let \(X\) be open.  Suppose real \(u\) and real vector field \(D u\) satisfy the project's local-square-integrability and weak-gradient predicates on \(X\), and suppose \(u=0\) almost everywhere on \(X\).  Then \(D u=0\) almost everywhere on \(X\).

Coordinatewise, the weak identity gives \(\int_X u\,\partial_i\phi=-\int_X(Du)_i\phi\); the left side is zero. The weak-gradient product-integrability clause, combined with a positive localized test, proves \((Du)_i\in L^1_{\mathrm{loc}}(X)\).  Theorem 1.2.5 then makes each component zero almost everywhere, and finite dimensionality combines the component conclusions.  The formal real test class is a conservative specialization of Hörmander's complex convention: real and imaginary parts recover all complex pairings.
\reviewlabel{Isabelle code}
\begin{lstlisting}
definition hormander_weak_gradient_locality_claim ::
  "'n::finite ucp_point set \<Rightarrow>
   ('n ucp_point \<Rightarrow> real) \<Rightarrow>
   ('n ucp_point \<Rightarrow> 'n ucp_point) \<Rightarrow> bool"
where
  "hormander_weak_gradient_locality_claim X u Du \<longleftrightarrow>
    ((open X \<and>
      ucp_locally_square_integrable_on X u \<and>
      ucp_locally_square_integrable_on X Du \<and>
      ucp_weak_gradient_on X u Du \<and>
      (AE x in lborel. x \<in> X \<longrightarrow> u x = 0))
     \<longrightarrow> (AE x in lborel. x \<in> X \<longrightarrow> Du x = 0))"
\end{lstlisting}
\reviewlabel{English mathematical translation of the Isabelle code}
If \(X\) is open, \(u\) and \(G\) are locally square integrable in the stipulated sense on \(X\), \(G\) satisfies the weak-gradient relation for \(u\), and, for ambient-a.e. \(x\),

\[
x\in X\Longrightarrow u(x)=0,
\]

then, for ambient-a.e. \(x\),

\[
x\in X\Longrightarrow G(x)=0.
\]

\dossierentry{\texttt{\detokenize{HORMANDER-REGULAR-DISTRIBUTION-INJECTIVITY}}}
\noindent\textit{Formal identifier:} \texttt{\detokenize{hormander_regular_distribution_injectivity_claim}}\par
\noindent\textit{Relationship to source:} direct.\par\smallskip
\reviewlabel{Source attribution and locator}
Hörmander \cite{Hormander1983}, Definitions 1.2.1--1.2.2 and Lemma 1.2.3, printed p. 14.
\reviewlabel{Complete source theorem statement (faithful mathematical restatement)}
These results define the smooth compactly supported tests on an open set, their support, and the localized smooth bumps used to separate a nonzero locally integrable function from the zero distribution.
\reviewlabel{Source attribution and locator}
Hörmander \cite{Hormander1983}, Theorem 1.2.5, printed p. 15.
\reviewlabel{Complete source theorem statement (faithful mathematical restatement)}
For an open \(X\subset\mathbb R^n\), if \(f,g\in L^1_{\mathrm{loc}}(X)\) satisfy \(\int_X f\phi=\int_Xg\phi\) for every \(\phi\in C_c^\infty(X)\), then \(f=g\) almost everywhere on \(X\).
\reviewlabel{Source attribution and locator}
Hörmander \cite{Hormander1983}, regular-distribution identification on printed p. 37.
\reviewlabel{Complete source theorem statement (faithful mathematical restatement)}
The source expresses the preceding result as injectivity of the regular-distribution embedding of \(L^1_{\mathrm{loc}}(X)\) modulo almost-everywhere equality.

\textbf{Direct formal specialization.}

For every open \(X\) and real \(u\in L^1(X)\), if \(u\phi\) is integrable and \(\int_Xu\phi=0\) for every project smooth compact test \(\phi\), then \(u=0\) almost everywhere on \(X\).  Global \(L^1(X)\) is stronger than the source's local integrability, the comparison function is specialized to zero, and the project test predicate is the globally zero-extended realization of \(C_c^\infty(X)\).  The explicit product-integrability premise prevents use of Isabelle's totalized integral outside the ordinary integrable case.
\reviewlabel{Isabelle code}
\begin{lstlisting}
definition hormander_regular_distribution_injectivity_claim ::
  "'n::finite itself \<Rightarrow> bool"
where
  "hormander_regular_distribution_injectivity_claim dimension_type \<longleftrightarrow>
    (\<forall>X :: 'n cgu_point set. \<forall>u.
      open X \<and> hormander_regular_distribution_test_zero_on X u
      \<longrightarrow>
      (AE x in restrict_space lborel X. u x = 0))"
\end{lstlisting}
\reviewlabel{English mathematical translation of the Isabelle code}
For every \(X\subseteq V_I\) and \(u:V_I\to\mathbb R\), if \(X\) is open, \(u\) is integrable on \(X\), and every \(\phi\in\mathcal D(X)\) has \(u\phi\) integrable with

\[
\int_Xu\phi\,dx=0,
\]

then

\[
u(x)=0
\quad\text{for a.e. }x
\text{ with respect to Lebesgue measure restricted to }X.
\]

\dossierentry{\texttt{\detokenize{EVANS-POINCARE-ELLIPTIC-GREEN}}}
\noindent\textit{Formal identifier:} \texttt{\detokenize{evans_poincare_elliptic_green_claim}}\par
\noindent\textit{Relationship to source:} derived.\par\smallskip
\reviewlabel{Source attribution and locator}
Evans \cite{Evans2010}, Section 5.6.1, Theorem 3 and its separately displayed particular consequence, printed pp. 279--280.
\reviewlabel{Complete source theorem statement (faithful mathematical restatement)}
If \(U\subset\mathbb R^n\) is bounded and open and \(u\in W^{1,p}_0(U)\), \(1\leq p<\infty\), then \(\|u\|_{L^p(U)}\leq C\|Du\|_{L^p(U)}\), with \(C\) depending on \(U,n,p\).  The interface uses only \(p=2\).
\reviewlabel{Source attribution and locator}
Evans \cite{Evans2010}, Section 5.9.1, Theorem 1, printed pp. 299--300.
\reviewlabel{Complete source theorem statement (faithful mathematical restatement)}
The negative Sobolev space \(H^{-1}(U)\) is the bounded dual of \(H^1_0(U)\).  Equivalently, every such functional has a representation by \(L^2\) functions \(f=f^0-\sum_{i=1}^n\partial_i f^i\), and conversely every such expression defines a bounded functional on \(H^1_0(U)\).
\reviewlabel{Source attribution and locator}
Evans \cite{Evans2010}, Sections 6.1.1--6.1.2, printed pp. 311--315.
\reviewlabel{Complete source theorem statement (faithful mathematical restatement)}
On a bounded open domain, a divergence-form operator with bounded measurable leading coefficients, uniformly elliptic almost everywhere (and symmetric in the symmetric case), gives the bilinear energy form on \(H^1_0(U)\).  A weak zero-Dirichlet solution is an element \(u\in H^1_0(U)\) satisfying that bilinear identity against every \(v\in H^1_0(U)\).
\reviewlabel{Source attribution and locator}
Evans \cite{Evans2010}, Section 6.2.1, Theorem 1 (Lax--Milgram), printed pp. 315--317.
\reviewlabel{Complete source theorem statement (faithful mathematical restatement)}
If \(H\) is a real Hilbert space and a bilinear form \(B:H\times H\to\mathbb R\) is bounded and coercive, \(|B[u,v]|\leq M\|u\|\|v\|\) and \(B[u,u]\geq\beta\|u\|^2\) with \(\beta>0\), then every \(f\in H^*\) has a unique \(u\in H\) satisfying \(B[u,v]=f(v)\) for all \(v\in H\), and \(\|u\|\leq\|f\|/\beta\).
\reviewlabel{Source attribution and locator}
Evans \cite{Evans2010}, Section 6.2.2, Theorem 2, printed pp. 317--318.
\reviewlabel{Complete source theorem statement (faithful mathematical restatement)}
Under the elliptic coefficient hypotheses, the weak bilinear form satisfies \(|B[u,v]|\leq C\|u\|_{H^1_0(U)}\|v\|_{H^1_0(U)}\) and a Gårding bound \(B[u,u]\geq\beta\|u\|_{H^1_0(U)}^2-\gamma\|u\|_{L^2(U)}^2\) for constants \(C,\beta>0\) and \(\gamma\geq0\).  Thus the leading elliptic part controls the gradient norm, while lower-order terms cost only an \(L^2\) term.
\reviewlabel{Source attribution and locator}
Evans \cite{Evans2010}, Section 6.2.2, Theorem 3, printed pp. 318--319.
\reviewlabel{Complete source theorem statement (faithful mathematical restatement)}
There is \(\gamma\geq0\) such that, for \(f\in L^2(U)\) and every \(\mu\geq\gamma\), the shifted zero-Dirichlet weak problem \(Lu+\mu u=f\) has a unique solution in \(H^1_0(U)\), with the corresponding energy estimate.  This theorem does not itself state the unshifted arbitrary-\(H^{-1}\) result below.
\reviewlabel{Exact derived mathematical result formalized}
Let \(U\neq\varnothing\) be bounded and open, and let \(a(x)\) be a chosen restricted-Borel measurable, pointwise bounded, pointwise symmetric matrix with one positive pointwise uniform ellipticity constant.  Then the pure leading energy form is bounded and strictly coercive on the project's \(H^1_0(U)\) quotient.  Moreover, one constant \(C\geq0\), depending only on \(U,a\), works for every bounded real linear source \(F\) with bound \(K\): there is a zero-boundary weak solution of \(-\operatorname{div}(a\nabla u)=F\), it satisfies \(\|u\|_{H^1(U)}\leq CK\), and it is unique modulo zero \(H^1\) distance.

The derivation takes \(H=H^1_0(U)\) and \(B[u,v]=\int_U(aDu)\cdot Dv\).  Coefficient boundedness gives continuity; ellipticity plus the \(p=2\) Poincaré inequality gives coercivity in the full project \(H^1\) norm; Section 5.9.1 identifies the source with an element of \(H^*\); Lax--Milgram gives the unshifted solution, uniqueness, and the uniform estimate.  The pointwise representative assumptions conservatively strengthen Evans's almost-everywhere coefficient assumptions.
\reviewlabel{Isabelle code}
\begin{lstlisting}
definition evans_poincare_elliptic_green_claim ::
  "'n::finite cgu_point set \<Rightarrow> 'n mclean_coefficient \<Rightarrow> bool"
where
  "evans_poincare_elliptic_green_claim U a \<longleftrightarrow>
    (0 < CARD('n) \<and> U \<noteq> {} \<and> open U \<and>
      evans_bounded_domain U \<and> evans_elliptic_coefficient_on U a)
    \<longrightarrow>
    (mclean_variational_form_on U a \<and>
      (\<exists>C. 0 \<le> C \<and>
        (\<forall>F K. mclean_source_bound_on U F K \<longrightarrow>
          (\<exists>u Du. mclean_green_solution_for U a F u Du \<and>
            mclean_h1_norm U u Du \<le> C * K \<and>
            (\<forall>v Dv. mclean_green_solution_for U a F v Dv
              \<longrightarrow>
              cgu_h1_squared_distance U u Du v Dv = 0)))))"
\end{lstlisting}
\reviewlabel{English mathematical translation of the Isabelle code}
Let \(U\subseteq V_I\) and let \(a\) be a matrix field.  If

\[
n>0,\qquad
U\ne\varnothing,\qquad
U\text{ is open and bounded},
\]

and \(a\) has the bounded symmetric uniformly elliptic property on \(U\), then:

\begin{enumerate}[start=1]
\item \(B_{a,U}\) has the bounded/coercive pair property;
\item there is one \(C\ge0\) such that, whenever \(\mathcal F\) is \(K\)-bounded-linear on zero-boundary pairs, a variational solution \((u,G)\) exists with
\end{enumerate}

\[
   \|(u,G)\|_U\le CK,
   \]

and every other variational solution \((v,H)\) satisfies

\[
   d_U^2\bigl((u,G),(v,H)\bigr)=0.
   \]

The second clause is existence plus zero-distance uniqueness, not literal equality of representatives.

\dossierentry{\texttt{\detokenize{KANG-YUN-LTU-LOCAL-BOUNDARY-METRIC}}}
\noindent\textit{Formal identifier:} \texttt{\detokenize{kang_yun_ltu_local_boundary_metric_claim_v6}}\par
\noindent\textit{Relationship to source:} derived.\par\smallskip
\reviewlabel{Source attribution and locator}
Kang--Yun \cite{KangYun2003}, weak energy pairing on printed p. 719 and localized-DN definition on printed p. 723.
\reviewlabel{Complete source theorem statement (faithful mathematical restatement)}
For an open connected measured boundary patch \(\Gamma\), the localized Dirichlet-to-Neumann map takes \(H^{1/2}\) boundary data supported in \(\Gamma\), solves the conductivity/metric equation, and restricts the Neumann response to \(\Gamma\).  Its weak form is the boundary energy pairing of an input solution with an arbitrary lift of the test trace.
\reviewlabel{Source attribution and locator}
Kang--Yun \cite{KangYun2003}, Theorem 1.3, printed p. 723.
\reviewlabel{Complete source theorem statement (faithful mathematical restatement)}
Let two uniformly elliptic Riemannian metrics have \(C^{m,p}\) regularity near a connected open boundary patch \(\Gamma\), with \(m\geq1\) and \(p>0\).  For every compact \(K\Subset\Gamma\), the local DN maps stably determine the metrics near \(K\), modulo a local diffeomorphism which fixes the boundary: after the appropriate pullback, the metric difference is bounded by a positive power of the local-DN operator-norm difference.  In particular, equality of the local maps gives equality up to that boundary-fixing gauge on \(K\).
\reviewlabel{Source attribution and locator}
Kang--Yun \cite{KangYun2003}, Lemmas 2.1 and 2.2, printed pp. 724--727.
\reviewlabel{Complete source theorem statement (faithful mathematical restatement)}
In boundary-normal coordinates and for a boundary point and tangential frequency, the authors construct high-frequency boundary packets supported in an arbitrarily small part of \(\Gamma\), extend them into a shrinking boundary box \(D_N\), and establish the packet, trace, and Sobolev norm bounds used in the asymptotic energy calculation.
\reviewlabel{Source attribution and locator}
Kang--Yun \cite{KangYun2003}, Section 3, equations (3.7)--(3.9) and the intervening Hardy inequality, printed pp. 728--731.
\reviewlabel{Complete source theorem statement (faithful mathematical restatement)}
Writing the exact solution as the explicit packet plus a correction, the source proves the local correction estimate (3.7), the artificial-top trace estimate (3.8), the displayed Hardy inequality before (3.9), and the weighted residual estimate (3.9).  In the notation used by the reviewed derivation, the local estimate is \(\|s_N\|_{H^1(D_N)}\leq C N^{1-|\alpha|/2}\), with the subsequent estimates making the correction terms negligible at the normalization used in the DN limit.
\reviewlabel{Source attribution and locator}
Kang--Yun \cite{KangYun2003}, Theorem 4.1 and equations (4.5)--(4.8), printed pp. 731--733.
\reviewlabel{Complete source theorem statement (faithful mathematical restatement)}
The normalized high-frequency local-DN energy limit at a chosen boundary point determines the tangential metric quadratic form there. Varying the tangential covector and polarizing determines the complete tangential bilinear form.
\reviewlabel{Source attribution and locator}
Kang--Yun \cite{KangYun2003}, boundary-normal-coordinate and diffeomorphism reduction, printed pp. 733--734.
\reviewlabel{Complete source theorem statement (faithful mathematical restatement)}
A smooth positive metric near the observed patch can be put into the boundary-normal form used by the packet computation.  The comparison of the two normal-coordinate charts is a local diffeomorphism whose restriction to the boundary is the identity; hence its differential is the identity on tangent vectors.
\reviewlabel{Source attribution and locator}
Lassas--Taylor--Uhlmann \cite{LassasTaylorUhlmann2003}, remote rough-boundary setting on printed p. 209 and local-symbol discussion in Section 2, printed pp. 210--211.
\reviewlabel{Complete source theorem statement (faithful mathematical restatement)}
In their analytic, Wiener-regular complete-manifold setting, boundary roughness is allowed away from the measured analytic patch, and the boundary symbol is obtained locally from the DN map.  This is corroboration of locality in a different setting, not a theorem proving the graph-domain corollary below.
\reviewlabel{Exact derived mathematical result formalized}
In dimension at least three, let \(U\) be a bounded connected graph-Lipschitz domain contained in the fixed subdivision model.  Let \(\gamma\) be a nonempty relatively open flat boundary patch with a full one-sided collar, disjoint from the subdivision skeleton.  Suppose each of two globally bounded uniformly elliptic piecewise-polynomial coefficients agrees on an open collar of \(\gamma\) with a smooth positive metric-density realization \(a_j=\sqrt{\det g_j}\,g_j^{-1}\).  If the completed real weak local-DN energy pairings agree for every supported input/test trace, with solutions existing for those data, then for every \(x\in\gamma\) and every pair of vectors tangent to the supporting hyperplane, \(g_1(x)(v,w)=g_2(x)(v,w)\).  No normal component or remote coefficient is identified.

The derivation localizes the Kang--Yun packet with a fixed smooth collar cutoff.  Cutoff derivatives occur a fixed positive normal distance from the boundary and are exponentially small; graph-domain trace lifting, Poincaré, and global ellipticity give the correction bound without remote smoothness. Four real bilinear pairings reconstruct the complex packet pairing.  The boundary-fixing gauge has identity tangential differential, so equality after pullback gives the stated covariant tangential equality.  Pointwise shrinking handles a disconnected formal patch.  Neither paper prints this exact graph-domain implication verbatim, and LTU supplies corroboration only.
\reviewlabel{Isabelle code}
\begin{lstlisting}
definition kang_yun_ltu_local_boundary_metric_claim_v6 ::
  "('n::finite, 'c::finite, 'h::finite) cgu_model \<Rightarrow>
    ('c \<Rightarrow> 'n cgu_matrix_polynomial) \<Rightarrow>
    ('c \<Rightarrow> 'n cgu_matrix_polynomial) \<Rightarrow>
    'n cgu_point set \<Rightarrow> 'n cgu_point set \<Rightarrow>
    'n cgu_point \<Rightarrow> real \<Rightarrow>
    'n lu_metric \<Rightarrow> 'n lu_metric \<Rightarrow>
    'n lu_metric \<Rightarrow> 'n lu_metric \<Rightarrow>
    'n lu_density \<Rightarrow> 'n lu_density \<Rightarrow> bool"
where
  "kang_yun_ltu_local_boundary_metric_claim_v6 M P1 P2 U gamma normal offset
      g1 g2 a1 a2 rho1 rho2 \<longleftrightarrow>
    (3 \<le> CARD('n) \<and>
     cgu_graph_subdivision_core_v7 M \<and>
     cgu_piecewise_polynomial_class M P1 \<and>
     cgu_piecewise_polynomial_class M P2 \<and>
     cgu_graph_lipschitz_domain_v5 U \<and>
     U \<subseteq> cgu_domain M \<and>
     cgu_flat_face_patch U normal offset gamma \<and>
     gamma \<inter> cgu_subdivision_edge_corner_set M = {} \<and>
     lu_smooth_metric_extension_near M P1 U gamma g1 a1 rho1 \<and>
     lu_smooth_metric_extension_near M P2 U gamma g2 a2 rho2 \<and>
     cgu_local_dn_equal M P1 P2 U gamma)
    \<longrightarrow> lu_induced_boundary_metrics_equal_on gamma normal g1 g2"
\end{lstlisting}
\reviewlabel{English mathematical translation of the Isabelle code}
Suppose

\begin{itemize}
\item \(n\ge3\);
\item \(M\) is an admissible subdivision;
\item \(P^{(1)},P^{(2)}\) are admissible coefficient families;
\item \(U\) is Lipschitz admissible and \(U\subseteq\Omega\);
\item \((U;\nu,b;\gamma)\) is an oriented flat boundary patch;
\item \(\gamma\cap\Sigma_M=\varnothing\);
\item packages \((g_s,a_s,\rho_s)\) realize \((M,P^{(s)},U,\gamma)\), \(s=1,2\);
\item the two families have the same stipulated local boundary-energy data on \((U,\gamma)\).
\end{itemize}

Then \(g_1\) and \(g_2\) have the same tangential bilinear pairing along \((\gamma,\nu)\).

\dossierentry{\texttt{\detokenize{MCLEAN-GRAPH-SUPPORTED-H1-ZERO}}}
\noindent\textit{Formal identifier:} \texttt{\detokenize{mclean_graph_supported_h1_zero_claim}}\par
\noindent\textit{Relationship to source:} direct.\par\smallskip
\reviewlabel{Source attribution and locator}
McLean \cite{McLean2000}, Definition 3.28, printed pp. 89--90.
\reviewlabel{Complete source theorem statement (faithful mathematical restatement)}
A Lipschitz domain has compact boundary and a finite cover by rigid Cartesian-coordinate neighborhoods in which the domain is a Lipschitz hypograph.  Replacing the local Lipschitz graph functions by continuous functions gives the book's \(C^0\)-domain class; every Definition 3.28 Lipschitz domain is therefore a \(C^0\) domain.
\reviewlabel{Source attribution and locator}
McLean \cite{McLean2000}, Theorem 3.29(ii), printed pp. 91--92.
\reviewlabel{Complete source theorem statement (faithful mathematical restatement)}
For a \(C^0\) domain \(\Omega\) and every real Sobolev order \(s\), the smooth compactly supported class \(\mathcal D(\Omega)\) is dense, in the ambient \(H^s(\mathbb R^n)\) norm, in the global Sobolev space consisting of distributions supported in \(\overline\Omega\).  Equivalently, the ambient compact-test closure equals that supported global Sobolev space.

\textbf{Direct formal specialization.}

At real scalar order \(s=1\), suppose \(U\) is a bounded graph-Lipschitz domain and \((u,Du)\) is a global project \(H^1(\mathbb R^n)\) pair whose potential and selected weak gradient vanish almost everywhere outside \(\overline U\).  Then \((u,Du)\) belongs to the project's sequential \(H^1_0(U)\) compact-test closure.  The graph-Lipschitz premise is stronger than the source's \(C^0\) hypothesis, the gradient-support premise is an additional compatible restriction, and source global convergence is stronger than the conclusion's restricted \(U\)-distance convergence.
\reviewlabel{Isabelle code}
\begin{lstlisting}
definition mclean_graph_supported_h1_zero_claim ::
  "'n::finite itself \<Rightarrow> bool"
where
  "mclean_graph_supported_h1_zero_claim dimension_type \<longleftrightarrow>
    0 < CARD('n) \<and>
    (\<forall>(U :: 'n cgu_point set) u Du.
      cgu_graph_lipschitz_domain_v5 U \<and>
      cgu_h1_pair_on UNIV u Du \<and>
      (AE x in lborel.
        x \<notin> closure U \<longrightarrow> u x = 0 \<and> Du x = 0)
      \<longrightarrow> cgu_h1_zero_pair_on U u Du)"
\end{lstlisting}
\reviewlabel{English mathematical translation of the Isabelle code}
This condition is the conjunction of:

\begin{enumerate}[start=1]
\item \(n>0\);
\item for every \(U\subseteq V_I\) and every pair \((u,G)\in\mathcal H^1_{\mathrm{pair}}(V_I)\), if \(U\) is Lipschitz admissible and, for ambient-a.e. \(x\),
\end{enumerate}

\[
   x\notin\overline U
   \Longrightarrow
   u(x)=0\ \text{and}\ G(x)=0,
   \]

then

\[
   (u,G)\in\mathcal H^1_{0,\mathrm{pair}}(U).
   \]

The support premise is almost-everywhere, not a literal compact-support assertion about the stored representatives.

\dossierentry{\texttt{\detokenize{MCLEAN-H1-ZERO-BIDUAL-REPRESENTATION}}}
\noindent\textit{Formal identifier:} \texttt{\detokenize{mclean_h1_zero_bidual_representation_claim_v2}}\par
\noindent\textit{Relationship to source:} derived.\par\smallskip
\reviewlabel{Source attribution and locator}
McLean \cite{McLean2000}, weak-derivative Sobolev construction and Hilbert conclusion, printed pp. 73--75.
\reviewlabel{Complete source theorem statement (faithful mathematical restatement)}
On every nonempty open \(U\), \(W^r_p(U)\) consists of \(L^p\) functions possessing all weak derivatives through order \(r\) in \(L^p\), with the derivative graph norm.  For \(p=2\), the sum of the \(L^2\) inner products of the weak derivatives induces that norm, and \(W^s(U)\) is a Hilbert space for every real \(s\geq0\).
\reviewlabel{Source attribution and locator}
McLean \cite{McLean2000}, zero-boundary closure notation, printed p. 77.
\reviewlabel{Complete source theorem statement (faithful mathematical restatement)}
The source defines the conventional zero-boundary space as the closure of \(\mathcal D(U)\) in the relevant Sobolev norm.  For this interface the cited page supplies nomenclature only; no equality with a distinct supported Bessel-potential space on an arbitrary open set is used.
\reviewlabel{Source attribution and locator}
McLean \cite{McLean2000}, Hilbert reflexivity statement, printed p. 79.
\reviewlabel{Complete source theorem statement (faithful mathematical restatement)}
Every Hilbert space is reflexive: the canonical isometric embedding \(J:H\to H^{**}\), \((Ju)(F)=F(u)\), is onto.  Thus every bounded linear functional on \(H^*\) is evaluation at some \(u\in H\).
\reviewlabel{Exact derived mathematical result formalized}
For every nonempty open \(U\), let \(H\) be the quotient Hilbert space represented by project \(H^1_0(U)\) potential/weak-gradient pairs.  There exist a set \(E\) and seminorm \(p\) such that \(E\) is exactly the set of total raw functions inducing bounded real linear functionals on \(H\), and \(p(F)\) is the least valid source bound, hence the operator seminorm.  Every carrier-relative real-linear \(g:E\to\mathbb R\) bounded by \(p\) is of the form \(g(F)=F(u,Du)\) for some project \(H^1_0(U)\) pair.

The project's potential-plus-gradient norm is the order-one graph norm, and the compact-test closure is a closed Hilbert subspace after quotienting zero distance.  Bounded raw sources descend to its continuous dual; the least bound is the dual norm.  A \(p\)-bounded \(g\) vanishes on the seminorm kernel, descends to the genuine bidual, and reflexivity represents it by evaluation. The raw carrier \(E\), least-bound plumbing, and displayed evaluation formula are standard consequences, not statements printed verbatim by McLean.
\reviewlabel{Isabelle code}
\begin{lstlisting}
definition mclean_h1_zero_bidual_representation_claim_v2 ::
  "'n::finite itself \<Rightarrow> bool"
where
  "mclean_h1_zero_bidual_representation_claim_v2 dimension_type \<longleftrightarrow>
    (\<forall>U :: 'n cgu_point set.
      U \<noteq> {} \<and> open U
      \<longrightarrow>
      (\<exists>E p. mclean_h1_zero_bidual_realization_on_v2 U E p))"
\end{lstlisting}
\reviewlabel{English mathematical translation of the Isabelle code}
For every \(U\subseteq V_I\), if \(U\ne\varnothing\) and \(U\) is open, then there exist \(E\) and \(p\) such that:

\begin{enumerate}[start=1]
\item \(E\) is precisely the total two-argument functionals that are \(K\)-bounded-linear on \(\mathcal H^1_{0,\mathrm{pair}}(U)\) for some \(K\), and \(E\) is closed under the stated vector operations;
\item \(p\) is the stipulated seminorm, \(p(\mathcal F)\) is an admissible bound for every \(\mathcal F\in E\), and it is no larger than every other admissible bound;
\item every \(p\)-dominated linear \(\ell:E\to\mathbb R\) has a representing pair \((u,G)\in\mathcal H^1_{0,\mathrm{pair}}(U)\) with
\end{enumerate}

\[
   \ell(\mathcal F)=\mathcal F(u,G)
   \qquad(\mathcal F\in E).
   \]

No quotient identifies total functionals that agree only on the zero-boundary pair class, and uniqueness of the representing pair is not asserted.

\dossierentry{\texttt{\detokenize{MCLEAN-H1-ZERO-EXTENSION-LOCALIZATION}}}
\noindent\textit{Formal identifier:} \texttt{\detokenize{mclean_h1_zero_extension_localization_claim}}\par
\noindent\textit{Relationship to source:} derived.\par\smallskip
\reviewlabel{Source attribution and locator}
McLean \cite{McLean2000}, smooth compact tests and the two compact-test closures, printed pp. 61, 65, and 77--78.
\reviewlabel{Complete source theorem statement (faithful mathematical restatement)}
The test class is \(\mathcal D(\Omega)=C_c^\infty(\Omega)\). \(H^s_0(\Omega)\) is its closure in the intrinsic/restriction \(H^s(\Omega)\) norm, whereas \(\widetilde H^s(\Omega)\) is its closure in the ambient \(H^s(\mathbb R^n)\) norm.  These two spaces are not identified on an arbitrary open set in this interface.
\reviewlabel{Source attribution and locator}
McLean \cite{McLean2000}, Theorem 3.6, printed p. 64.
\reviewlabel{Complete source theorem statement (faithful mathematical restatement)}
Given a closed set and an open neighborhood (equivalently, an arbitrarily small prescribed neighborhood), there is a smooth cutoff taking values between zero and one, equal to one on the closed set and zero outside the prescribed neighborhood. For a compact set contained in an open set, the cutoff may be chosen compactly supported in that open set.
\reviewlabel{Source attribution and locator}
McLean \cite{McLean2000}, Theorem 3.16, printed p. 80.
\reviewlabel{Complete source theorem statement (faithful mathematical restatement)}
For every nonnegative integer \(m\), the global weak-derivative space \(W^m(\mathbb R^n)\) and the Bessel-potential space \(H^m(\mathbb R^n)\) are the same set with equivalent norms.
\reviewlabel{Source attribution and locator}
McLean \cite{McLean2000}, Theorem 3.20, printed p. 83.
\reviewlabel{Complete source theorem statement (faithful mathematical restatement)}
Multiplication by a smooth compactly supported cutoff is a bounded linear operator on \(H^s(\Omega)\) and, separately, on \(\widetilde H^s(\Omega)\).  The theorem does not state a multiplier clause for \(H^s_0(\Omega)\) under that name.
\reviewlabel{Exact derived mathematical result formalized}
Two order-one representative statements are formalized.  First, if \(\varnothing\neq U\subseteq W\) are open and \((u,Du)\) is in the project \(H^1_0(U)\) closure, then the zero extensions of both representatives form a project \(H^1_0(W)\) pair.  Second, if nonempty open \(V\subseteq W\), compact \(K\subseteq V\), and a project \(H^1_0(W)\) pair vanishes almost everywhere together with its gradient off \(K\), then the same representatives form a project \(H^1_0(V)\) pair.

For extension, each compact test in \(U\) is also a compact test in \(W\), and zero extension turns intrinsic order-one convergence into global \(W^1\) convergence; Theorem 3.16 supplies the ambient \(H^1\) bridge.  For localization, choose the Theorem 3.6 cutoff equal to one near \(K\); Theorem 3.20 gives continuity on the ambient compact-test closure, and the product equals the original pair because both representatives vanish off \(K\).  If \(K=\varnothing\), the conclusion is the zero class.  No arbitrary-domain identity \(H^1_0=\widetilde H^1\) is used.
\reviewlabel{Isabelle code}
\begin{lstlisting}
definition mclean_h1_zero_extension_localization_claim ::
  "'n::finite itself \<Rightarrow> bool"
where
  "mclean_h1_zero_extension_localization_claim dimension_type \<longleftrightarrow>
    ((\<forall>(U :: 'n cgu_point set) W u Du.
        U \<noteq> {} \<and> open U \<and> open W \<and> U \<subseteq> W \<and>
        cgu_h1_zero_pair_on U u Du
        \<longrightarrow>
        cgu_h1_zero_pair_on W
          (mclean_zero_extension_on U u)
          (mclean_zero_extension_on U Du)) \<and>
     (\<forall>(W :: 'n cgu_point set) V K u Du.
        W \<noteq> {} \<and> V \<noteq> {} \<and> open W \<and> open V \<and> V \<subseteq> W \<and>
        compact K \<and> K \<subseteq> V \<and>
        cgu_h1_zero_pair_on W u Du \<and>
        (AE x in restrict_space lborel W.
          x \<notin> K \<longrightarrow> u x = 0 \<and> Du x = 0)
        \<longrightarrow> cgu_h1_zero_pair_on V u Du))"
\end{lstlisting}
\reviewlabel{English mathematical translation of the Isabelle code}
Both of the following universally quantified implications hold.

\begin{enumerate}[start=1]
\item If \(U\ne\varnothing\), \(U,W\) are open, \(U\subseteq W\), and \((u,G)\in\mathcal H^1_{0,\mathrm{pair}}(U)\), then
\end{enumerate}

\[
   \bigl(
   \operatorname{Ext}_Uu,
   \operatorname{Ext}_UG
   \bigr)
   \in
   \mathcal H^1_{0,\mathrm{pair}}(W).
   \]

\begin{enumerate}[start=2]
\item If \(W,V\ne\varnothing\), \(W,V\) are open, \(V\subseteq W\), \(K\subseteq V\) is compact, \((u,G)\in\mathcal H^1_{0,\mathrm{pair}}(W)\), and, for a.e. \(x\) with respect to Lebesgue measure restricted to \(W\),
\end{enumerate}

\[
   x\notin K
   \Longrightarrow
   u(x)=0\ \text{and}\ G(x)=0,
   \]

then

\[
   (u,G)\in\mathcal H^1_{0,\mathrm{pair}}(V).
   \]

\dossierentry{\texttt{\detokenize{MCLEAN-HYPERPLANE-OFFSET-SOURCE}}}
\noindent\textit{Formal identifier:} \texttt{\detokenize{mclean_hyperplane_offset_source_claim}}\par
\noindent\textit{Relationship to source:} derived.\par\smallskip
\reviewlabel{Source attribution and locator}
McLean \cite{McLean2000}, Theorems 3.37--3.38, printed pp. 102--104.
\reviewlabel{Complete source theorem statement (faithful mathematical restatement)}
At the order used here, the trace is a bounded surjection \(H^1\to H^{1/2}\) on a Lipschitz boundary, has a bounded right inverse, and agrees with ordinary restriction for smooth functions.  By duality, \(H^{-1/2}\) boundary distributions act continuously on \(H^1\) traces. These statements provide the trace and Fourier-estimate setting for the flat-hyperplane specialization.
\reviewlabel{Source attribution and locator}
McLean \cite{McLean2000}, Lemma 3.39, printed pp. 104--105.
\reviewlabel{Complete source theorem statement (faithful mathematical restatement)}
A distribution supported on the hyperplane \(\{x_n=0\}\) has a finite representation by tangential distributions tensored with normal derivatives of the Dirac mass, \(\sum_j v_j\otimes D_n^j\delta_0\).  For membership in ambient \(H^s\), each occurring order satisfies \(s+j<-1/2\) and the tangential coefficient has order \(v_j\in H^{s+j+1/2}(\mathbb R^{n-1})\).
\reviewlabel{Source attribution and locator}
McLean \cite{McLean2000}, equation (3.33), printed pp. 104--105.
\reviewlabel{Complete source theorem statement (faithful mathematical restatement)}
For an admissible term,

\[
 \|v\otimes D_n^j\delta_0\|_{H^s(\mathbb R^n)}^2
   = C_{s,j}\|v\|_{H^{s+j+1/2}(\mathbb R^{n-1})}^2,
\]

with the finite positive constant given by the normal-frequency integral. In particular, at \(s=-1,j=0\), every \(v\in H^{-1/2}(\mathbb R^{n-1})\) defines an ambient \(H^{-1}\) source.
\reviewlabel{Exact derived mathematical result formalized}
Let the ambient dimension be one larger than a positive tangential dimension. Fix an isometric affine hyperplane chart, a unit normal, a smooth compactly supported tangential density \(g\), and \(\delta>0\) such that the full swept support cylinder for offsets \(0\leq t\leq\delta\) lies in an open set \(\Omega\).  Then there are bounded quotient-compatible sources \(S_t\) and one \(C\geq0\) such that, on every smooth test \(\phi\),

\[
 S_t(\phi)=\int g(y)\phi(\operatorname{chart}(y)+t\nu)\,dy,
\]

with \(\|S_t\|_{H^{-1}(\Omega)}\leq C\) and \(\|S_t-S_0\|_{H^{-1}(\Omega)}\leq C\sqrt t\).

Translation in the normal variable multiplies the Fourier transform by \(e^{-2\pi i\xi_nt}\).  In the integral underlying (3.33), bounding \(|e^{-2\pi i\xi_nt}-1|\) by \(\min(2,C|\xi_n|t)\) and splitting at \(|\xi_n|=1/t\) gives a squared \(H^{-1}\) bound of order \(t\), hence the displayed \(\sqrt t\) norm bound. That difference estimate is an elementary consequence of (3.33), not a numbered theorem.  The formal source is an abstract bounded extension on the \(H^1_0\) quotient; only its action on smooth tests is given pointwise.
\reviewlabel{Isabelle code}
\begin{lstlisting}
definition mclean_hyperplane_offset_source_claim ::
  "('m::finite itself) \<Rightarrow> ('n::finite itself) \<Rightarrow> bool"
where
  "mclean_hyperplane_offset_source_claim tangential_type ambient_type \<longleftrightarrow>
    (1 \<le> CARD('m) \<and> CARD('n) = Suc (CARD('m))) \<longrightarrow>
    (\<forall>(Omega :: 'n cgu_point set)
      (chart :: 'm cgu_point \<Rightarrow> 'n cgu_point) normal g delta.
      open Omega \<and> cgu_test_function_on UNIV g \<and>
      mclean_flat_chart chart normal \<and> 0 < delta \<and>
      mclean_flat_offset_cylinder_inside Omega chart normal g delta
      \<longrightarrow>
      (\<exists>source_at C. 0 \<le> C \<and>
        (\<forall>t. 0 \<le> t \<and> t \<le> delta \<longrightarrow>
          mclean_source_extends_flat_offset_action
            Omega chart normal t g (source_at t) \<and>
          mclean_source_bound_on Omega
            (source_at t) C \<and>
          mclean_source_bound_on Omega
            (\<lambda>u Du.
              source_at t u Du - source_at 0 u Du)
            (C * sqrt t))))"
\end{lstlisting}
\reviewlabel{English mathematical translation of the Isabelle code}
If

\[
|J|\ge1,
\qquad
|I|=|J|+1,
\]

then the following holds.  Let \(\Omega\subseteq V_I\) be open, let \((\chi,\nu)\) be a flat isometric chart, let \(g\in\mathcal D(V_J)\), let \(\delta>0\), and assume the support-tube condition.  There exist a family \((\mathcal F_t)_{t\in\mathbb R}\) and one \(C\ge0\) such that, for every \(0\le t\le\delta\),

\[
\mathcal F_t(\phi,\nabla\phi)
=
\int_{V_J}g(x)\phi(\chi(x)+t\nu)\,dx
\qquad(\phi\in\mathcal D(\Omega)),
\]

\(\mathcal F_t\) is \(C\)-bounded-linear on \(\mathcal H^1_{0,\mathrm{pair}}(\Omega)\), and the total functional

\[
(u,G)\longmapsto
\mathcal F_t(u,G)-\mathcal F_0(u,G)
\]

is \(C\sqrt t\)-bounded-linear there.  The same family and the same \(C\) serve every \(t\in[0,\delta]\).

\dossierentry{\texttt{\detokenize{MCLEAN-SMOOTH-GENERATOR-AMBIENT-LIFT}}}
\noindent\textit{Formal identifier:} \texttt{\detokenize{mclean_smooth_generator_ambient_lift_claim_v2}}\par
\noindent\textit{Relationship to source:} derived.\par\smallskip
\reviewlabel{Source attribution and locator}
McLean \cite{McLean2000}, Theorem 3.6, printed p. 64.
\reviewlabel{Complete source theorem statement (faithful mathematical restatement)}
A closed set contained in a prescribed open neighborhood admits a smooth cutoff equal to one on the closed set and supported in that neighborhood; for a compact set in an open set the cutoff can be compactly supported there.
\reviewlabel{Source attribution and locator}
McLean \cite{McLean2000}, Theorems 3.37--3.38, printed pp. 102--104.
\reviewlabel{Complete source theorem statement (faithful mathematical restatement)}
At order one on a Lipschitz domain, the trace \(H^1(U)\to H^{1/2}(\partial U)\) is bounded and surjective, has a continuous right inverse, and agrees with classical boundary restriction on smooth functions.
\reviewlabel{Source attribution and locator}
McLean \cite{McLean2000}, Theorem 3.40, printed pp. 105--106.
\reviewlabel{Complete source theorem statement (faithful mathematical restatement)}
On a \(C^{k-1,1}\) domain, \(H^s_0(U)=H^s(U)\) for \(0\leq s\leq1/2\); for \(1/2<s\leq k\), the zero-boundary space is characterized by vanishing of the applicable boundary traces.  At the used specialization \(s=k=1\), the boundary hypothesis is \(C^{0,1}\) (Lipschitz), and \(H^1_0(U)=\ker(\gamma:H^1(U)\to H^{1/2}(\partial U))\).
\reviewlabel{Exact derived mathematical result formalized}
Let \(\Omega_1\) and \(U\subseteq\Omega_1\) be graph-Lipschitz domains in dimension at least two, let \(\gamma\subseteq\Omega_1\cap\partial U\), and fix one globally smooth compactly supported boundary generator \(p\) whose boundary support lies in \(\gamma\).  Then there is a project \(H^1_0(\Omega_1)\) pair \((z,Dz)\) whose quotient trace on \(U\) equals the trace of \((p,\nabla p)\).

Let \(K\) be the compact closure of the nonzero boundary support of \(p\). If \(K\neq\varnothing\), choose a generator-dependent cutoff \(\chi\in C_c^\infty(\Omega_1)\), \(\chi=1\) near \(K\), and set \(z=\chi p\).  The difference \((1-\chi)p\) has zero boundary trace on \(U\), so Theorem 3.40 puts it in \(H^1_0(U)\).  If \(K=\varnothing\), use the zero ambient pair and apply the same kernel theorem directly.  The witness is chosen after the single generator; no linear, bounded, or simultaneous lift operator is asserted.
\reviewlabel{Isabelle code}
\begin{lstlisting}
definition mclean_smooth_generator_ambient_lift_claim_v2 ::
  "'n::finite itself \<Rightarrow> 'c::finite itself \<Rightarrow>
    'h::finite itself \<Rightarrow> bool"
where
  "mclean_smooth_generator_ambient_lift_claim_v2 ambient_type cell_type
      hyperplane_type \<longleftrightarrow>
    2 \<le> CARD('n) \<longrightarrow>
    (\<forall>(M :: ('n, 'c, 'h) cgu_model)
      (Omega1 :: 'n cgu_point set) U gamma p.
      cgu_graph_lipschitz_domain_v5 Omega1 \<and>
      cgu_graph_lipschitz_domain_v5 U \<and>
      U \<subseteq> Omega1 \<and>
      gamma \<subseteq> Omega1 \<inter> frontier U \<and>
      cgu_boundary_test_on M U gamma p
      \<longrightarrow>
      (\<exists>z Dz.
        cgu_h1_zero_pair_on Omega1 z Dz \<and>
        mclean_same_trace_on U p (cgu_classical_gradient p) z Dz))"
\end{lstlisting}
\reviewlabel{English mathematical translation of the Isabelle code}
If \(n\ge2\), then for every finite cell and hyperplane types, every geometric datum \(M\), and all \(\Omega_1,U,\gamma,p\), the following implication holds.  If

\begin{itemize}
\item \(\Omega_1\) and \(U\) are Lipschitz admissible;
\item \(U\subseteq\Omega_1\);
\item \(\gamma\subseteq\Omega_1\cap\partial U\);
\item \(p\) is boundary-supported relative to \((M,U,\gamma)\);
\end{itemize}

then there exists \((z,Z)\in\mathcal H^1_{0,\mathrm{pair}}(\Omega_1)\) such that

\[
\mathcal R_U\bigl((p,\nabla p),(z,Z)\bigr).
\]

\dossierentry{\texttt{\detokenize{MCLEAN-SMOOTH-TEST-H1-ZERO}}}
\noindent\textit{Formal identifier:} \texttt{\detokenize{mclean_smooth_test_h1_zero_claim}}\par
\noindent\textit{Relationship to source:} derived.\par\smallskip
\reviewlabel{Source attribution and locator}
McLean \cite{McLean2000}, test class, derivative convention, and first-order Sobolev definition, printed pp. 65, 68, and 73--74.
\reviewlabel{Complete source theorem statement (faithful mathematical restatement)}
\(\mathcal D(U)=C_c^\infty(U)\).  Distributional derivatives use the integration-by-parts sign, and \(W^1_2(U)\) consists of square-integrable functions whose first distributional derivatives are square-integrable.  A smooth compactly supported function and its classical first derivatives are square-integrable, and its distributional and classical derivatives agree; Hence every test belongs to \(W^1_2(U)\) with its classical gradient.
\reviewlabel{Source attribution and locator}
McLean \cite{McLean2000}, zero-boundary closure definition, printed p. 77.
\reviewlabel{Complete source theorem statement (faithful mathematical restatement)}
The zero-boundary Sobolev space is the norm closure of \(\mathcal D(U)\).  In particular every generator is in that closure, witnessed by the constant sequence.
\reviewlabel{Source attribution and locator}
McLean \cite{McLean2000}, Theorem 3.16, printed p. 80.
\reviewlabel{Complete source theorem statement (faithful mathematical restatement)}
On \(\mathbb R^n\), the integer-order weak-derivative and Bessel-potential spaces coincide with equivalent norms.  The theorem is used only for the globally smooth compactly supported representative; no equality of full spaces on an arbitrary open set is asserted.
\reviewlabel{Exact derived mathematical result formalized}
For every open \(U\) and every globally smooth real \(\phi\) with compact closed nonzero support contained in \(U\), the pair \((\phi,\nabla\phi)\) is a project \(H^1_0(U)\) pair.  Smoothness and compact support give measurability and square integrability; integration by parts identifies the classical gradient with the weak gradient; the constant sequence \(\phi\) proves membership in the project's compact-test closure. McLean states the definitions for nonempty open sets.  When \(U\) is empty, the formal support condition forces \(\phi=0\), so the conclusion is an elementary zero-function extension, not a source attribution.
\reviewlabel{Isabelle code}
\begin{lstlisting}
definition mclean_smooth_test_h1_zero_claim ::
  "'n::finite itself \<Rightarrow> bool"
where
  "mclean_smooth_test_h1_zero_claim dimension_type \<longleftrightarrow>
    (\<forall>(U :: 'n cgu_point set) phi.
      open U \<longrightarrow>
      cgu_test_function_on U phi \<longrightarrow>
      cgu_h1_zero_pair_on U phi (cgu_classical_gradient phi))"
\end{lstlisting}
\reviewlabel{English mathematical translation of the Isabelle code}
For every \(U\subseteq V_I\) and every \(\phi\), the literal right-associated implication is

\[
U\text{ open}
\Longrightarrow
\left(
\phi\in\mathcal D(U)
\Longrightarrow
(\phi,\nabla\phi)
\in\mathcal H^1_{0,\mathrm{pair}}(U)
\right).
\]

\dossierentry{\texttt{\detokenize{MCLEAN-TRACE-DIRICHLET-GREEN}}}
\noindent\textit{Formal identifier:} \texttt{\detokenize{mclean_trace_dirichlet_green_claim_v3}}\par
\noindent\textit{Relationship to source:} derived.\par\smallskip
\reviewlabel{Source attribution and locator}
McLean \cite{McLean2000}, Definition 3.28, printed pp. 89--90.
\reviewlabel{Complete source theorem statement (faithful mathematical restatement)}
The relevant \(C^{0,1}\) geometry is a domain whose compact boundary is locally a Lipschitz graph after a rigid change of Cartesian coordinates.  The project uses the bounded connected subclass.
\reviewlabel{Source attribution and locator}
McLean \cite{McLean2000}, Theorem 3.37, printed p. 102.
\reviewlabel{Complete source theorem statement (faithful mathematical restatement)}
On a Lipschitz domain, at order one the trace \(\gamma:H^1(U)\to H^{1/2}(\partial U)\) is bounded and surjective and has a continuous right inverse.  Consequently the boundary space is isomorphic, with equivalent norm, to the quotient of \(H^1(U)\) by the trace kernel.
\reviewlabel{Source attribution and locator}
McLean \cite{McLean2000}, Theorem 3.40, printed pp. 105--106.
\reviewlabel{Complete source theorem statement (faithful mathematical restatement)}
On a \(C^{k-1,1}\) domain the zero-boundary space equals all of \(H^s\) for \(0\leq s\leq1/2\), while for \(1/2<s\leq k\) it is characterized by vanishing of the relevant boundary traces.  At \(s=k=1\), used here, the domain is Lipschitz \(C^{0,1}\) and \(H^1_0(U)=\ker\gamma\).
\reviewlabel{Source attribution and locator}
McLean \cite{McLean2000}, Lemma 4.3, printed pp. 116--117.
\reviewlabel{Complete source theorem statement (faithful mathematical restatement)}
Once the right-hand side of a weak elliptic equation is fixed, a weak solution has a unique conormal boundary distribution in \(H^{-1/2}(\partial U)\) satisfying the variational Green formula.  Its norm is bounded in terms of the solution and source norms; the proof uses a continuous trace right inverse.
\reviewlabel{Source attribution and locator}
McLean \cite{McLean2000}, Theorem 4.4, printed p. 118.
\reviewlabel{Complete source theorem statement (faithful mathematical restatement)}
The first Green identity expresses the variational form as the interior source pairing plus the conormal/trace boundary pairing.  Applying it in both orders gives the second Green identity.  For a homogeneous solution, its energy against an arbitrary \(H^1\) lift therefore depends only on the lift's boundary trace.
\reviewlabel{Source attribution and locator}
McLean \cite{McLean2000}, Theorem 4.10, printed pp. 128--130, with the pure Dirichlet specialization on printed p. 131.
\reviewlabel{Complete source theorem statement (faithful mathematical restatement)}
The mixed variational boundary problem on a bounded Lipschitz domain satisfies a Fredholm alternative: in the trivial-kernel branch every admissible datum has a unique solution; otherwise solvability is subject to the stated orthogonality conditions against the adjoint homogeneous kernel.  For pure Dirichlet data the energy space is \(V=H^1_0(U)\) and \(V^*=H^{-1}(U)\), so arbitrary bounded dual sources are allowed.  The formal strict-coercivity premise forces the trivial-kernel branch.
\reviewlabel{Source attribution and locator}
McLean \cite{McLean2000}, equations (4.35)--(4.38), printed p. 145.
\reviewlabel{Complete source theorem statement (faithful mathematical restatement)}
The solution operators furnished by the variational problem are bounded, and the Steklov--Poincaré/Dirichlet-to-Neumann operator obtained from the conormal pairing is a bounded operator between the trace space and its dual.
\reviewlabel{Exact derived mathematical result formalized}
Let \(U\) be a bounded graph-Lipschitz domain in dimension at least two, and let \(a\) define an integrable bilinear energy form which is bounded on project \(H^1(U)\) pairs and strictly coercive on project \(H^1_0(U)\).  Then:

\begin{enumerate}[start=1]
\item every project quotient trace has a weak homogeneous Dirichlet solution, unique modulo zero \(H^1\) distance;
\item one constant bounds the homogeneous solution's energy pairing by the product of the two quotient trace norms, so the pairing is a bounded DN form; and
\item one constant gives, for every bounded real-linear source on \(H^1_0(U)\), a zero-boundary Green solution with an \(H^1\)-norm estimate proportional to the source bound, unique modulo zero distance.
\end{enumerate}

The project trace is the quotient \(H^1(U)/H^1_0(U)\), with norm the infimum of the norms of its lifts.  Theorems 3.37 and 3.40 identify this with the source trace space; Lemma 4.3 and Theorem 4.4 make homogeneous energy depend only on quotient traces; strict coercivity places Theorem 4.10 in its unique pure-Dirichlet branch and yields the solution bounds.  The graph-domain premise is a conservative bounded subclass of McLean's \(C^{0,1}\) geometry.
\reviewlabel{Isabelle code}
\begin{lstlisting}
definition mclean_trace_dirichlet_green_claim_v3 ::
  "'n::finite cgu_point set \<Rightarrow> 'n mclean_coefficient \<Rightarrow> bool"
where
  "mclean_trace_dirichlet_green_claim_v3 U a \<longleftrightarrow>
    (2 \<le> CARD('n) \<and> cgu_graph_lipschitz_domain_v5 U \<and>
      mclean_variational_form_on U a)
    \<longrightarrow>
    ((\<forall>f Df. cgu_h1_pair_on U f Df \<longrightarrow>
      (\<exists>u Du. mclean_dirichlet_solution_for U a f Df u Du \<and>
        (\<forall>v Dv. mclean_dirichlet_solution_for U a f Df v Dv
          \<longrightarrow> cgu_h1_squared_distance U u Du v Dv = 0))) \<and>
     (\<exists>C. 0 \<le> C \<and>
      (\<forall>f Df u Du h Dh.
        mclean_dirichlet_solution_for U a f Df u Du \<and>
        cgu_h1_pair_on U h Dh
        \<longrightarrow>
        mclean_energy_integrable U a Du Dh \<and>
        abs (mclean_energy U a Du Dh) \<le>
          C * mclean_trace_norm U f Df * mclean_trace_norm U h Dh)) \<and>
     (\<exists>C. 0 \<le> C \<and>
      (\<forall>F K. mclean_source_bound_on U F K \<longrightarrow>
        (\<exists>u Du. mclean_green_solution_for U a F u Du \<and>
          mclean_h1_norm U u Du \<le> C * K \<and>
          (\<forall>v Dv. mclean_green_solution_for U a F v Dv
            \<longrightarrow> cgu_h1_squared_distance U u Du v Dv = 0)))))"
\end{lstlisting}
\reviewlabel{English mathematical translation of the Isabelle code}
If

\[
\begin{gathered}
n\ge2,\qquad U\text{ is Lipschitz admissible},\\
B_{a,U}\text{ has the bounded/coercive pair property},
\end{gathered}
\]

then all three conclusions hold.

\begin{enumerate}[start=1]
\item Every \((f,F)\in\mathcal H^1_{\mathrm{pair}}(U)\) has an \(a\)-Dirichlet pair \((u,G)\), and every other such pair \((v,H)\) satisfies
\end{enumerate}

\[
   d_U^2\bigl((u,G),(v,H)\bigr)=0.
   \]

\begin{enumerate}[start=2]
\item There exists one \(C_2\ge0\) such that, whenever \((u,G)\) is an \(a\)-Dirichlet pair for \((f,F)\) and \((h,H)\in\mathcal H^1_{\mathrm{pair}}(U)\), the energy integrand is integrable and
\end{enumerate}

\[
   |B_{a,U}(G,H)|
   \le
   C_2\,q_U(f,F)\,q_U(h,H).
   \]

\begin{enumerate}[start=3]
\item There exists a separately quantified \(C_3\ge0\) such that every \(K\)-bounded-linear total functional \(\mathcal F\) has a variational solution \((u,G)\) with
\end{enumerate}

\[
   \|(u,G)\|_U\le C_3K,
   \]

and every other variational solution \((v,H)\) satisfies

\[
   d_U^2\bigl((u,G),(v,H)\bigr)=0.
   \]

The constants \(C_2\) and \(C_3\) are independent existential witnesses.  Both uniqueness clauses assert zero distance rather than literal equality.

\end{document}